\documentclass{amsart}
\usepackage{graphicx}

\newtheorem{theorem}{Theorem}[section]

\newtheorem{lemma}[theorem]{Lemma}
\theoremstyle{definition}
\newtheorem{definition}[theorem]{Definition}

\theoremstyle{remark}

\numberwithin{equation}{section}

\begin{document}

\title[The spectra of Jacobi operators for CMC tori of revolution 
in $\mathbb{S}^3$]{CMC tori of revolution in $\mathbb{S}^3$: additional 
data on the spectra of their Jacobi operators}

\author{Wayne Rossman}
\address{Department of Mathematics, Kobe University,
Kobe 657-8501, Japan}
\email{wayne@math.kobe-u.ac.jp}

\author{Nahid Sultana}
\address{Department of Mathematics, Kobe University,
Kobe 657-8501, Japan}
\email{nahid@math.kobe-u.ac.jp}


\date{\today}


\begin{abstract}
We prove a theorem about elliptic operators with symmetric 
potential functions, defined on a function space over a closed loop. The result is 
similar to a known result for a function space on an interval with Dirichlet boundary 
conditions. These theorems provide accurate numerical methods for finding the spectra 
of those operators over either type of function space.  
As an application, we numerically compute the Morse index of constant mean curvature 
tori of revolution in the unit $3$-sphere $\mathbb{S}^3$, confirming 
that every such torus has Morse index at least five, and showing that other known 
lower bounds for this Morse index are close to optimal.  
\end{abstract}

\maketitle

\section{Introduction}
\label{Introduction}

Our goal is to study the Morse index of constant mean curvature (CMC) tori 
of revolution in the spherical $3$-space $\mathbb{S}^3$, where 
the Morse index is the number of negative eigenvalues of the Jacobi 
operators of those surfaces.  The central tool 
we use is a result about the number of nodes of eigenfunctions 
of those Jacobi operators.  The result, proven with the standard 
Sturm comparison technique in ordinary differential equations and 
closely related to classically known results, is proven here 
before being applied to the index of CMC surfaces of revolution in 
$\mathbb{S}^3$.  So let us start by considering an operator of the form
\[\mathcal L=- \tfrac{d^2}{dx^2}-V \; , \;\;\; \text{i.e.} \;\; 
\mathcal L(f)=- \tfrac{d^2}{dx^2}f-V \cdot f \; , \]
on function spaces $\mathcal F_p$ over a closed loop or $\mathcal F_0$ over an 
interval $[0,a]$ with Dirichlet boundary conditions:
\begin{eqnarray}\nonumber
\mathcal F_p =& \mathcal F_p(a)=&\{f: \mathbb{R} \stackrel {C^{\infty}}{\longrightarrow} \mathbb{R} \; |\; f(x)=f(x+a) \}\;, \\ \nonumber
\mathcal F_0 =& \mathcal F_0(a)=&\{f: [0,a]\stackrel {C^{\infty}}{\longrightarrow} \mathbb{R}\; |\; f(0)=f(a)=0 \}\;, \;\;a>0\;.
\end{eqnarray}
We assume the potential function $V=V(x)$  is real-valued and real-analytic on the 
closed interval $[0,a]$, and $V \in \mathcal F_p$ when the function space $\mathcal F_p$ is used.  However, we do not assume $V$ is in 
$\mathcal F_0$ when the function space $\mathcal F_0$ is used, that is, we do not assume 
$V(0)$ and $V(a)$ are zero.  

The eigenvalue problem is to find $\lambda \in \mathbb{R}$ and $f\in 
\mathcal F_p$ (or $\mathcal F_0$) that solve the second-order ordinary differential 
equation (ODE)
\begin{equation}\label{ode for the eigenvalue problem}
\mathcal L(f)=\lambda f\;.
\end{equation}

The operator $-\mathcal L$ is elliptic and it is well-known (\cite{berard}, \cite{berger}, \cite{mikhlin}, \cite{urakawa}) that the eigenvalues of $\mathcal L$ are real and form a discrete sequence 
\[\lambda_1<\lambda_2 \leq \lambda_3 \leq ...\uparrow +\infty\]
(each considered with multiplicity 1) whose first eigenvalue $\lambda_1$ is simple. The eigenvalues form a discrete spectrum, and corresponding eigenfunctions 
\[f_1,f_2,f_3,... \;\;\text{in}\;\;\mathcal F_p \;\;\text{or in}\;\;\mathcal F_0 \; 
, \;\;\; \mathcal{L}(f_j)=\lambda_j f_j \; , \]
can be chosen to form an orthonormal basis with respect to the standard $L^2$ norm 
on $\mathcal{F}_p$ or $\mathcal{F}_0$ over $[0,a]$.  

Let \[ \Sigma_p=\mathbb{R} /(x \sim x+a)\;\;\; \text{and} \;\;\; \Sigma_0=[0,a]\] denote 
the domains of the functions in $\mathcal F_p$ and $\mathcal F_0$, respectively. The 
{\em nodes} of an eigenfunction $f \in \mathcal F_p$ (or $\mathcal F_0$) are those points of $\Sigma_p$ (or $\Sigma_0$) at which $f$ vanishes. 
When $f$ is not identically zero, the fact that $\mathcal{L}$ is second-order and 
linear implies all zeros of $f$ are isolated and of lowest order, i.e. if $f(0)=0$, 
then $(\tfrac{d}{dx}f)(0) \neq 0$.  

We have the following two theorems, the second of which uses a symmetry condition on 
$V$. The first theorem is well-known and can be proven using Sturm comparison and 
Courant's nodal domain theorem (see \cite{courantHil}, \cite{Died}, \cite{Hart}, 
\cite{Hille}, for example):  

\begin{theorem}\label{initial conditions for Dirichlet case}
Consider the operator $\mathcal L$ on the function space $\mathcal F_0$ of 
$C^\infty$ functions 
over $\Sigma_0$ with Dirichlet boundary conditions. Then all eigenspaces are 
$1$-dimensional, and to find a nonzero solution $f \in \mathcal{F}_0$ 
of $\mathcal L(f)=\lambda f$ for some eigenvalue $\lambda$, without loss of 
generality we may assume:
\[f(0)=0 \;\;,\;\; (\tfrac{d}{dx} f)(0)=1 \;.\]
Furthermore, any eigenfunction associated 
to the $j$'th eigenvalue $\lambda_j$ of $\mathcal L$ has exactly $j+1$ nodes.
\end{theorem}

The following theorem can be similarly proven, but is a bit more complicated, because 
in this case the eigenvalues are not always simple. We will prove 
Theorem \ref{initial conditions for closed case} here (and in the process also 
prove Theorem \ref{initial conditions for Dirichlet case}).  The conclusions about the 
initial conditions in these two 
theorems are quite trivial; it is the conclusions about the number of nodes of the 
eigenfunctions that are of the most interest to us.  

\begin{theorem}\label{initial conditions for closed case}
Consider the operator $\mathcal L$ on the function space 
$\mathcal F_p$ of $C^\infty$ periodic functions over 
$\Sigma_p$. Suppose the real-analytic function $V \in \mathcal F_p$ has 
the symmetry
\begin{equation}\label{Vsymmetry} 
   V(x)=V(-x) \;\; \forall x\in \mathbb{R} \; . 
\end{equation}
Let $\lambda_1<\lambda_2 \leq \lambda_3 \leq ...\uparrow +\infty$ be the 
spectrum of $\mathcal L$ with a corresponding basis $f_1, f_2,f_3,... \in 
\mathcal{F}_p$ of eigenfunctions. Then the eigenspaces are each at 
most $2$-dimensional, and to find a basis 
for the eigenspace associated to $\lambda_j$, we may assume:
\begin{itemize}
\item When the eigenspace for $\lambda_j$ is $1$-dimensional, we may take 
$f_j$ so that one of 
\[ f_j(0)=1 \;\;,\;\; (\tfrac{d}{dx} f_{j})(0)=0 \;\;\; \;\;\;\; \text{or} \;\;\;\;\; \;\; f_j(0)=0 \;\;,\;\; (\tfrac{d}{dx} f_{j})(0)=1 \;\;\;\; \; \text{holds.} \]
\item  When the eigenspace for $\lambda_j$ is $2$-dimensional, and $\lambda_j=\lambda_{j+1}$, we may take 
\[f_j(0)=1 \;\;,\;\; (\tfrac{d}{dx} f_{j})(0)=0 \;\;\;\;\;\;\;  {\text and} 
\;\;\;\;\;\;\;f_{j+1}(0)=0 \;\;,\;\; (\tfrac{d}{dx} f_{j+1})(0)=1 \;.\]
\end{itemize}
Furthermore, any eigenfunction in $\mathcal{F}_p$ associated to 
$\lambda_j$ has exactly $j$ nodes if $j$ is even, and $j-1$ nodes otherwise.  
\end{theorem}

After proving these results in Section \ref{Proofs of the Theorems}, 
we will see in Section \ref{Computation of the spectra} that 
Theorem \ref{initial conditions for closed case} 
gives a method to numerically compute the spectra of the operator $\mathcal L$. Then, 
in Section \ref{Application to CMC surfaces of revolution in S3}, we apply that method 
to study the index of CMC surfaces of revolution in the 
round $3$-sphere.   

In \cite{rossman1}, a method was given for computing
the eigenvalues for the Jacobi operator of a Wente torus, involving
the Rayleigh-Ritz method and restricting to finite dimensional subspaces
of function spaces defined over tori.  Then in \cite{rossman2},
both this method and a second more direct method were given for computing
the first eigenvalue of the Jacobi operator of a Delaunay surface with
respect to periodic functions, and the second
method depended on Delaunay surfaces being surfaces of revolution.
It was argued in \cite{rossman2} that, although the second method was
clearly the simpler of the two, the first method was still of value
because it could compute any eigenvalue of the Jacobi operator, while the
second method computed only the first eigenvalue.  However, via 
Theorem \ref{initial conditions for closed case}, 
the second method in \cite{rossman2} in fact extends to a method that gives any
eigenvalue and hence is both simpler and equally as robust as the
first method.  Additionally, this extended 
second method involves only using any standard ODE
solver, such as the Euler algorithm or the Runge-Kutta algorithm,
and so has only as much numerical error as those algorithms have,
whereas the first method involves restrictions to finite dimension
subspaces for which the numerical error cannot be easily estimated and 
appears to be very much larger than for the extended 
second method.  (This can be seen by comparing the respective errors of the 
two methods in cases where the spectra are explicitly known.)

Certainly the first method was necessary in \cite{rossman1}, because
Wente surfaces are not surfaces of revolution.  
But for the above reasons, the method we give here is in every
way superior to the methods found in \cite{rossman1} and \cite{rossman2}, in the 
case of CMC surfaces of revolution.

\begin{figure}[phbt]
  \centering
  \includegraphics[scale=0.22]{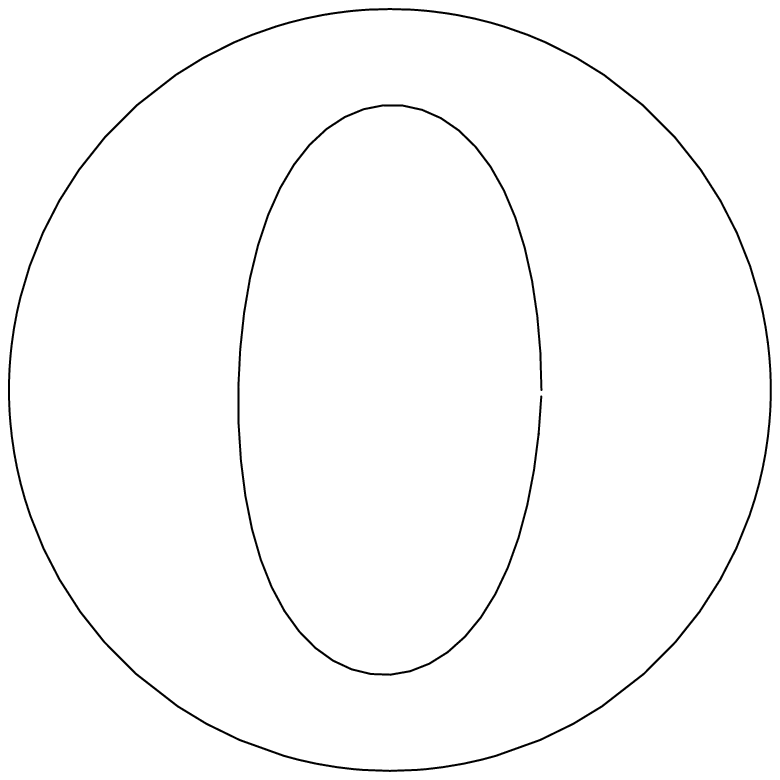}
  \includegraphics[scale=0.22]{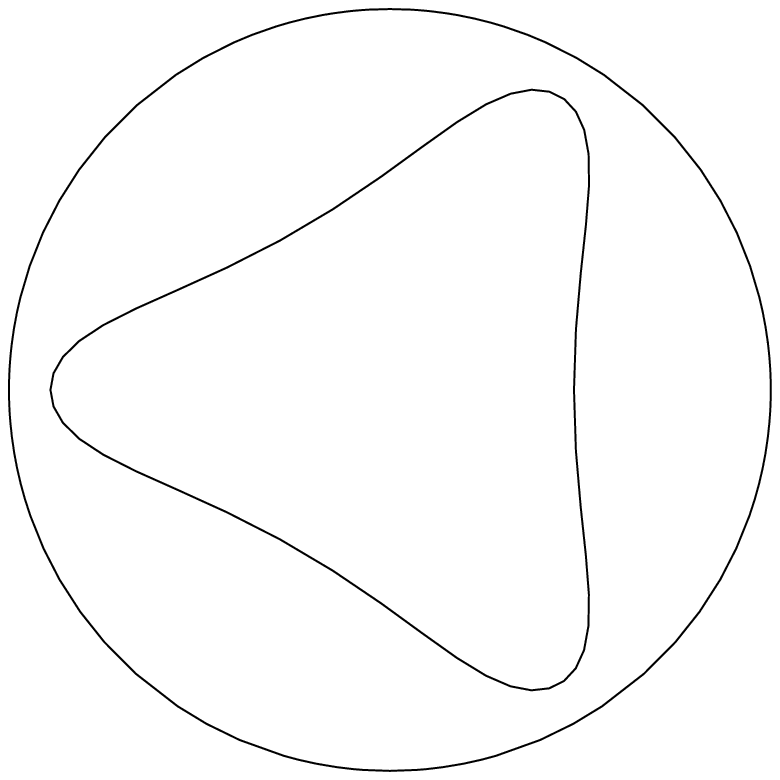}
  \includegraphics[scale=0.22]{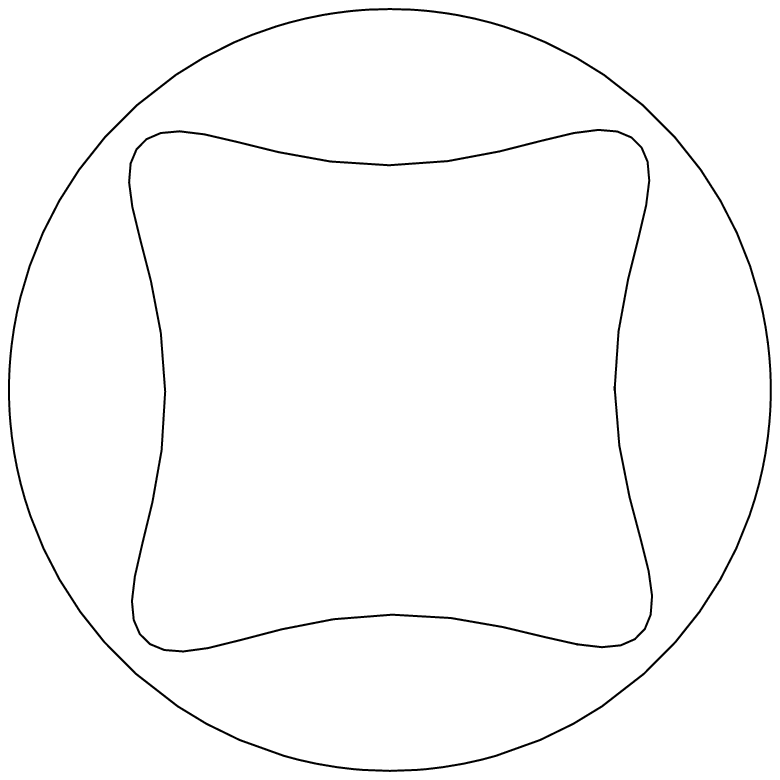}
  \includegraphics[scale=0.22]{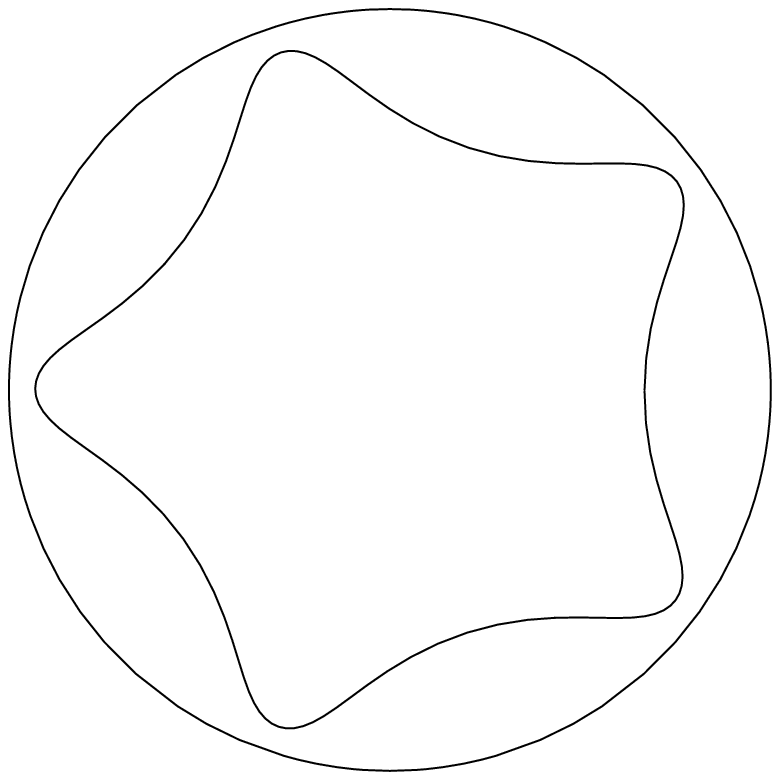}
  \includegraphics[scale=0.22]{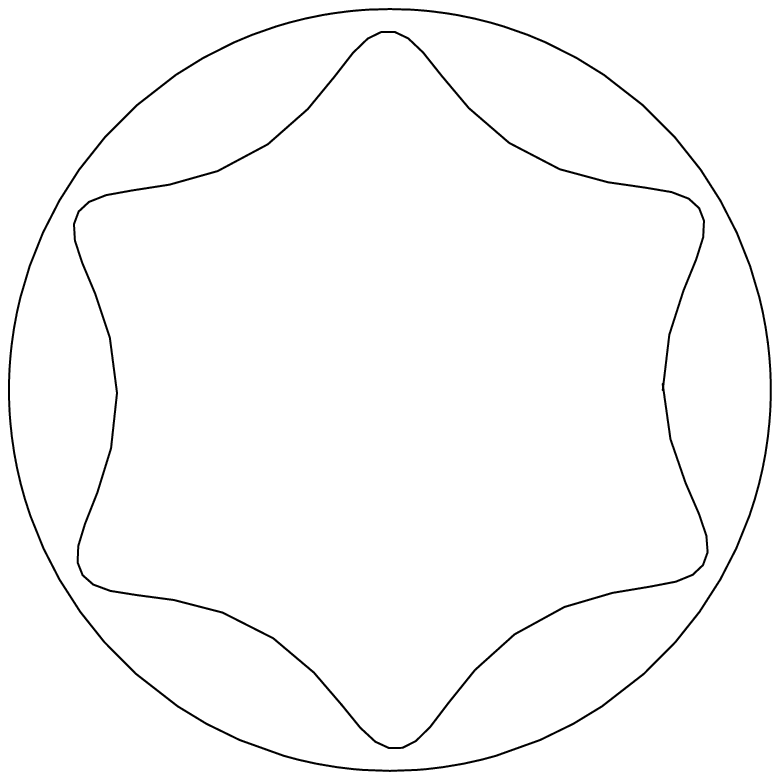}
  \includegraphics[scale=0.35]{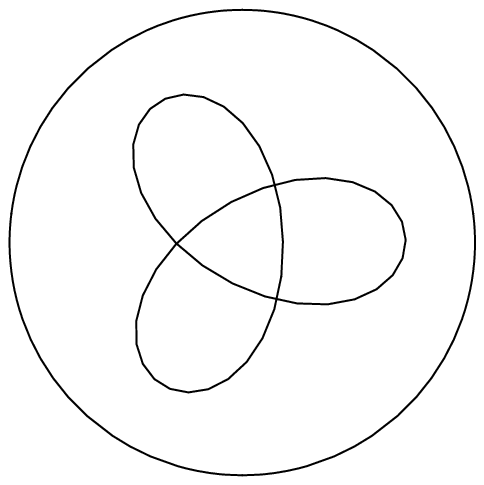}
  \includegraphics[scale=0.22]{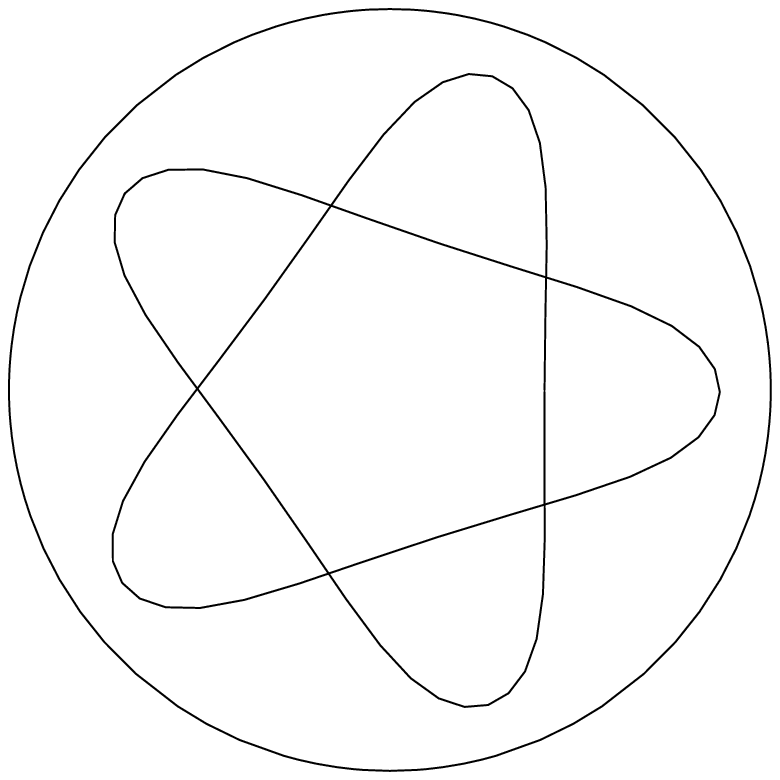}
  \includegraphics[scale=0.22]{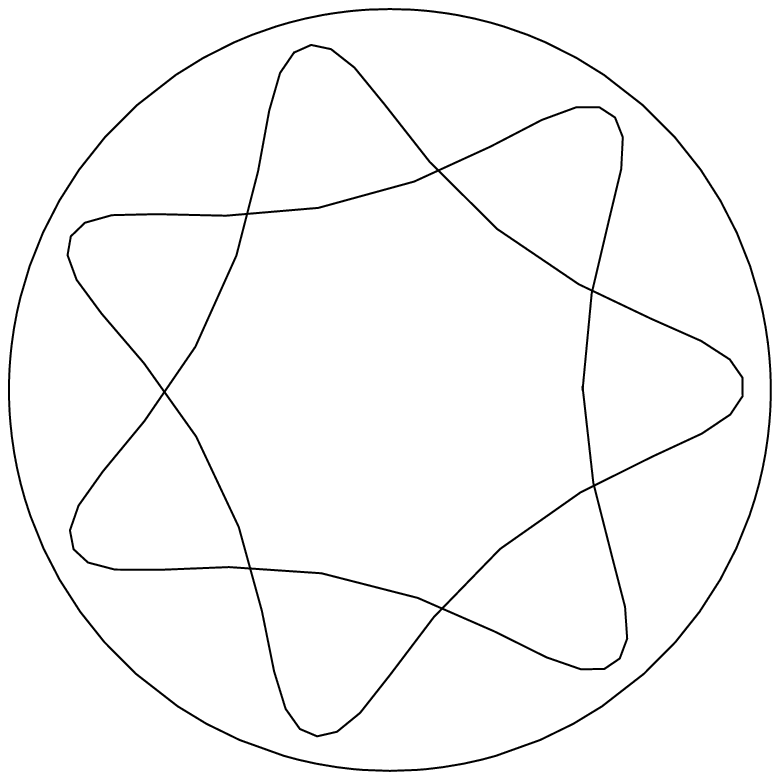}
    \includegraphics[scale=0.22]{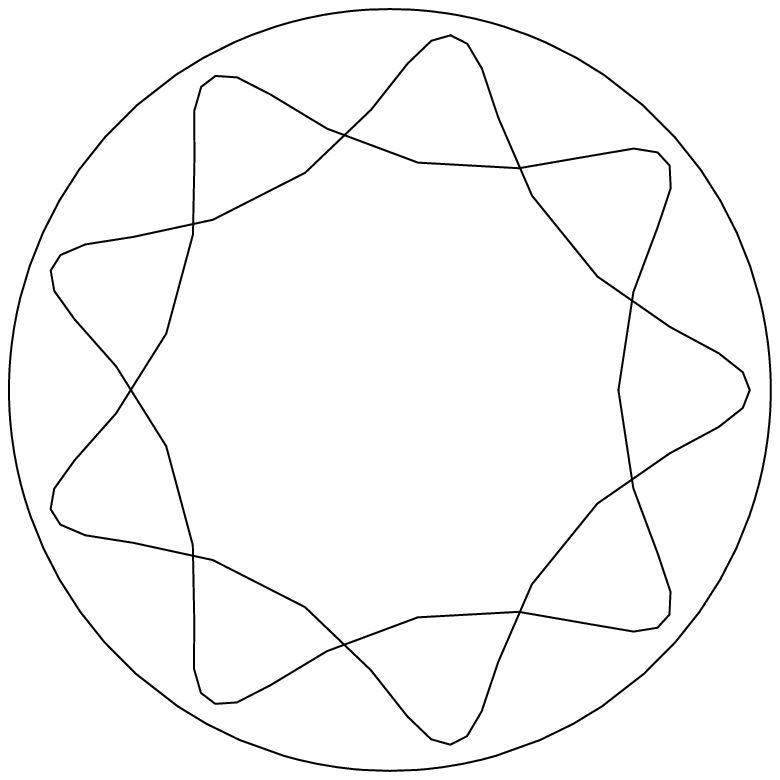}
      \includegraphics[scale=0.35]{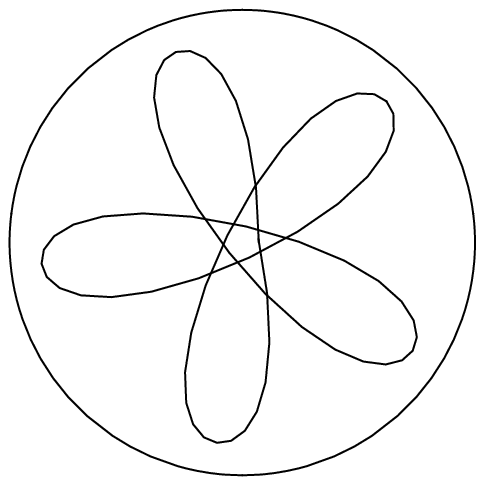}
   \includegraphics[scale=0.22]{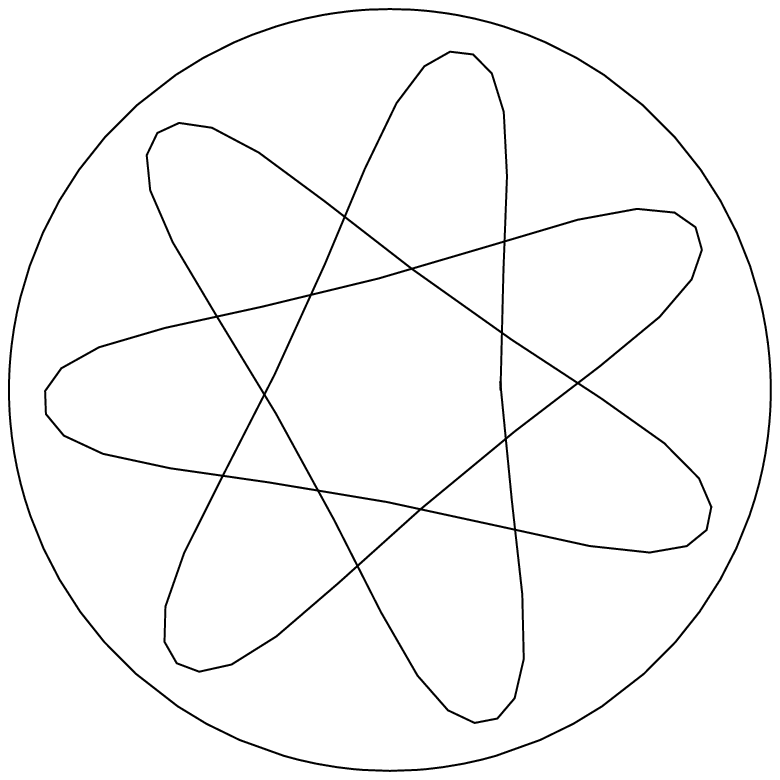}
   \includegraphics[scale=0.22]{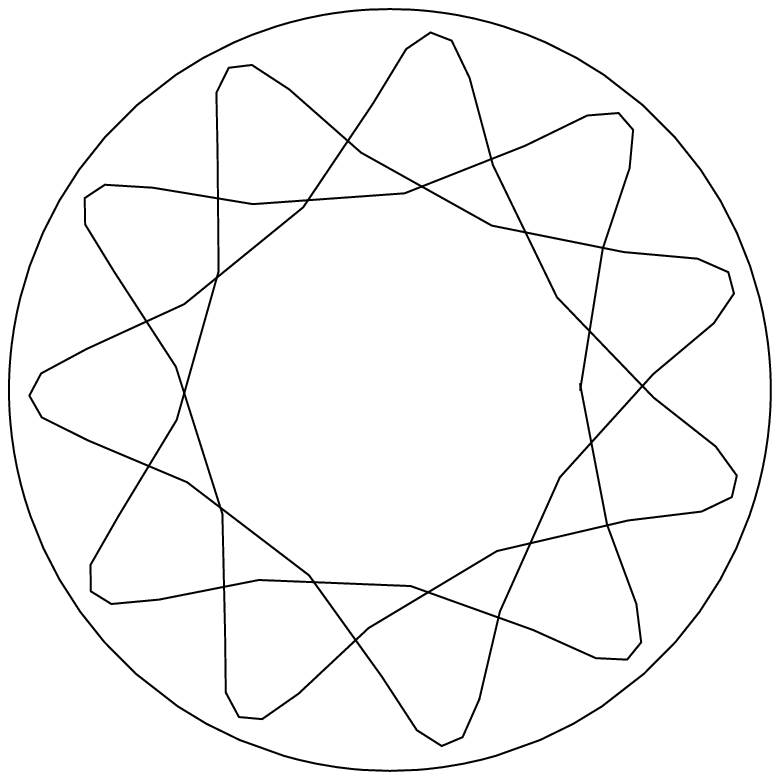}
     \includegraphics[scale=0.35]{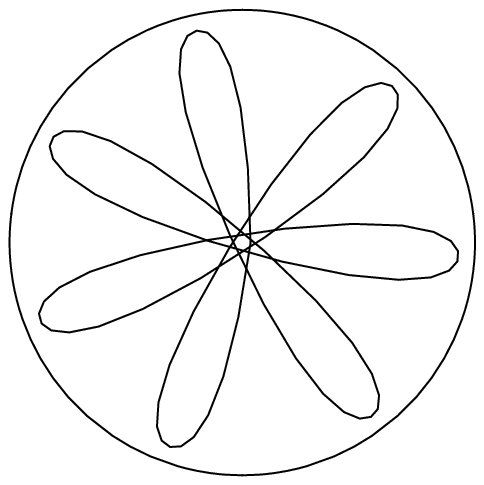}
  \includegraphics[scale=0.22]{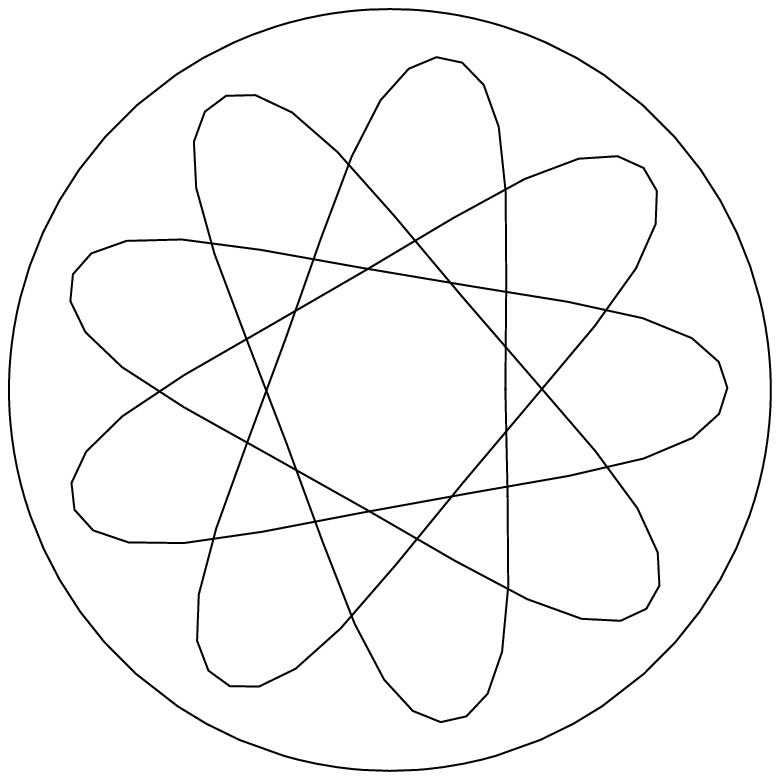}
  \includegraphics[scale=0.22]{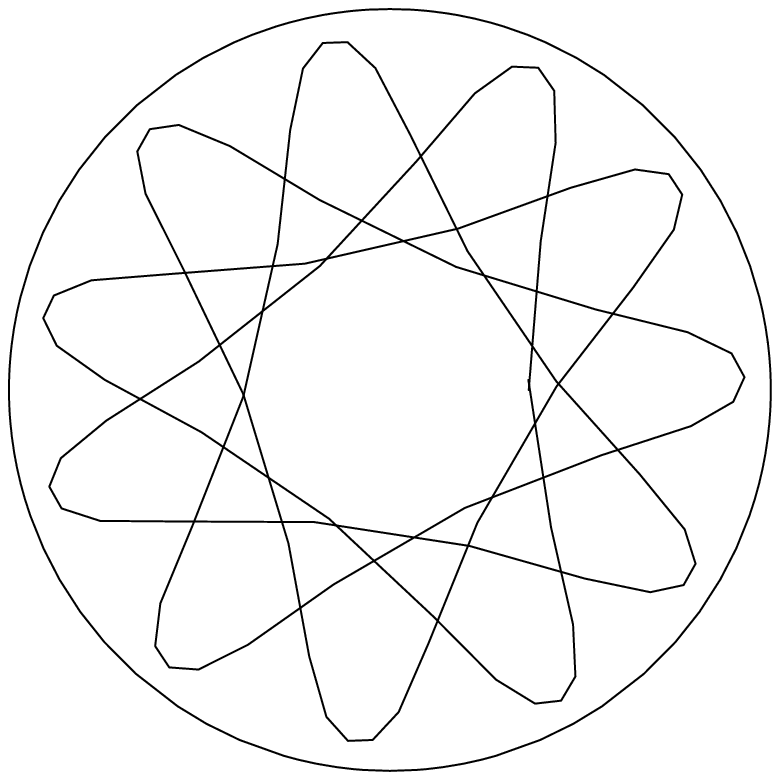}
    \includegraphics[scale=0.35]{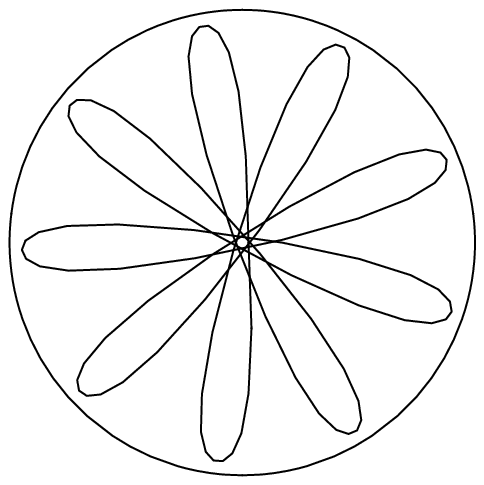}
   \includegraphics[scale=0.22]{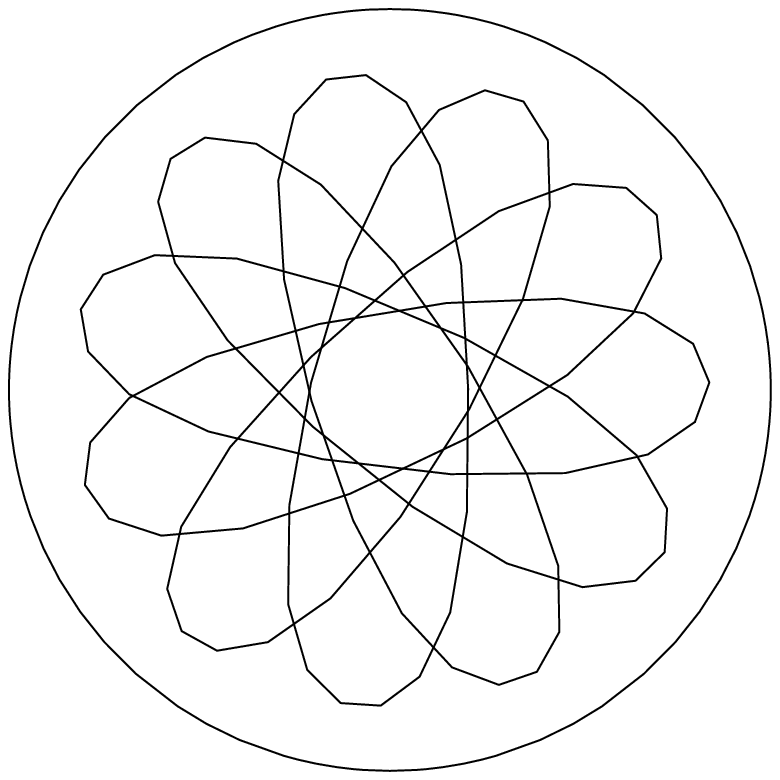}
   \caption{Profile curves for surfaces that we label $U_i$ for $i=1,2,...,17$ 
         ($U_1$,...,$U_9$ from left to right in the upper row, 
         $U_{10}$,...,$U_{17}$ from left to right in the lower row). 
           These are profile curves of 
           CMC tori of revolution, shown in totally geodesic 
           hemispheres having the rotation axis as boundary.  
           The images are stereographic projections from $\mathbb{S}^3$ to 
           $\mathbb{R}^3 \cup \{ \infty \}$.  The outer circle is the rotation axis, 
           with profile curve inside.  All of these surfaces are unduloidal, in the 
           sense that the projections of these curves to the nearest points in the 
           rotation axis are everywhere continuous local injections.}
  \label{figure1}
\end{figure}

A surface of revolution in the unit 3-sphere $\mathbb{S}^3$ is generated by revolving a given planar curve about a geodesic line in the geodesic plane containing this given 
curve. The given curve is called the profile curve and the geodesic line is called the axis of revolution. The profile curves of non-spherical non-flat 
CMC surfaces of revolution in $\mathbb{S}^3$ will periodically have minimal and 
maximal distances to the axis of revolution 
\cite{hsiang}. We call the points of minimal distance 
the {\em necks}, and the points of maximal distance 
the {\em bulges}. In general, when these surfaces close to become compact surfaces 
without boundary, they are of the following 3 types: 
\begin{itemize}
\item {\em round spheres}, every point of which is the same distance 
from a fixed point (the center), 
\item {\em flat} CMC tori, every point of which is the same distance 
from a closed geodesic (the axis of revolution), 
\item {\em non-flat} CMC tori, where 
the distances from the axis of revolution to the necks and bulges are not equal. 
\end{itemize} 

Because these surfaces are closed, the number of negative eigenvalues of their 
Jacobi operators, counted with multiplicity and called the {\em Morse index}, is 
finite.  The Morse index is of interest because it is a measure of the degree of 
instability of the surface.  
In the first two cases above, the Morse index is easily explicitly computed 
\cite{WN}, being $1$ for the first case (this is closely related to the 
fact that spheres are stable \cite{bdcjost}) and always at least $5$ for the second 
case.  Regarding the third case, the authors proved the following in \cite{WN}: 

\begin{theorem}\label{previousresult} 
Let $\mathcal{S}$ be a non-flat closed CMC torus of revolution in 
$\mathbb{S}^3$, with $k$ bulges and $k$ necks.  Let $w$ denote the wrapping 
number of the projection of a profile curve of $\mathcal{S}$ to the axis 
circle of revolution.  Then: 
\begin{itemize}
\item $\mathcal{S}$ has index at least $\max(5,2k+1)$.  
\item If $\mathcal{S}$ is nodoidal with $k \geq 2$, then 
$\mathcal{S}$ has index at least $\max (11,2k+5)$.  
\item If $\mathcal{S}$ is unduloidal with $w \geq 2$, 
then $\mathcal{S}$ has index at least $\max (6w-1,2k+4w-3)$.  
\end{itemize}
\end{theorem}

\begin{figure}[phbt]
  \centering
  \includegraphics[scale=0.22]{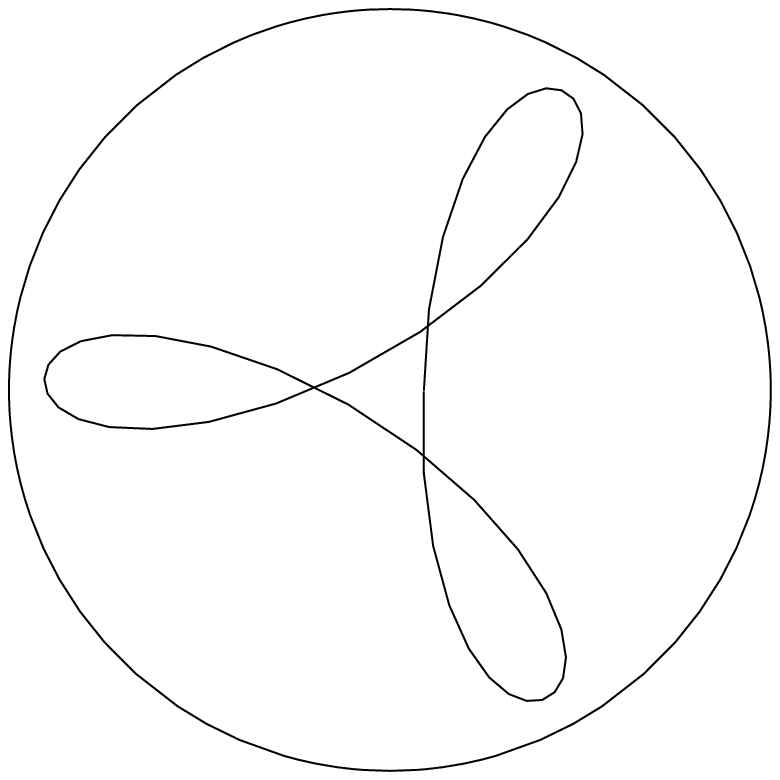}
  \includegraphics[scale=0.31]{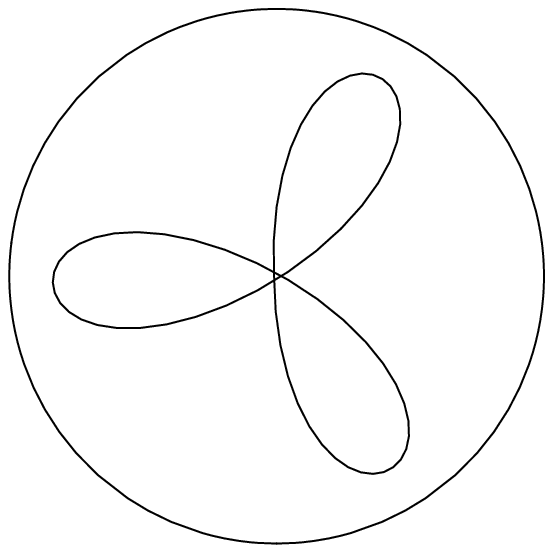}
  \includegraphics[scale=0.22]{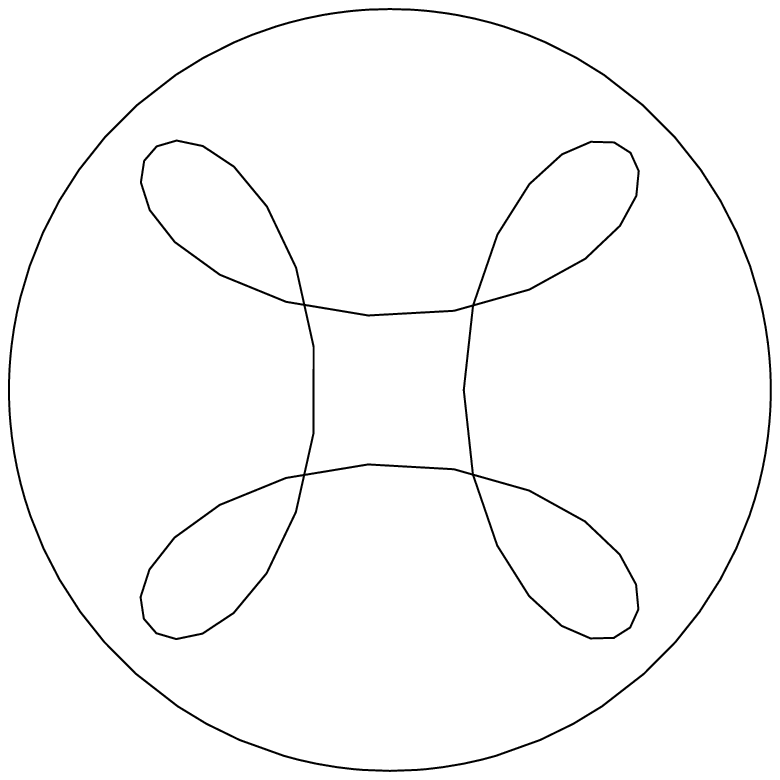}
  \includegraphics[scale=0.22]{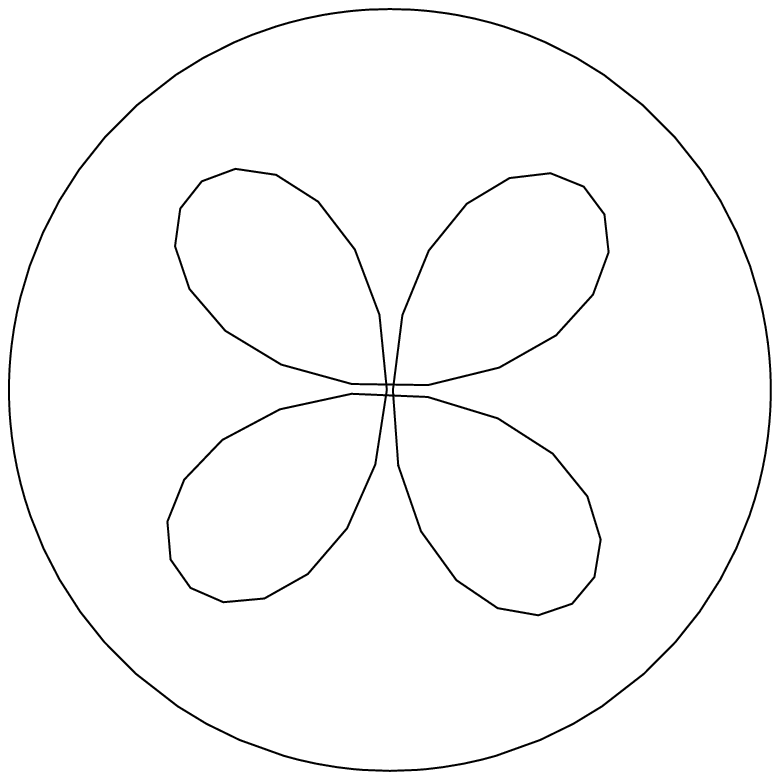}
  \includegraphics[scale=0.22]{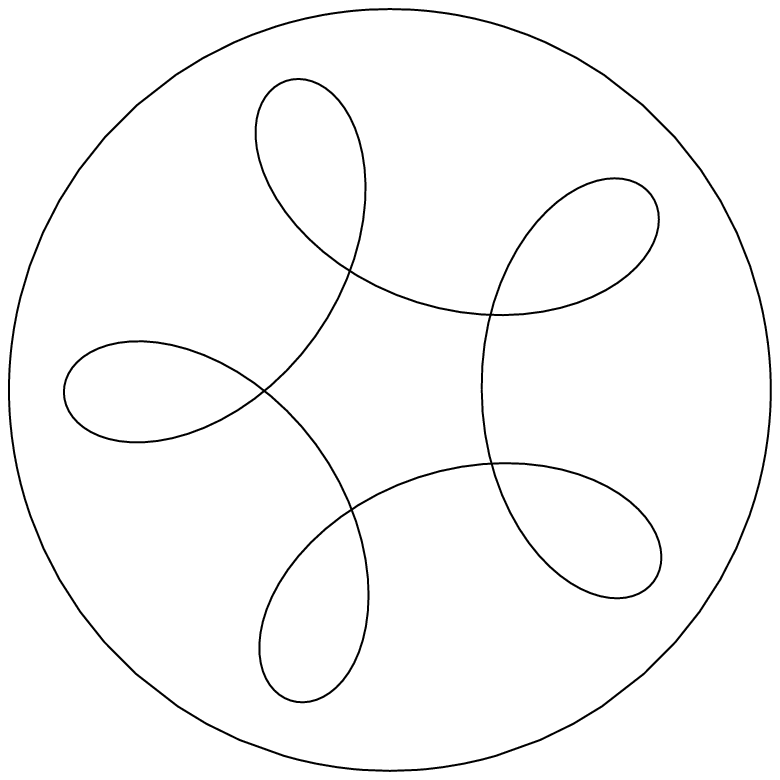}
  \includegraphics[scale=0.22]{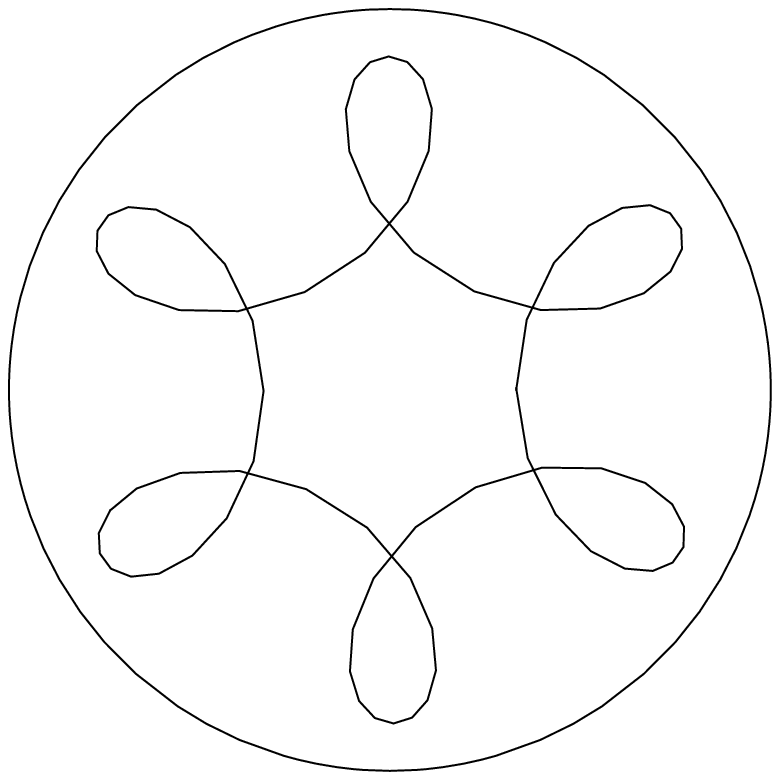}
  \includegraphics[scale=0.22]{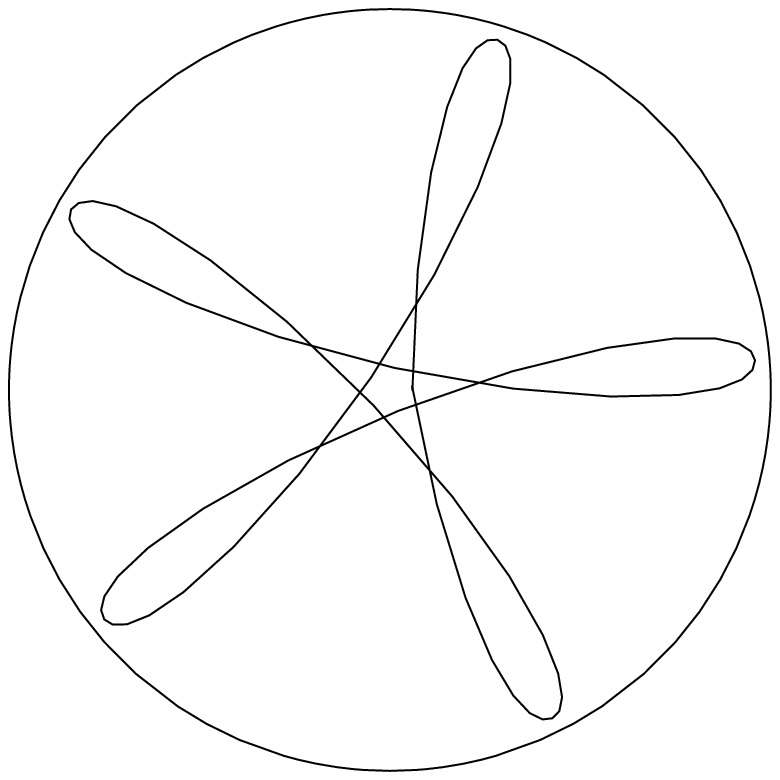}
  \includegraphics[scale=0.22]{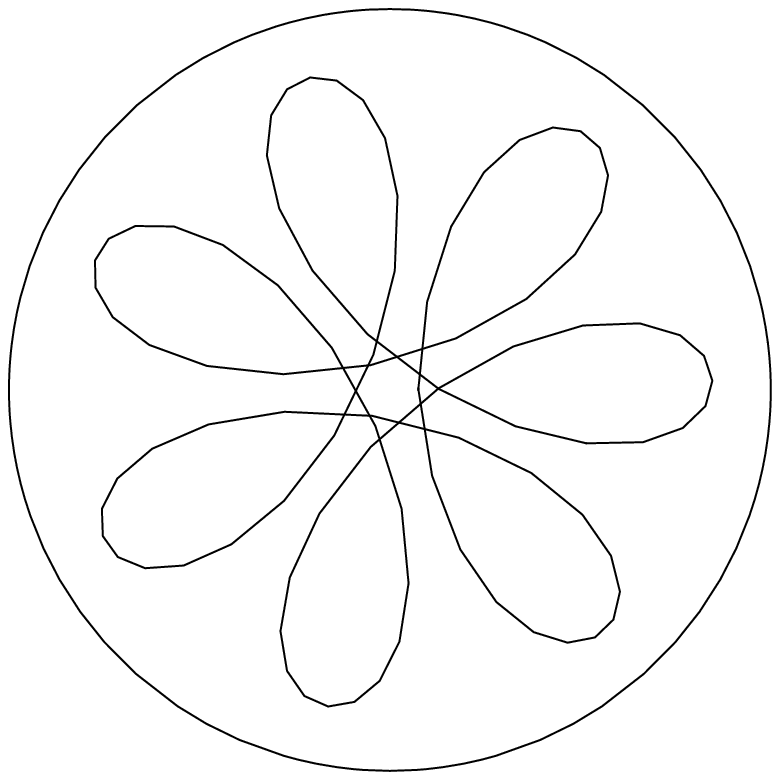}
  \includegraphics[scale=0.22]{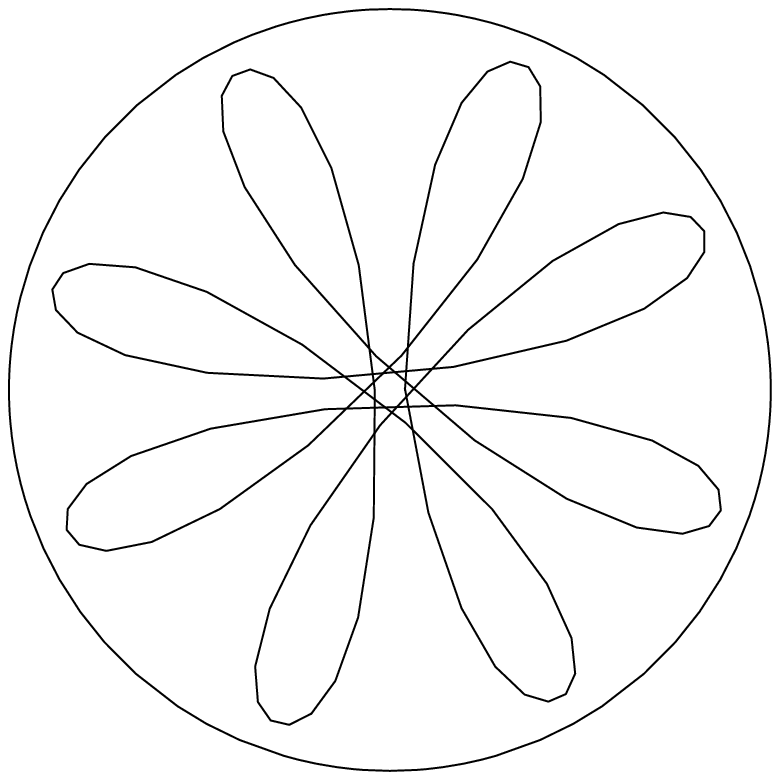}
  \includegraphics[scale=0.22]{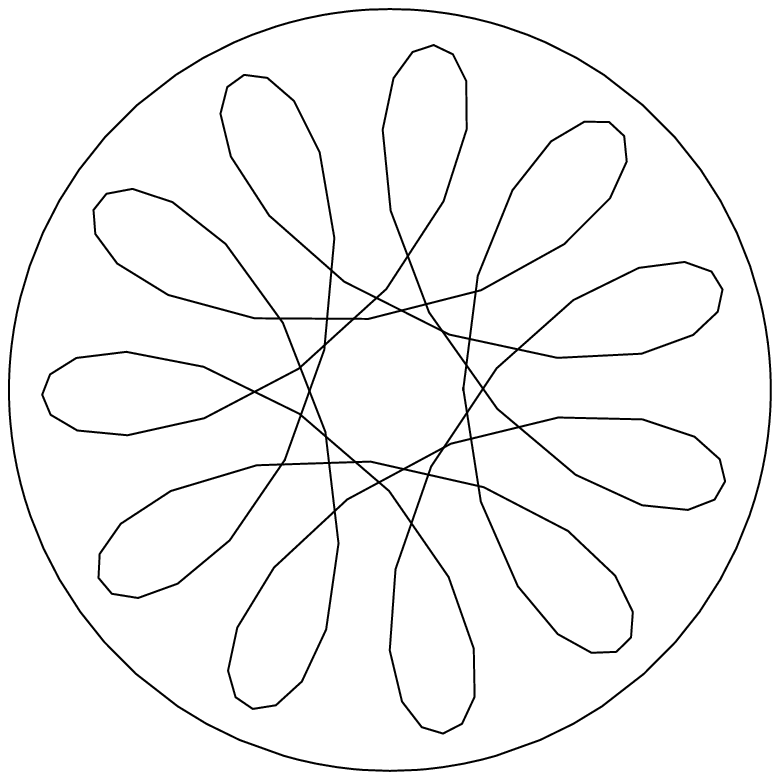}
  \includegraphics[scale=0.22]{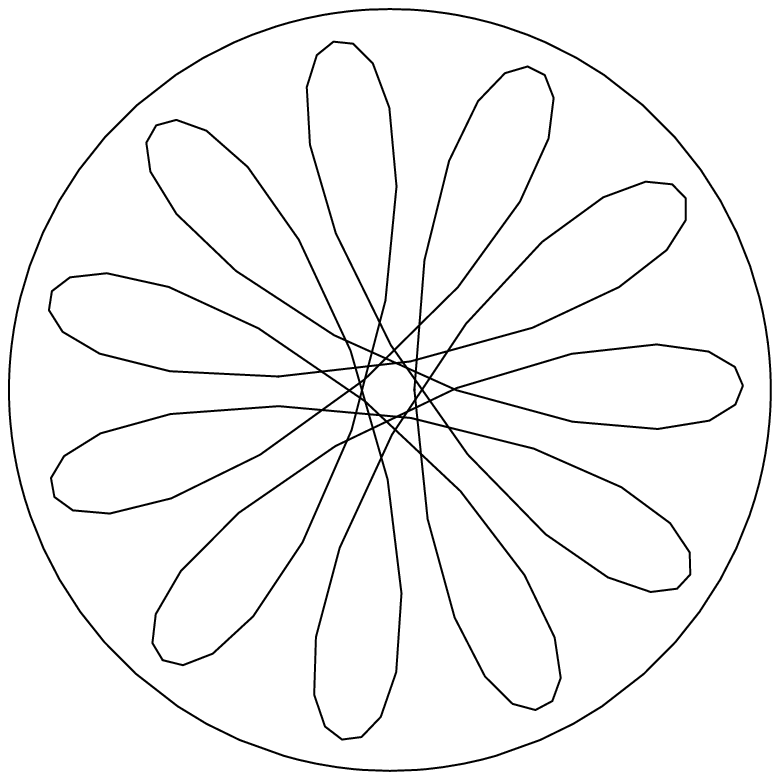}
  \caption{Profile curves for surfaces that we label $N_i$ for $i=1,2,...,11$ 
  ($N_1$,...,$N_9$ from left to right in the upper row, 
         $N_{10}$, $N_{11}$ from left to right in the lower row). 
         All of these surfaces are nodoidal, i.e. they are not unduloidal.}
  \label{figure2}
\end{figure}

The numerical results here show the lower bounds in the above theorem are 
very close to the true value for the Morse index in the case of unduloids.  
For example, using Table \ref{table2} and Lemma \ref{lm index}, 
the numerically computed index of the unduloid $U_1$ 
(resp. $U_2$, $U_3$, ..., $U_{17}$)
is $6$ (resp. $8$, $10$, $12$, $14$, $12$, $16$, $20$, $24$, $20$, 
$24$, $32$, $28$, $32$, $36$, $36$, $44$), 
while the above theorem gives the lower bound $5$ 
(resp. $7$, $9$, $11$, $13$, $11$, $15$, $19$, $23$, $19$, 
$23$, $31$, $27$, $31$, $35$, $35$, $43$)
for the index.  In all cases, the lower bound in Theorem \ref{previousresult} 
differs from the numerically computed value for the index by only $1$, thus 
the lower bound is quite sharp.  

\begin{figure}[phbt]
  \centering
  \includegraphics[scale=0.25]{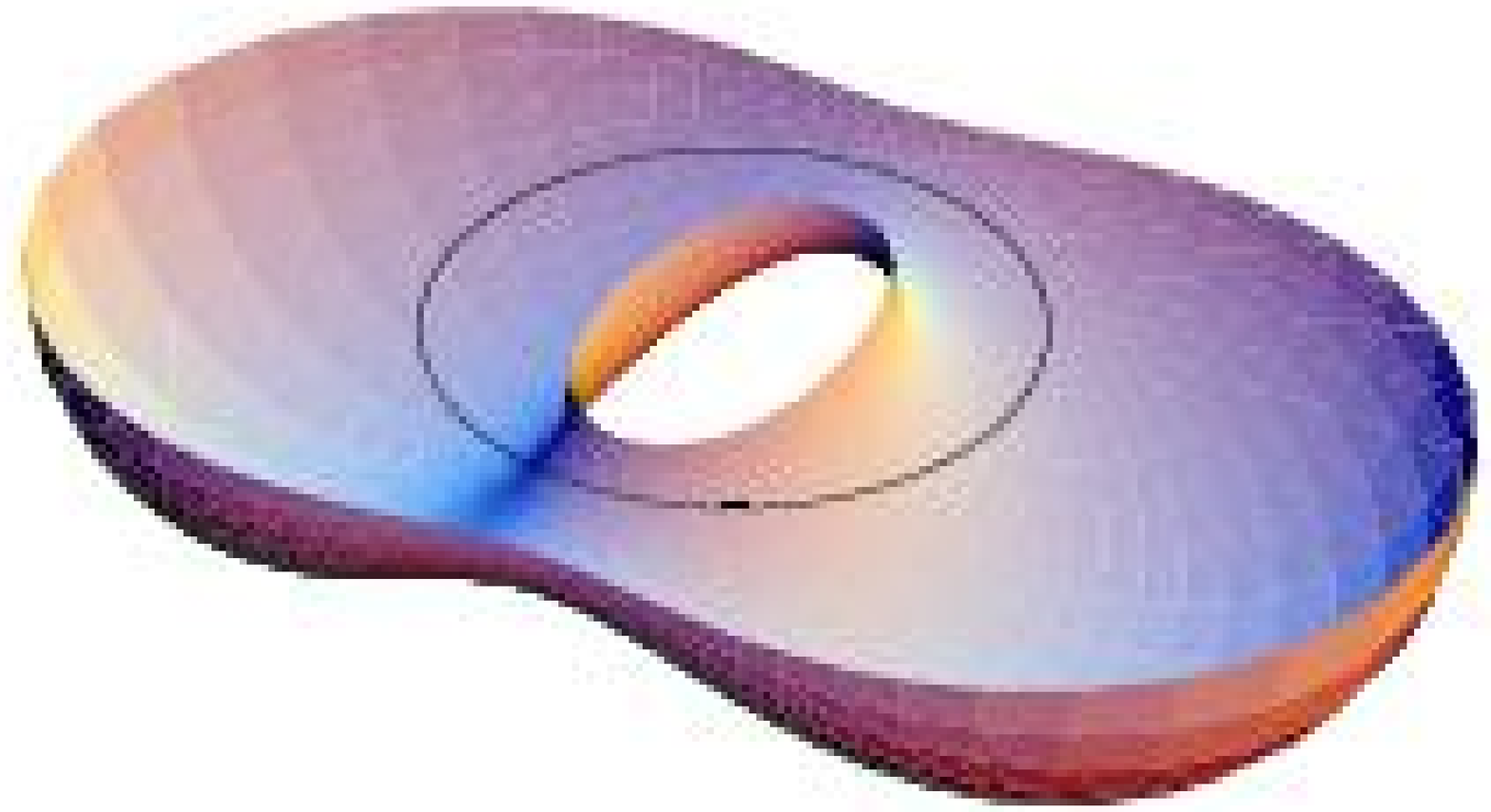}
  \includegraphics[scale=0.2]{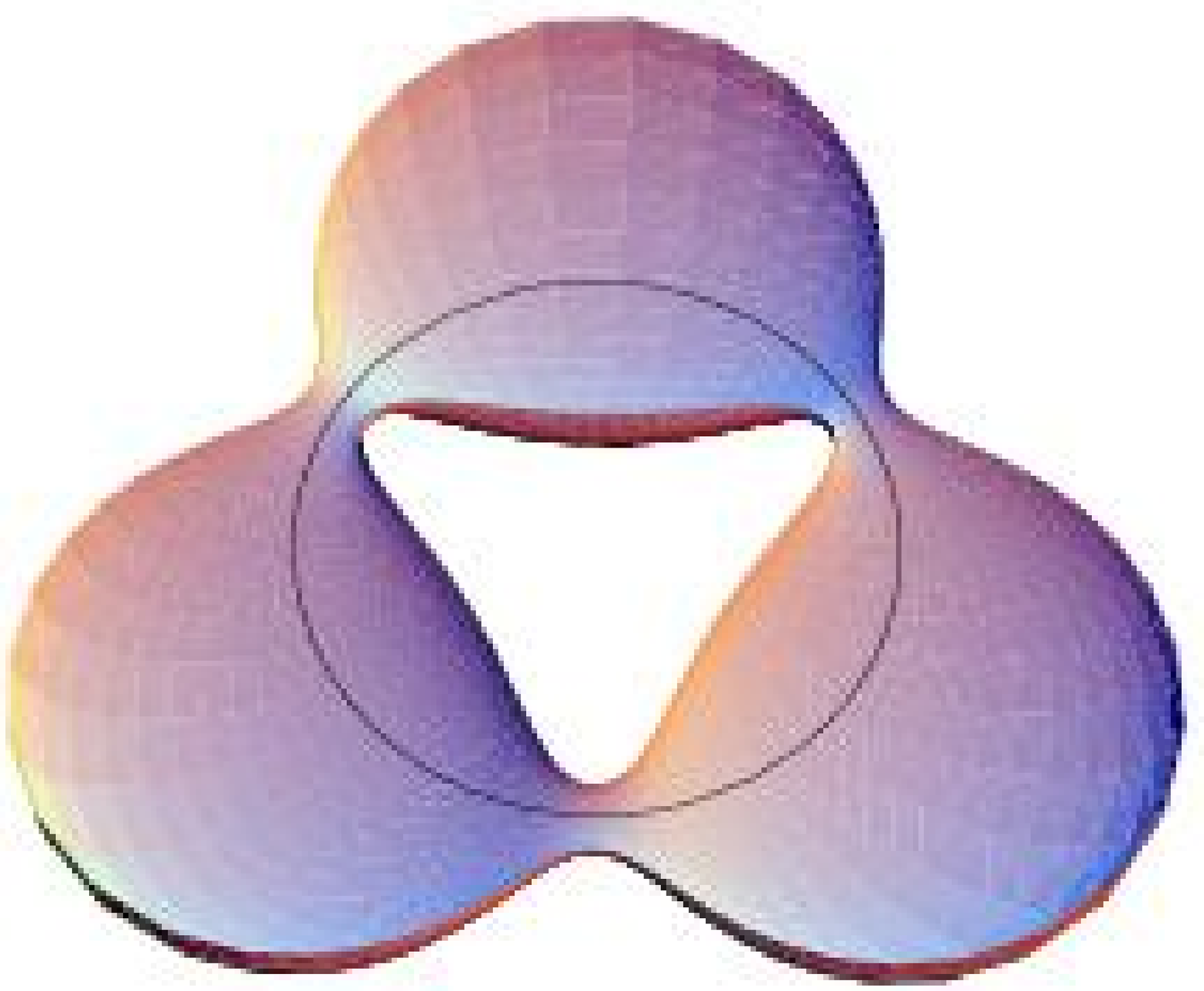}
  \includegraphics[scale=0.2]{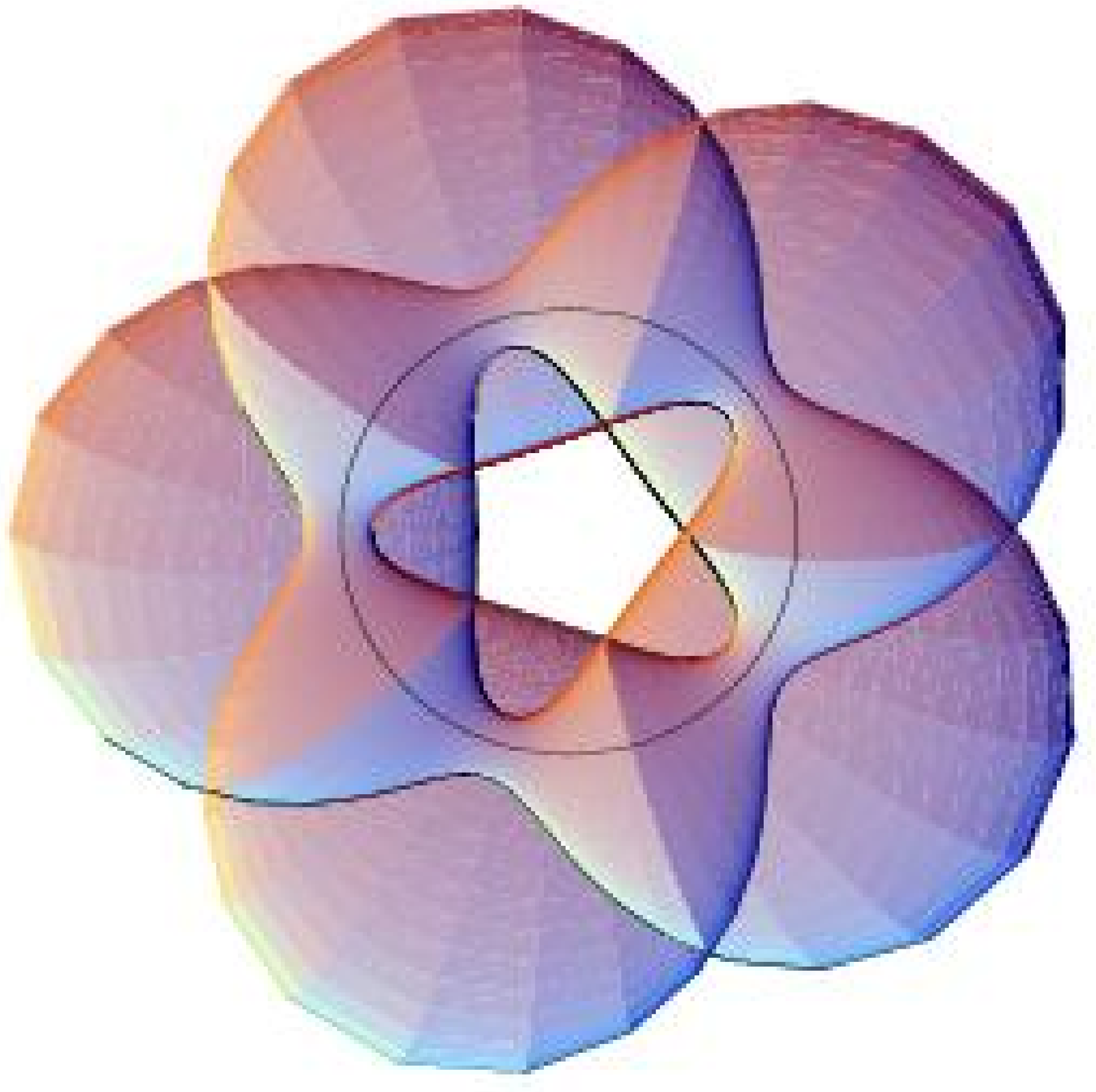}
  \caption{Half of each of the surfaces $U_1$, $U_2$, $U_{7}$.}
  \label{fullpics1}
\end{figure}

The lower bounds in Theorem \ref{previousresult} are 
not as sharp in the case of nodoids, but 
still are greater than half of the numerically computed value for all of the 
surfaces shown in Figure \ref{figure2}.  The numerically computed index of the 
nodoid $N_1$ (resp. $N_2$, $N_3$, ..., $N_{11}$)
is $12$ (resp. $12$, $18$, $18$, $24$, $30$, $20$, $32$, $34$, $52$, $48$), 
while the above theorem gives the lower bound $11$ 
(resp. $11$, $13$, $13$, $15$, $17$, $15$, $19$, $21$, $27$, 
$27$) for the index.  

\section{\bf Proofs of the Theorems \ref{initial conditions for Dirichlet case} and 
\ref{initial conditions for closed case}.}
\label{Proofs of the Theorems}
To prove Theorems \ref{initial conditions for Dirichlet case} and 
\ref{initial conditions for closed case}, we give a series of lemmas.  We first note 
that: 
\begin{itemize}
\item For each $\lambda$, it is easily shown that the space of solutions of 
\eqref{ode for the eigenvalue problem} 
amongst functions $f:\mathbb{R} \to \mathbb{R}$ is $2$-dimensional.  
\item Consider the eigenvalue problem \eqref{ode for the eigenvalue problem} over the 
function space $\mathcal F_0$ on the interval $\Sigma_0$ with 
Dirichlet boundary conditions. Suppose $\hat f_1$ and $\hat f_2$ are 
two linearly independent eigenfunctions corresponding to some eigenvalue $\lambda$. 
Noting that $(\tfrac{d}{dx} \hat f_1)(0)$ and $(\tfrac{d}{dx} \hat f_2)(0)$ are 
both nonzero, 
take the linear combination $\hat f_3= (\tfrac{d}{dx} \hat f_2)(0) \cdot \hat f_1 - 
(\tfrac{d}{dx} \hat f_1)(a) \cdot \hat f_2$.  Then $\hat f_3(0)=
(\tfrac{d}{dx} \hat f_3)(0)=0$, and it follows that $f_3$ is identically zero, 
contradicting the linear independence of $f_1$ and $f_2$.  
Hence the eigenvalues are simple.  Hence the eigenvalues are always simple for the 
Dirichlet eigenvalue problem.  Furthermore, multiplying by a scalar factor if necessary, 
we may assume the initial conditions for an eigenfunction $f$ is $f(0)=0$ 
and $(\tfrac{d}{dx} f)(0)=1$. 
\end{itemize}

For the closed eigenvalue problem \eqref{ode for the eigenvalue problem} 
with $f \in \mathcal F_p$, the 
eigenspace associated to any eigenvalue $\lambda$ is either $1$ or $2$ dimensional, 
and we have the following lemma regarding the initial conditions to find a basis for 
the eigenspace: 

\begin{figure}[phbt]
  \centering
  \includegraphics[scale=0.2]{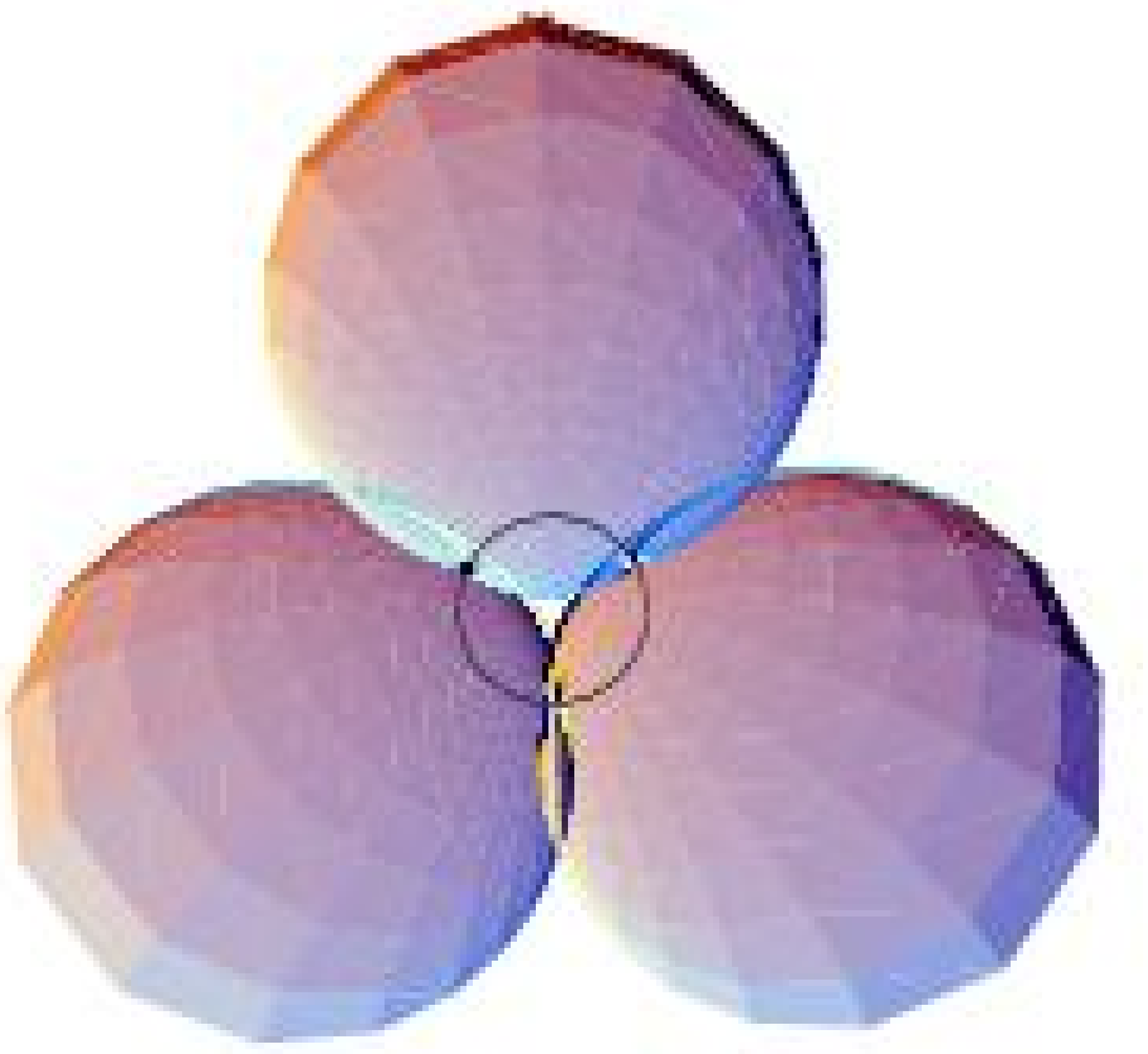}
  \includegraphics[scale=0.2]{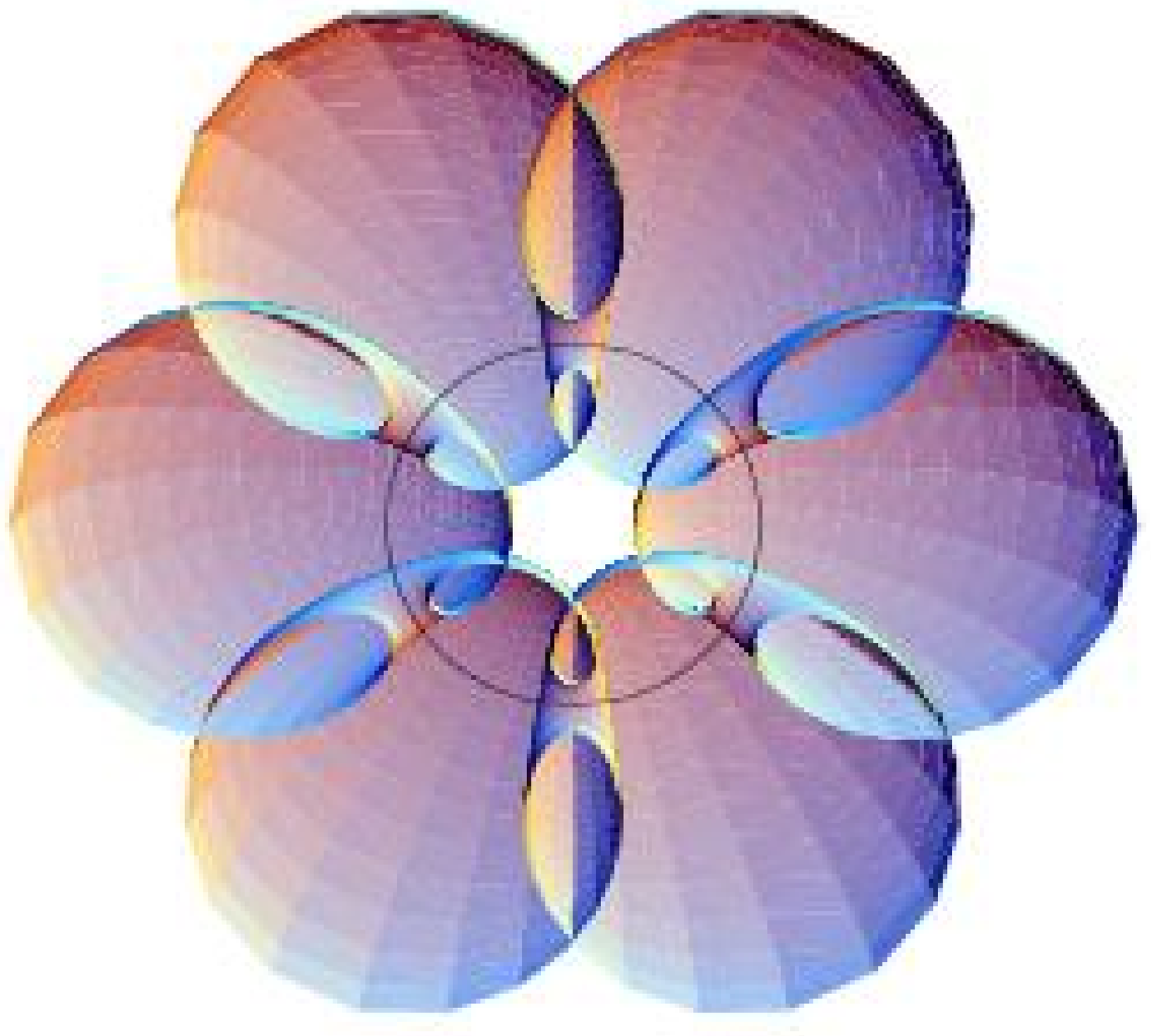}
  \caption{Half of each of the surfaces $N_1$, $N_6$.}
  \label{fullpics2}
\end{figure}

\begin{lemma}\label{initial condition lemma for closed case}
Suppose $V \in \mathcal F_p$ has the symmetry \eqref{Vsymmetry}. 
Let $\lambda_1<\lambda_2 \leq \lambda_3 \leq ...\uparrow +\infty$ be the spectrum of 
$\mathcal L$ over the space $\mathcal F_p$ with a corresponding basis 
$f_1, f_2,,f_3,...$ of eigenfunctions. Then the eigenspaces are each at most 
$2$-dimensional, and to find a basis for $\{ f \in \mathcal F_p \, | \, 
\mathcal L(f)=\lambda_j f\}$ for some eigenvalue $\lambda_j$, we may assume: 
\begin{itemize}
\item When the eigenspace for $\lambda_j$ is $1$-dimensional, we may take a single 
eigenfunction $f_j \in \mathcal F_p$ such that 
\[ \text{either} \;\; f_j(0)=1 \;\;,\;\; (\tfrac{d}{dx} f_{j})(0)=0 \;\; 
\text{or} \;\; f_j(0)=0 \;\;,\;\; (\tfrac{d}{dx} f_{j})(0)=1 \; . \]
\item  When the eigenspace for $\lambda_j=\lambda_{j+1}$ is $2$-dimensional, 
we may take two eigenfunctions $f_j, f_{j+1} \in \mathcal F_p$ such that 
\[f_j(0)=1 \;\;,\;\; (\tfrac{d}{dx} f_{j})(0)=0 \;\; 
{\text and}\;\;f_{j+1}(0)=0 \;\;,\;\; (\tfrac{d}{dx} f_{j+1})(0)=1 \;.\]
\end{itemize}
\end{lemma}
\begin{proof}
First we consider the case of a $1$-dimensional eigenspace. Let $f_j \in 
\mathcal F_p$ be a basis element of this eigenspace. If $f_j$ has neither the 
symmetry $f_j(x)=f_j(-x)$ nor $f_{j}(-x)=-f_j(x)$, then $f_j(x)$ and $f_j(-x)$ would 
be two linearly independent eigenfunctions with eigenvalue $\lambda_j$, a 
contradiction. Hence $f_j(x)=f_j(-x)$ or $f_j(-x)=-f_j(x)$ for all $x\in \mathbb{R}$, 
and so $f_j(0)=0$ or $(\tfrac{d}{dx} f_j)(0)=0$. Furthermore, because 
multiplying $f_j$ by a real constant still gives a solution to 
\eqref{ode for the eigenvalue problem} with $\lambda=\lambda_j$, 
we may assume either $f_j(0)=1$ or $(\tfrac{d}{dx} f_j)(0)=1$.  
Hence the first part of the lemma is shown.  

For the case of a $2$-dimensional eigenspace, any $C^\infty$ solution $f: \mathbb{R} 
\to \mathbb{R}$ to \eqref{ode for the eigenvalue problem} with $\lambda=\lambda_j$ 
lies in $\mathcal F_p$, hence we can choose a basis $f_j, f_{j+1} \in \mathcal F_p$ 
with the initial conditions as in the second part of the lemma.  
\end{proof}

The following lemma is known as Courant's nodal domain theorem, and the 
proof, which applies in our setting with either function space $\mathcal{F}_p$ 
or $\mathcal{F}_0$, can be found in \cite{c} (see also \cite{rossman3}).

\begin{lemma}\label{Courants nodal domain theorem}
{\em (Courant's nodal domain theorem.)} The number of nodes of any 
eigenfunction for \eqref{ode for the eigenvalue problem} in $\mathcal F_p$ (resp. 
$\mathcal F_0$) associated to the $j$'th eigenvalue $\lambda_j$ is at most $j$ 
(resp. $j+1$).  
\end{lemma}

Lemma \ref{Courants nodal domain theorem} may be strengthened using 
Sturm comparison, as we will see in the course of proving Theorems 
\ref{initial conditions for Dirichlet case} and 
\ref{initial conditions for closed case}. 

The following lemma is a slight generalization of a result in \cite{Hille}: 

\begin{lemma}\label{sturm comparison}
Consider the following two equations
\begin{equation}\label{first equation}
\tfrac{d^2}{dx^2}f+(V+\lambda) f=0 \; , \;\;\; \tfrac{d^2}{dx^2}\hat f+
(V+\hat \lambda) \hat f=0 
\end{equation}
with $V$ as in \eqref{ode for the eigenvalue problem} and $\lambda < 
\hat \lambda$. Suppose that the first equation in \eqref{first equation} has a 
solution $f(x) \not\equiv 0$ having two consecutive zeros at $x=\xi_1$ and $x=\xi_2$, 
with $\xi_1 <\xi_2$. Let $\hat f(x)$ be a solution of the second equation in 
\eqref{first equation}, then 
$\hat f(x)$ has at least one zero $x=\xi_3$ with $\xi_1 < \xi_3 < \xi_2$.  
\end{lemma}
\begin{proof}
Multiplying the first equation in \eqref{first equation} by 
$\hat f$ and the second equation 
in \eqref{first equation} by $f$, then subtracting the first expression from the 
second and integrating, we have 
\begin{equation}\label{comparison identity}
\int_{\xi_1}^{\xi_2} (\hat \lambda-\lambda) f\hat f \; dx = 
\left[(\tfrac{d}{dx} f) \hat f-(\tfrac{d}{dx} \hat f) f \right]_{\xi_1}^{\xi_2} 
= (\tfrac{d}{dx} f)(\xi_2) \hat f(\xi_2)-(\tfrac{d}{dx} f)(\xi_1) \hat f(\xi_1) \;, 
\end{equation}
as here $f(\xi_1)=f(\xi_2)=0$. Multiplying by the scalar $-1$ if necessary, 
we may assume $f(x)>0$ for $x \in (\xi_1,\xi_2)$, so 
$(\tfrac{d}{dx} f)(\xi_1) >0$ and $(\tfrac{d}{dx} f)(\xi_2)< 0$. 
If $\hat f(x)$ is positive everywhere in $(\xi_1,\xi_2)$, then 
$\int_{\xi_1}^{\xi_2} (\hat \lambda - \lambda) f\hat f \, dx>0$ and 
$(\tfrac{d}{dx} f)(\xi_2) \hat f(\xi_2)-(\tfrac{d}{dx} f)(\xi_1) \hat f(\xi_1) \leq 0$, 
contradicting \eqref{comparison identity}. Similarly, $\hat f(x)$ cannot be 
negative everywhere in $(\xi_1,\xi_2)$.  
\end{proof}

\begin{lemma}\label{sturm eigenvalue problem}
Consider the eigenvalue problem \eqref{ode for the eigenvalue problem} on 
$\mathcal F_0$ over the interval $\Sigma_0$ with Dirichlet boundary conditions, 
and with corresponding spectrum $\lambda_1<\lambda_2<...$ of simple eigenvalues.
Then any eigenfunction associated with $\lambda_j$ has exactly $j+1$ nodes. 
\end{lemma}
\begin{proof}
Denote a nonzero eigenfunction corresponding to eigenvalue $\lambda_j$ by $f_j$. 
Lemma \ref{Courants nodal domain theorem} 
implies $f_1$ has exactly two nodes (at $x=0$ and $x=a$).  
Assume $f_j$ has exactly $j+1$ nodal domains and let us prove $f_{j+1}$ 
has exactly $j+2$ nodes. From \eqref{ode for the eigenvalue problem} we have
$\tfrac{d^2}{dx^2} f_j+V f_j+\lambda_j f_j=0$ and 
$\tfrac{d^2}{dx^2} f_{j+1}+V f_{j+1}+\lambda_{j+1} f_{j+1}=0$.  
Let $\xi_1,\xi_2,...,\xi_{j-1}$ be the zeros of $f_j$ in the interval $(0,a)$. Since 
$\lambda_{j+1}>\lambda_{j}$, applying Lemma \ref{sturm comparison}, we conclude that 
$f_{j+1}$ must vanish in each 
the intervals $(0,\xi_1), (\xi_1,\xi_2),...,(\xi_{j-1},a)$ and hence that it has 
at least $j+2$ nodes.  Lemma \ref{Courants nodal domain theorem} 
implies it has exactly $j+2$ nodes.  
\end{proof}

The following lemma is proven in \cite{Hille}: 

\begin{lemma}\label{zeros of independent solutions}
Let $f$ and $\hat f$ be two linearly independent solutions of 
Equation \eqref{ode for the eigenvalue problem} for the same $\lambda$, and 
suppose that $f$ has two consecutive zeros $\xi_1$ and $\xi_2$ such that 
$\xi_1 < \xi_2$, then $\hat f$ has one and only one zero in $(\xi_1,\xi_2)$. 
\end{lemma}
\begin{proof}
We may assume $f$ is positive for all $x \in (\xi_1,\xi_2)$, then we have 
$(\tfrac{d}{dx} f)(\xi_1)>0$ and $(\tfrac{d}{dx} f)(\xi_2)<0$. Because 
$f$ and $\hat f$ are independent, $\hat f(\xi_k) \neq 0$ for $k=1,2$.  
Here $\tfrac{d}{dx} ((\tfrac{d}{dx} \hat f) f-(\tfrac{d}{dx} f) \hat f)=0$ for all $x$, 
so $(\tfrac{d}{dx} f)(\xi_1)\hat f(\xi_1) 
= (\tfrac{d}{dx} f)(\xi_2)\hat f(\xi_2)$.  
Hence $\hat f$ cannot keep a constant sign throughout the interval 
$(\xi_1,\xi_2)$, i.e. $\hat f$ has at least one zero in $(\xi_1,\xi_2)$.

Now suppose $\eta_1$ and $\eta_2$ are two zeros of $\hat f$ in $(\xi_1,\xi_2)$. If we 
interchange the roles of $f$ and $\hat f$ in the above argument, 
we conclude that $f$ has at least one zero in $(\eta_1,\eta_2)$, a 
contradiction.  Hence $\hat f$ has exactly one zero in $(\xi_1,\xi_2)$.  
\end{proof}

\begin{lemma}\label{nodal domains of independent eigenfunctions}
Any two eigenfunctions of \eqref{ode for the eigenvalue problem} in $\mathcal F_p$ 
associated with equal eigenvalues have the same number of nodes. 
\end{lemma}
\begin{proof}
Let $f$ and $\hat f$ be two eigenfunctions 
associated with $\lambda_j=\lambda_{j+1}$ in the spectrum of $\mathcal L$ over 
$\mathcal F_p$.  If $f$ and $\hat f$ are 
linearly dependent, then the lemma clearly holds, so we assume they are linearly 
independent.  

Suppose $f$ has $k$ nodes $\xi_1,\xi_2,...,\xi_k \in [0,a)$.  
Then $\hat f(\xi_\ell) \neq 0$ for $\ell=1,...,k$, and by Lemma 
\ref{zeros of independent solutions}, $\hat f$ has a unique node in each of 
$(\xi_1,\xi_2)$, $(\xi_2,\xi_3)$, ..., $(\xi_{k-1},\xi_k)$ and 
$(\xi_{k},\xi_1+a)$.  Hence $\hat f$ has exactly $k$ nodes. 
\end{proof}

\begin{lemma}\label{nodal domains of eigenfunctions associated unequal eigenvalues}
Take $\lambda_j$ as in Theorem \eqref{initial conditions for closed case}.  
Let $f_j$ and $f_{k}$ in $\mathcal F_p$ be two eigenfunctions of $\mathcal L$ 
corresponding to eigenvalues $\lambda_j$ and $\lambda_{k}$ with $\lambda_j < 
\lambda_{k}$, and with either of the initial conditions as in Lemma 
\eqref{initial condition lemma for closed case}.  Let $n_j$ and $n_k$ denote the 
number of nodes in $\Sigma_p$ of $f_j$ and $f_k$, respectively. If $f_j$ and 
$f_k$ have the same initial conditions, resp. different initial conditions, 
then $n_k \geq n_j+2$, resp. $n_k \geq n_j$. 
\end{lemma}
\begin{proof}
Suppose $f_j(0)=f_k(0)=0$, $(\tfrac{d}{dx} f_j)(0)=(\tfrac{d}{dx} f_k)(0)=1$.  So 
$f_j$ has $n_j-1$ nodes between $x=0$ and $x=a$. Then by Lemma \ref{sturm comparison}, 
$f_{k}$ has at least $n_j$ nodes in the open interval $(0,a)$, and so 
$n_k > n_j$. Since $n_j$ and $n_k$ are both even, $n_k \geq n_j+2$.  

Now suppose $f_j(0)=f_k(0)=1$, $(\tfrac{d}{dx} f_j)(0)=(\tfrac{d}{dx} f_k)(0)=0$, 
then $f_j$ has $n_j$ nodes $\xi_1,...,\xi_j$ in the open interval $(0,a)$. Also, 
$f_{j}$ and $f_k$ have the symmetry $f_j(x)=f_j(-x)$ and 
$f_k(x)=f_k(-x)$ for all $x \in [0,a]$, by the 
symmetry \eqref{Vsymmetry}.  By Lemma \ref{sturm comparison}, $f_{k}$ has a node 
in each interval $(\xi_\ell,\xi_{\ell+1})$ for $\ell=1,...,n_j-1$.  Also, it has a node 
in $(-\xi_1,\xi_1)$, so by the above symmetry, it has at least two nodes in 
$(-\xi_1,\xi_1)$, implying $n_k > n_j$ and so $n_k \geq n_j+2$.

If $f_j$ and $f_{k}$ have different initial conditions, then 
Lemma \ref{sturm comparison} immediately implies $n_{k} \geq n_j$.  
\end{proof}

\begin{lemma}\label{comparison of three eigenfunctions}
Take $\lambda_j$ as in Theorem \eqref{initial conditions for closed case}.  
Let $f_{j-1}$, $f_j$ and $f_{j+1}$ in $\mathcal F_p$ be three consecutive 
eigenfunctions associated with $\lambda_{j-1}$, $\lambda_j$ and $\lambda_{j+1}$, 
respectively, each with either of the initial conditions as in Lemma 
\eqref{initial condition lemma for closed case}.  
Let $n_{j-1}$, $n_j$ and $n_{j+1}$ denote the number of nodes 
of $f_{j-1}$, $f_j$ and $f_{j+1}$ respectively.  Then $n_{j+1} \geq n_{j-1}+2$. 
\end{lemma}
\begin{proof}
Here each of the eigenfunctions  $f_{j-1}$, $f_j$ and $f_{j+1}$ has either of the 
two initial conditions given in Lemma 
\eqref{initial condition lemma for closed case}. Thus two of these functions will 
have the same initial conditions, hence by Lemma 
\ref{nodal domains of eigenfunctions associated unequal eigenvalues} we have 
at least one of $n_{j} \geq n_{j-1}+2$ or $n_{j+1} \geq n_{j-1}+2$ or 
$n_{j+1} \geq n_{j}+2$.   Lemma 
\ref{nodal domains of eigenfunctions associated unequal eigenvalues} also 
implies $n_{j-1} \leq n_j \leq n_{j+1}$, hence the result is shown.  
\end{proof}

All of the lemmas in this section (Section \ref{Proofs of the Theorems}) 
immediately imply Theorems \ref{initial conditions for Dirichlet case} and 
\ref{initial conditions for closed case}. 

\section{\bf Computation of the spectrum of $\mathcal L$ over $\mathcal{F}_p$ with 
symmetric $V$}
\label{Computation of the spectra}

A numerical method for computing the spectrum of $\mathcal L$ on the 
function space $\mathcal F_p$ is as follows:
\begin{enumerate}
\item The eigenfunctions are in $\mathcal{F}_p$ and so are periodic, and 
the real-analytic $V \in \mathcal F_p$ is assumed to have the symmetry 
\eqref{Vsymmetry}.  Theorem \ref{initial conditions for closed case} implies we can 
numerically solve \eqref{ode for the eigenvalue problem} for 
$f$ with the initial conditions just 
either $f(0)=1$, $(\tfrac{d}{dx} f)(0)=0$ or $f(0)=0$, 
$(\tfrac{d}{dx} f)(0)=1$ by a numerical ODE solver, and search for the values 
of $\lambda$ that give periodic solutions $f$, i.e. give $f \in \mathcal{F}_p$.  
Such values of $\lambda$ are amongst the $\lambda_j$. 

\item By Theorem \ref{initial conditions for closed case}, we know the eigenspaces are 
at most $2$-dimensional.  If, for some $\lambda = \lambda_j$, one 
of the two types of initial conditions in Theorem 
\ref{initial conditions for closed case} for $f$ 
gives a solution $f \in \mathcal{F}_p$ and the other does not, 
then the eigenspace of $\lambda_j$ is $1$-dimensional; if both types of 
initial conditions give solutions $f \in \mathcal{F}_p$, 
then the eigenspace of $\lambda_j$ is $2$-dimensional.  

\item From Theorem \ref{initial conditions for closed case}, we know that any 
eigenfunction $f$ associated to $\lambda_j$ has exactly $j$ nodes when $j$ is even, 
and $j-1$ nodes otherwise.  So the value of $j$ is determined simply by counting the 
number of nodes of $f$.  Because we can determine $j$, we will 
know when we have found all $\lambda_j \leq M$ for any given $M \in \mathbb{R}$.  
\end{enumerate}

\section{\bf Application to CMC surfaces of revolution in $\mathbb{S}^3$}
\label{Application to CMC surfaces of revolution in S3}

\begin{table}[htb]
{\tiny 
\begin{center}
\begin{tabular}{c|l }
\text{surface} & nonpositive eigenvalues $\lambda_{1,0}, \lambda_{2,0}, 
... , \lambda_{k,0}$ (where $\lambda_{k,0}=0$ and $\lambda_{k+1,0}>0$), 
for the operator $\hat{\mathcal L}$  \\
& \\ \hline
$U_{1}$ & -1.28, -1,  -1, -0.25, 0 \\ \hline 
$U_{2}$ &  -1.08, -1, -1, -0.76,  -0.76,  -0.51, 0\\ \hline 
$U_{3}$ & -1.04, -1, -1, -0.87, -0.87, -0.67, -0.67, -0.52, 0 \\ \hline 
$U_{4}$ & -1.03, -1, -1, -0.91, -0.91, -0.78, -0.78, -0.6, -0.6, -0.48, 0 \\ \hline 
$U_{5}$ & -1.02, -1, -1, -0.94, -0.94, -0.84, -0.84, -0.714, 
-0.714, -0.57,-0.57,-0.48, 0 \\ \hline 
$U_{6}$ & -1.64, -1.48, -1.48, -1, -1,-0.36, 0   \\ \hline 
$U_{7}$ & -1.13, -1.1, -1.1, -1, -1, -0.84, -0.84, -0.64, -0.64, -0.5, 0  \\ \hline 
$U_{8}$ & -1.05, -1.04, -1.04, -1, -1, -0.94, -0.94, -0.86, -0.86, -0.77, -0.77, -0.68, 
-0.68, -0.64, 0  \\ \hline 
$U_{9}$ & -1.029, -1.022, -1.022, -1, -1, -0.96, -0.96, -0.92, -0.92, -0.86, -0.86, 
-0.79, -0.79, -0.72, -0.72, -0.66, -0.66, -0.63, 0   \\ \hline 
$U_{10}$ & -1.28, -1.25, -1.25, -1.15, -1.15, -1, -1, -0.83, -0.83, -0.72, 0  \\ \hline 
$U_{11}$& -1.11, -1.1, -1.1, -1.06, -1.06, -1, -1, -0.92, -0.92, -0.83, 
-0.83, -0.75, -0.75, -0.71, 0   \\ \hline 
$U_{12}$ &-1.04, -1.03, -1.03, -1.02, -1.02, -1, -1, -.97, -.97, -0.94, -0.94, 
-0.896, -0.896, 
-0.85, -0.85, -0.8, -0.8, -0.76, -0.76, -0.73, -0.73, -0.72, 0 \\ \hline
$U_{13}$ & -1.14, -1.13, -1.13, -1.1, -1.1, -1.06, -1.06, -1, -1, -0.94, -0.94, -0.88, 
-0.88, -0.86, 0  \\ \hline 
$U_{14}$ & -1.14, -1.13, -1.13, -1.1, -1.1, -1.06, -1.06, -1, -1, -0.93, -0.93, 
-0.84, -0.84, -0.75, -0.75, -0.67, -0.67, -0.63, 0 \\ \hline
$U_{15}$ & -1.07, -1.06, -1.06, -1.05, -1.05, -1.03, -1.03, -1, -1, -0.96, -0.96, 
-0.92, -0.92, -0.87, -0.87, -0.83, -0.83, 
-0.78, -0.78, -0.75, -0.75, -0.74, 0 \\ \hline 
$U_{16}$ & -1.11, -1.1, -1.1, -1.09, -1.09, -1.07, -1.07, -1.04, -1.04, -1, -1, 
-0.96,-0.96, -0.93, -0.93, -0.9, -0.9, -0.89, 0 \\ \hline
$U_{17}$ & -1.26, -1.25, -1.25, -1.23, -1.23, -1.19, -1.19, -1.14, -1.14, 
-1.08, -1.08, -1, -1, -0.91, -0.91, -0.81, -0.81, -0.71, -0.71, -0.62, -0.62, 
-0.59, 0 \\ \hline
$N_{1}$ & -1.26, -1.19, -1.19, -1, -1, -0.85, 0\\ \hline 
$N_{2}$ & -1.42, -1.31, -1.31, -1, -1, -0.696, 0\\ \hline
$N_{3}$ &  -1.43, -1.37, -1.37, -1.22, -1.22, -1, -1, -0.85, 0 \\ \hline 
$N_{4}$ & -1.85, -1.76, -1.76, -1.47, -1.47, -1, -1, -0.55, 0 \\ \hline
$N_{5}$ & -1.67, -1.62, -1.62, -1.49, -1.49, -1.27, -1.27, -1, -1, -0.83, 0\\ \hline 
$N_{6}$ & -1.7, -1.67, -1.67, -1.58, -1.58, -1.43, -1.43, -1.22, -1.22, 
-1, -1, -0.88, 0 \\ \hline
$N_{7}$ & -1.09, -1.08, -1.08, -1.05, -1.05, -1, -1, -0.95, -0.95, -0.93, 0 \\ \hline
$N_{8}$ & -1.47, -1.45, -1.45, -1.39, -1.39, -1.29, -1.29, -1.16, -1.16, -1, -1, 
-0.84, -0.84, -0.75, 0   \\ \hline 
$N_{9}$ & -1.18, -1.176, -1.176, -1.15, -1.15, -1.11, -1.11, -1.06, -1.06, 
-1, -1, -0.94, -0.94, -0.89, -0.89, -0.87, 0 \\ \hline
$N_{10}$ & -1.31, -1.3, -1.3, -1.29, -1.29, -1.26, -1.26, -1.22, -1.22, -1.18, 
-1.18, -1.12, -1.12, -1.06, -1.06, -1, -1, -0.94, -0.94, -0.9, -0.9, -0.89, 0 \\ \hline
$N_{11}$ & -1.19, -1.18, -1.18, -1.17, -1.17, -1.15, -1.15, -1.12, -1.12, -1.09, 
-1.09, -1.05, -1.05, -1, -1, -0.95, -0.95, -0.91, -0.91, -0.89, -0.89, -0.88, 0
\end{tabular}
\end{center}
}
\caption{Numerical estimates for the nonpositive eigenvalues of the operator 
$\hat{\mathcal L}$ for the specific examples of 
CMC non-flat tori of revolution shown in Figures \ref{figure1} and 
\ref{figure2}.  The values are rounded off to the nearest hundredth 
or thousandth.}\label{table2}
\end{table}

As an application of the numerical approach described in Section 
\ref{Computation of the spectra}, we consider CMC tori of revolution in the 
unit $3$-sphere 
$\mathbb{S}^3 \subset \mathbb{R}^4$ 
and compute the spectra  of their Jacobi operators.  This 
gives us a numerical evaluation of the Morse index of these surfaces.

Let $\mathcal S(x,y):\mathcal T=\{(x,y) \in \mathbb{R}^2|(x,y) \equiv (x+a,y) \equiv 
(x, y+2\pi)\} \rightarrow \mathbb{S}^3$ be a conformal immersion from the torus 
$\mathcal T$ to $\mathbb{S}^3$, 
with mean curvature $H$ and Gauss curvature $K$. When $H$ is constant, $\mathcal S$ is 
critical for a variation problem whose associated  {\em Jacobi operator}  is
\[ - \Delta - \hat V \;\; \text{with}\;\; \hat V=4+4H^2-2K\;,\]
where $\Delta$ is the Laplace-Beltrami operator of the induced metric 
$ds^2=g(dx^2+dy^2)$ for some smooth function $g=g(x,y)$ ($g$ is in fact 
real-analytic in the application here).  
We take $\mathcal S$ to be a non-flat CMC torus of revolution. 

Let us define \[ \mathcal L = - g \Delta - g \hat V = 
\tfrac{\partial^2}{\partial x^2}-\tfrac{\partial^2}{\partial y^2}- V \; , \] 
where $V=g \hat V$.  Then the eigenvalues of $\mathcal L$ form a discrete sequence whose 
corresponding eigenfunctions can be chosen to form an orthonormal basis for the 
$L^2$ norm over $\mathcal T$ with respect to the Euclidean 
metric $dx^2+dy^2$.  Let \[ \lambda_1<\lambda_2 \leq \lambda_3 
\leq ...\uparrow +\infty \] be the spectrum of $\mathcal L$.  

By using Rayleigh quotient characterizations for eigenvalues 
it can be shown that $\mathcal L$ and $- \Delta - \hat V$ will give the same number 
of negative eigenvalues (counted with multiplicity), although these two 
operators will have different eigenvalues. Hence we can use either 
$\mathcal L$ or $- \Delta - \hat V$ to find the Morse index of the surface $\mathcal S$: 

\begin{definition}\label{def index}
The {\em Morse index} $\text{Ind}(\mathcal{S})$ of $\mathcal{S}$ is the sum of 
multiplicities of the negative eigenvalues of $- \Delta - \hat V$ with function space 
the smooth functions from $\mathcal T$ to $\mathbb{R}$.  Equivalently, 
it is the sum of the multiplicities of the negative eigenvalues of $\mathcal L$ 
over the same function space.  
\end{definition}

  \begin{table}[htb]
{ 
    \begin{center}
    \begin{tabular}{c|c|c|c|c|c|c|c|c|c}
      & & & & & & & & & numerical   \\
      \text{surf-} & & & & & & & & & 
      value  \\
     $\text{ace}$ & $s$ & $t$ & a & $k$ & $w$ & 
$\mathcal{B}_1$ & $\mathcal{B}_2$ & $\mathcal{B}_3$ & for \\
  $\mathcal{S}$ & & & & &  & & & & Ind($\mathcal S$) \\\hline
$U_1$ & $0.4078$ & $0.1583$  &$11.7053$& $2$ & $1$& $0$ &$1$ &$1$& $6$ \\ \hline 
$U_2$ & $0.4392$ & $0.0812$  &$21.1215$& $3$ & $1$& $0$ &$1$ &$3$& $8$ \\ \hline 
$U_3$ & $0.4352$ & $0.0758$  &$28.9593$& $4$ & $1$& $0$ &$1$ &$5$& $10$\\ \hline 
$U_4$ & $0.4275$ & $0.0796$  &$36.0835$& $5$ & $1$& $0$ & $1$&$7$& $12$ \\ \hline 
$U_5$ & $0.4259$ & $0.0789$  &$43.5185$& $6$ & $1$& $0$ &$1$ &$9$& $14$ \\ \hline
$U_6$ & $0.4703$ & $0.1697$  &$15.6572$& $3$ & $2$& $0$ &$3$ &$1$& $12$ \\ \hline  
$U_7$ & $0.4431$ & $0.0881$  &$34.0978$& $5$ & $2$& $0$    &$3$ &$5$& $16$\\ \hline 
$U_8$ & $0.4561$ & $0.0559$  &$53.6192$& $7$ & $2$& $0$  &$3$ &$9$& $20$ \\ \hline
$U_9$ & $0.4526$ & $0.0545$  &$69.8309$& $9$ & $2$& $0$  &$3$ &$13$& $24$ \\ \hline 
$U_{10}$ &$0.4949$ &$0.0707$ &$33.7818$ & $5$ & $3$& $0$  &$5$ &$3$& $20$ \\ \hline
$U_{11}$&$0.4738$ & $0.0528$  &$53.0235$& $7$ & $3$ & $0$&$5$ &$7$& $24$\\ \hline 
$U_{12}$&$0.4667$ & $0.0426$  &$89.2538$& $11$& $3$ & $0$&$5$ &$15$& $32$\\ \hline 
$U_{13}$& $0.4987$& $0.0354$ &$56.6566$& $7$ & $4$ & $0$  &$7$ &$5$& $28$ \\ \hline 
$U_{14}$&$0.4659$ &$0.0675$  &$64.3347$& $9$ & $4$ & $0$  &$7$ &$9$& $32$ \\ \hline 
$U_{15}$&$0.47302$&$0.0438$ &$87.7273$& $11$ & $4$ & $0$  &$7$ &$13$& $36$\\ \hline 
$U_{16}$& $0.4993$& $0.0269$ &$77.7075$ & $9$ & $5$ & $0$  &$9$ &$7$& $36$\\ \hline 
$U_{17}$&$0.4719$ &$0.0893$  &$71.551$& $11$ & $6$ & $0$  &$11$ &$9$& $44$\\ \hline 
$N_1$& $0.5112$ & $-0.0502$ &$21.7946$ & $3$ & $1$&  $0$  &$3$ &$1$& $12$\\ \hline 
$N_2$& $0.5061$ & $-0.089$  &$18.6334$ & $3$ & $1$&  $0$  &$3$ &$1$& $12$ \\ \hline
$N_3$& $0.5292$ & $-0.068$  &$26.0688$ & $4$ & $1$& $0$   &$5$ &$1$& $18$ \\ \hline 
$N_4$& $0.5257$ & $-0.155$  &$20.1429$ & $4$ & $1$&  $0$  &$5$ &$1$& $18$ \\ \hline
$N_5$& $0.5501$ & $-0.095$  &$28.6743$ & $5$ & $1$& $0$   &$7$ &$1$& $24$ \\ \hline 
$N_6$& $0.56002$& $-0.092$  &$34.3367$ & $6$ & $1$& $0$   &$9$ &$1$& $30$ \\ \hline 
$N_7$ & $0.5027$& $-0.0197$ &$46.0084$ & $5$ & $2$& $0$   &$5$ &$3$& $20$\\ \hline  
$N_8$& $0.5199$ & $-0.087$  &$43.0143$ & $7$ & $2$&  $0$  &$9$ &$3$& $32$ \\ \hline 
$N_9$& $0.5047$ & $-0.039$  &$62.452$  & $8$ & $3$&  $0$  &$9$ &$5$& $34$\\ \hline 
$N_{10}$&$0.5211$&$-0.051$  &$78.4551$& $11$ & $3$& $0$  &$15$ &$5$& $52$ \\ \hline 
$N_{11}$&$0.5064$& $-0.039$ &$86.0723$ & $11$ & $4$&  $0$  &$13$ &$7$& $48$
\end{tabular}
\end{center}
}
\caption{Here $\mathcal{B}_1$ is the number of eigenvalues less than $-4$, 
$\mathcal{B}_2$ is the number of eigenvalues in $[-4,-1)$, and $\mathcal{B}_3$ is the 
number of eigenvalues in $(-1,0)$, all counted with 
multiplicity.}
\label{table1}
\end{table}

\begin{figure}[phbt]
  \centering
  \includegraphics[scale=0.4]{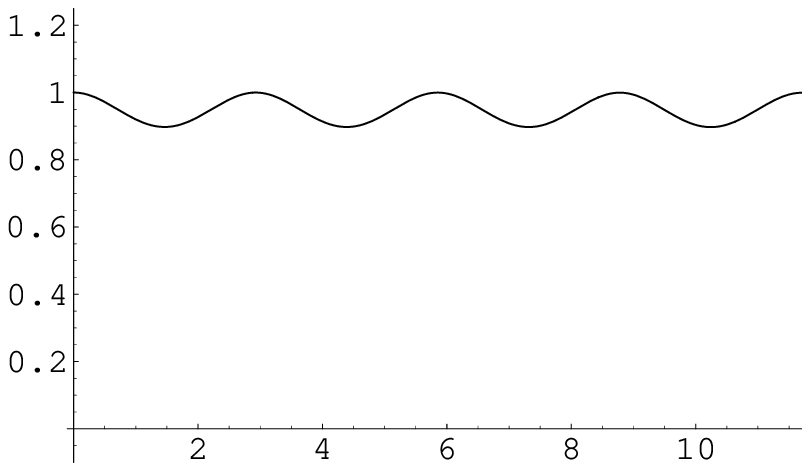}
  \includegraphics[scale=0.9]{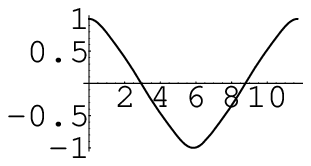}
  \includegraphics[scale=0.9]{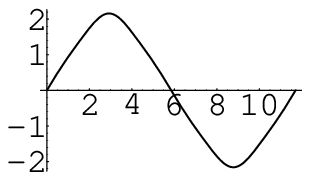}
  \includegraphics[scale=0.9]{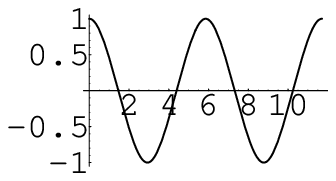}
  \includegraphics[scale=0.9]{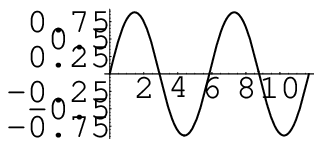}
 \caption{Eigenfunctions associated to the eigenvalues $\lambda_{1,0}$, ..., 
           $\lambda_{4,0}$, $\lambda_{5,0}=0$ of the surface $U_{1}$.}
  \label{A}
\end{figure}

\begin{figure}[phbt]
  \centering
  \includegraphics[scale=0.38]{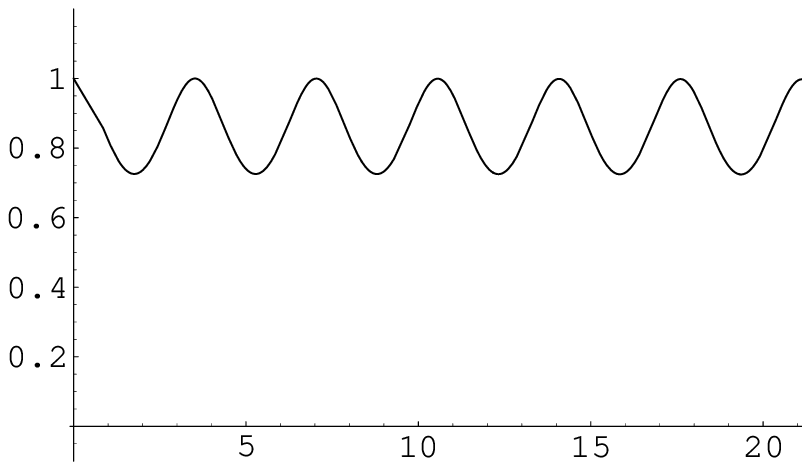}
  \includegraphics[scale=0.5]{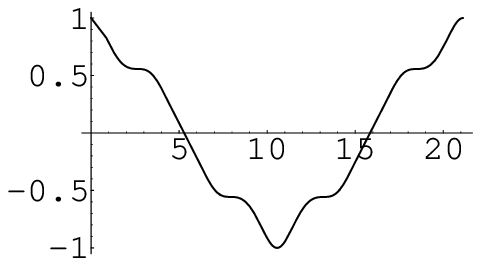}
  \includegraphics[scale=0.5]{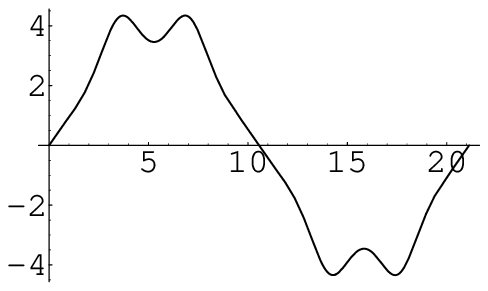}
  \includegraphics[scale=0.5]{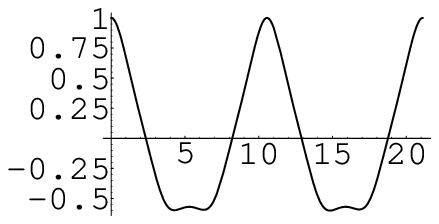}
  \includegraphics[scale=0.5]{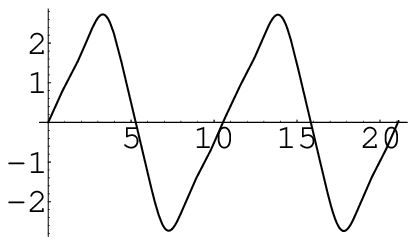}
  \includegraphics[scale=0.5]{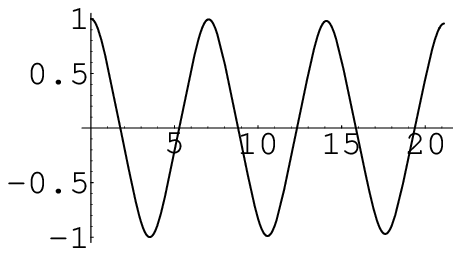}
  \includegraphics[scale=0.5]{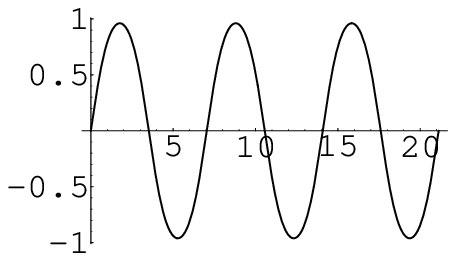}
 \caption{Eigenfunctions associated to the eigenvalues $\lambda_{1,0}$, ..., 
           $\lambda_{6,0}$, $\lambda_{7,0}=0$ of the surface $U_{2}$.}
  \label{B}
\end{figure}

\begin{figure}[phbt]
  \centering
  \includegraphics[scale=0.32]{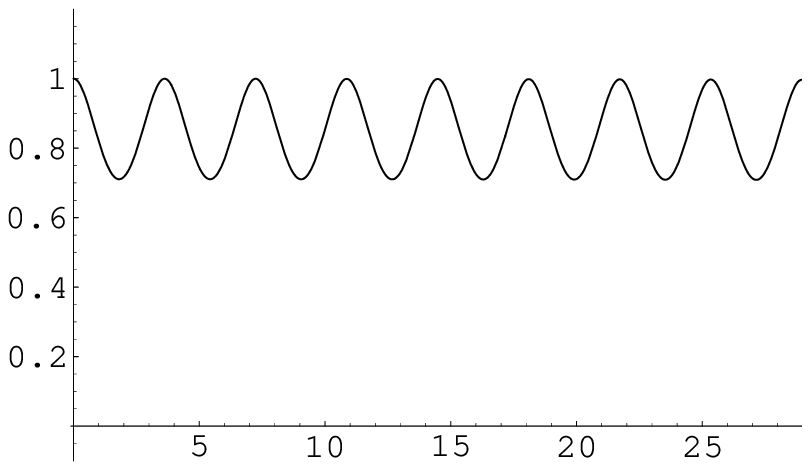}
  \includegraphics[scale=0.65]{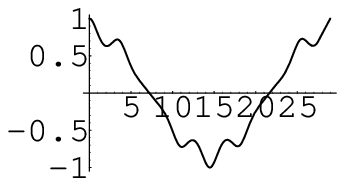}
  \includegraphics[scale=0.65]{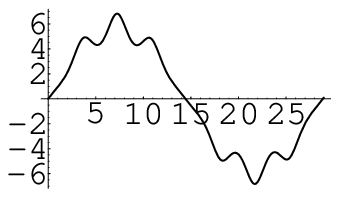}
  \includegraphics[scale=0.65]{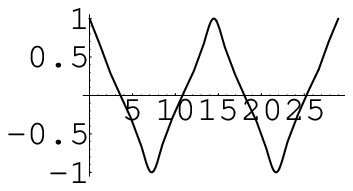} 
  \includegraphics[scale=0.65]{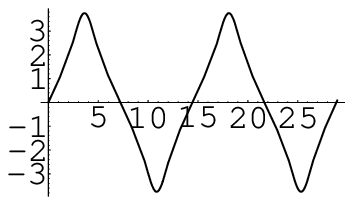}
  \includegraphics[scale=0.65]{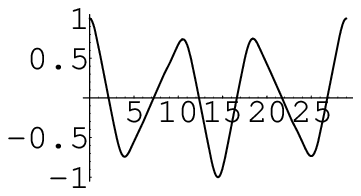}
  \includegraphics[scale=0.65]{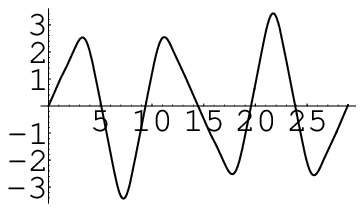}
  \includegraphics[scale=0.65]{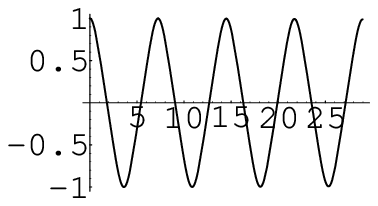}
  \includegraphics[scale=0.65]{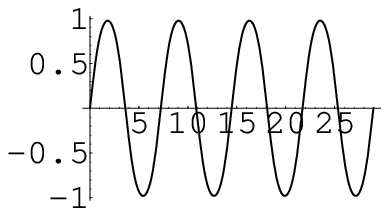}
 \caption{Eigenfunctions associated to the eigenvalues $\lambda_{1,0}$, ..., 
           $\lambda_{8,0}$, $\lambda_{9,0}=0$ of the surface $U_{3}$.}
  \label{C}
\end{figure}

\begin{figure}[phbt]
  \centering
  \includegraphics[scale=0.42]{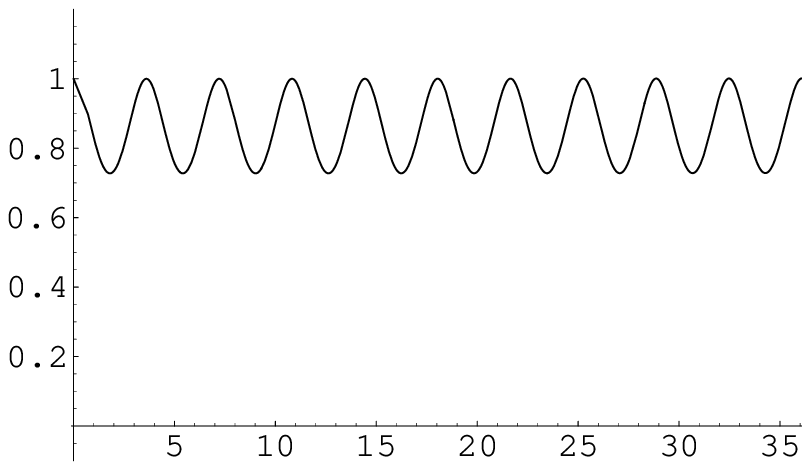}
  \includegraphics[scale=0.85]{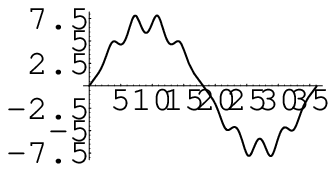}
  \includegraphics[scale=0.85]{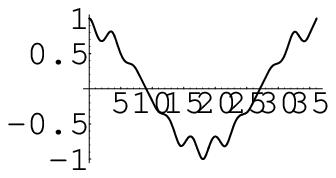}
  \includegraphics[scale=0.85]{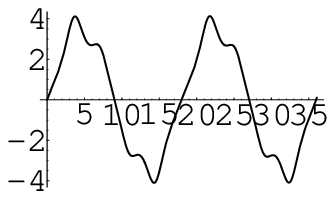}
     \includegraphics[scale=0.85]{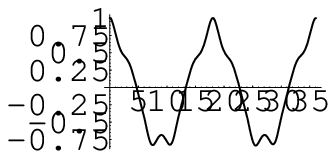}
  \includegraphics[scale=0.85]{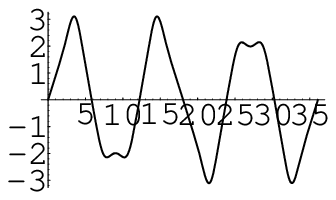}
  \includegraphics[scale=0.85]{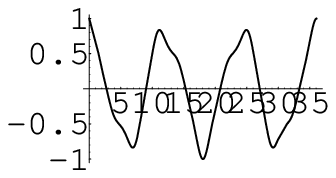}
  \includegraphics[scale=0.85]{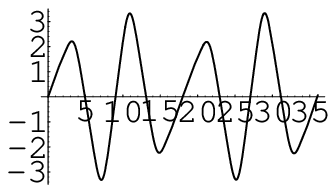}
     \includegraphics[scale=0.9]{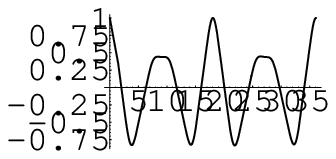}
   \includegraphics[scale=0.9]{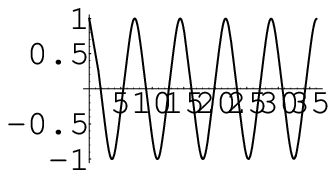}
  \includegraphics[scale=0.9]{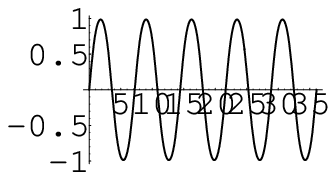}
 \caption{Eigenfunctions associated to the eigenvalues $\lambda_{1,0}$, ..., 
           $\lambda_{10,0}$, $\lambda_{11,0}=0$ of the surface $U_{4}$.}
  \label{D}
\end{figure}

\begin{figure}[phbt]
  \centering
  \includegraphics[scale=0.35]{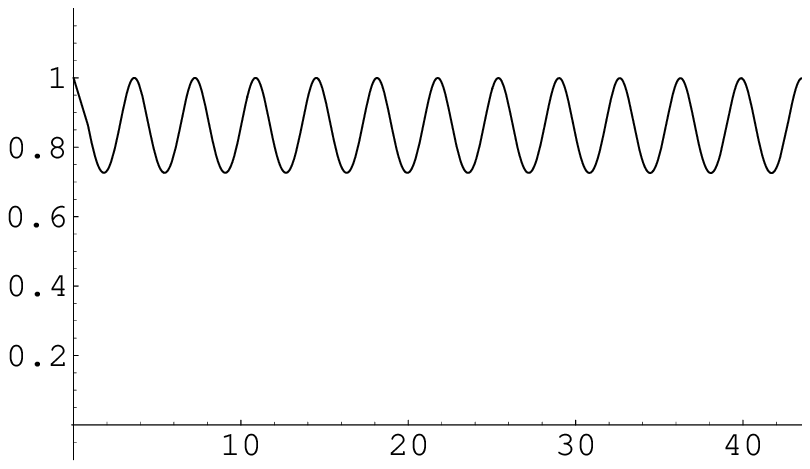}
  \includegraphics[scale=0.35]{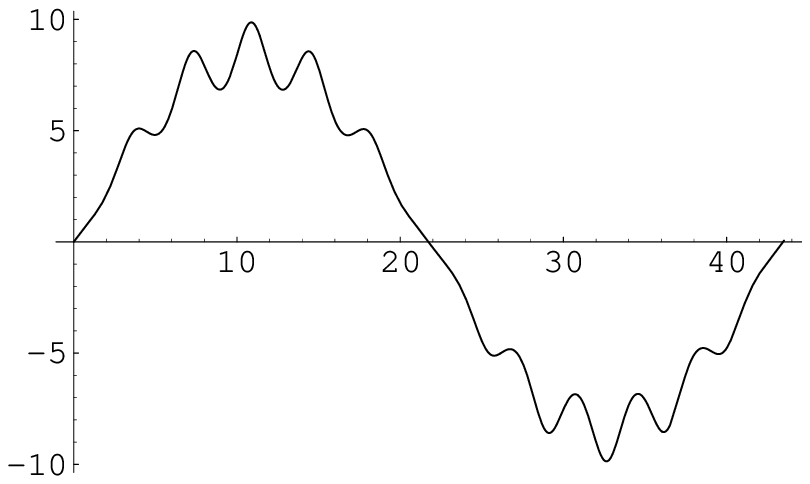}
  \includegraphics[scale=0.35]{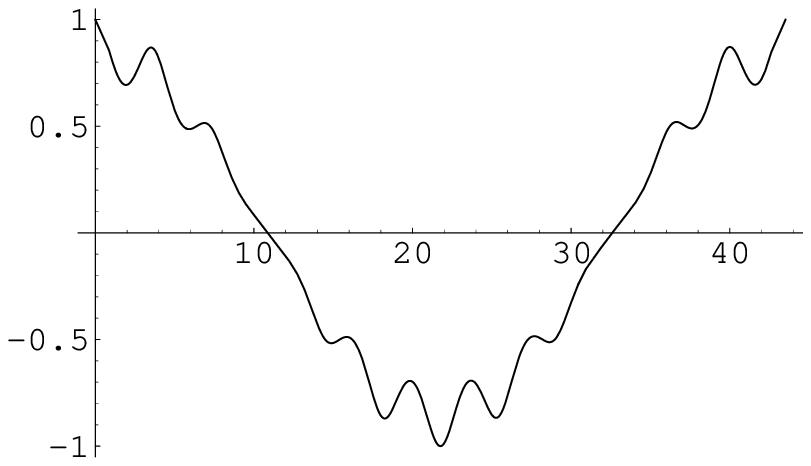}
  \includegraphics[scale=0.35]{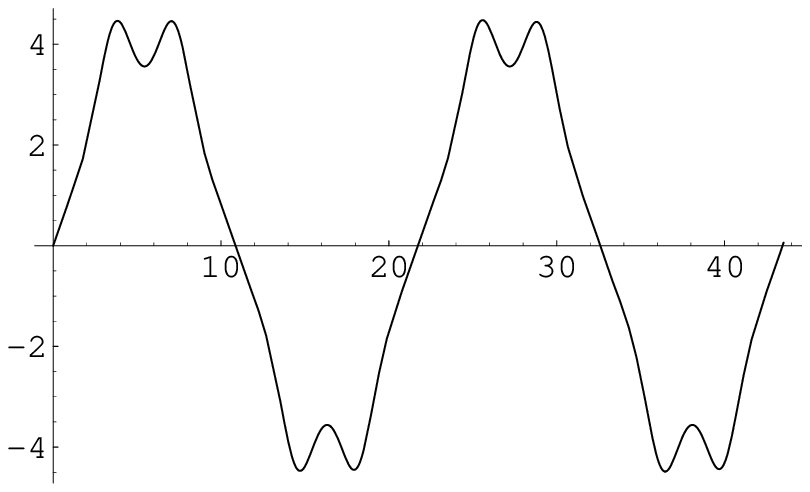}
    \includegraphics[scale=0.35]{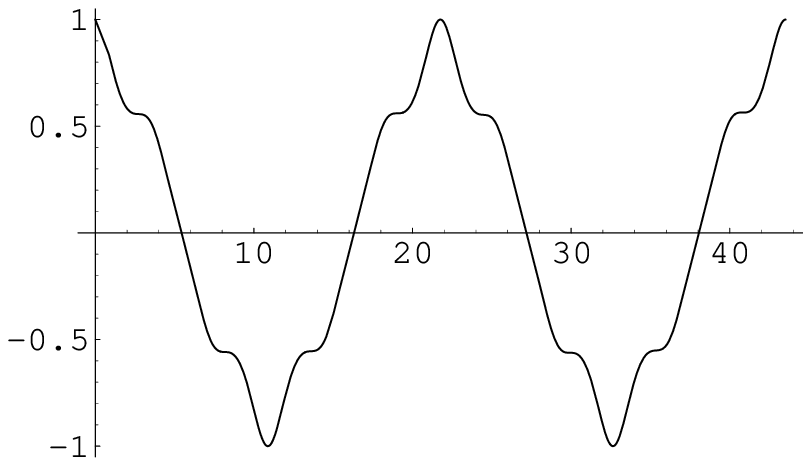}
  \includegraphics[scale=0.35]{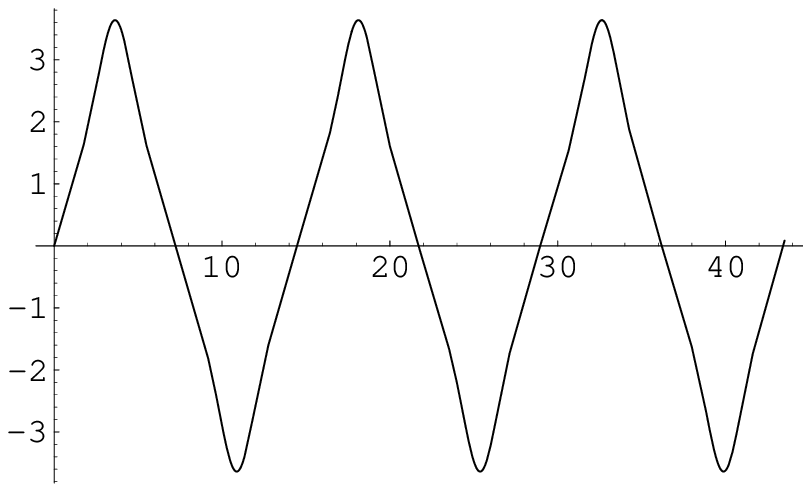}
  \includegraphics[scale=0.35]{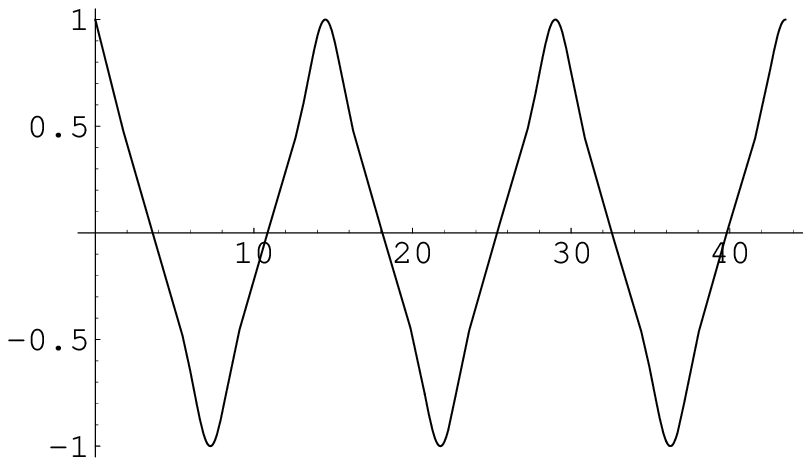}
  \includegraphics[scale=0.35]{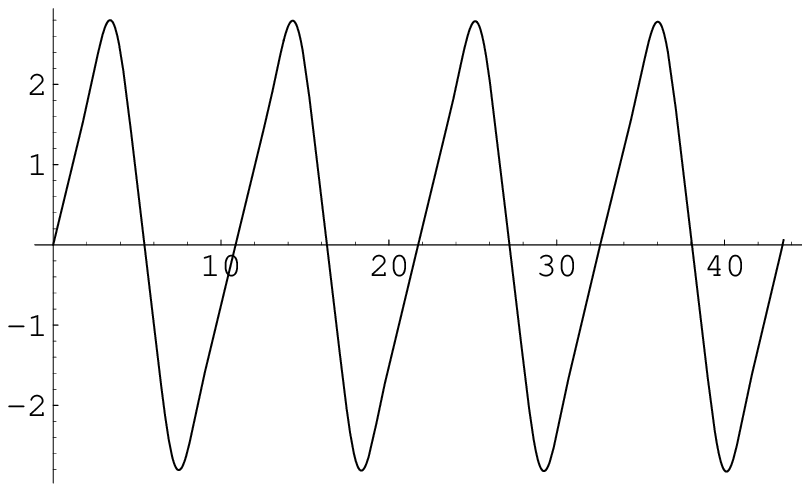}
   \includegraphics[scale=0.35]{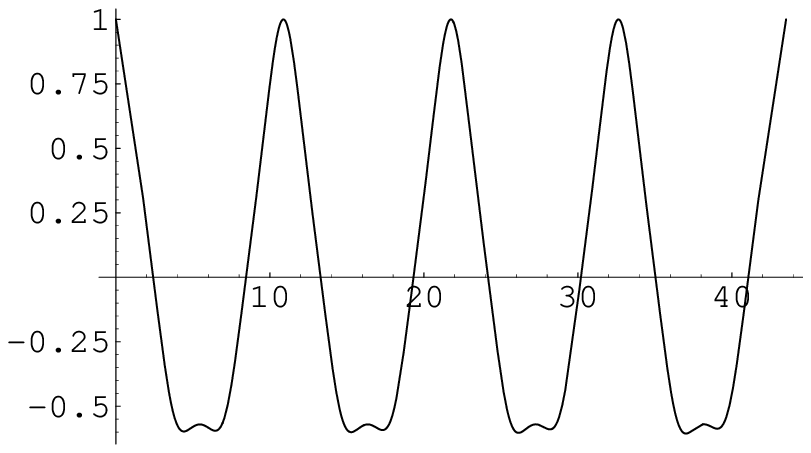}
   \includegraphics[scale=0.35]{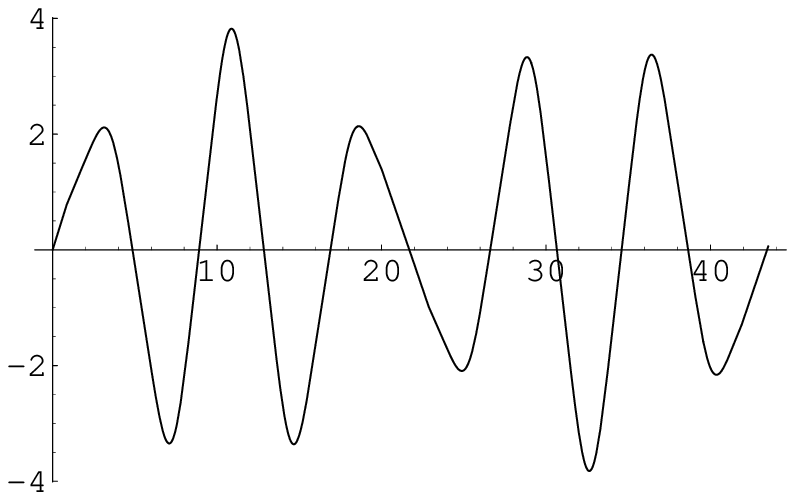}
  \includegraphics[scale=0.35]{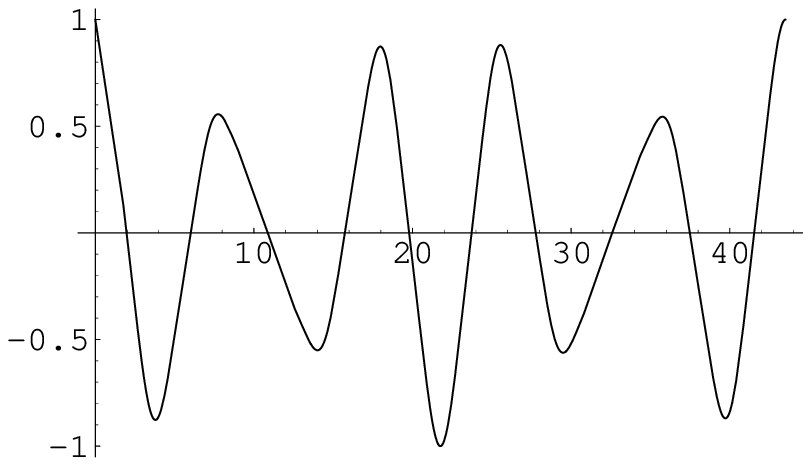}
   \includegraphics[scale=0.35]{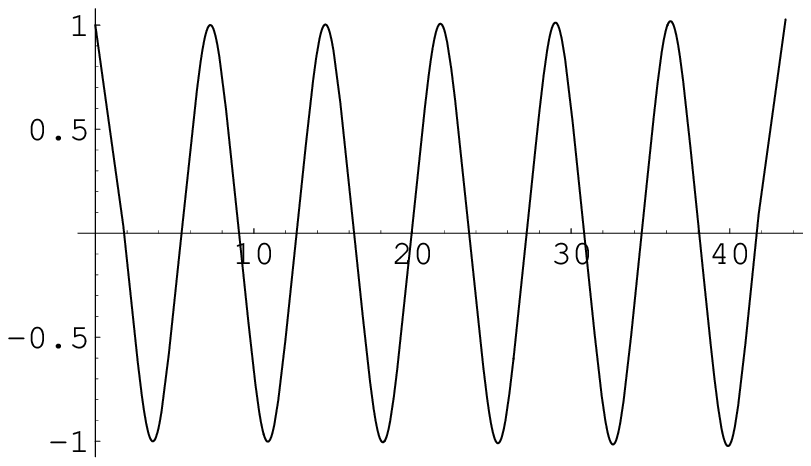}
      \includegraphics[scale=0.35]{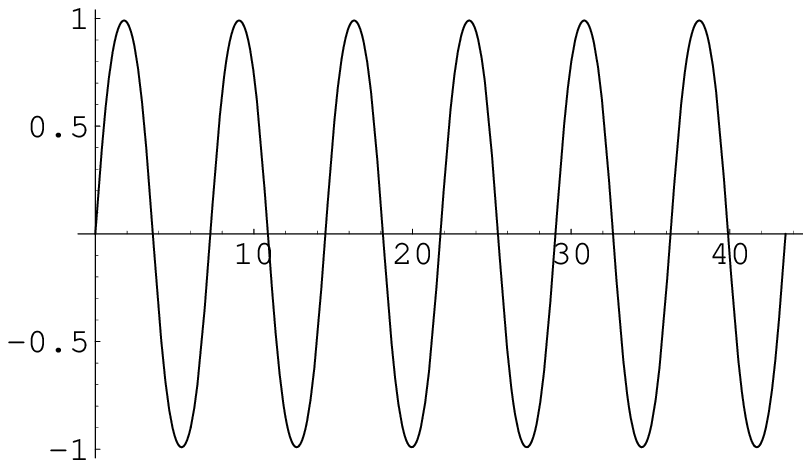}
 \caption{Eigenfunctions associated to the eigenvalues $\lambda_{1,0}$, ..., 
           $\lambda_{12,0}$, $\lambda_{13,0}=0$ of the surface $U_{5}$.}
  \label{E}
\end{figure}

\begin{figure}[phbt]
  \centering
  \includegraphics[scale=0.3]{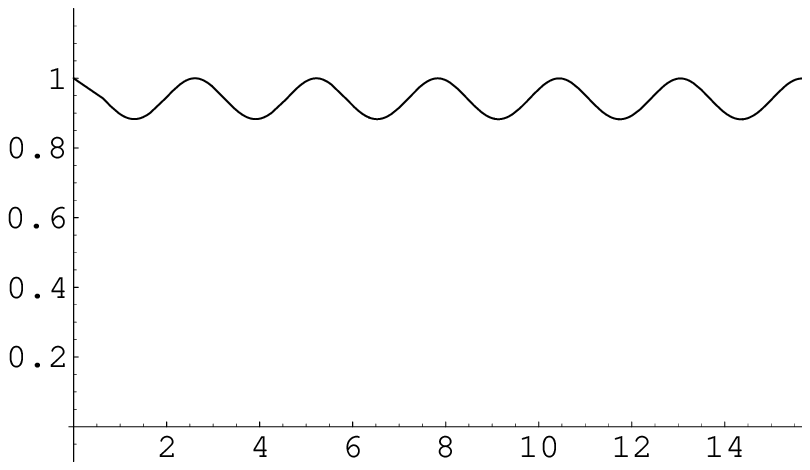}
  \includegraphics[scale=0.3]{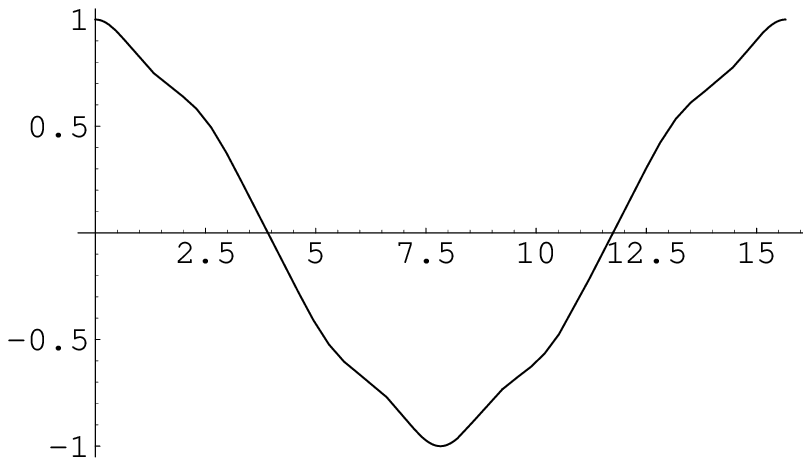}
  \includegraphics[scale=0.3]{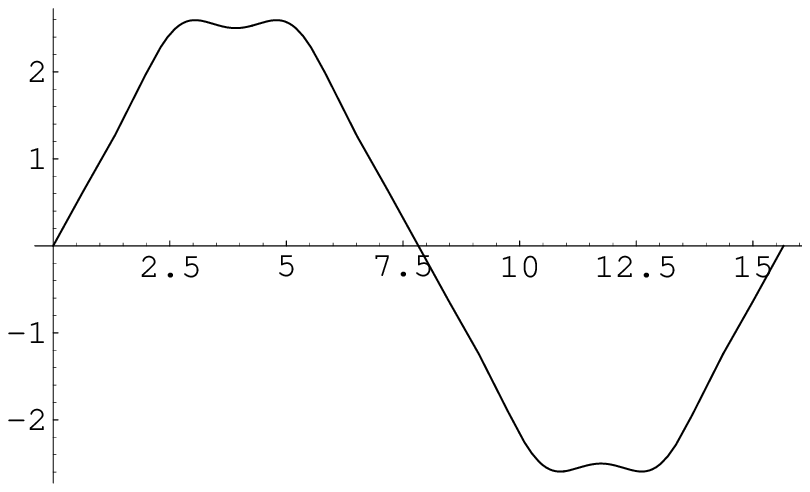}
  \includegraphics[scale=0.3]{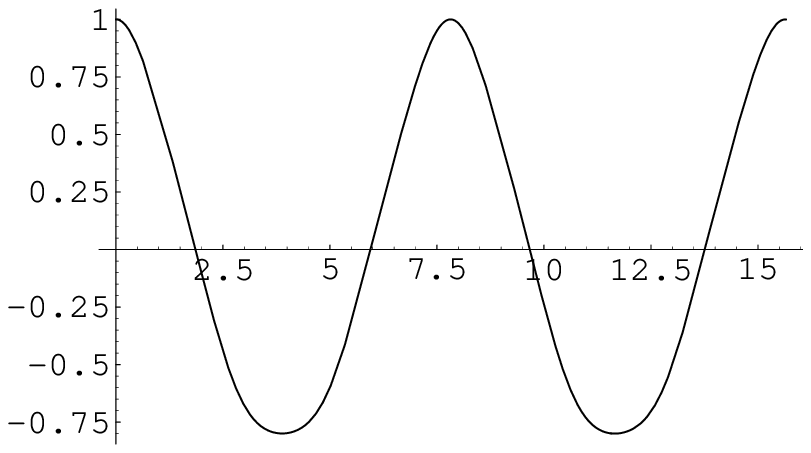}
    \includegraphics[scale=0.3]{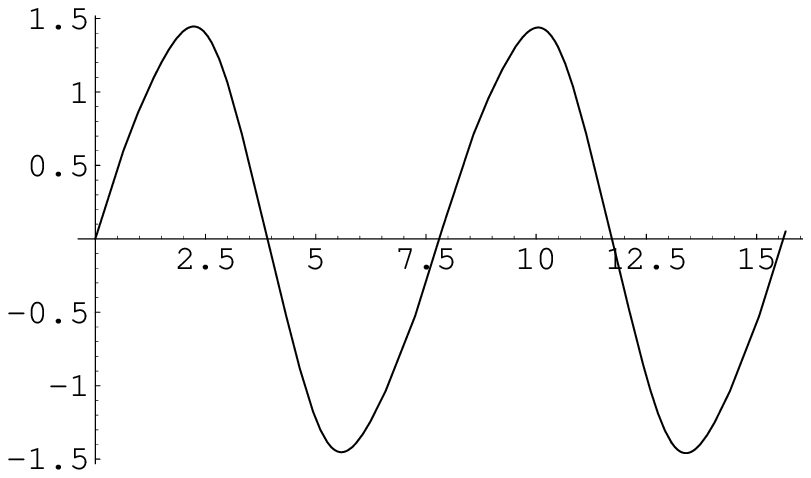}
  \includegraphics[scale=0.3]{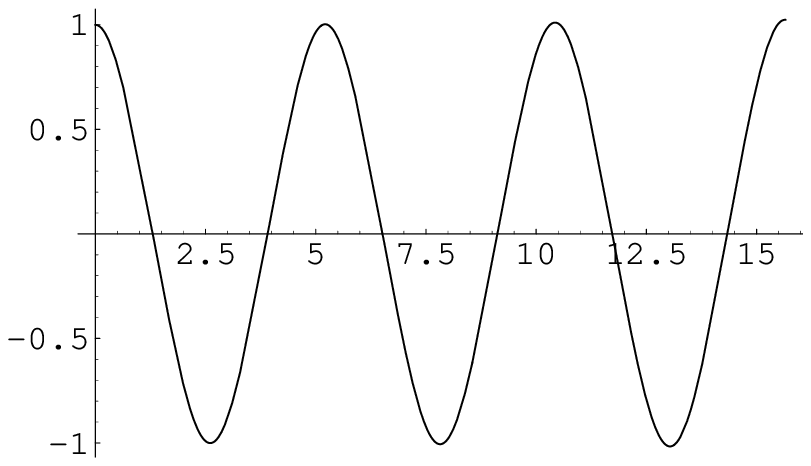}
  \includegraphics[scale=0.3]{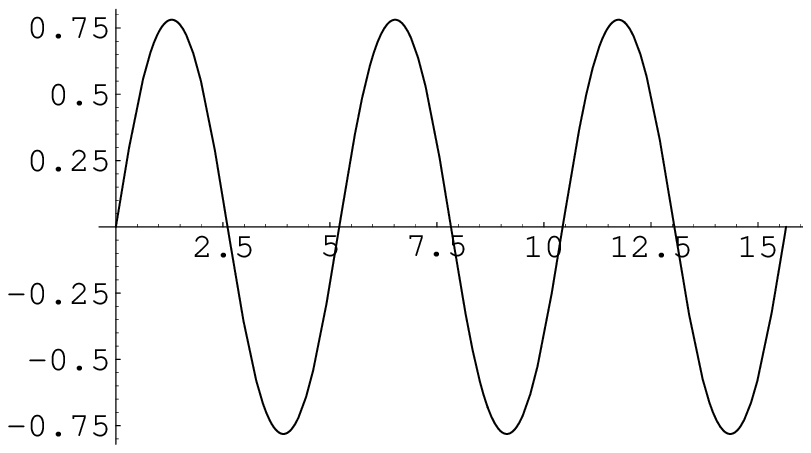}
    \caption{Eigenfunctions associated to the eigenvalues $\lambda_{1,0}$, ..., 
           $\lambda_{6,0}$, $\lambda_{7,0}=0$ of the surface $U_{6}$.}
  \label{F}
\end{figure}

\begin{figure}[phbt]
  \centering
  \includegraphics[scale=0.43]{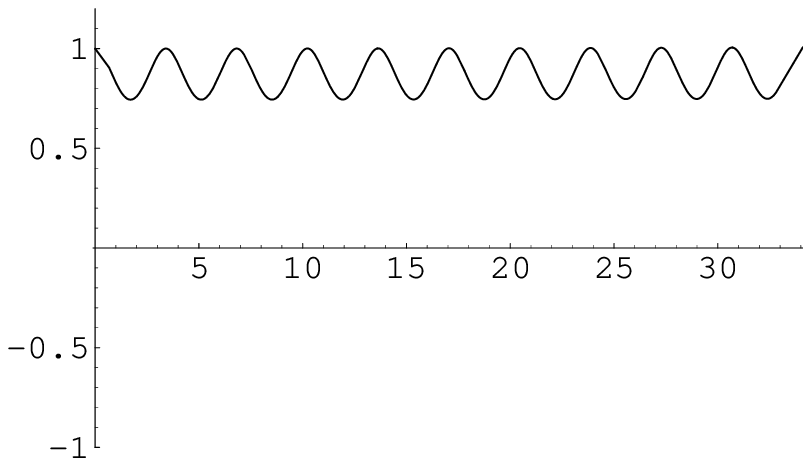}
  \includegraphics[scale=0.85]{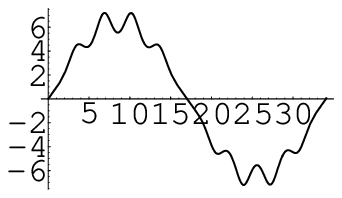}
  \includegraphics[scale=0.85]{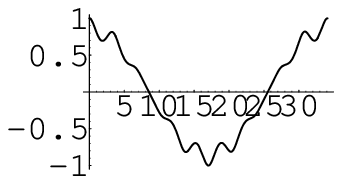}
  \includegraphics[scale=0.85]{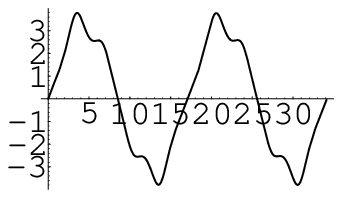}
    \includegraphics[scale=0.85]{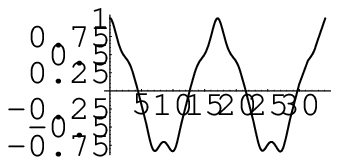}
  \includegraphics[scale=0.85]{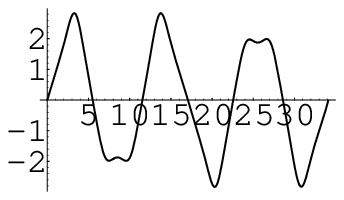}
  \includegraphics[scale=0.85]{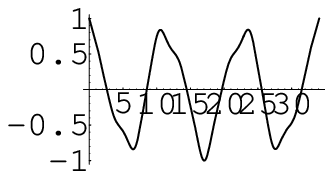}
  \includegraphics[scale=0.85]{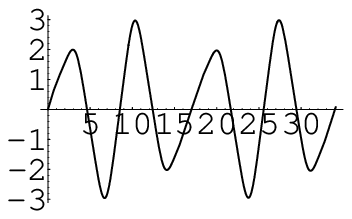}
   \includegraphics[scale=0.9]{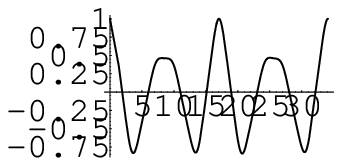}
   \includegraphics[scale=0.9]{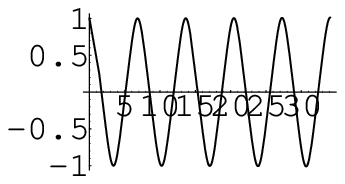}
  \includegraphics[scale=0.9]{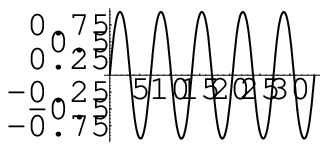}
 \caption{Eigenfunctions associated to the eigenvalues $\lambda_{1,0}$, ..., 
           $\lambda_{10,0}$, $\lambda_{11,0}=0$ of the surface $U_{7}$.}
  \label{G}
\end{figure}

A $C^\infty$ function $f=f(x,y):\mathcal T \rightarrow \mathbb{R}$ can 
be decomposed into a series of spherical harmonics as follows: 
\begin{equation}\label{harmonics} 
f=\sum_{j=0}^{\infty} u_{j,1}(x) \cos(jy)+u_{j,2}(x) \sin(jy) \; , 
\end{equation}
where $u_{j,1}, u_{j,2} \in \mathcal F_p$ for $\Sigma_p$ with the given $a>0$.  
The operator $\hat{\mathcal L}$ on the function space $\mathcal F_p$ is defined by 
\begin{equation*} 
\hat{\mathcal L} = -\tfrac{d^2}{dx^2} -V \;, \end{equation*} 
and the spectrum \[ \lambda_{1,0} < \lambda_{2,0} \leq 
\lambda_{3,0} \leq ...\uparrow +\infty \] of $\hat{\mathcal L}$ has all the 
analogous properties as those of the spectrum for $\mathcal L$.  
Furthermore, by uniqueness of the spherical harmonics decomposition, 
$f$ is an eigenfunction of $\mathcal L$ for the eigenvalue $\lambda$ if and only if 
each $u_{j,k}$, $k=1,2$, is an eigenfunction of $\hat{\mathcal L}$ for the 
eigenvalue $\lambda-j^2$.  And if $f$ is not identically zero, then some 
$u_{j,k}$ will also be not identically zero.  Thus we can say: 
\begin{itemize}
\item $\lambda$ is an eigenvalue for the operator $\mathcal L$ if and only if 
$\lambda-n^2$ is an eigenvalue for the operator $\hat{\mathcal L}$ for some $n \in 
\mathbb{N} \cup \{ 0 \}$.  
\item For any eigenvalue $\lambda_{j,0} < -n^2$, $n \in \mathbb{N} \cap [2,\infty)$, of 
$\hat{\mathcal L}$, with associated eigenfunction $f_j \in \mathcal F_p$, 
the eigenvalues of $\mathcal L$ associated to the eigenfunctions 
$f_j \cdot \cos (k y)$ and $f_j \cdot \sin (k y)$, for integers $k \leq [0,n]$, 
will be negative.  
\end{itemize}
Furthermore, we can conclude the following:
\begin{lemma}\label{lm index}
We have 
\[ \text{Ind}(\mathcal S)=\sum_{j \in \mathbb{N}} 
\lambda_{j,0} \cdot  \ell(\lambda_{j,0}) \; , \;\;\; 
\text{where} \;\;\; 
\ell (\lambda)=\left\{ \begin{array}{ll} 0 & \text{if} \;\; 
\lambda \geq 0 \; , \\ 2 i-1 & \text{if} \;\; \lambda \in [-i^2,-(i-1)^2) 
\;\; \text{for}\;\; i\in  \mathbb{N} \; . \end{array}\right. \]
\end{lemma}

\begin{proof}
Let $E(\lambda)$, resp. $E_0(\lambda)$, denote the eigenspace of solutions 
of $\mathcal L(f)=\lambda f$ for smooth $f: \mathcal T \to \mathbb{R}$, 
resp. $\hat{\mathcal L}(f)=\lambda f$ for $f \in \mathcal F_p$.  Then 
$\dim E(\lambda) = 0$, resp. $\dim E_0(\lambda)=0$ whenever $\lambda$ is 
not an eigenvalue of $\mathcal L$, resp. $\hat{\mathcal L}$, and is a positive 
integer otherwise.  Then, by the uniqueness of the spherical harmonics decomposition, 
\[ \sum_{\lambda < 0} \dim E(\lambda) = \sum_{\lambda < 0} \left( \dim E_0(\lambda) + 
2 \sum_{j \geq 1} \dim E_0(\lambda-j^2) \right) = \left( \sum_{\lambda < 0} 
\dim E_0(\lambda) \right) + 
2 \sum_{j \geq 1} \sum_{\lambda < -j^2} \dim E_0(\lambda) \; . \] 
\end{proof}

Here $\mathcal S$ is a CMC surface of revolution, so, following \cite{WN}, we can 
consider 
\[\mathcal L = -\tfrac{\partial^2}{ \partial x^2} -\tfrac{\partial^2}{ \partial y^2}-2 
v^2 - 32 s^2 t^2 v^{-2} \;\;\; \text{and} \;\;\; 
\hat{\mathcal L} = -\tfrac{d^2}{ dx^2}-2 v^2 - 32 s^2 t^2 v^{-2} \; ,
\]  
where $s \in \mathbb{R}^{+}$, $t \in (-s,s) \setminus \{ 0 \}$, 
$\gamma \in (0,\pi/4]$, and $s$, $t$, $\gamma$ (we note that 
$\cot(2 \gamma)$ is the mean curvature 
of $\mathcal S$) satisfy the 
conditions $(s+t)^2-4 s t \sin^2 \gamma=1/4 $ and 
\[st \in \left( -(16 \sin^2 \gamma)^{-1},0 \right) \cup 
\left( 0,(16 \cos^2 \gamma)^{-1} \right) \;,\]
and $v=v(x)$ is the elliptic function 
\[v=\frac{2 t}{\text{dn}_{\tau} (2sx)}\;\;\text{with period}\;\; 
x_0=\frac{1}{s} \int_0^1 \frac{d\varrho}{\sqrt{(1-\varrho^2)(1- \tau^2 \varrho^2)}}\;, 
\;\;\; \text{where} \;\; \tau=\sqrt{1-t^2/s^2} \; . \]  

When $st > 0$, we have unduloidal surfaces.
When $st<0$, we have either nodoidal or unduloidal surfaces (see \cite{WN}).  

Using the method in Section \ref{Computation of the spectra}, we can 
numerically compute the negative eigenvalues of $\hat{\mathcal L}$, and can then apply 
Lemma \ref{lm index} to find Ind($\mathcal S$).  
We do this for the CMC tori of revolution shown in Figures \ref{figure1} and 
\ref{figure2}.  In \cite{WN}, it is shown that $0$ is an eigenvalue of 
$\hat{\mathcal L}$, and $-1$ is an eigenvalue of $\hat{\mathcal L}$ with multiplicity 
$2$.  Since $\ell(\lambda)$ is discontinuous at $\lambda = 0$ and 
$\lambda = -1$, it is crucial 
to know that both $0$ and $-1$ are eigenvalues of $\hat{\mathcal L}$ in order to 
determine Ind($\mathcal S$).  Furthermore, as the eigenvalue $-1$ has multiplicity 
$2$ and $\lambda_{1,0}$ must be simple, we have $\lambda_{1,0}<-1$.  (In 
the numerical experiments here, we find that $0$ is always a simple eigenvalue.)

\begin{figure}[phbt]
  \centering
  \includegraphics[scale=0.4]{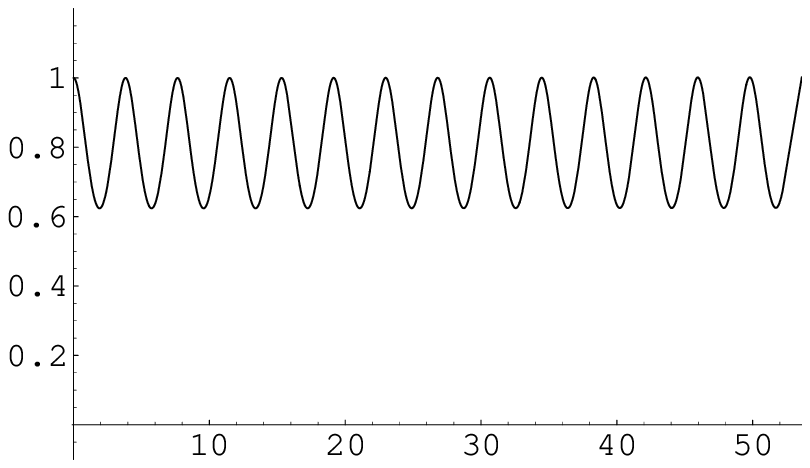}
  \includegraphics[scale=0.8]{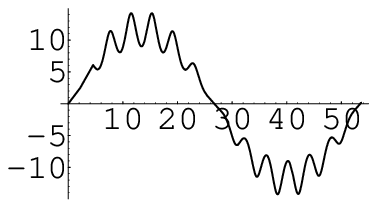}
  \includegraphics[scale=0.8]{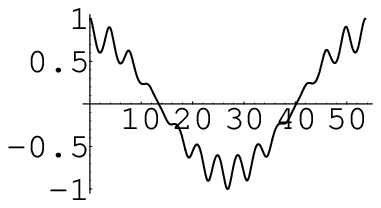}
  \includegraphics[scale=0.8]{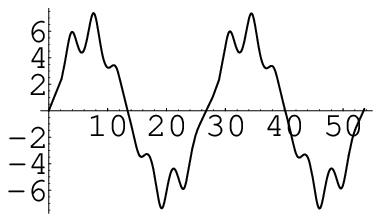}
  \includegraphics[scale=0.8]{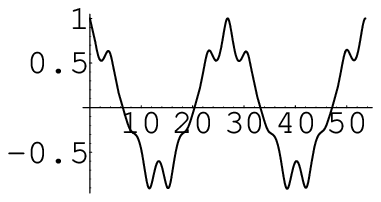}
  \includegraphics[scale=0.8]{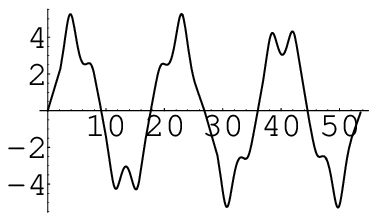}
  \includegraphics[scale=0.8]{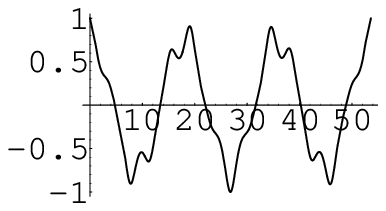}
  \includegraphics[scale=0.8]{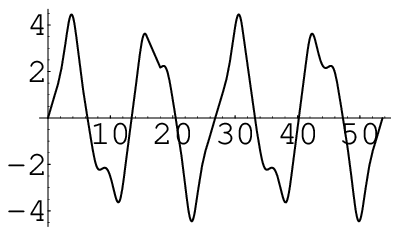}
  \includegraphics[scale=0.8]{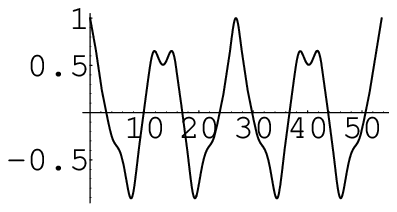}
   \includegraphics[scale=0.8]{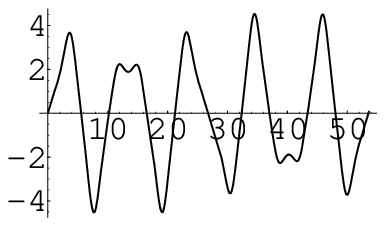}
  \includegraphics[scale=0.8]{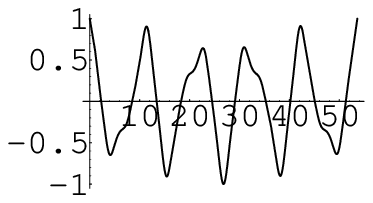}
   \includegraphics[scale=0.8]{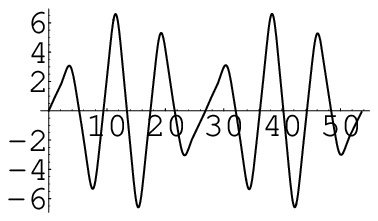}
      \includegraphics[scale=0.85]{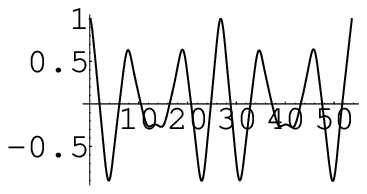}
     \includegraphics[scale=0.85]{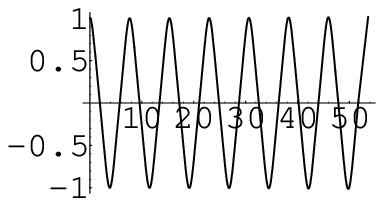}
      \includegraphics[scale=0.85]{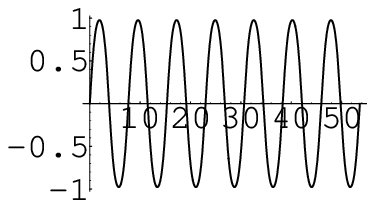}
 \caption{Eigenfunctions associated to the eigenvalues $\lambda_{1,0}$, ..., 
           $\lambda_{14,0}$, $\lambda_{15,0}=0$ of the surface $U_{8}$.}
  \label{H}
\end{figure}

\begin{figure}[phbt]
  \centering
  \includegraphics[scale=0.35]{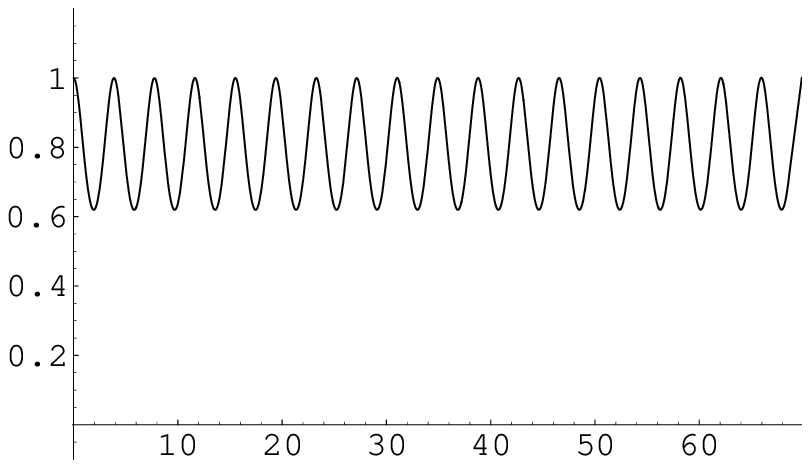}
  \includegraphics[scale=0.35]{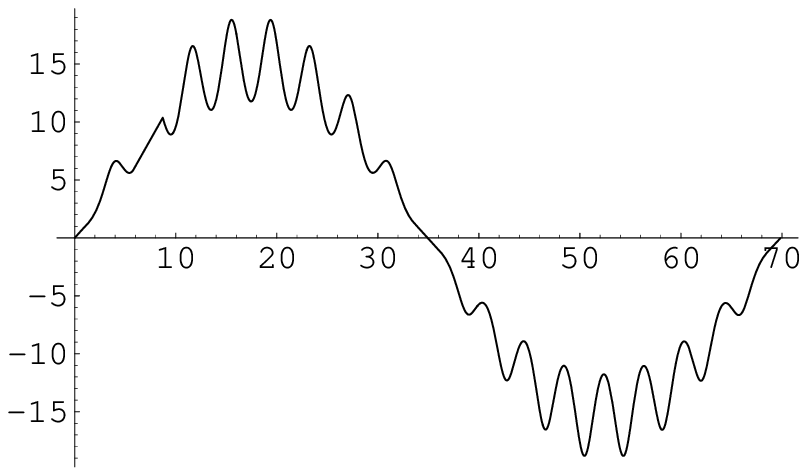}
  \includegraphics[scale=0.35]{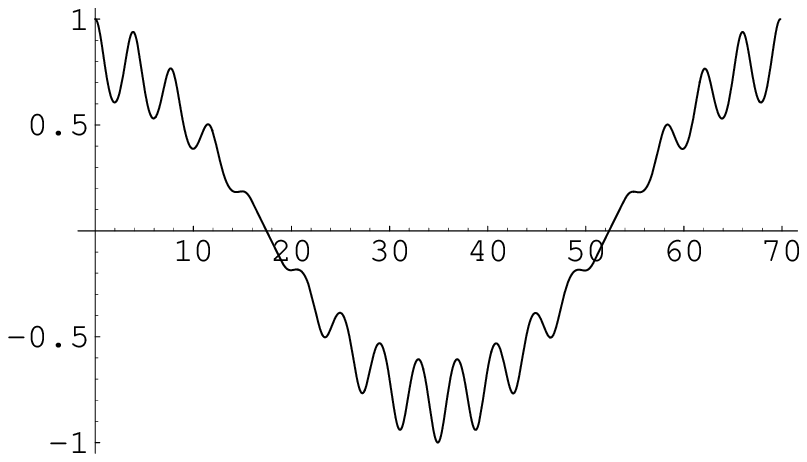}
  \includegraphics[scale=0.35]{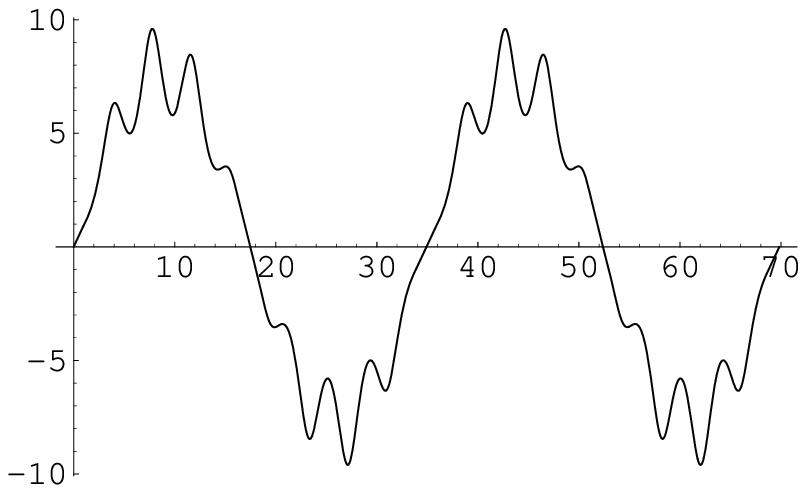}
    \includegraphics[scale=0.35]{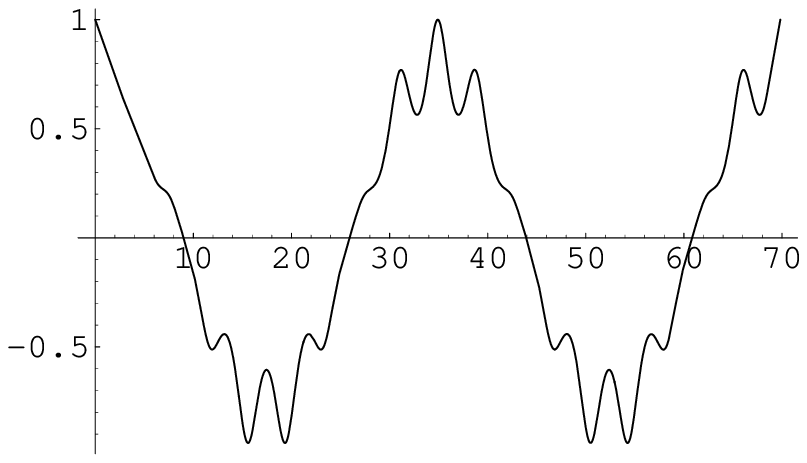}
  \includegraphics[scale=0.35]{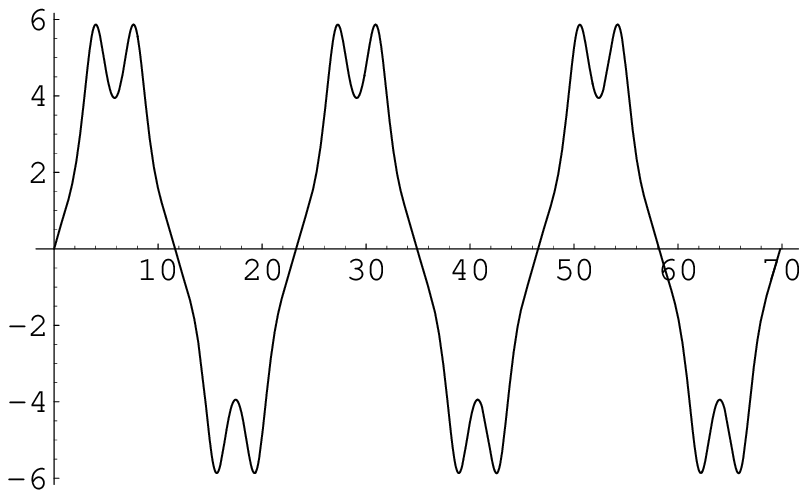}
  \includegraphics[scale=0.35]{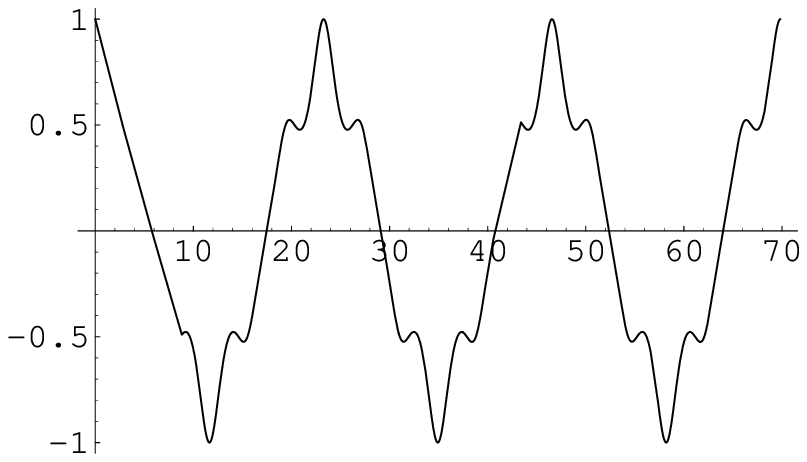}
  \includegraphics[scale=0.35]{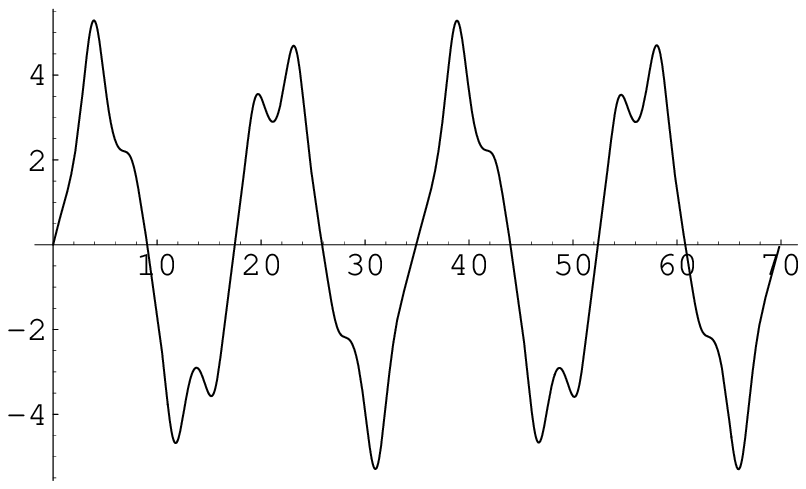}
   \includegraphics[scale=0.35]{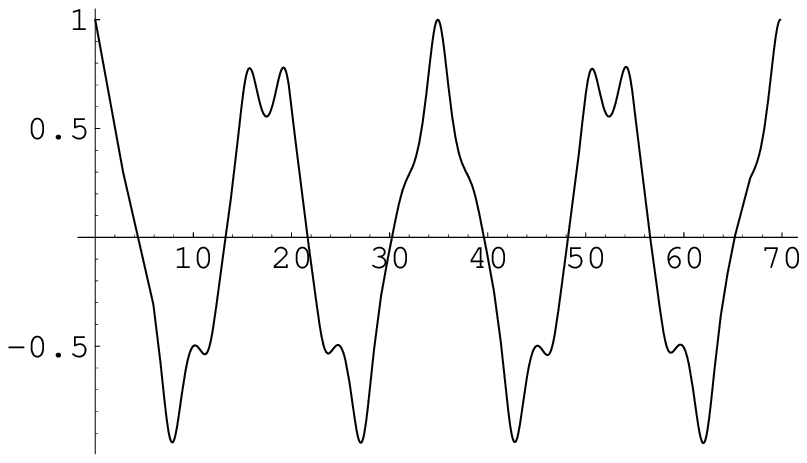}
   \includegraphics[scale=0.35]{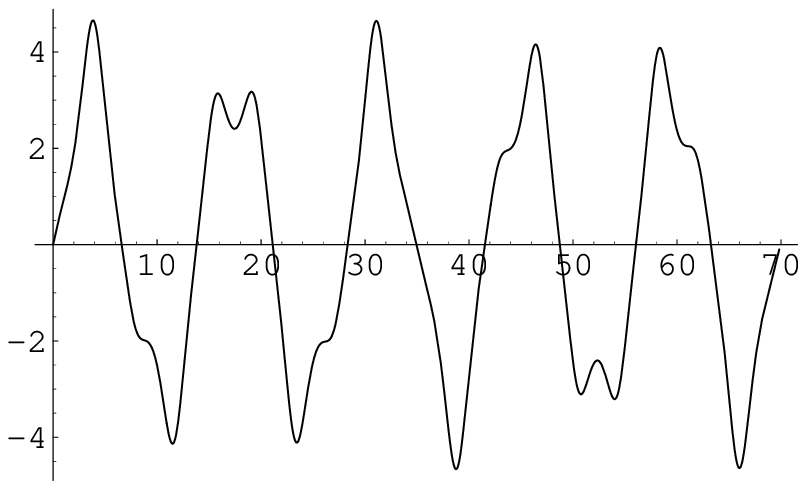}
  \includegraphics[scale=0.35]{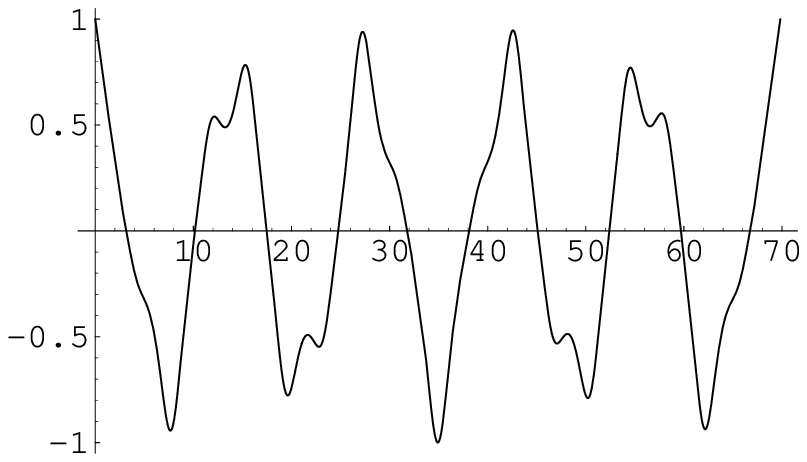}
   \includegraphics[scale=0.35]{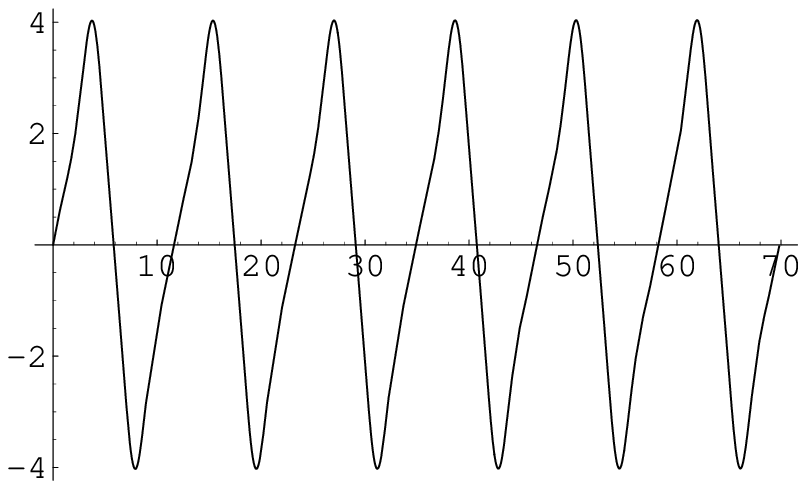}  
     \includegraphics[scale=0.35]{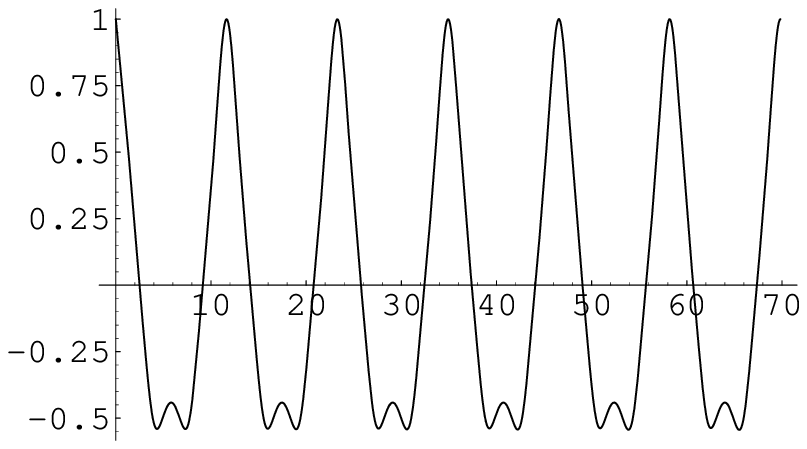}
     \includegraphics[scale=0.35]{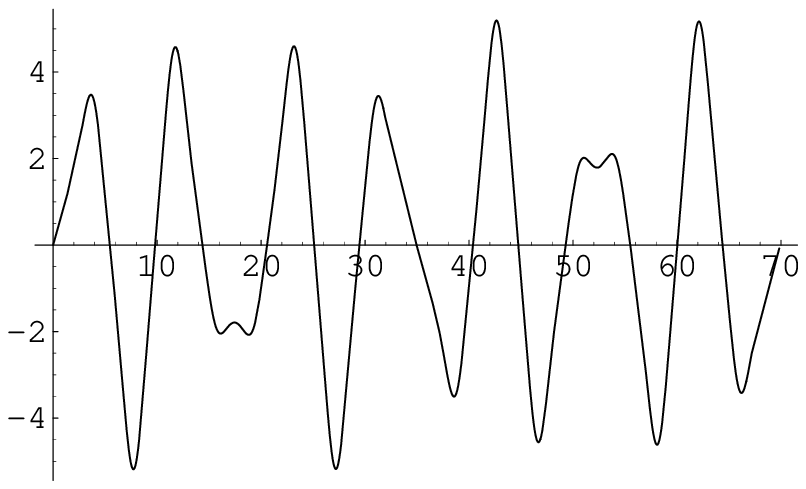}
      \includegraphics[scale=0.35]{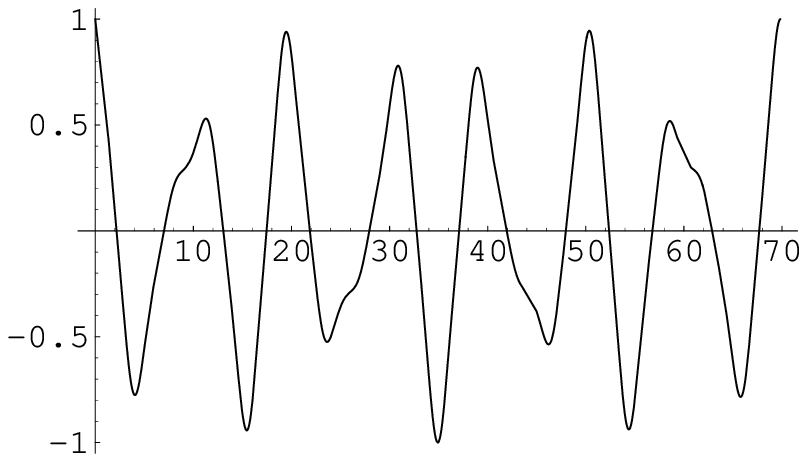}
             \includegraphics[scale=0.35]{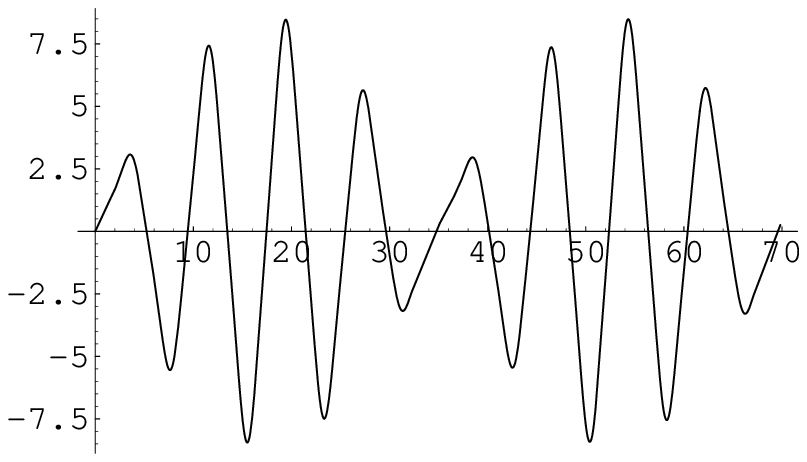}    
      \includegraphics[scale=0.35]{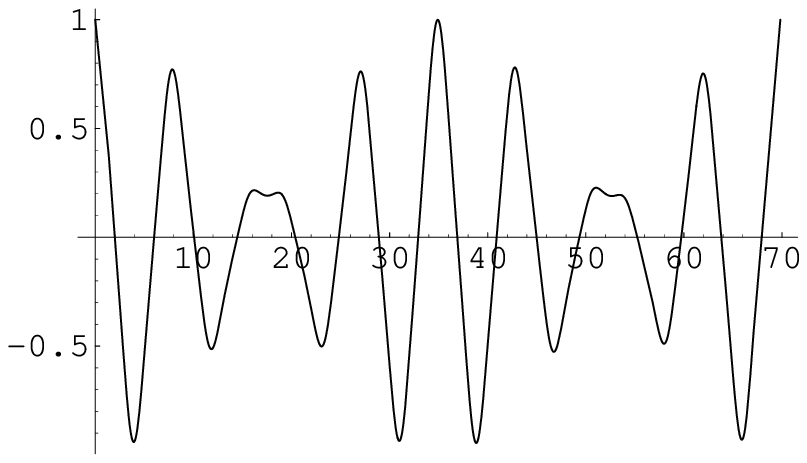}
      \includegraphics[scale=0.35]{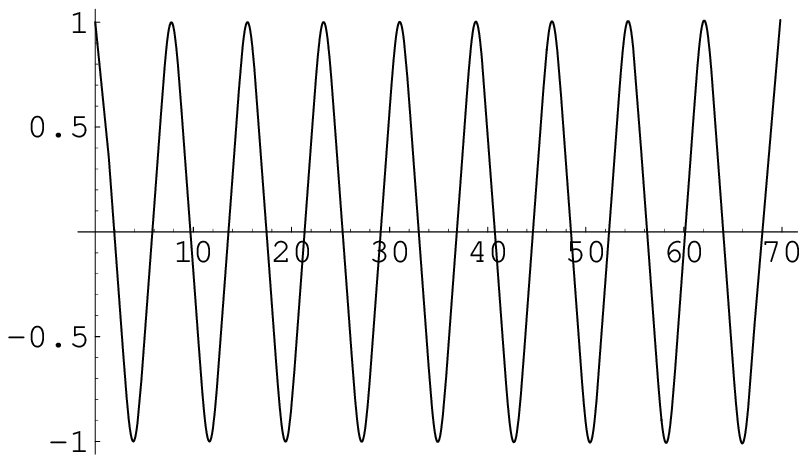}
              \includegraphics[scale=0.35]{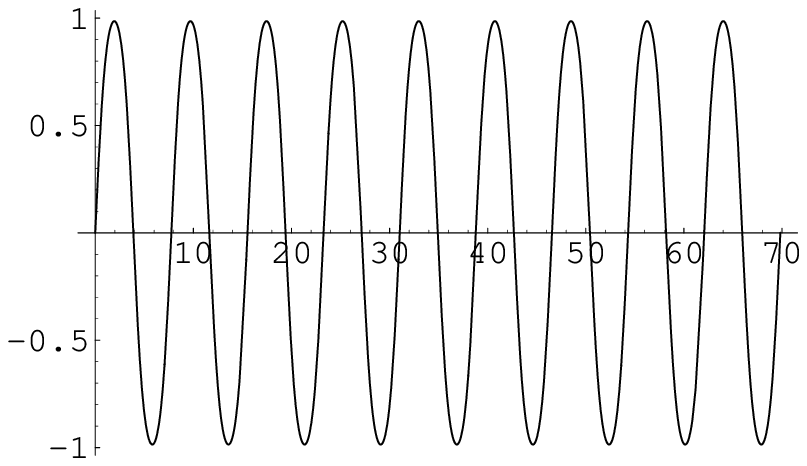}
 \caption{Eigenfunctions associated to the eigenvalues $\lambda_{1,0}$, ..., 
           $\lambda_{18,0}$, $\lambda_{19,0}=0$ of the surface $U_{9}$.}
  \label{I}
\end{figure}

\begin{figure}[phbt]
  \centering
  \includegraphics[scale=0.3]{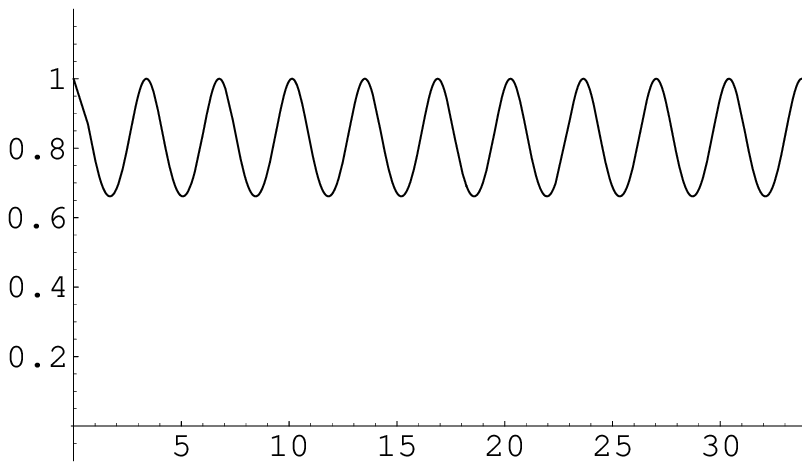}
  \includegraphics[scale=0.3]{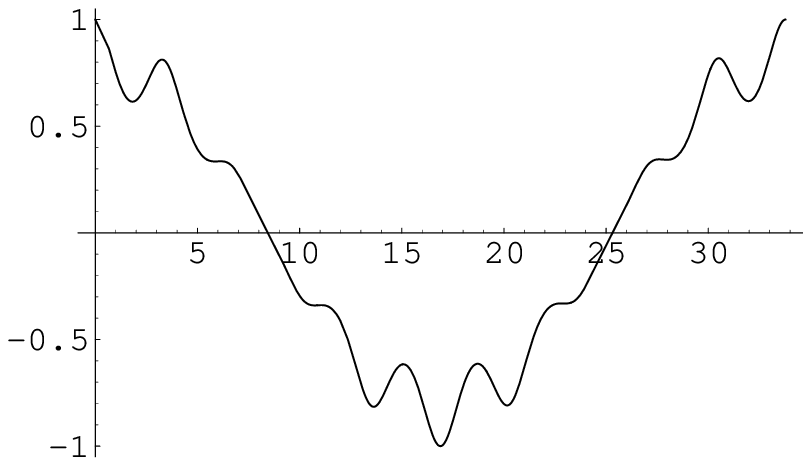}
  \includegraphics[scale=0.3]{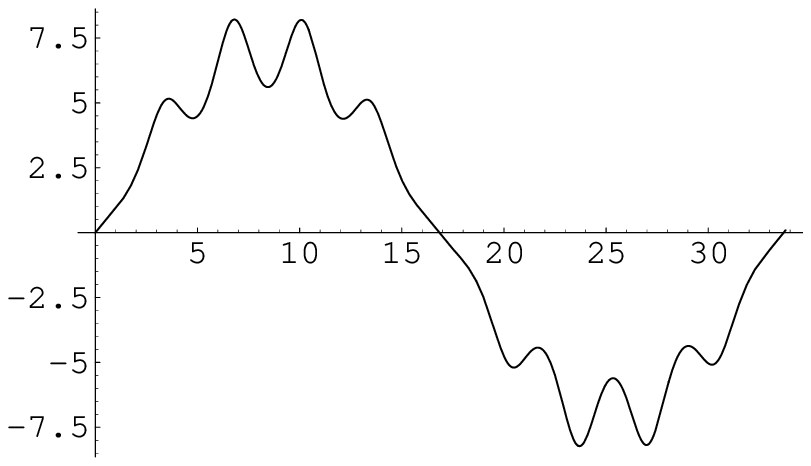}
  \includegraphics[scale=0.3]{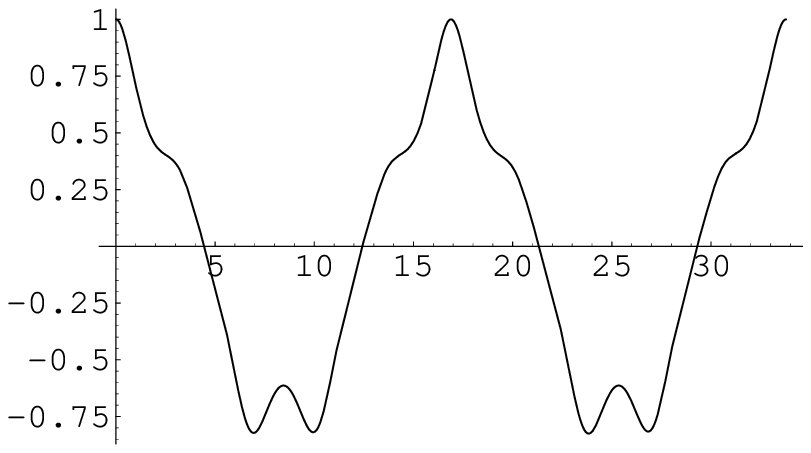}
   \includegraphics[scale=0.3]{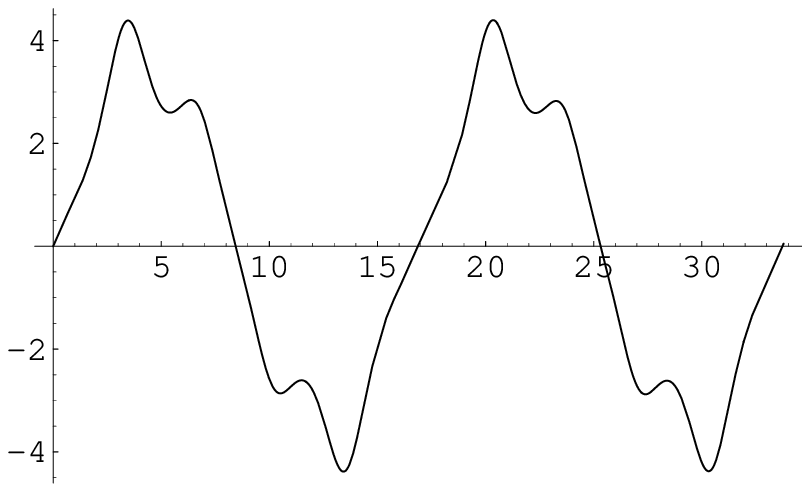}
  \includegraphics[scale=0.3]{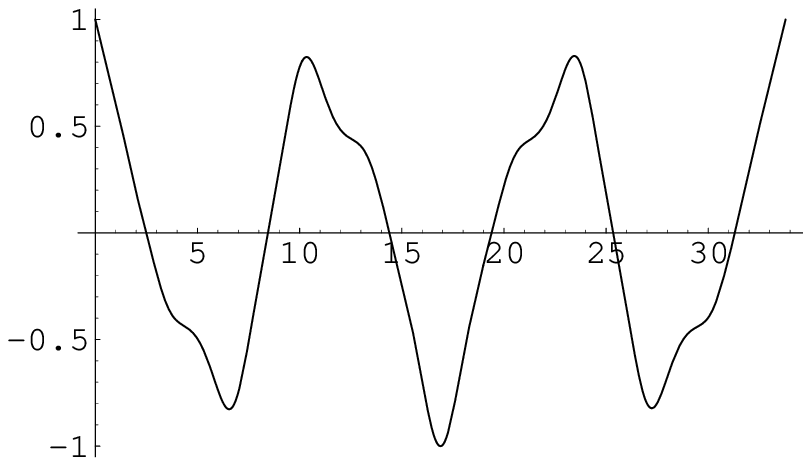}
  \includegraphics[scale=0.3]{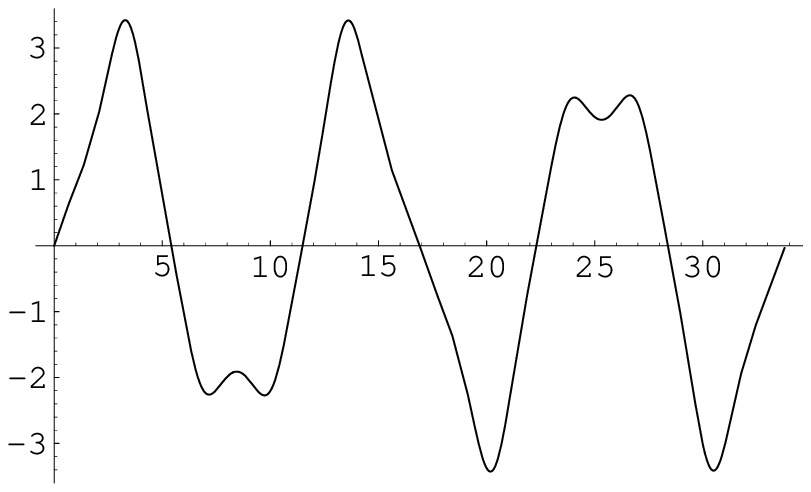}
  \includegraphics[scale=0.3]{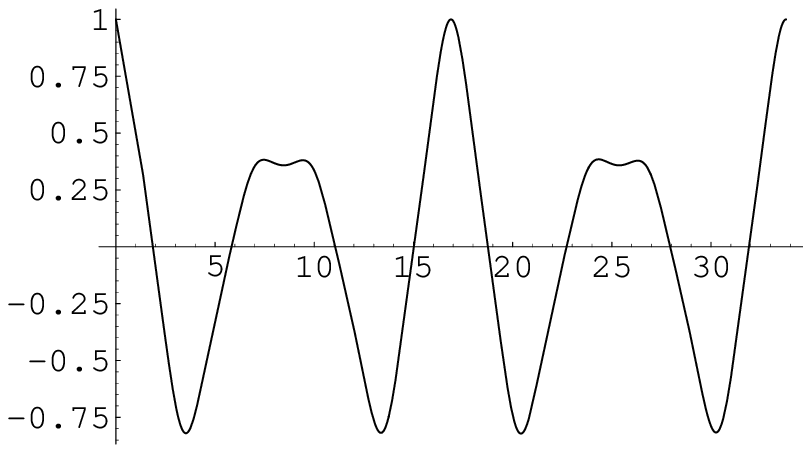}
     \includegraphics[scale=0.3]{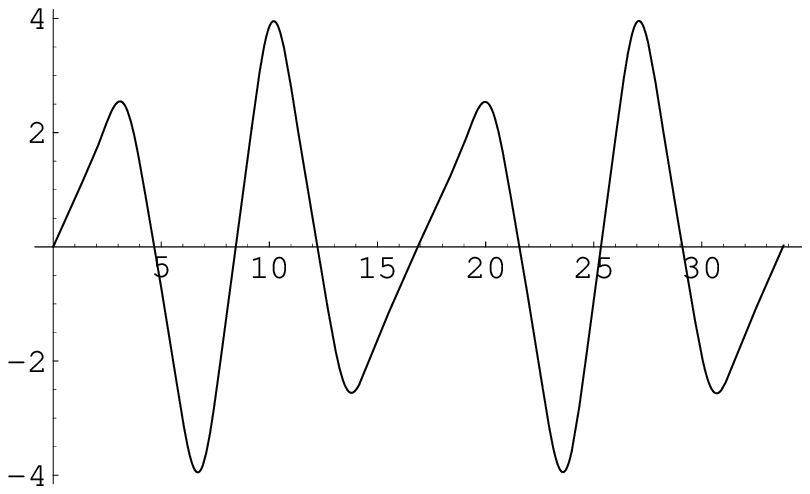}
   \includegraphics[scale=0.3]{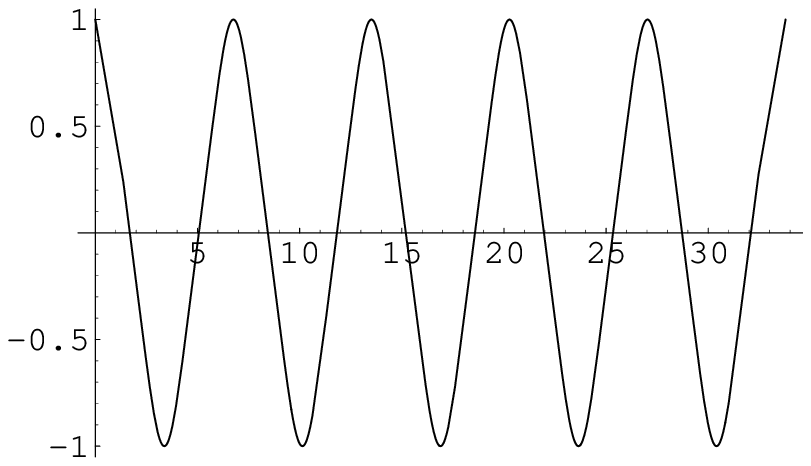}
  \includegraphics[scale=0.3]{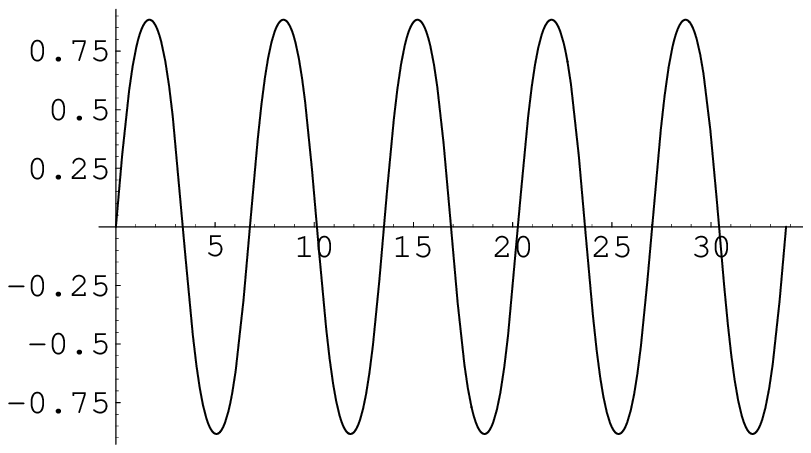}
     \caption{Eigenfunctions associated to the eigenvalues $\lambda_{1,0}$, ..., 
           $\lambda_{10,0}$, $\lambda_{11,0}=0$ of the surface $U_{10}$.}
  \label{J}
\end{figure}

\begin{figure}[phbt]
  \centering
  \includegraphics[scale=0.4]{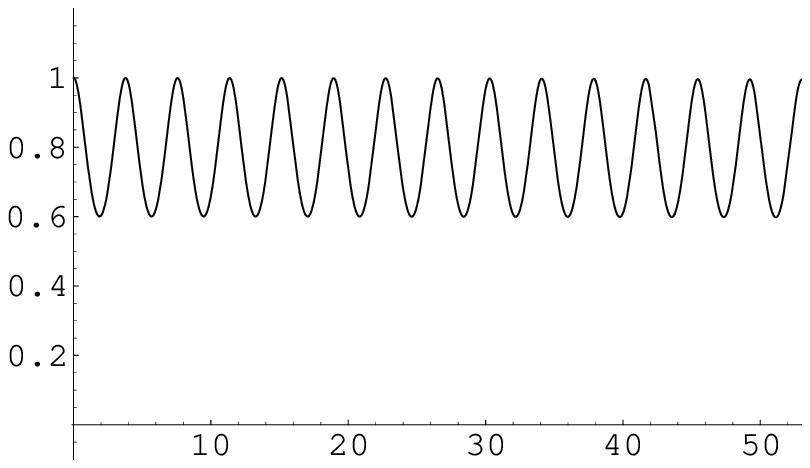}
  \includegraphics[scale=0.7]{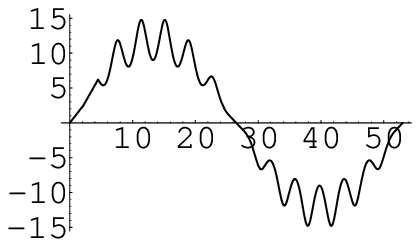}
  \includegraphics[scale=0.7]{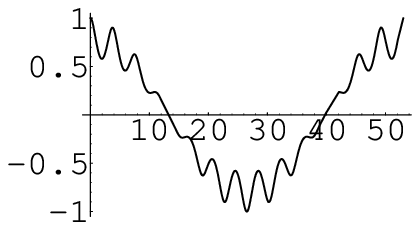}
  \includegraphics[scale=0.7]{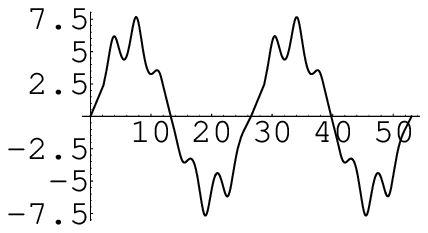}
      \includegraphics[scale=0.75]{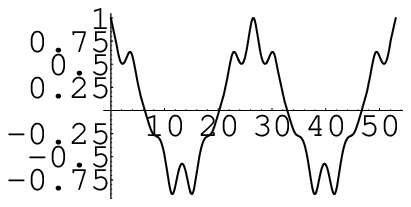}
  \includegraphics[scale=0.75]{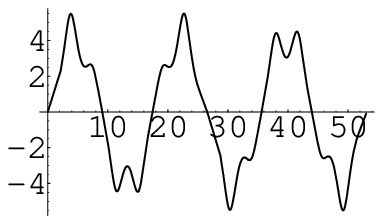}
  \includegraphics[scale=0.75]{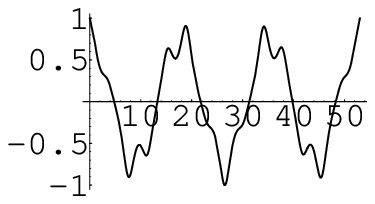}
  \includegraphics[scale=0.75]{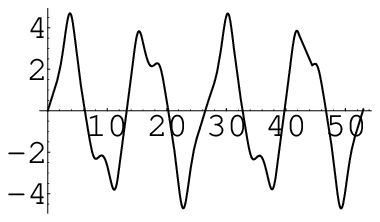}
  \includegraphics[scale=0.75]{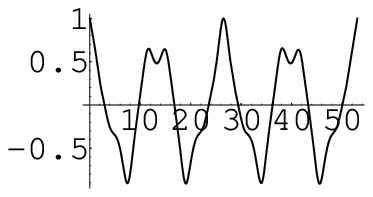}
   \includegraphics[scale=0.75]{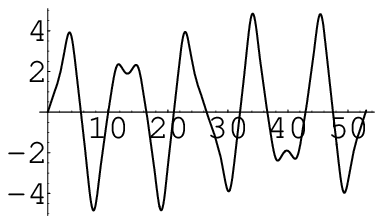}
  \includegraphics[scale=0.75]{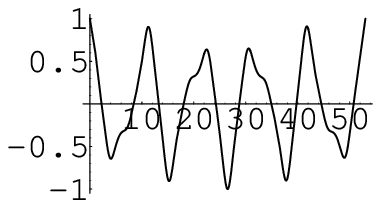}
   \includegraphics[scale=0.75]{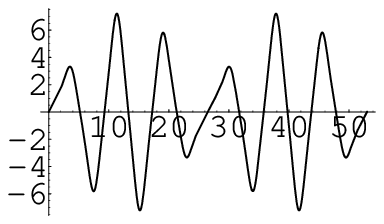}
       \includegraphics[scale=0.8]{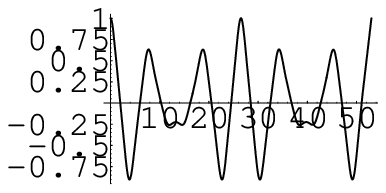}
     \includegraphics[scale=0.8]{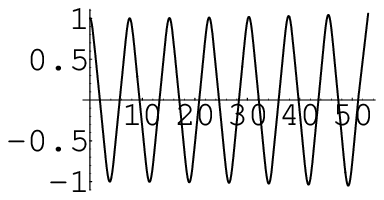}
      \includegraphics[scale=0.8]{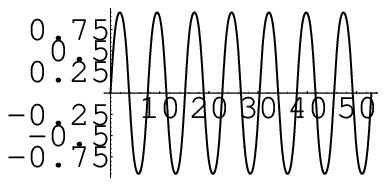}
 \caption{Eigenfunctions associated to the eigenvalues $\lambda_{1,0}$, ..., 
           $\lambda_{14,0}$, $\lambda_{15,0}=0$ of the surface $U_{11}$.}
  \label{K}
\end{figure}

\begin{figure}[phbt]
  \centering
  \includegraphics[scale=0.35]{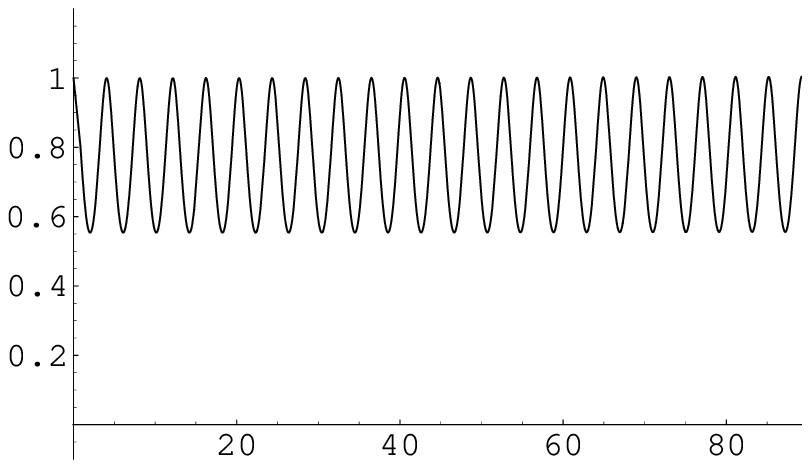}
  \includegraphics[scale=0.35]{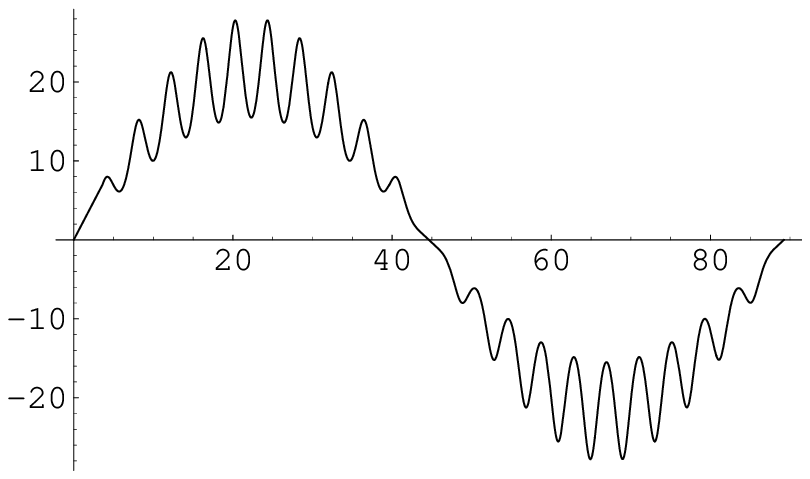}
  \includegraphics[scale=0.35]{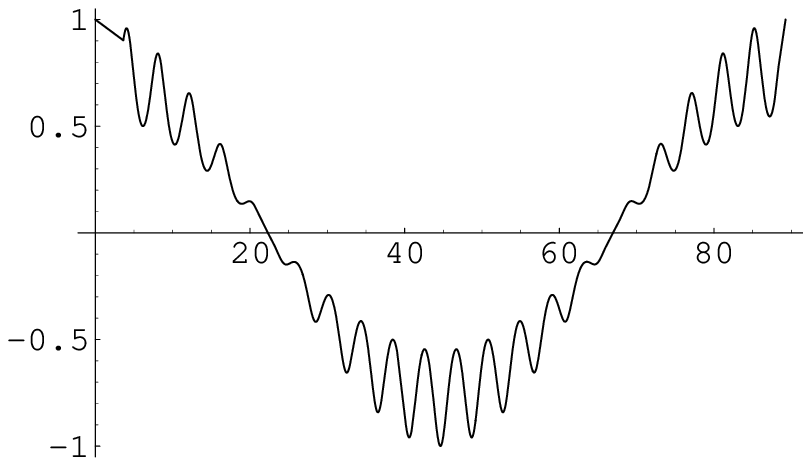}
  \includegraphics[scale=0.35]{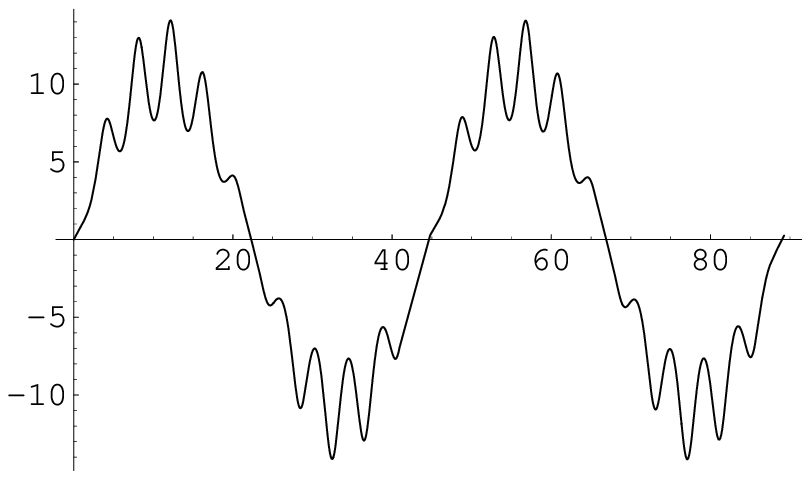}
      \includegraphics[scale=0.35]{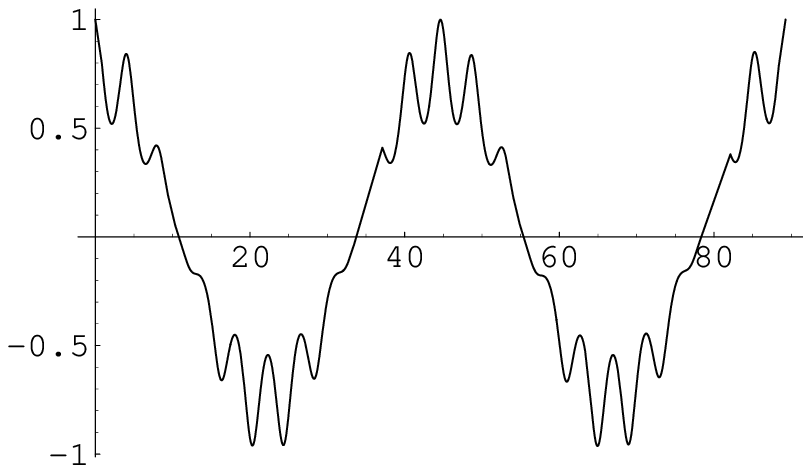}
  \includegraphics[scale=0.35]{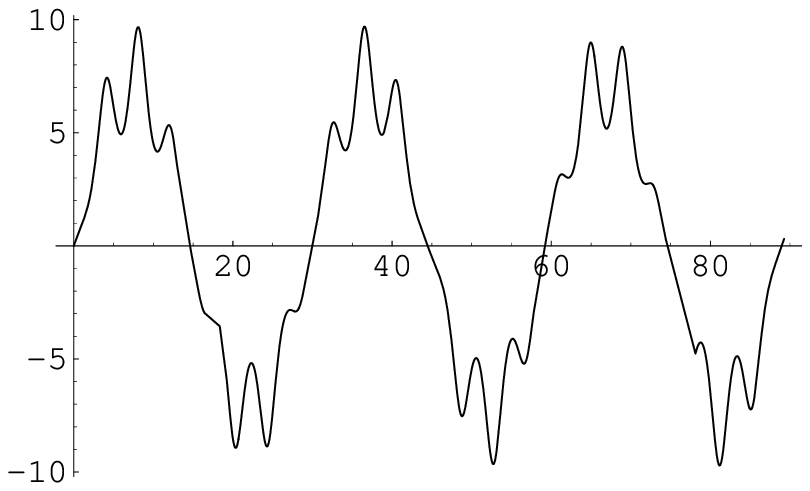}
  \includegraphics[scale=0.35]{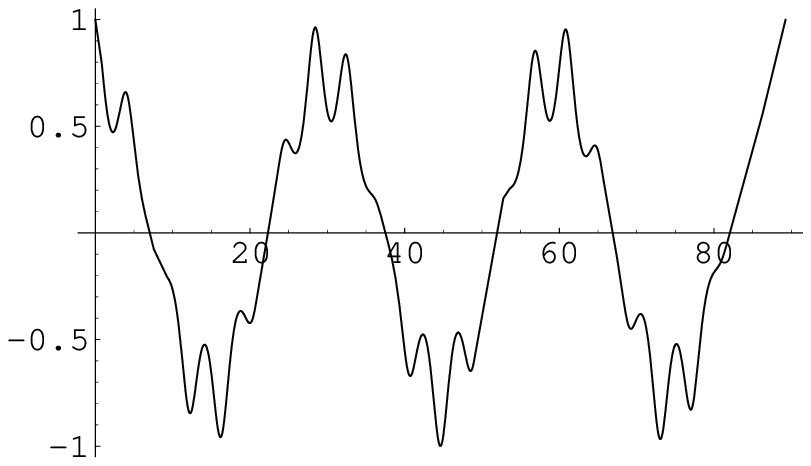}
  \includegraphics[scale=0.35]{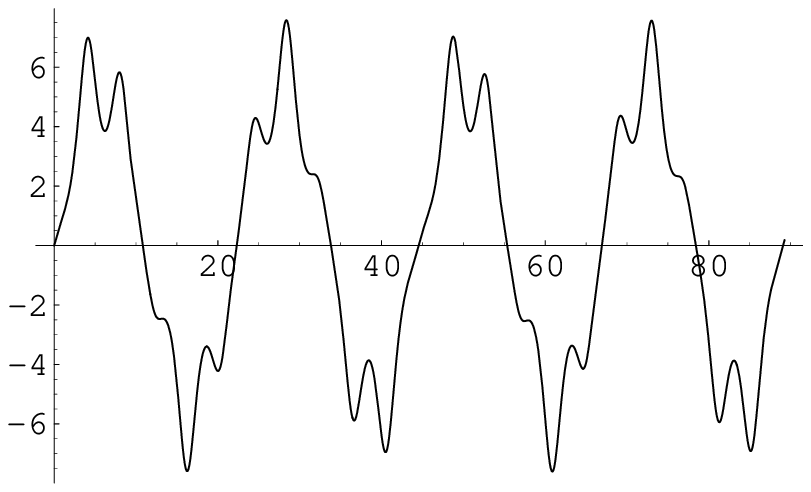}
     \includegraphics[scale=0.35]{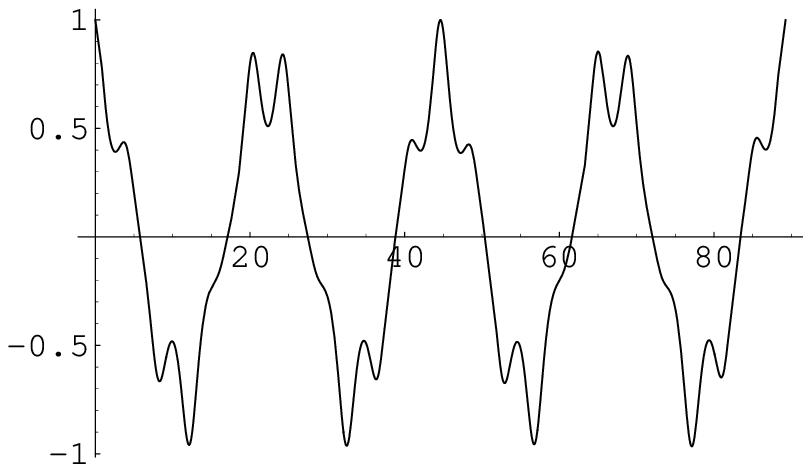}
  \includegraphics[scale=0.35]{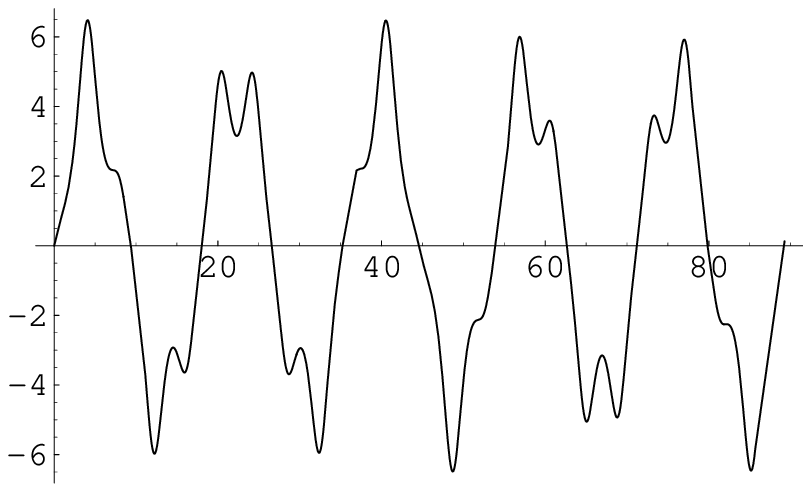}
  \includegraphics[scale=0.35]{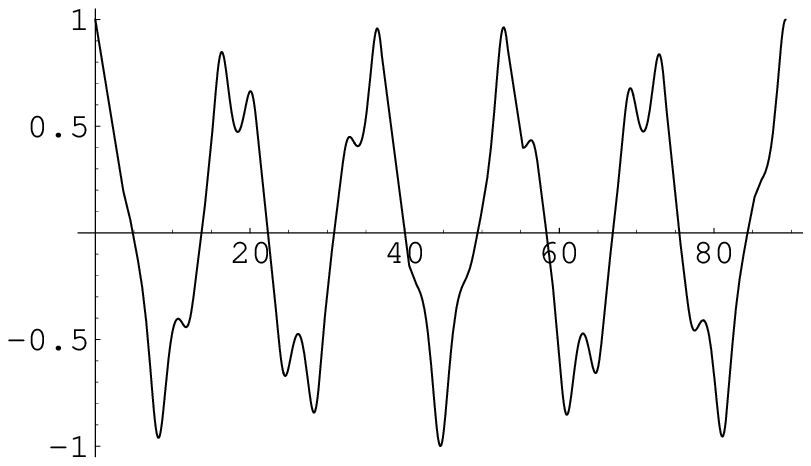}
  \includegraphics[scale=0.35]{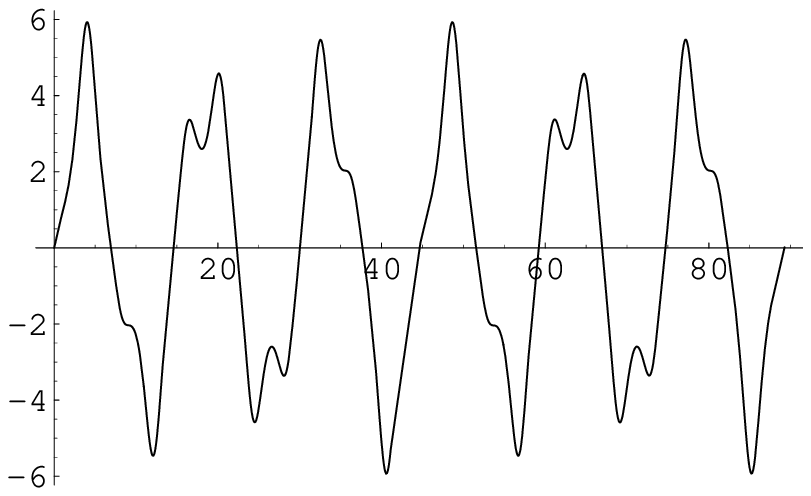}
     \includegraphics[scale=0.35]{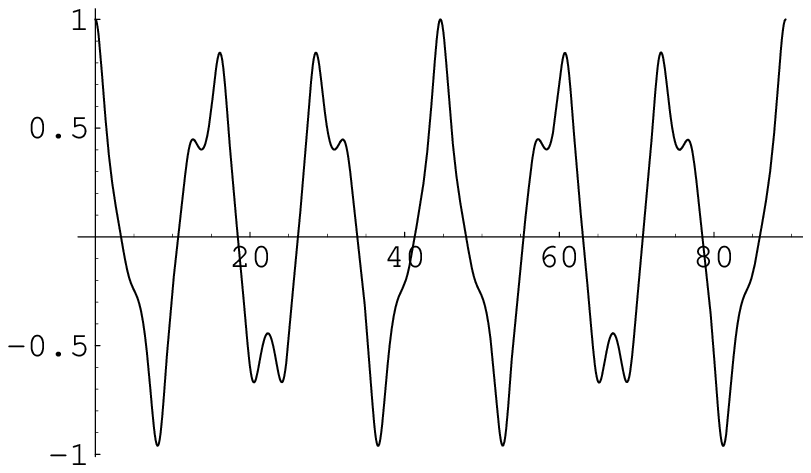}
  \includegraphics[scale=0.35]{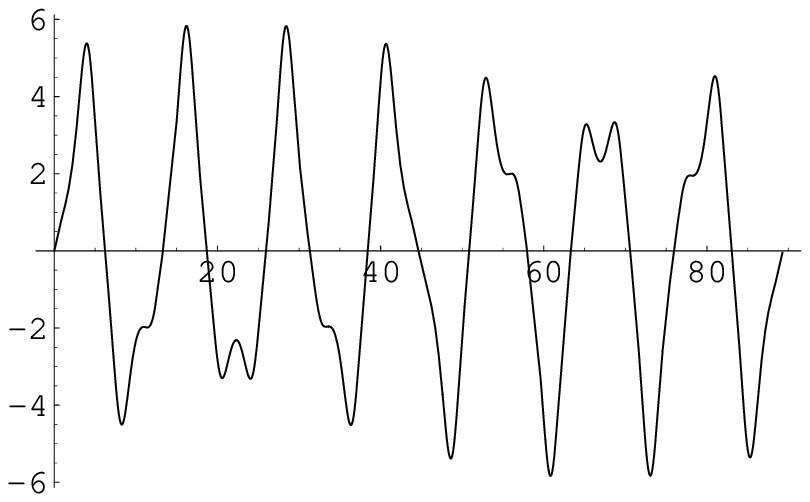}
  \includegraphics[scale=0.35]{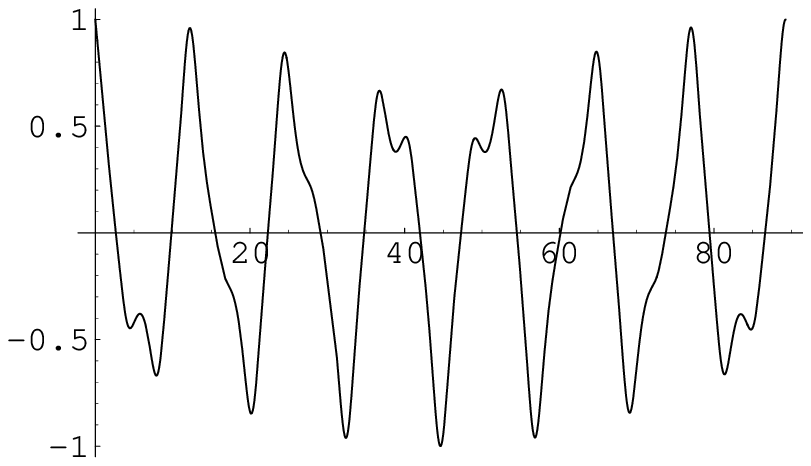}
  \includegraphics[scale=0.35]{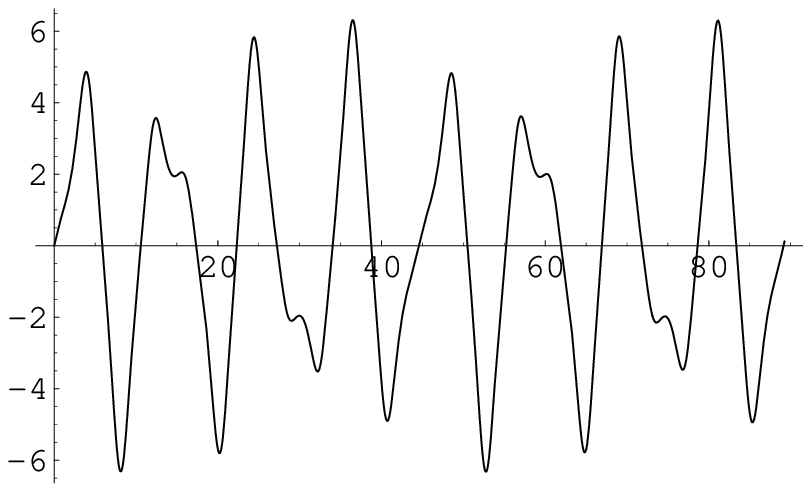}
  \includegraphics[scale=0.35]{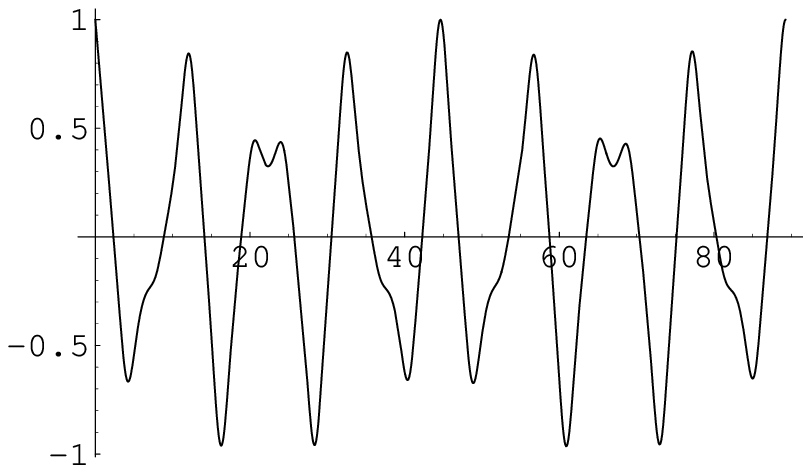}
  \includegraphics[scale=0.35]{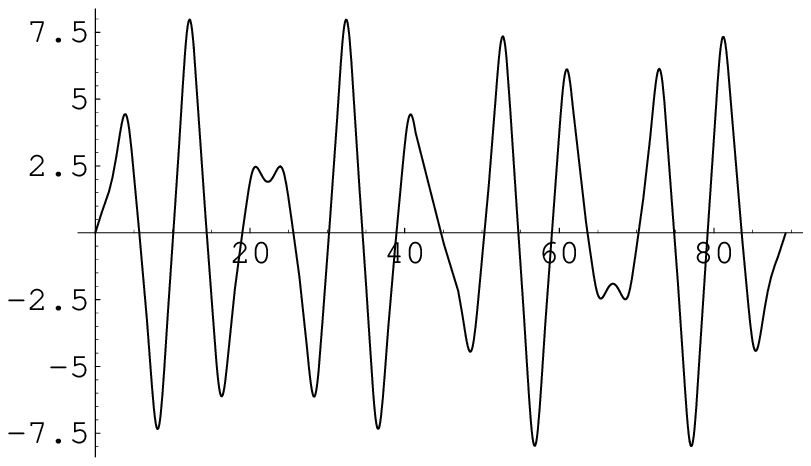}
  \includegraphics[scale=0.35]{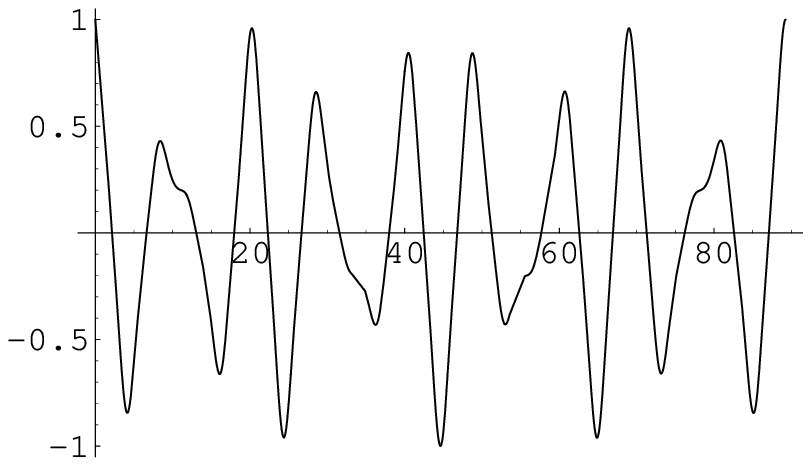}
   \includegraphics[scale=0.35]{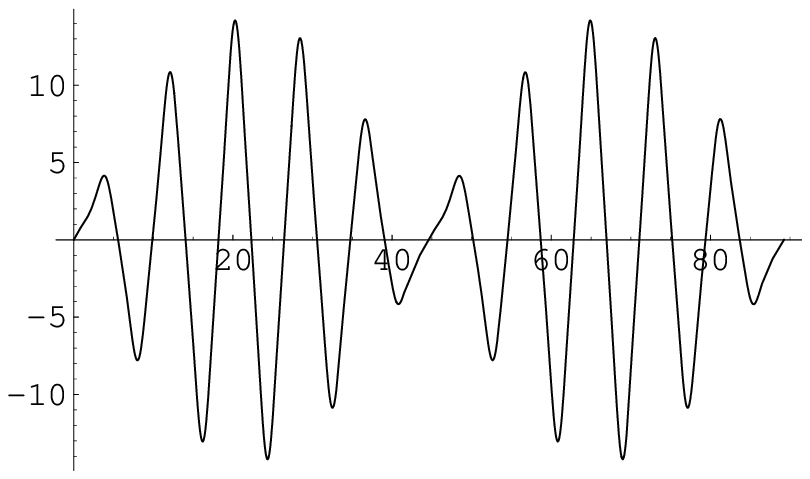}
      \includegraphics[scale=0.35]{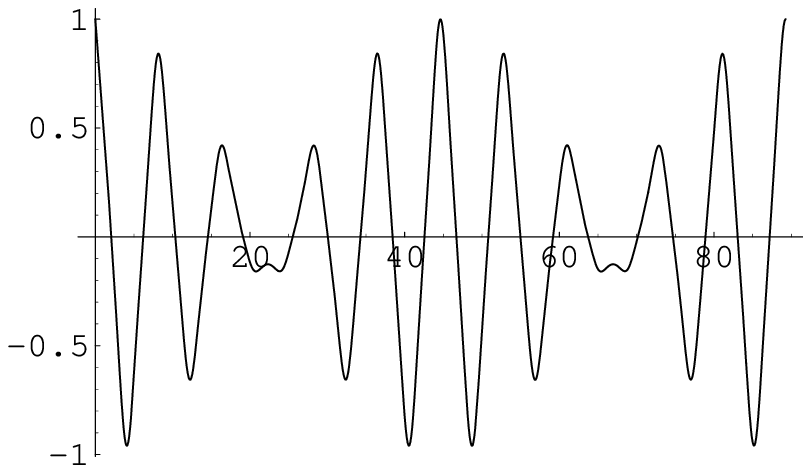}
    \includegraphics[scale=0.35]{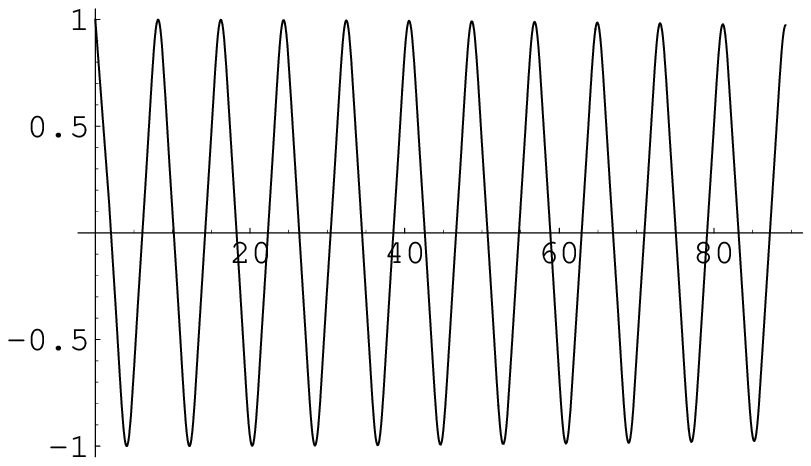}
     \includegraphics[scale=0.35]{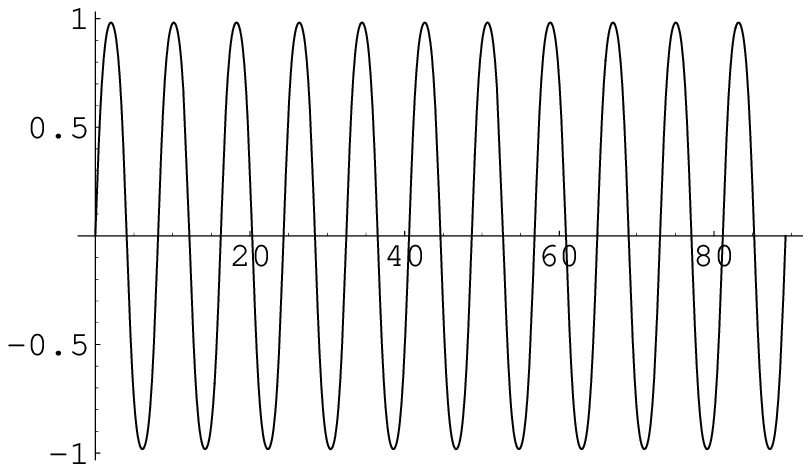}
 \caption{Eigenfunctions associated to the eigenvalues $\lambda_{1,0}$, ..., 
           $\lambda_{22,0}$, $\lambda_{23,0}=0$ of the surface $U_{12}$.}
  \label{L}
\end{figure}

\begin{figure}[phbt]
  \centering
  \includegraphics[scale=0.35]{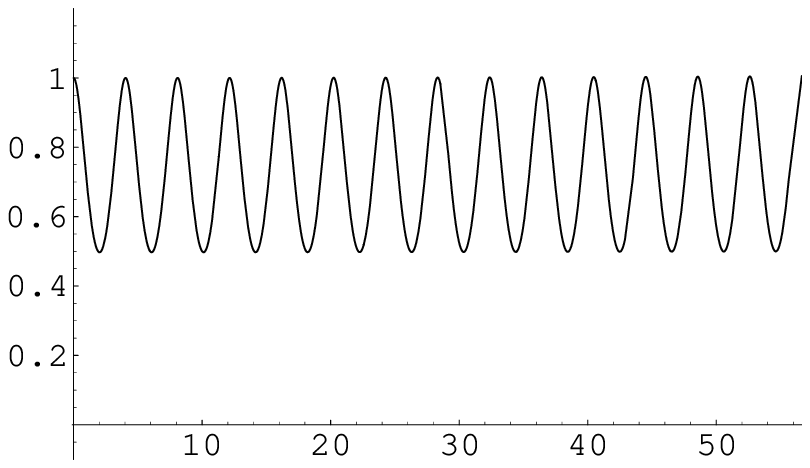}
  \includegraphics[scale=0.35]{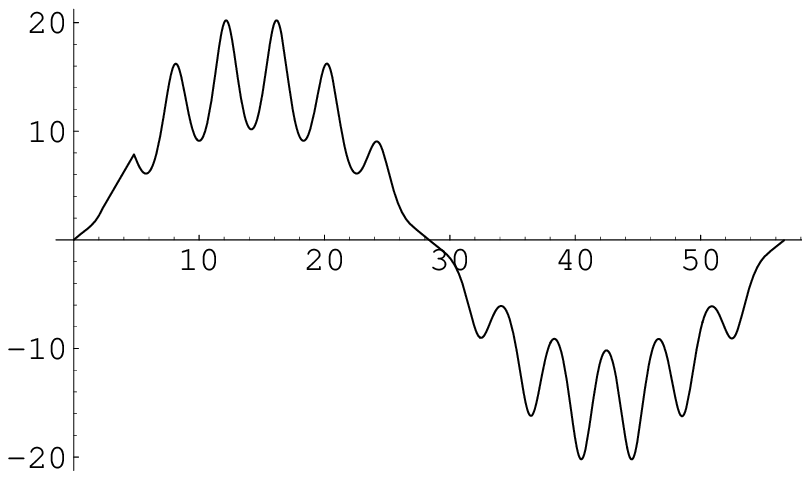}
  \includegraphics[scale=0.35]{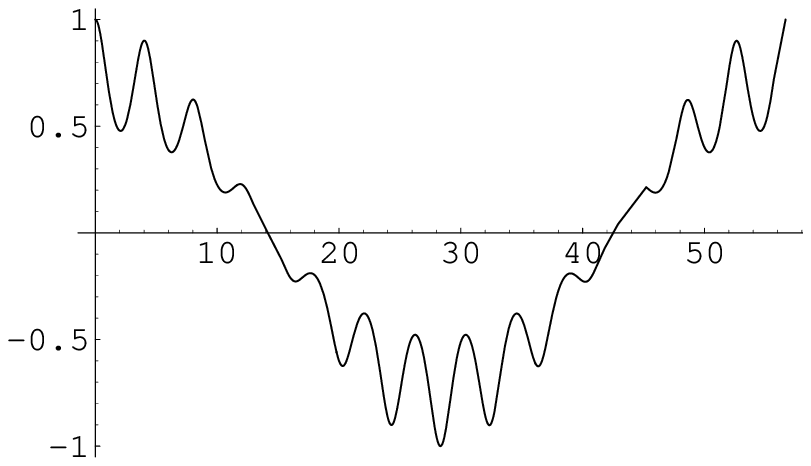}
  \includegraphics[scale=0.35]{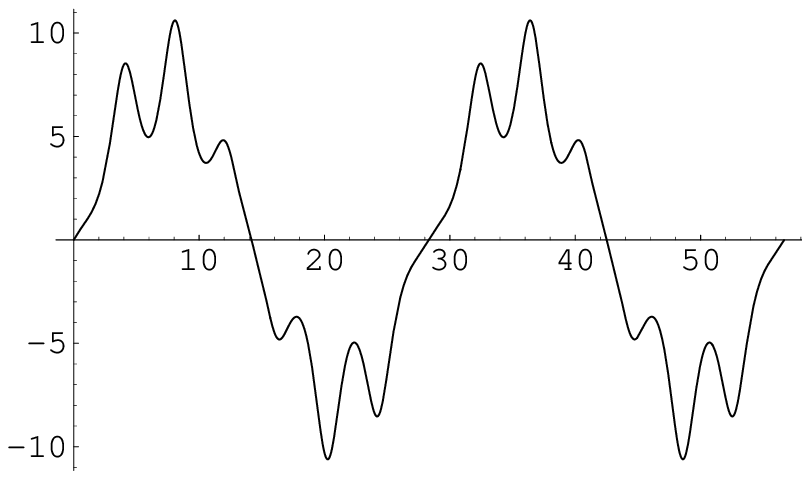}
  \includegraphics[scale=0.35]{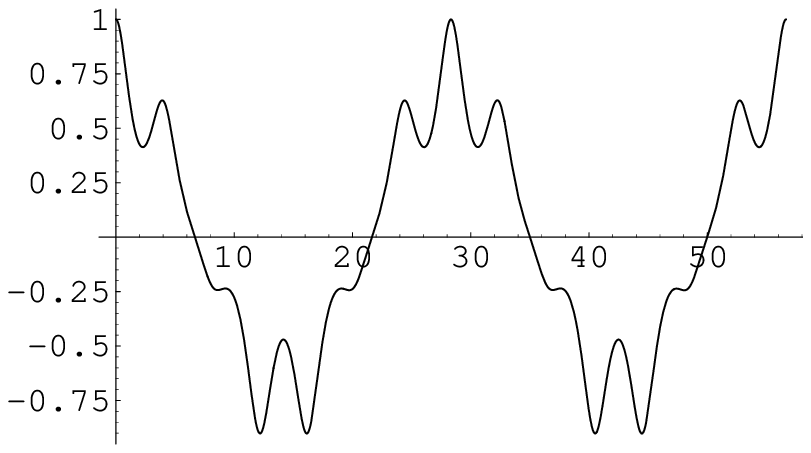}
  \includegraphics[scale=0.35]{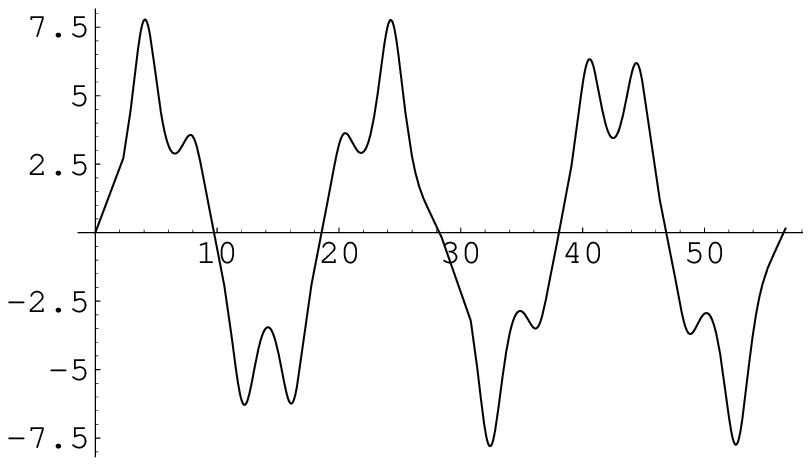}
  \includegraphics[scale=0.35]{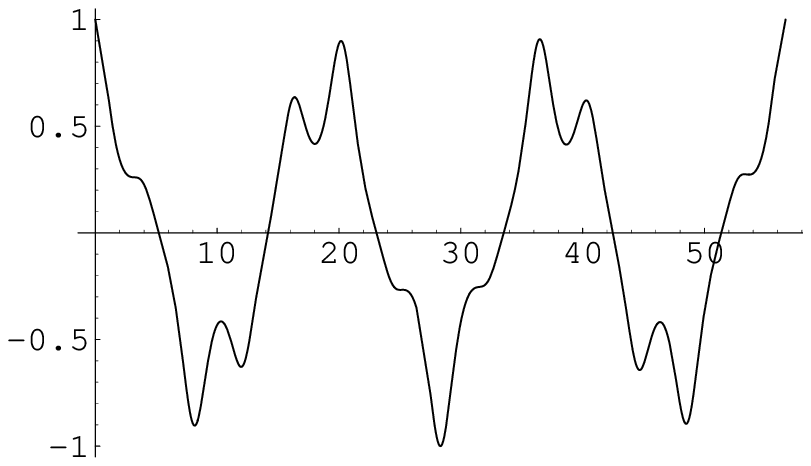}
  \includegraphics[scale=0.35]{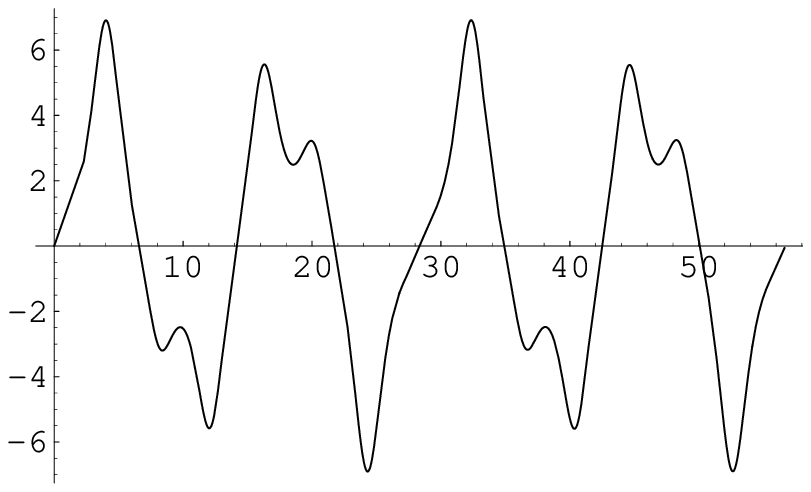}
   \includegraphics[scale=0.35]{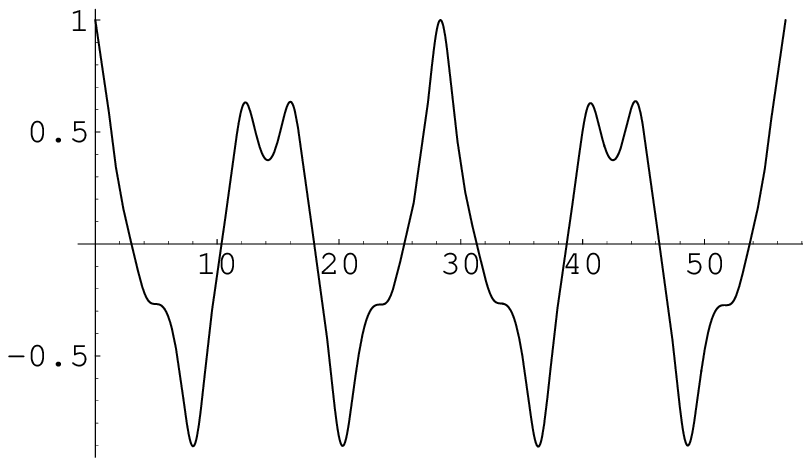}
   \includegraphics[scale=0.35]{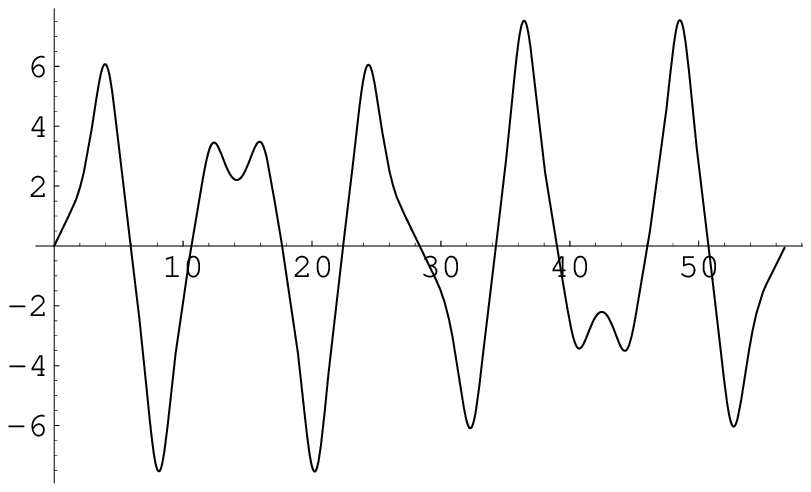}
  \includegraphics[scale=0.35]{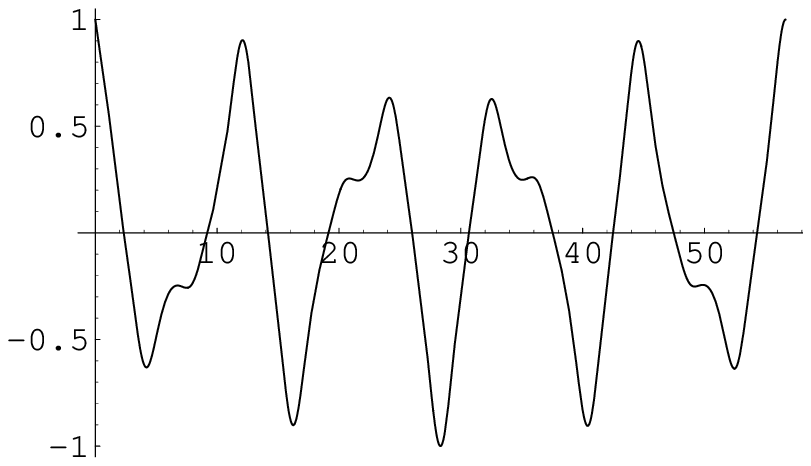}
  \includegraphics[scale=0.35]{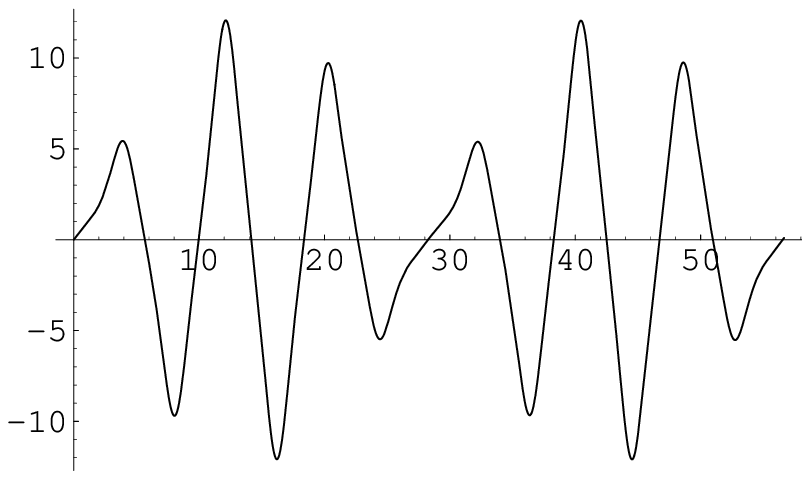}  
   \includegraphics[scale=0.35]{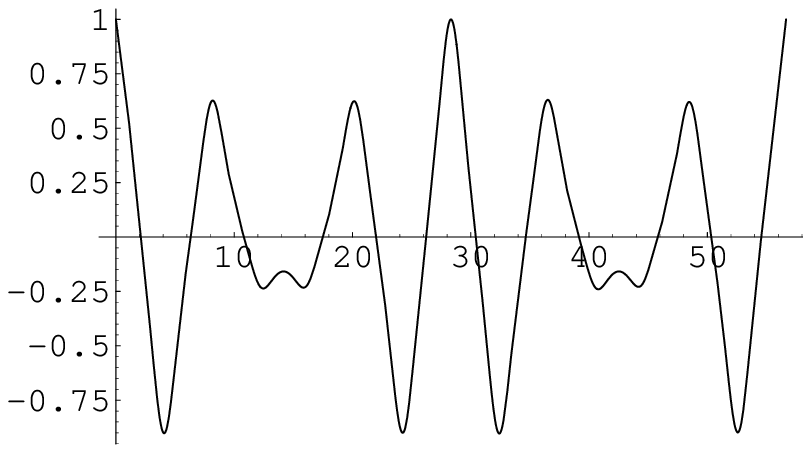}
  \includegraphics[scale=0.35]{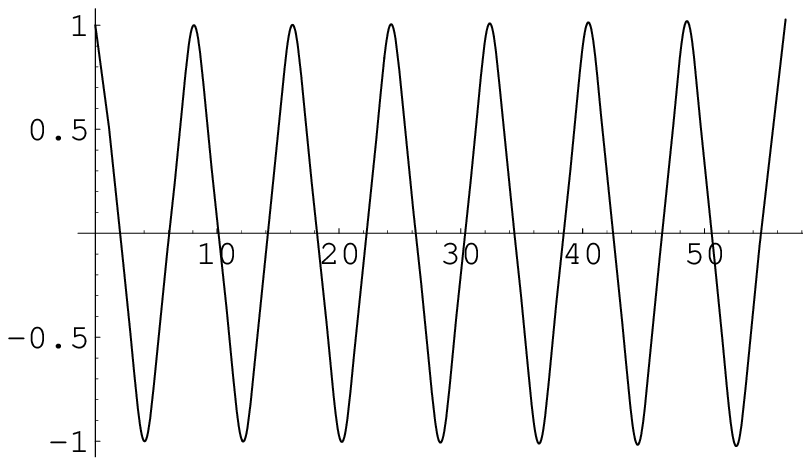}
  \includegraphics[scale=0.35]{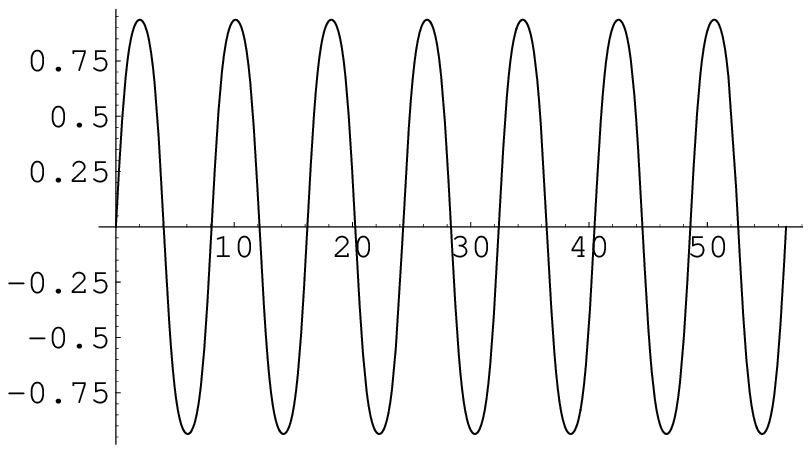}
     \caption{Eigenfunctions associated to the eigenvalues $\lambda_{1,0}$, ..., 
           $\lambda_{14,0}$, $\lambda_{15,0}=0$ of the surface $U_{13}$.}
  \label{M}
\end{figure}

\begin{figure}[phbt]
  \centering
  \includegraphics[scale=0.35]{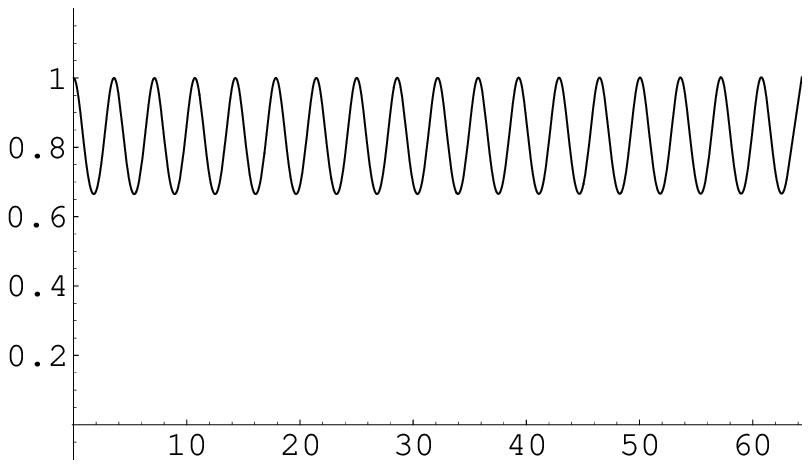}
  \includegraphics[scale=0.35]{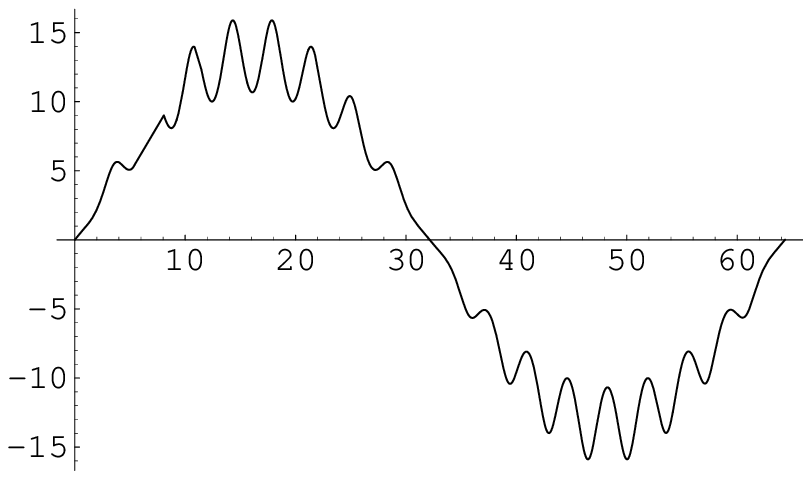}
  \includegraphics[scale=0.35]{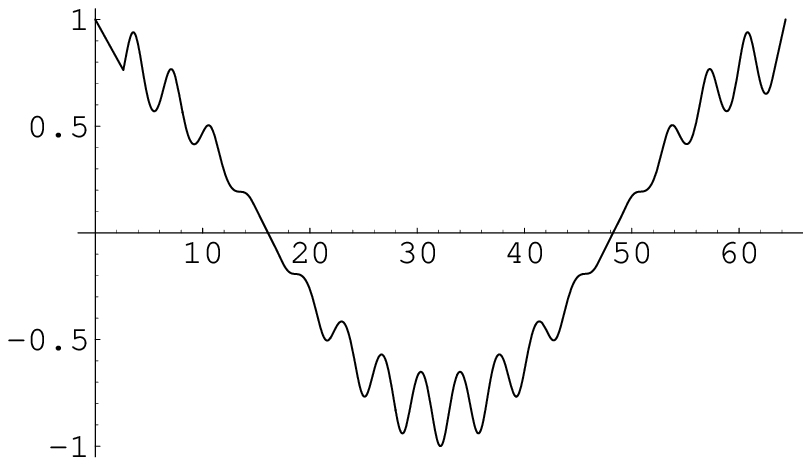}
  \includegraphics[scale=0.35]{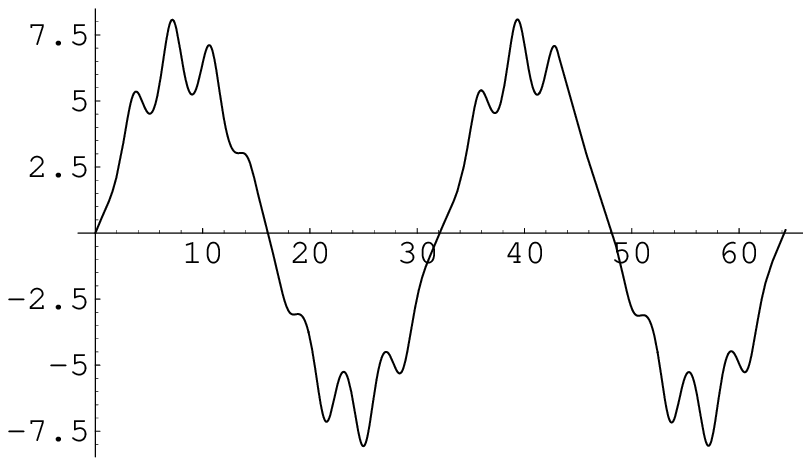}
      \includegraphics[scale=0.35]{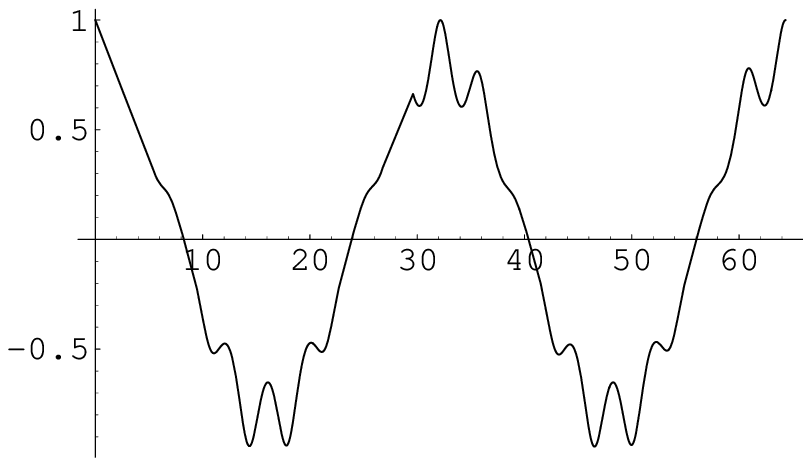}
  \includegraphics[scale=0.35]{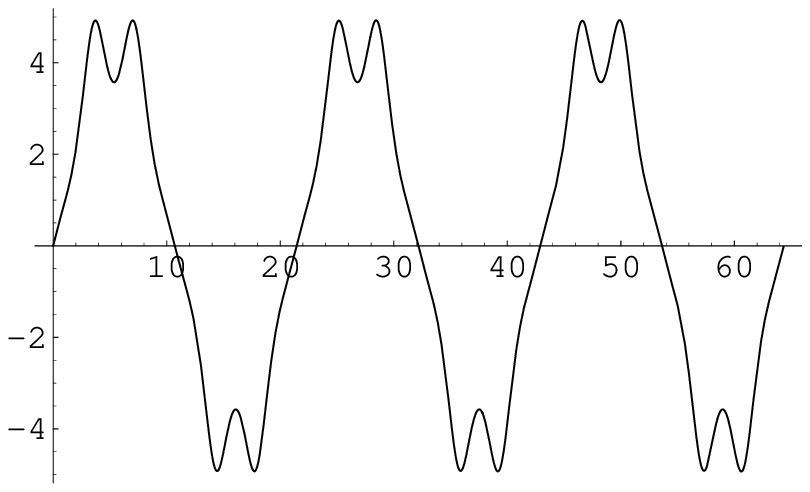}
  \includegraphics[scale=0.35]{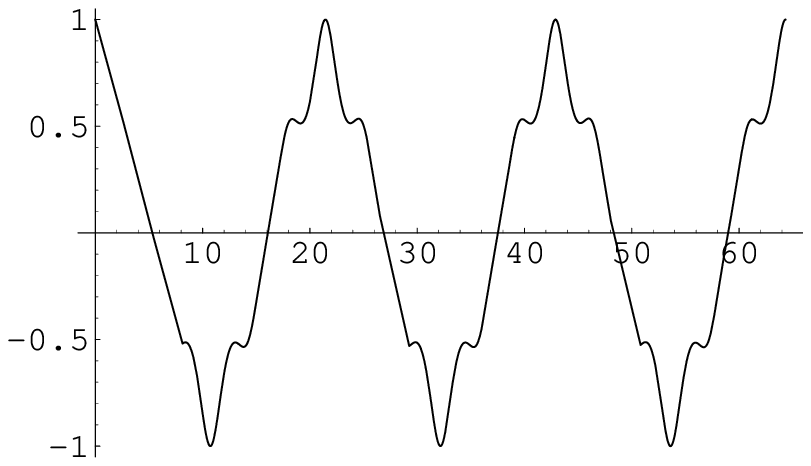}
  \includegraphics[scale=0.35]{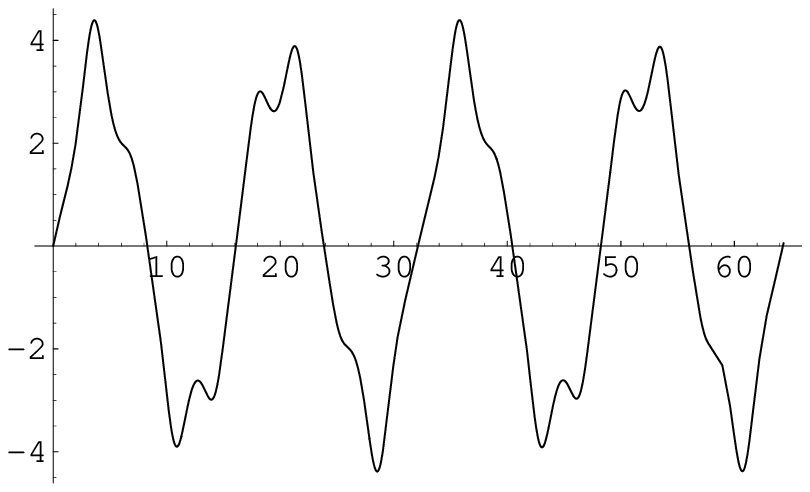}
    \includegraphics[scale=0.35]{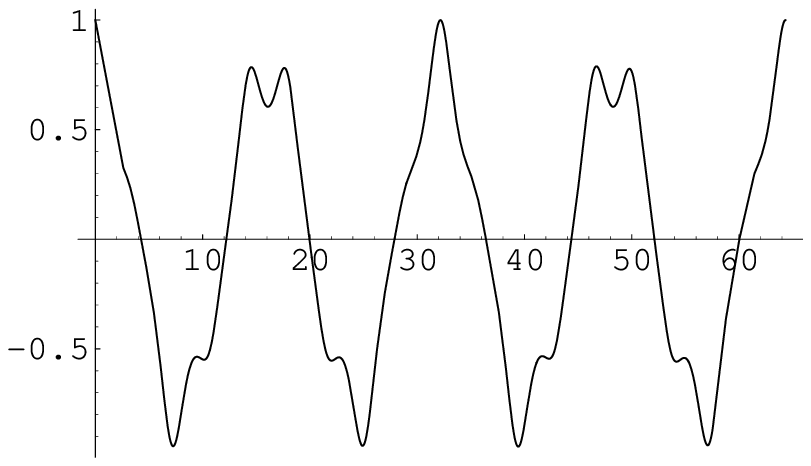}
  \includegraphics[scale=0.35]{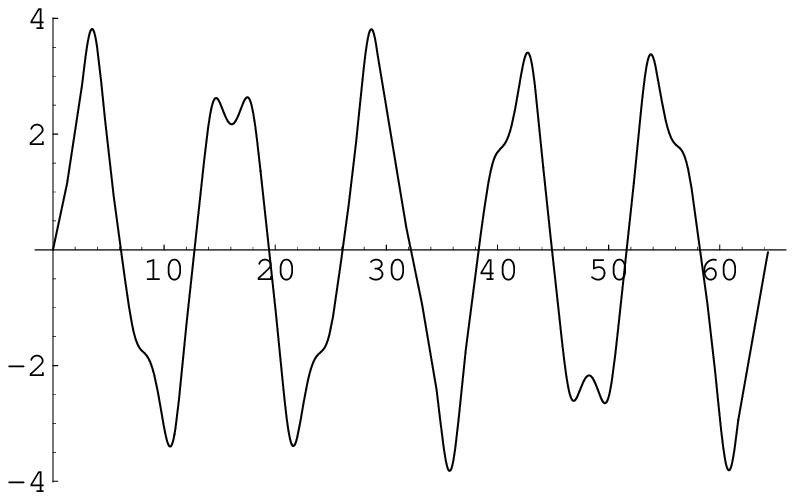}
  \includegraphics[scale=0.35]{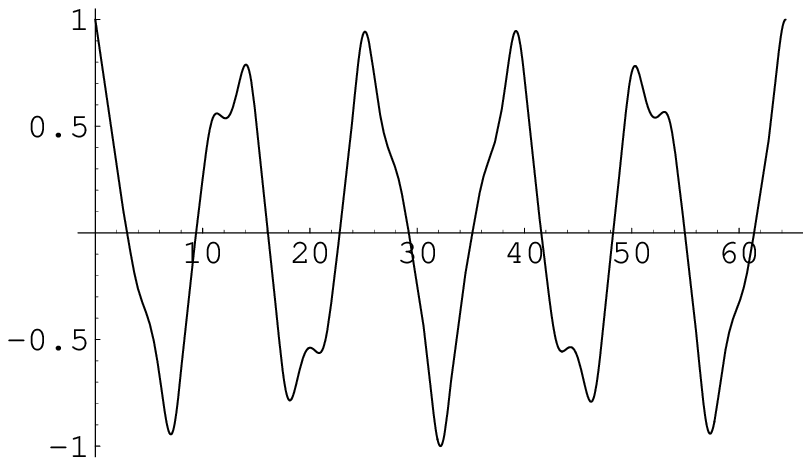}
  \includegraphics[scale=0.35]{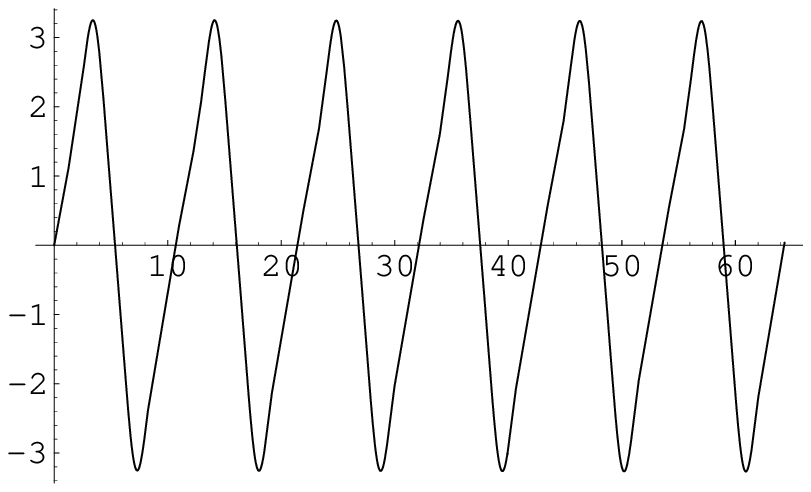} 
  \includegraphics[scale=0.35]{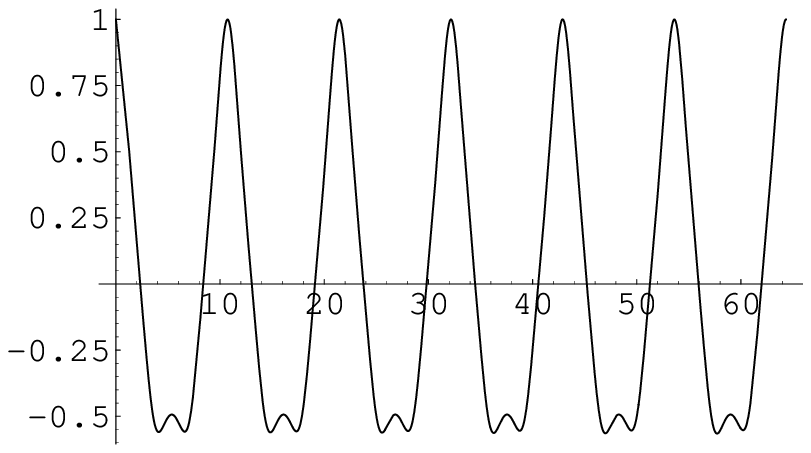}
  \includegraphics[scale=0.35]{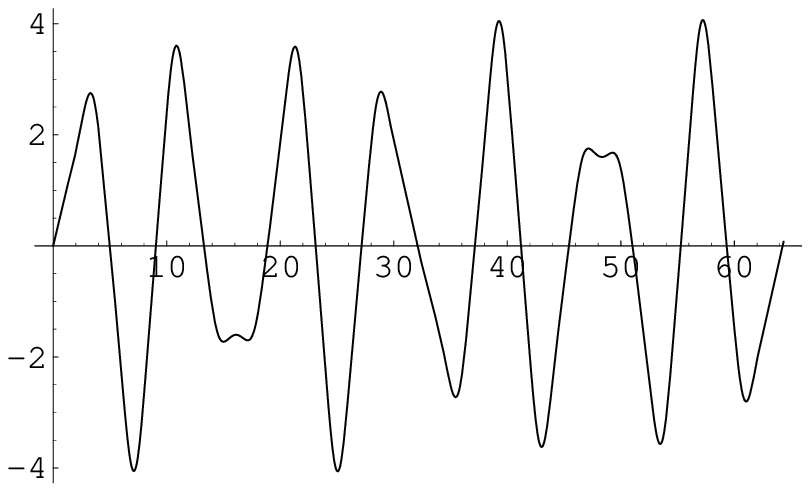}
  \includegraphics[scale=0.35]{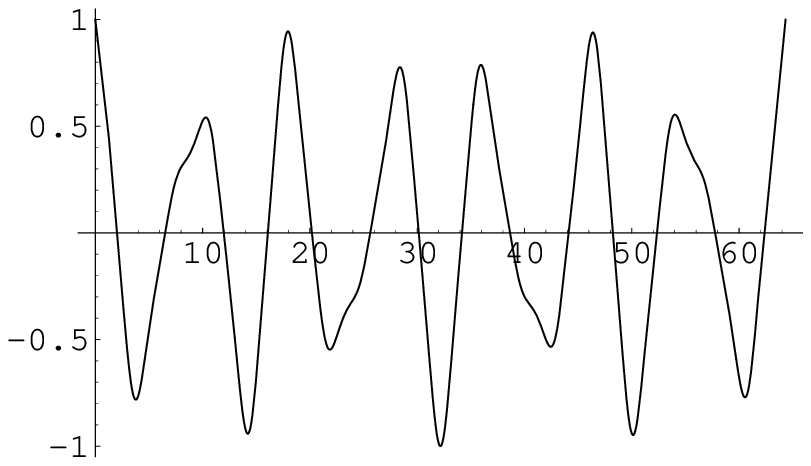}
  \includegraphics[scale=0.35]{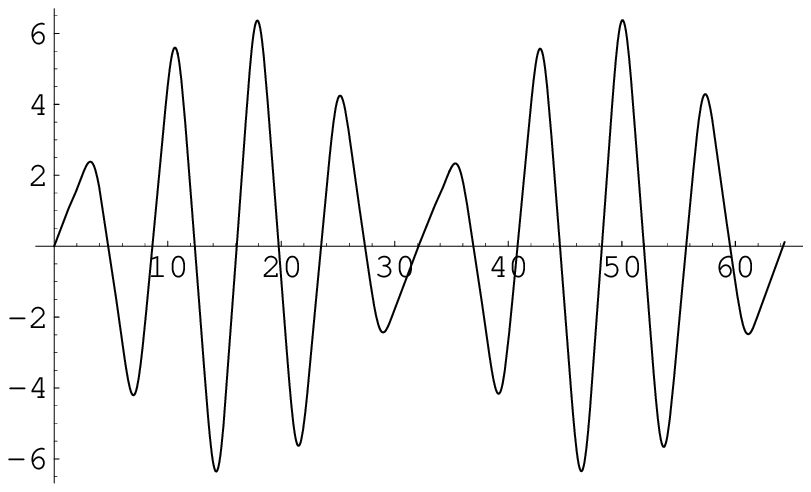}
      \includegraphics[scale=0.35]{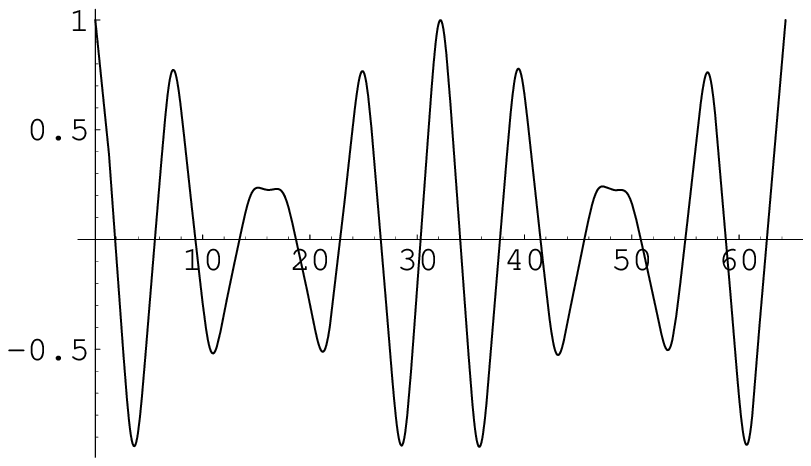}
  \includegraphics[scale=0.35]{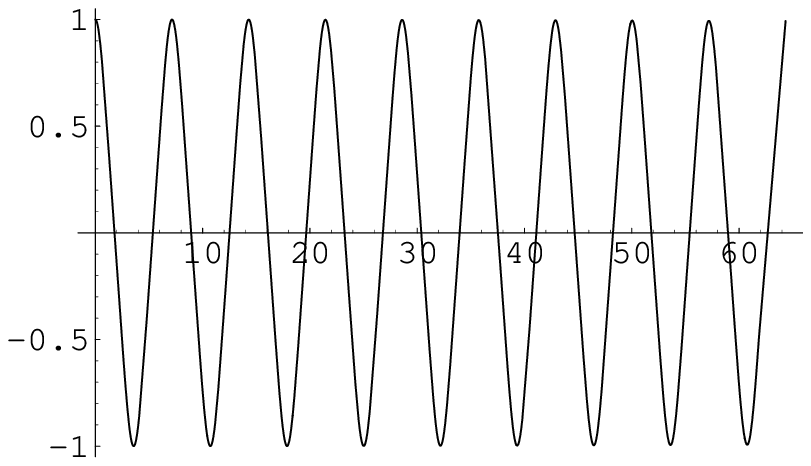}
  \includegraphics[scale=0.35]{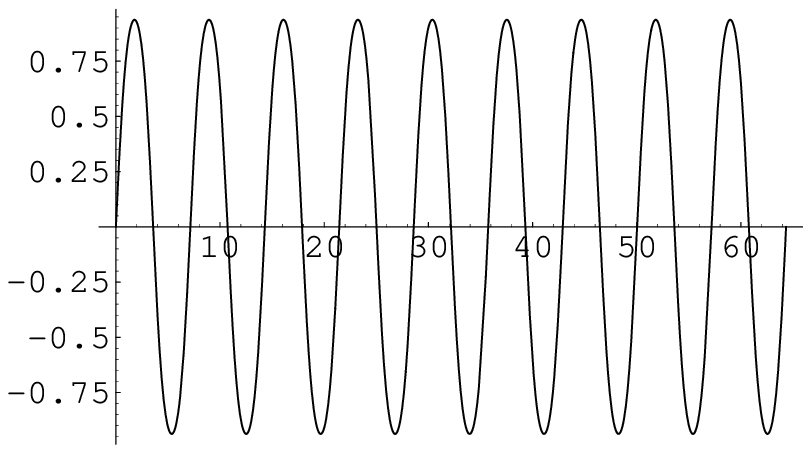}
    \caption{Eigenfunctions associated to the eigenvalues $\lambda_{1,0}$, ..., 
           $\lambda_{18,0}$, $\lambda_{19,0}=0$ of the surface $U_{14}$.}
  \label{P}
\end{figure}

\begin{figure}[phbt]
  \centering
  \includegraphics[scale=0.35]{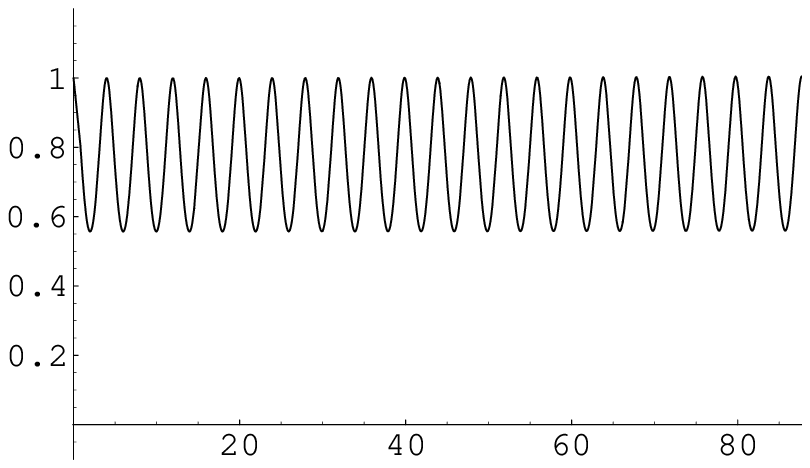}
  \includegraphics[scale=0.35]{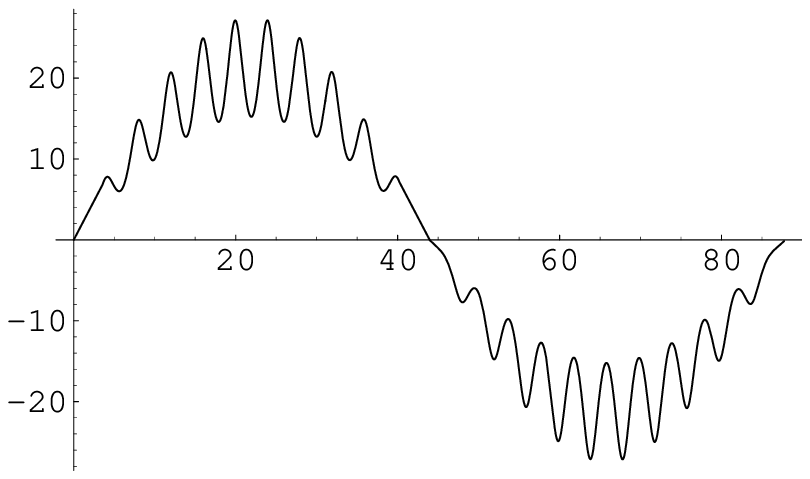}
  \includegraphics[scale=0.35]{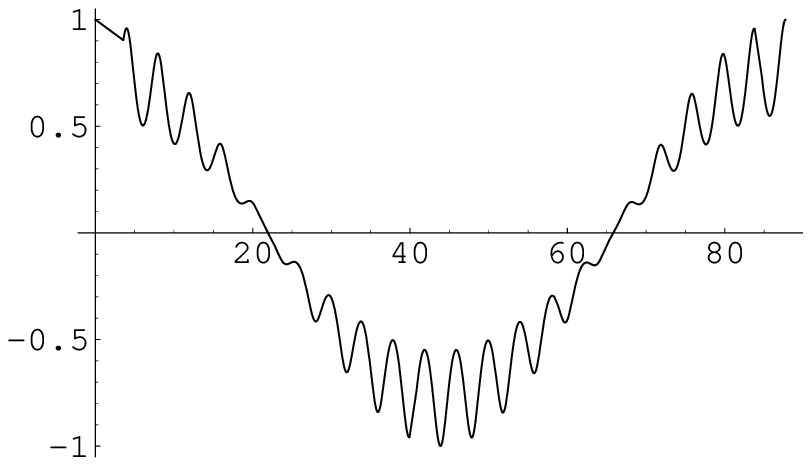}
  \includegraphics[scale=0.35]{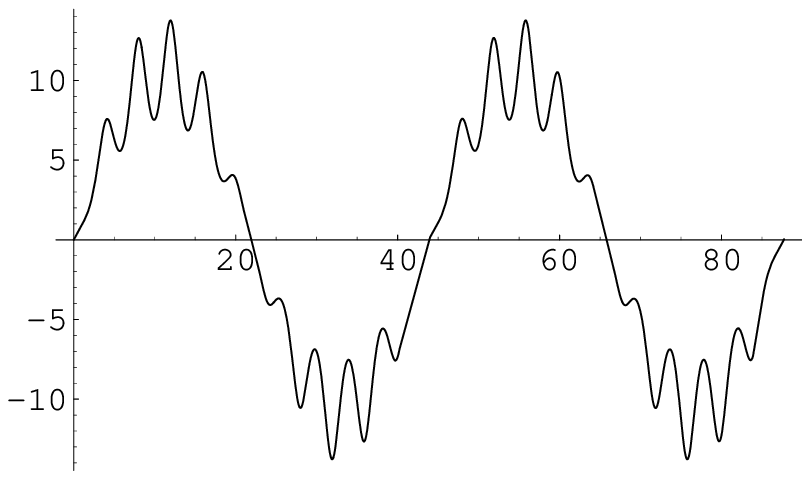}
   \includegraphics[scale=0.35]{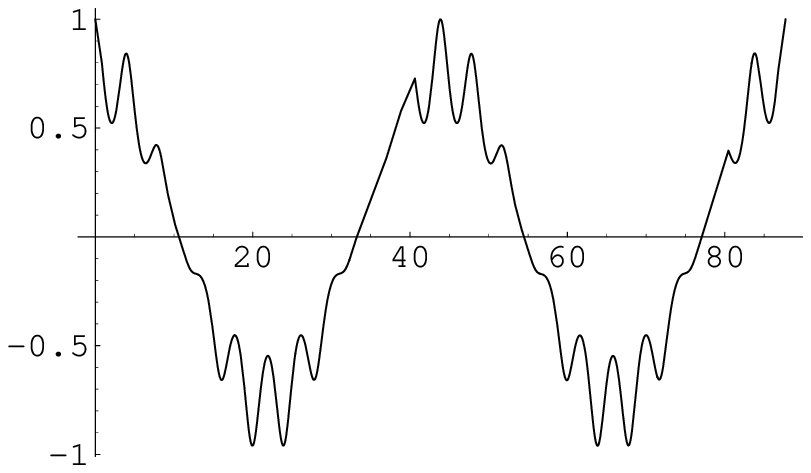}
  \includegraphics[scale=0.35]{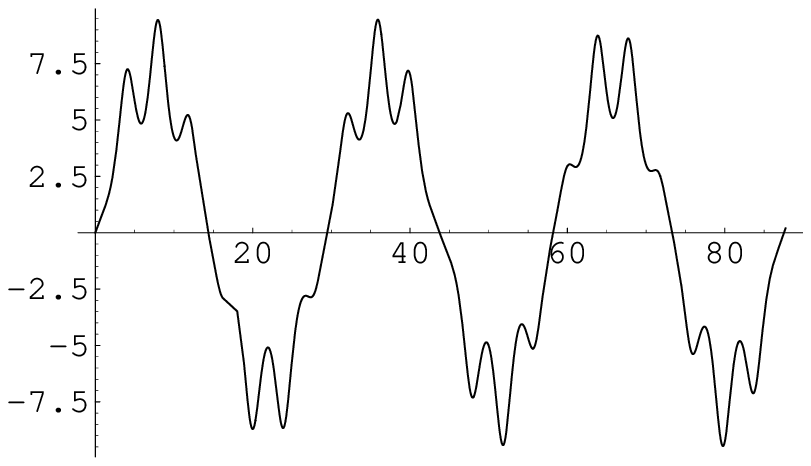}
  \includegraphics[scale=0.35]{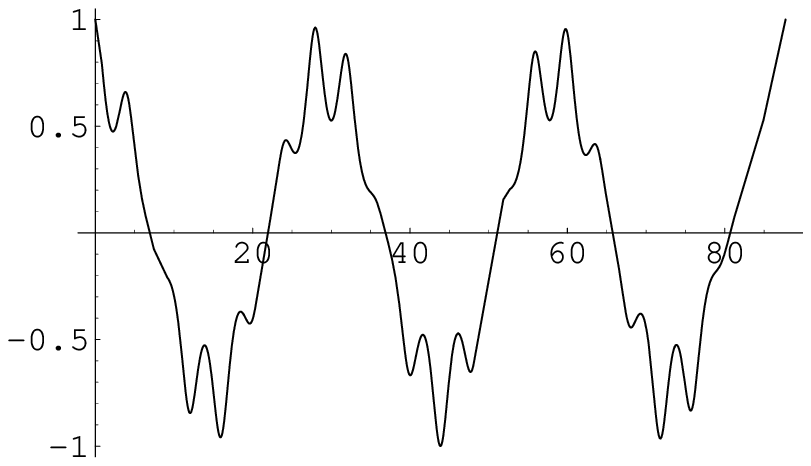}
  \includegraphics[scale=0.35]{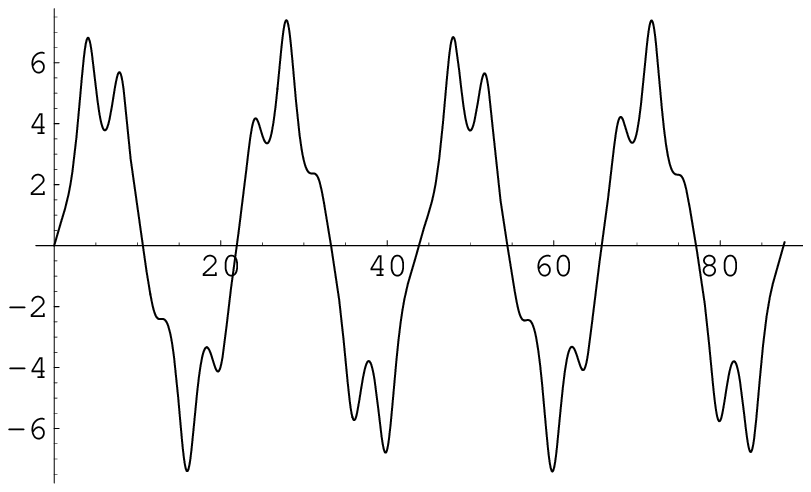}
    \includegraphics[scale=0.35]{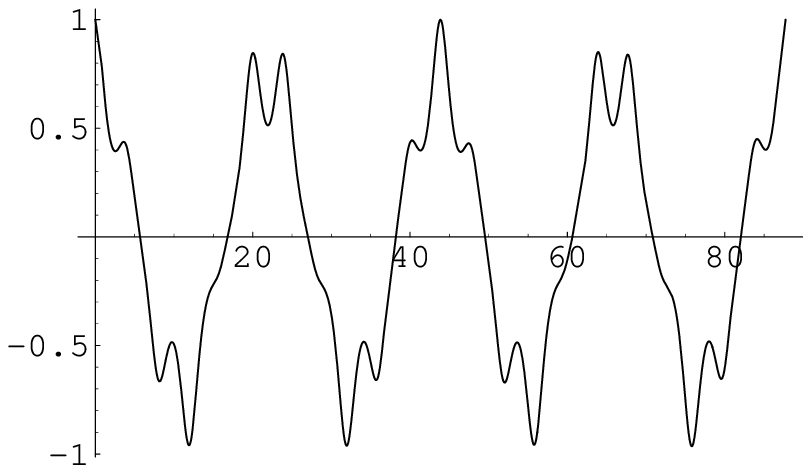}
  \includegraphics[scale=0.35]{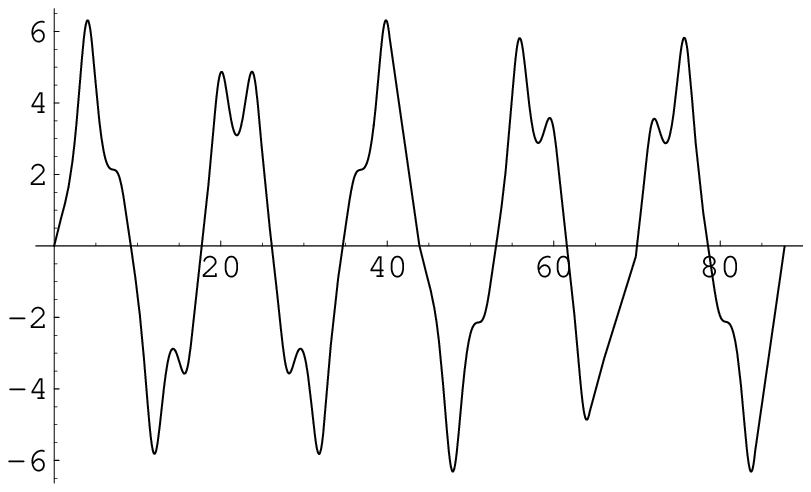}
  \includegraphics[scale=0.35]{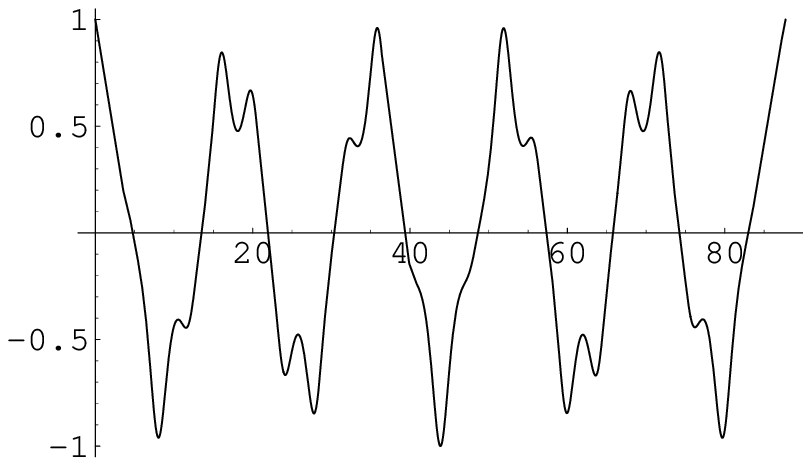}
  \includegraphics[scale=0.35]{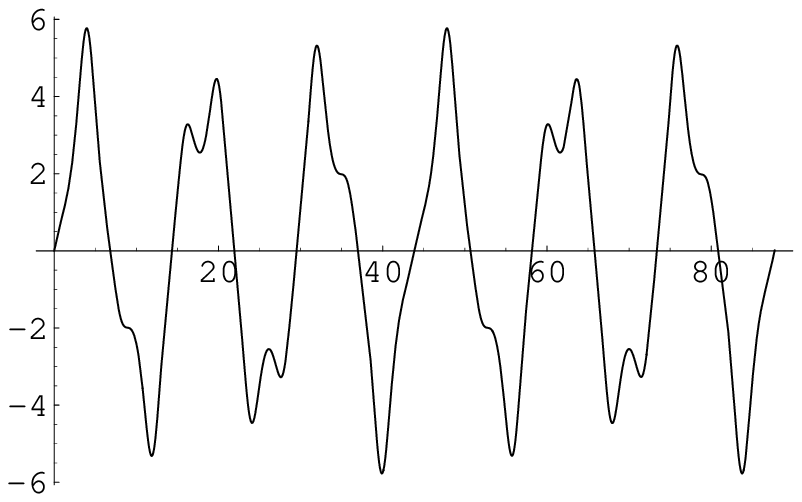} 
  \includegraphics[scale=0.35]{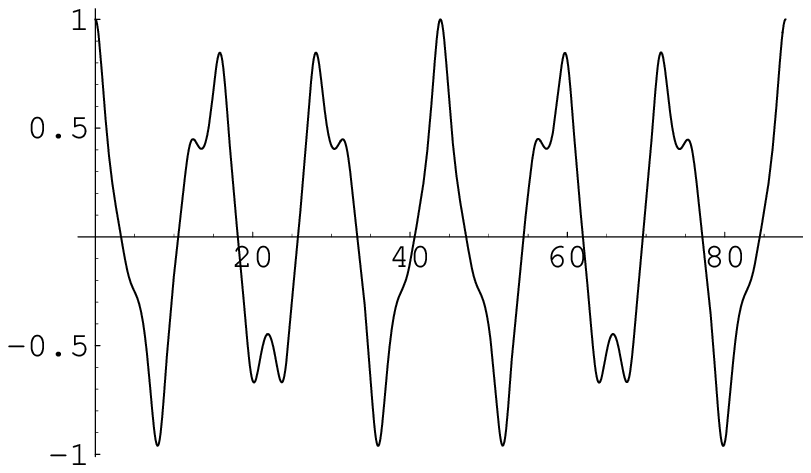}
  \includegraphics[scale=0.35]{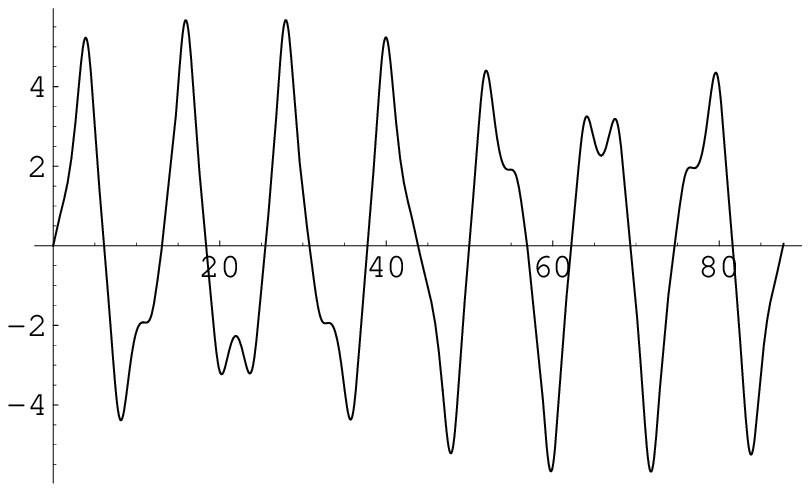}
  \includegraphics[scale=0.35]{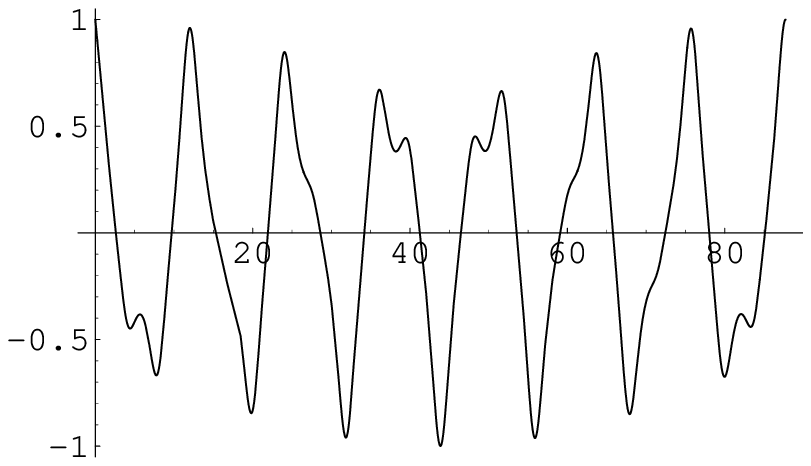}
  \includegraphics[scale=0.35]{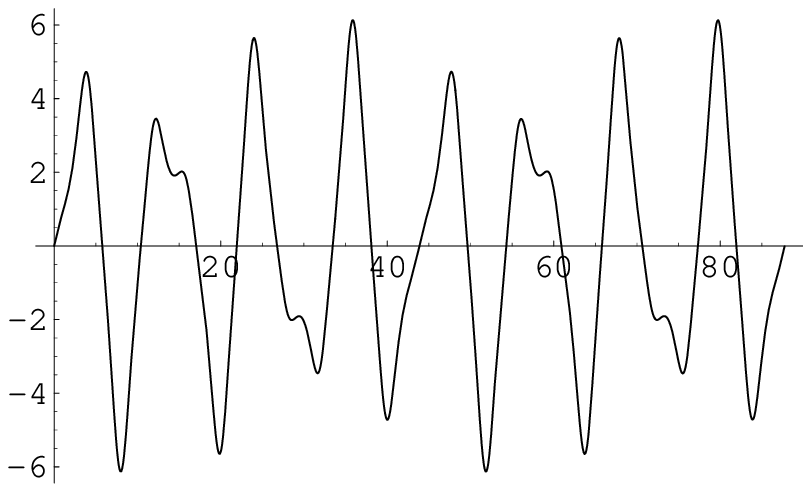}
   \includegraphics[scale=0.35]{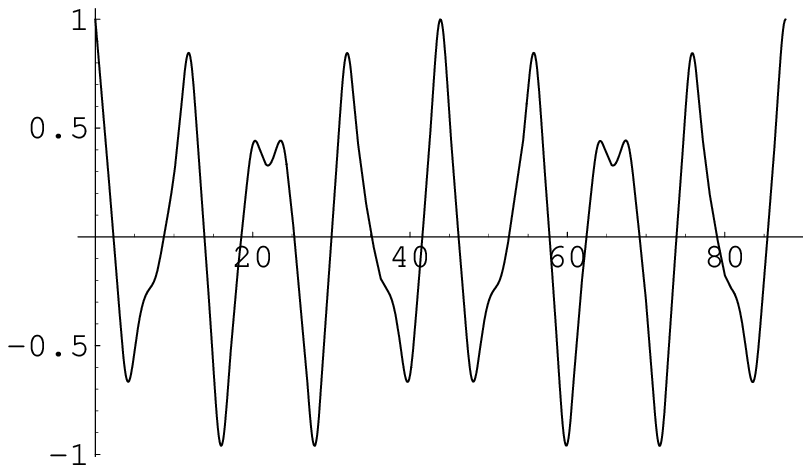}
  \includegraphics[scale=0.35]{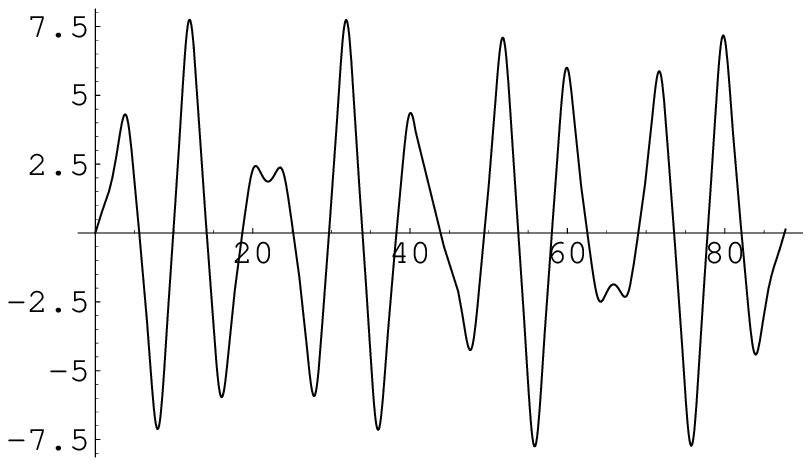}
  \includegraphics[scale=0.35]{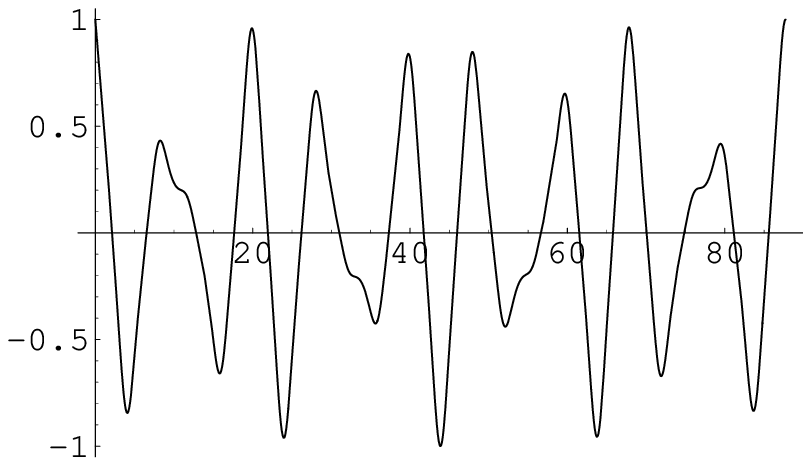}
   \includegraphics[scale=0.35]{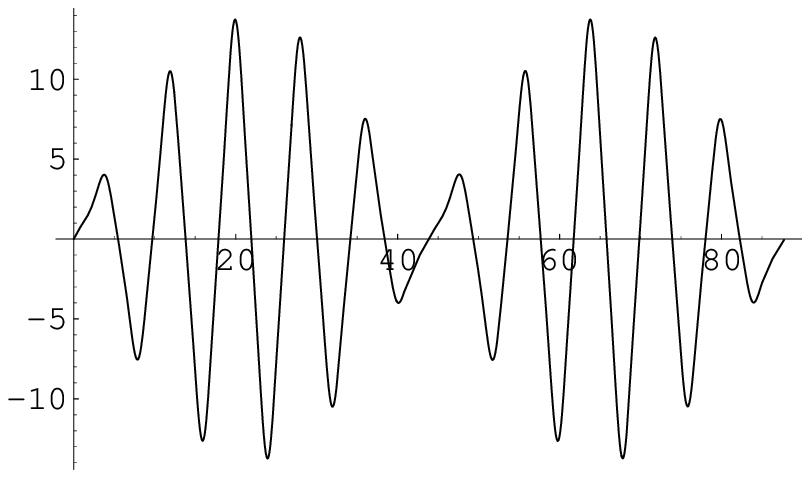}
      \includegraphics[scale=0.35]{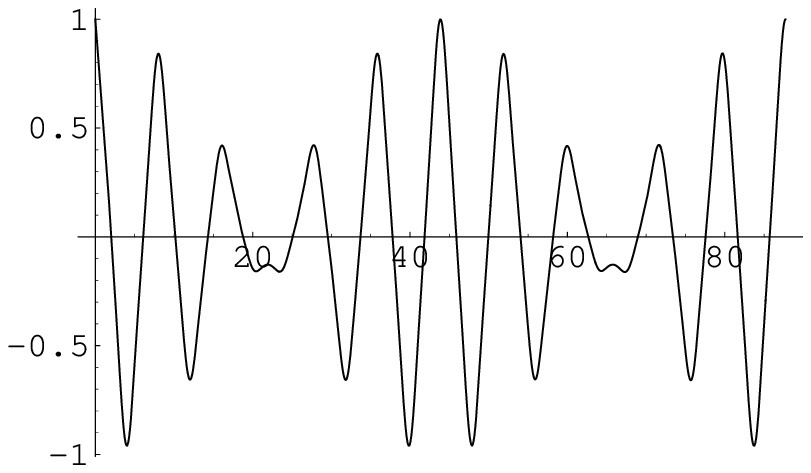}
    \includegraphics[scale=0.35]{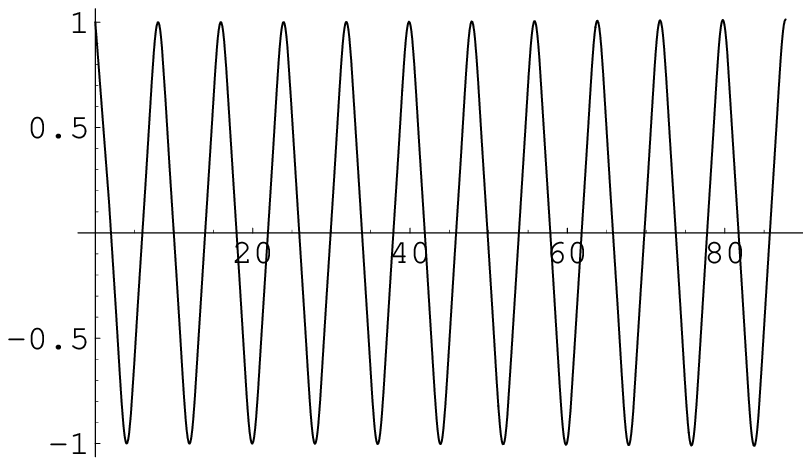}
     \includegraphics[scale=0.35]{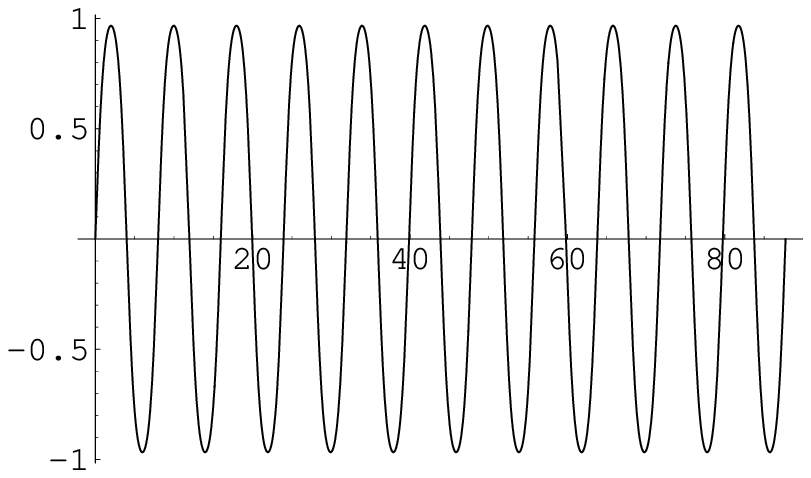}
 \caption{Eigenfunctions associated to the eigenvalues $\lambda_{1,0}$, ..., 
           $\lambda_{22,0}$, $\lambda_{23,0}=0$ of the surface $U_{15}$.}
  \label{O}
\end{figure}

\begin{figure}[phbt]
  \centering
  \includegraphics[scale=0.35]{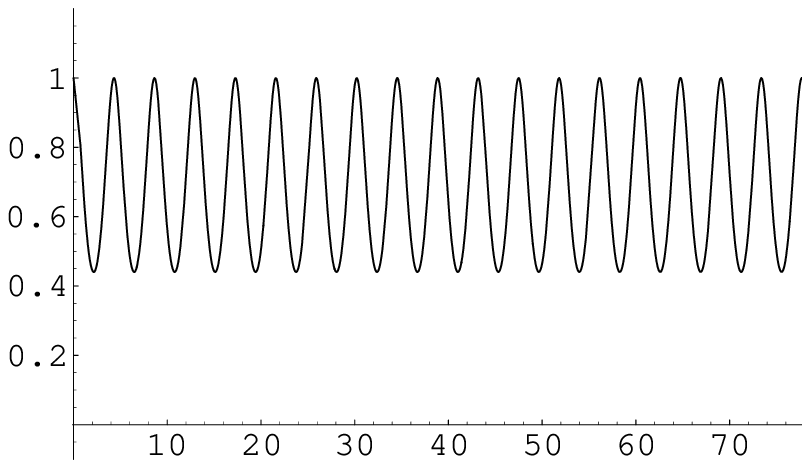}
  \includegraphics[scale=0.35]{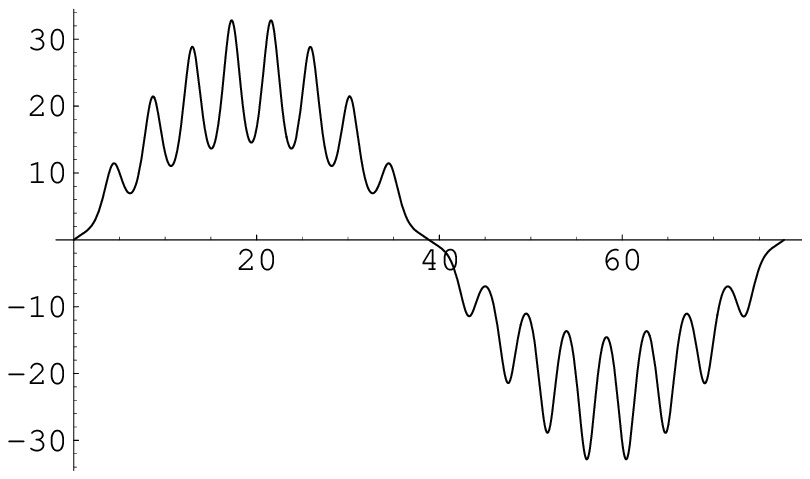}
  \includegraphics[scale=0.35]{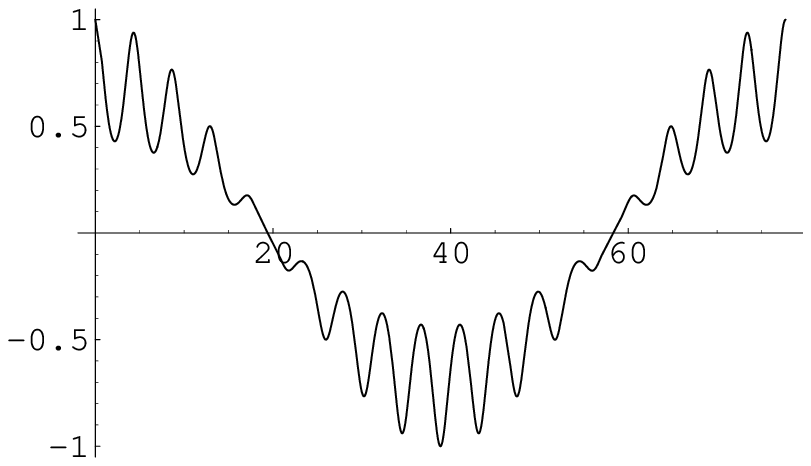}
  \includegraphics[scale=0.35]{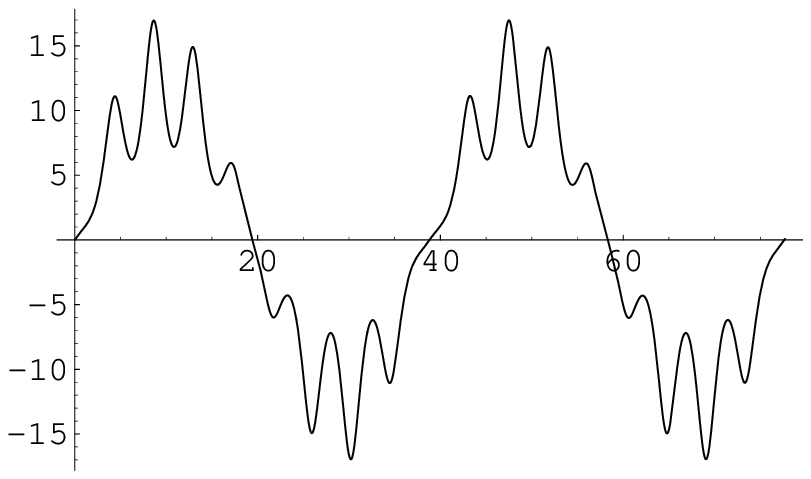}
     \includegraphics[scale=0.35]{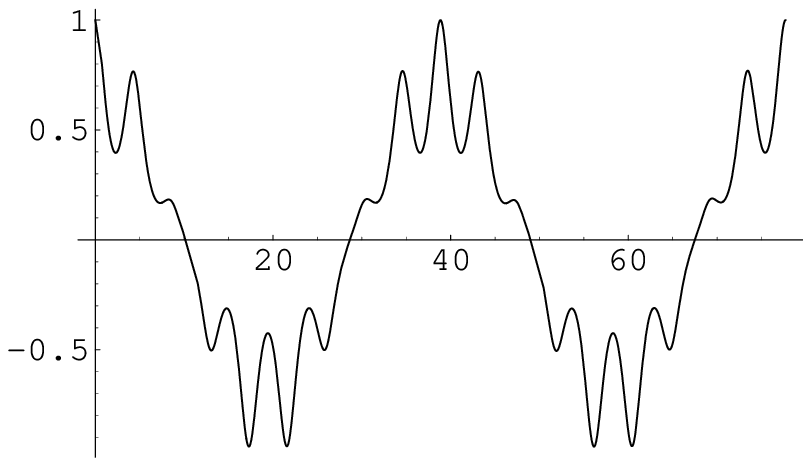}
  \includegraphics[scale=0.35]{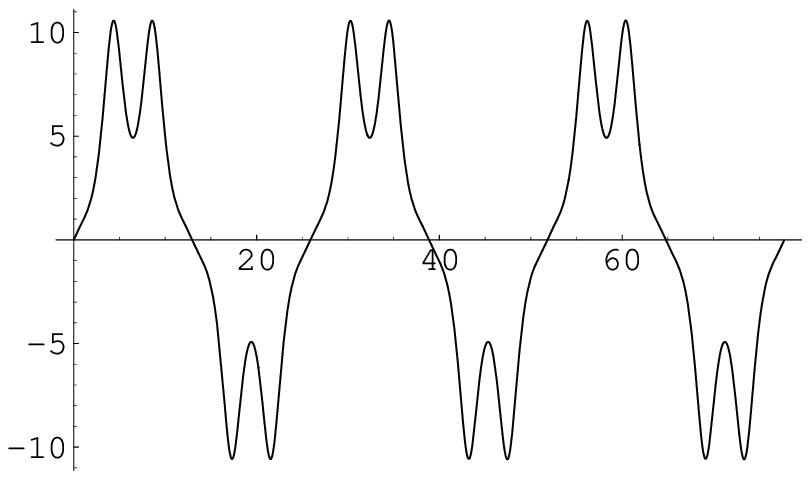}
  \includegraphics[scale=0.35]{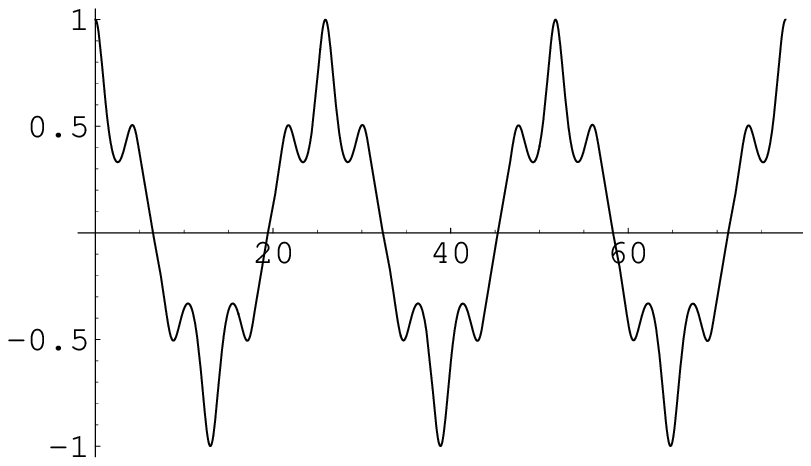}
  \includegraphics[scale=0.35]{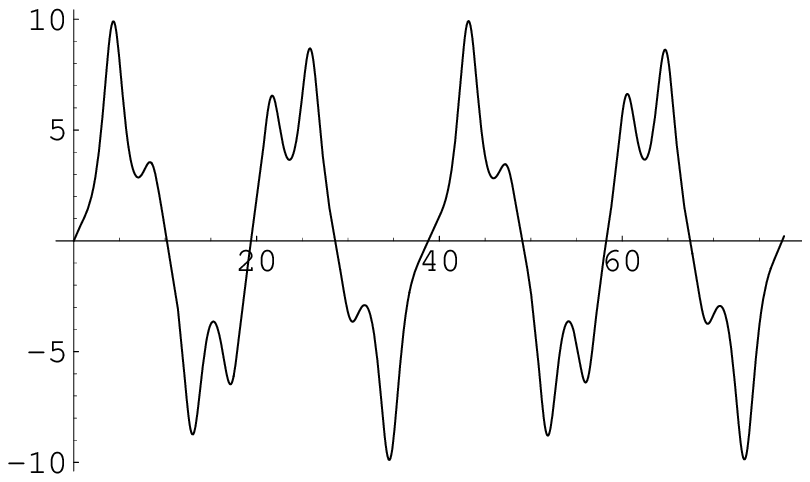}
   \includegraphics[scale=0.35]{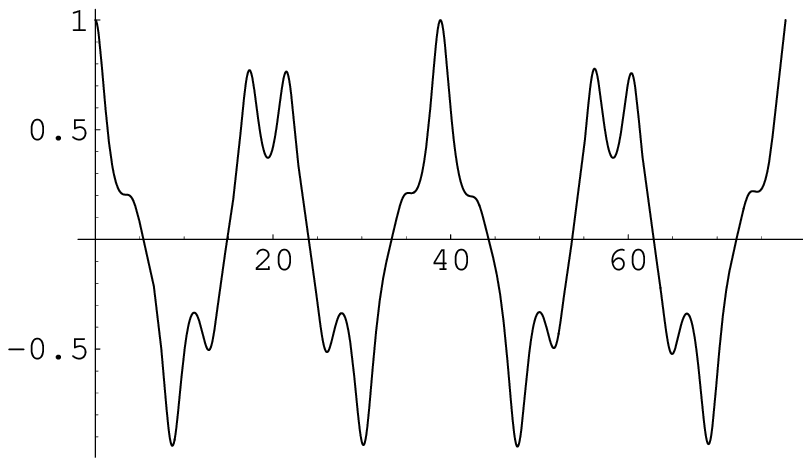}
  \includegraphics[scale=0.35]{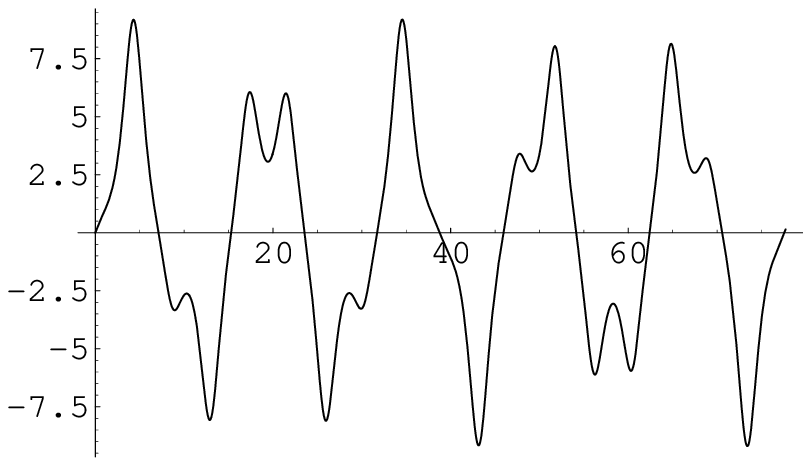}
  \includegraphics[scale=0.35]{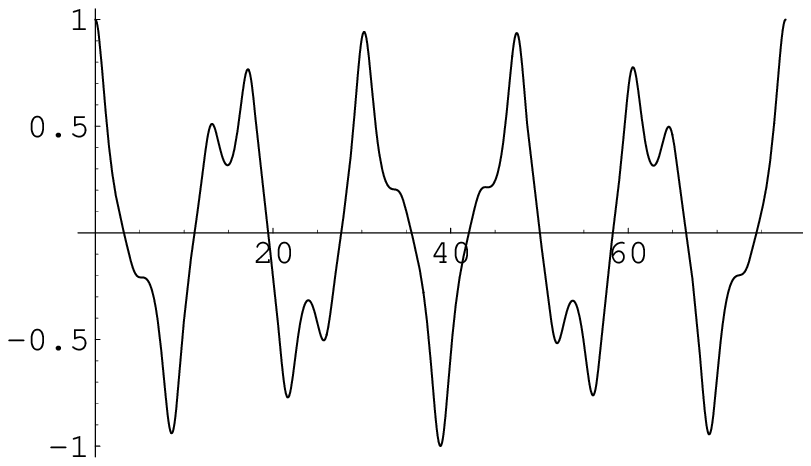}
  \includegraphics[scale=0.35]{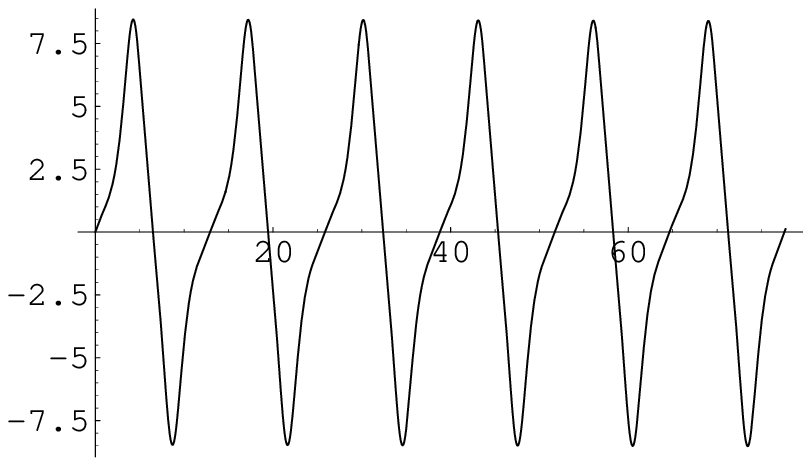}
    \includegraphics[scale=0.35]{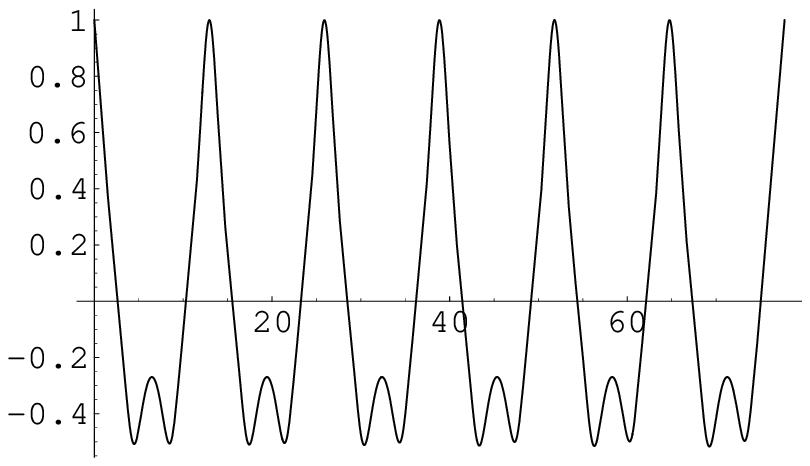}
  \includegraphics[scale=0.35]{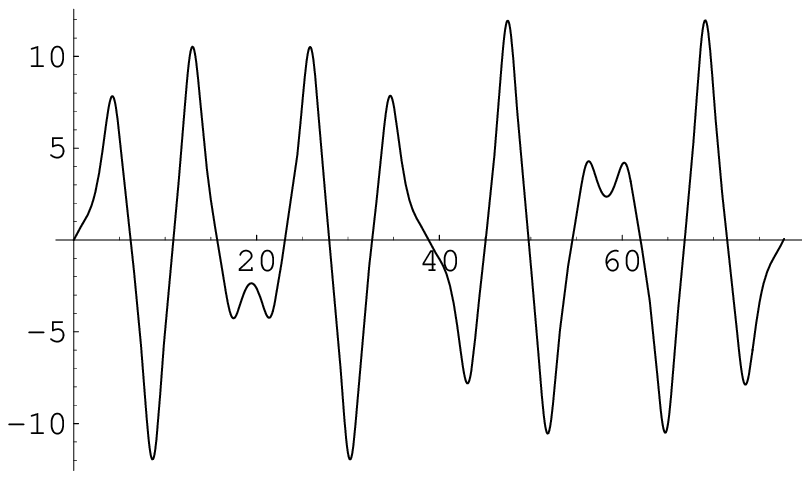}
  \includegraphics[scale=0.35]{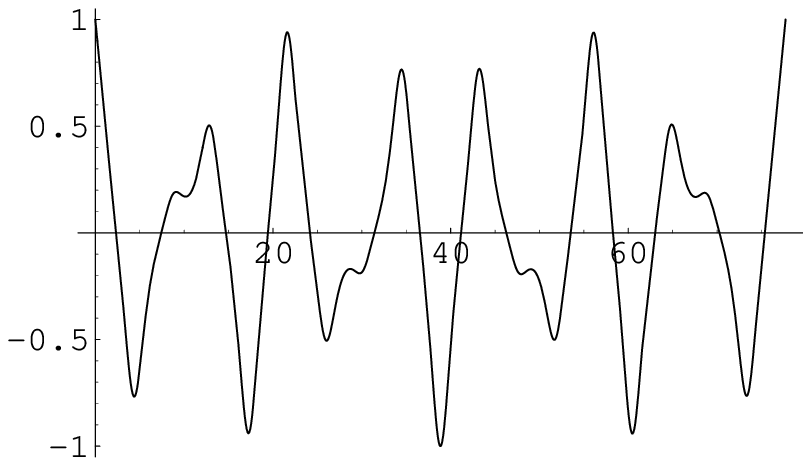}
  \includegraphics[scale=0.35]{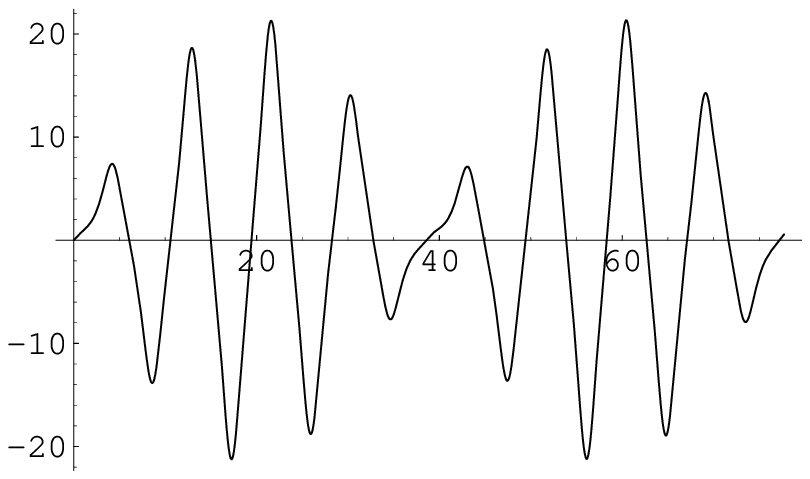}
    \includegraphics[scale=0.35]{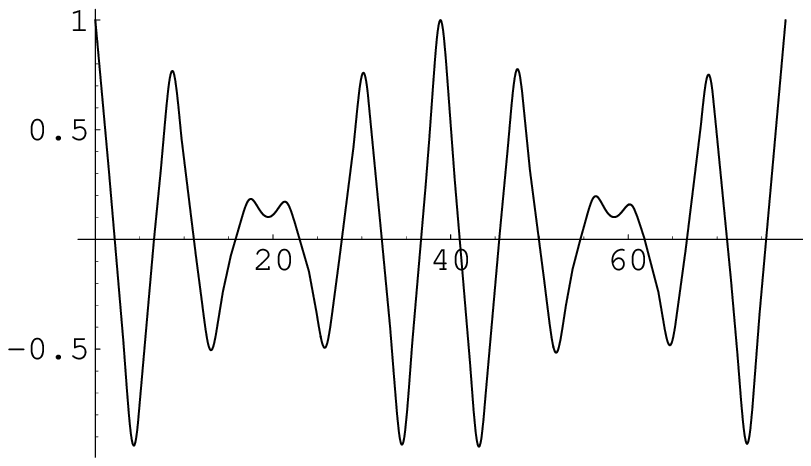}
  \includegraphics[scale=0.35]{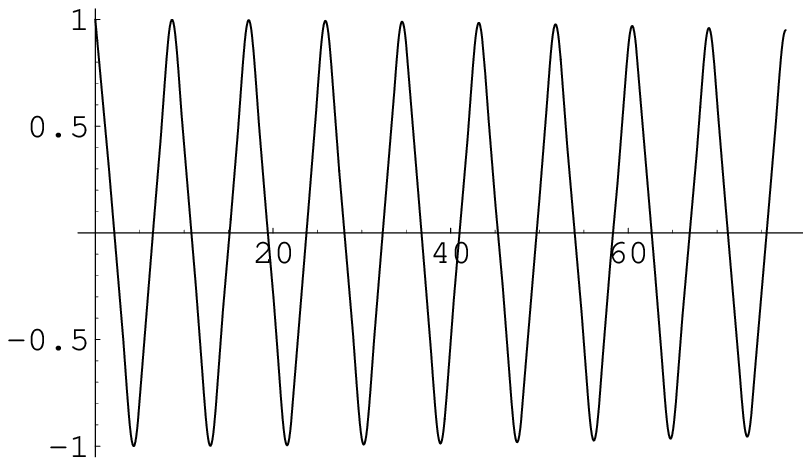}
  \includegraphics[scale=0.35]{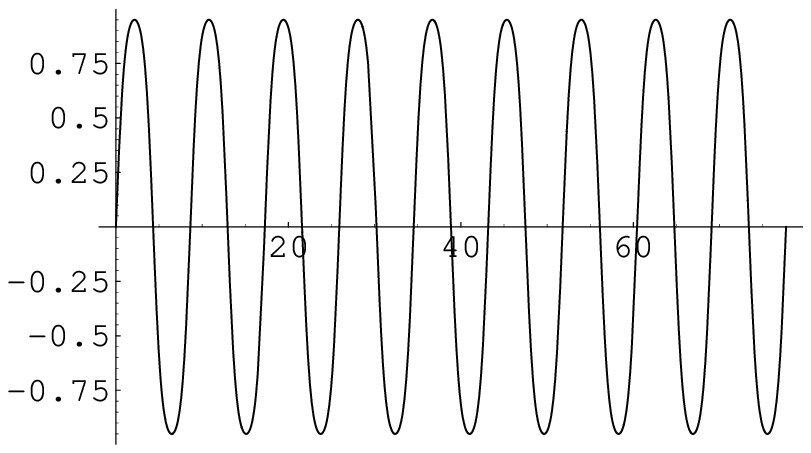}
    \caption{Eigenfunctions associated to the eigenvalues $\lambda_{1,0}$, ..., 
           $\lambda_{18,0}$, $\lambda_{19,0}=0$ of the surface $U_{16}$.}
  \label{Q}
\end{figure}

\begin{figure}[phbt]
  \centering
  \includegraphics[scale=0.35]{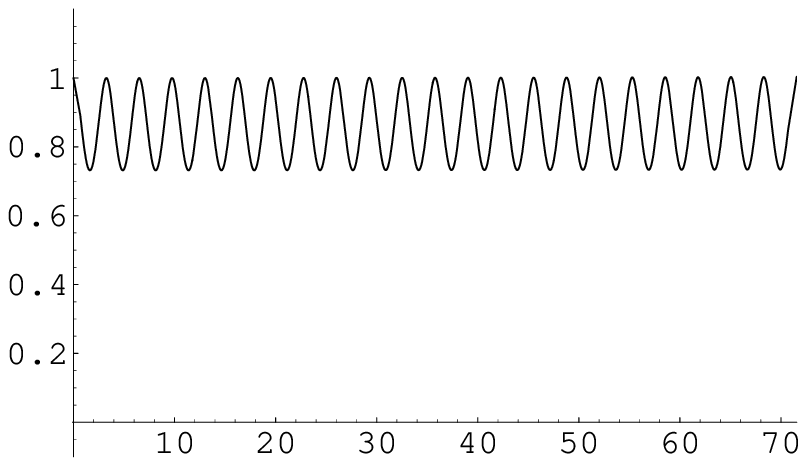}
  \includegraphics[scale=0.35]{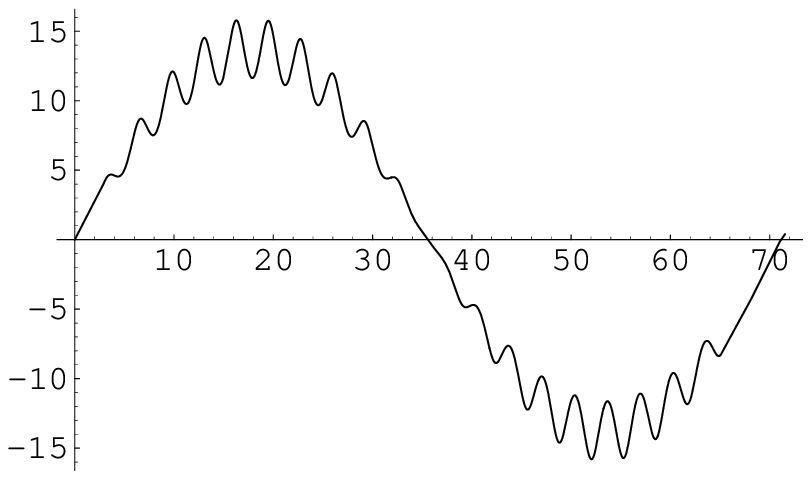}
  \includegraphics[scale=0.35]{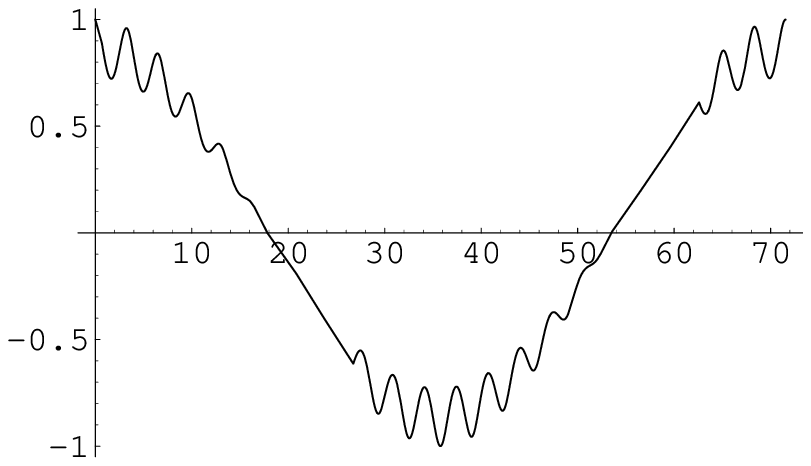}
  \includegraphics[scale=0.35]{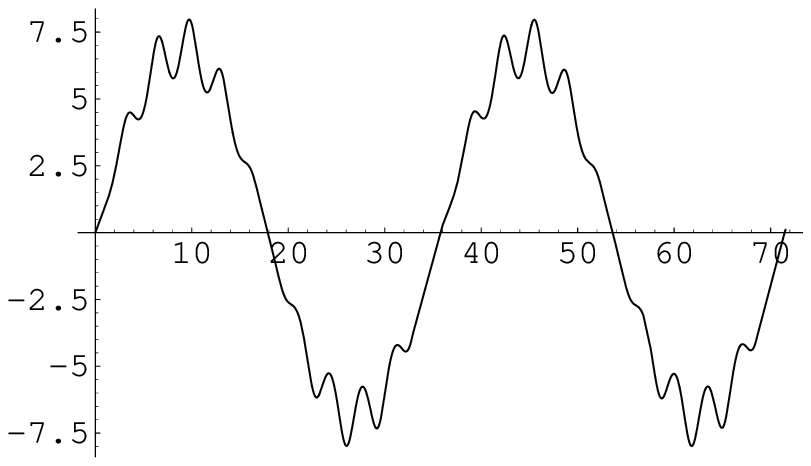} 
     \includegraphics[scale=0.35]{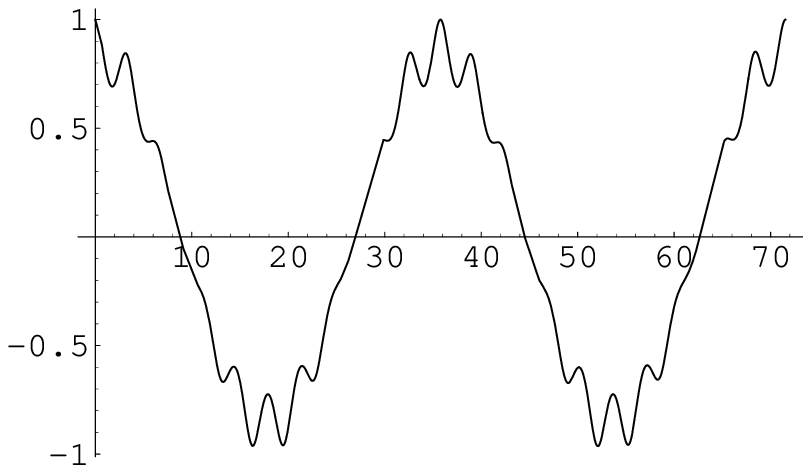}
  \includegraphics[scale=0.35]{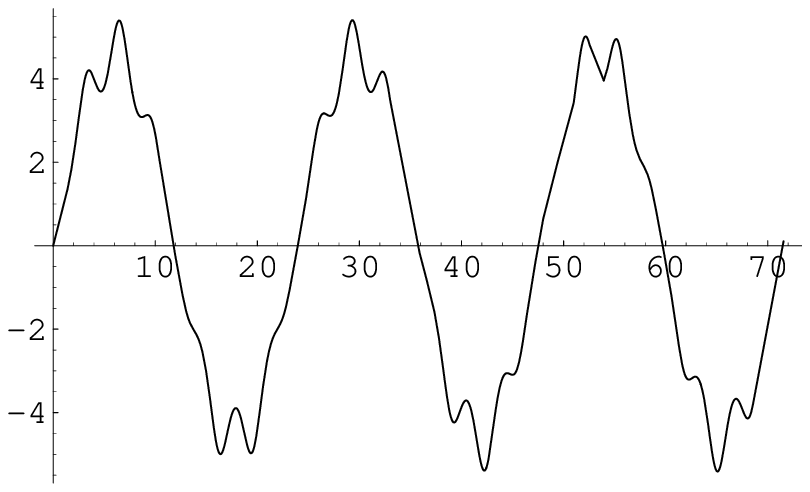}
  \includegraphics[scale=0.35]{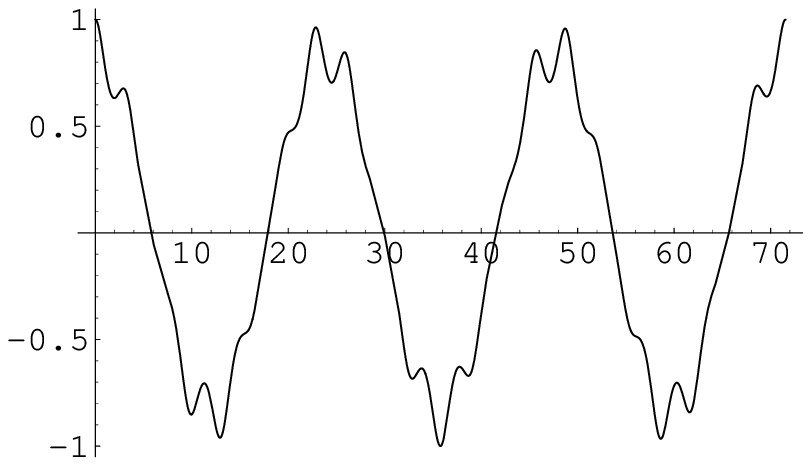}
  \includegraphics[scale=0.35]{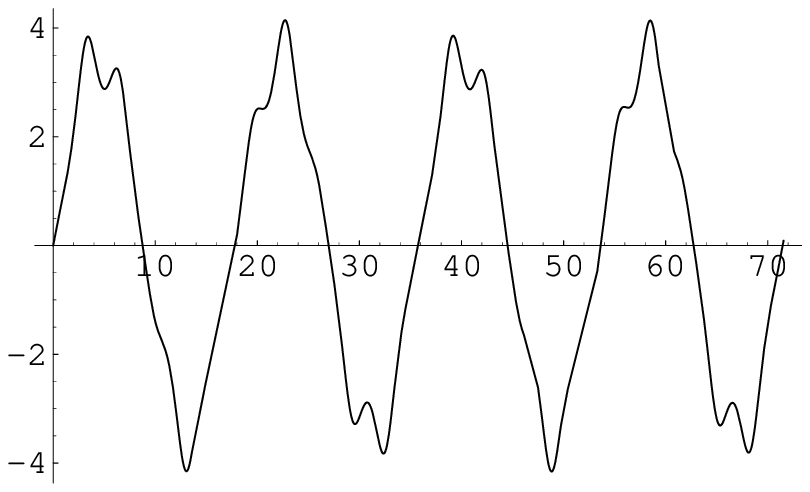}
   \includegraphics[scale=0.35]{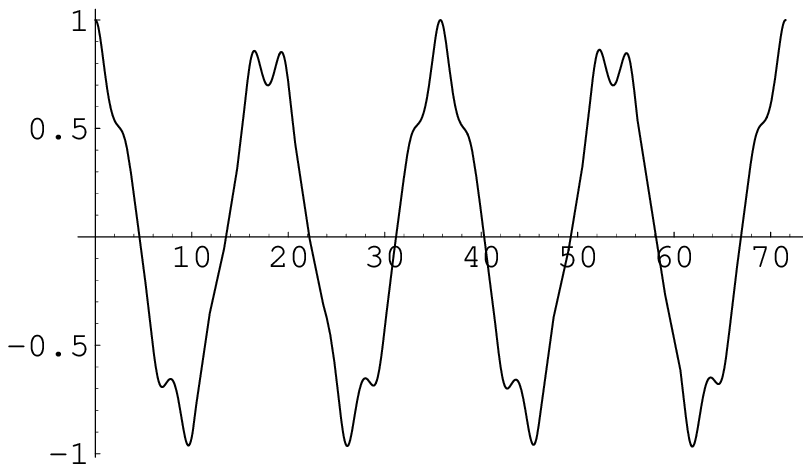}
  \includegraphics[scale=0.35]{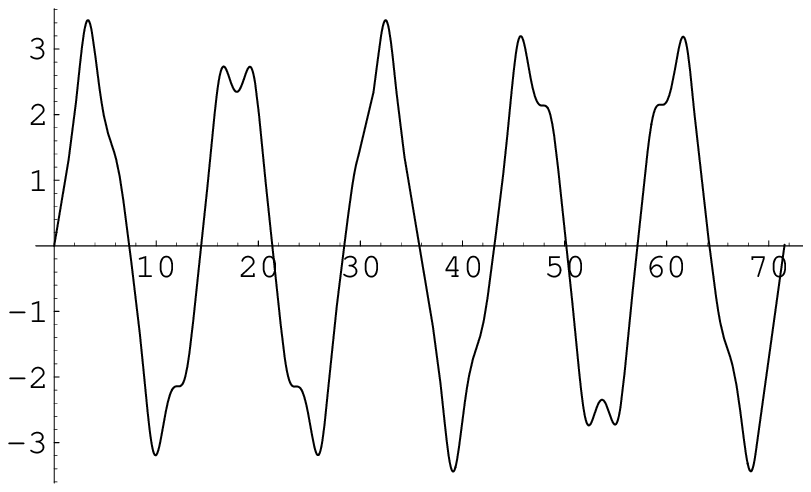}
  \includegraphics[scale=0.35]{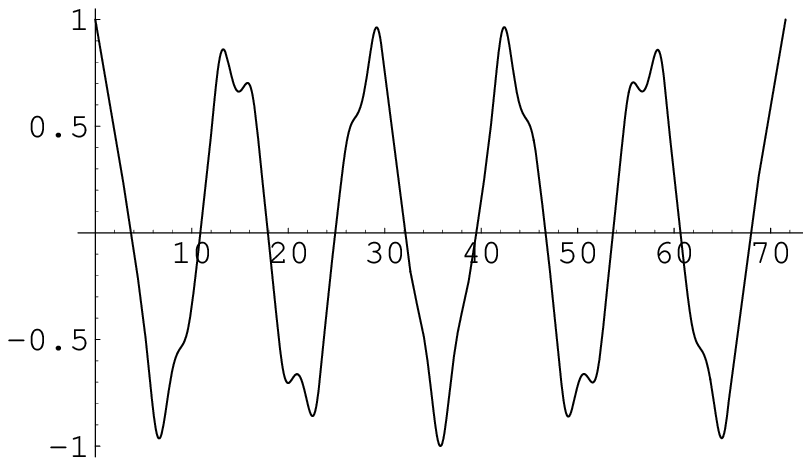}
  \includegraphics[scale=0.35]{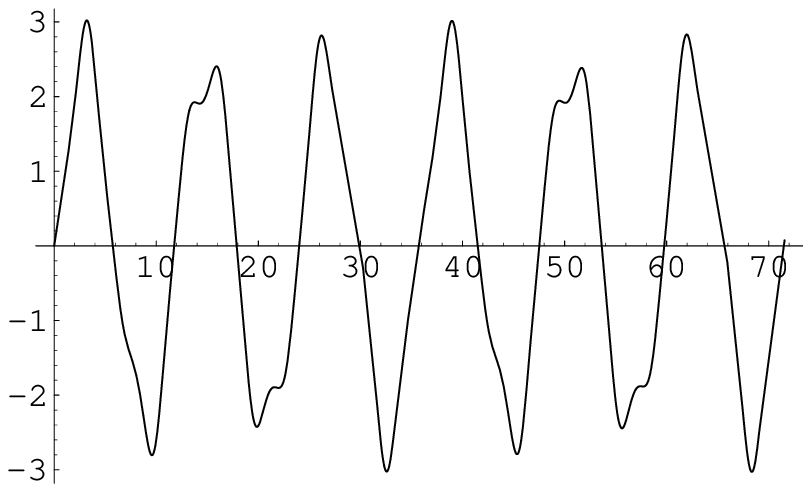}  
    \includegraphics[scale=0.35]{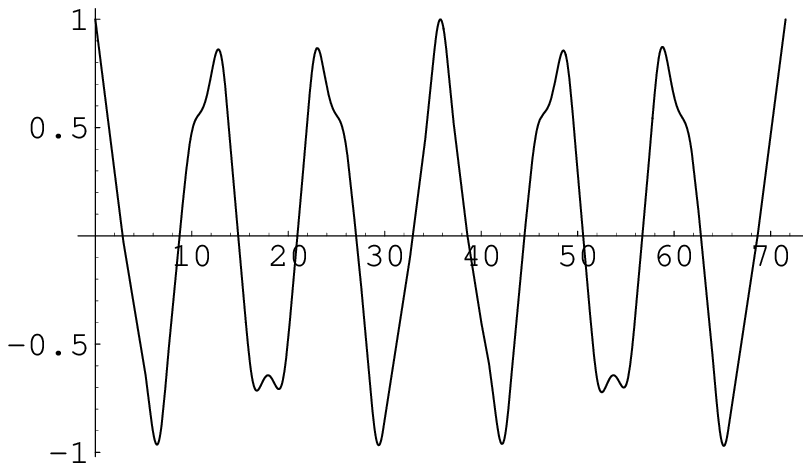}
  \includegraphics[scale=0.35]{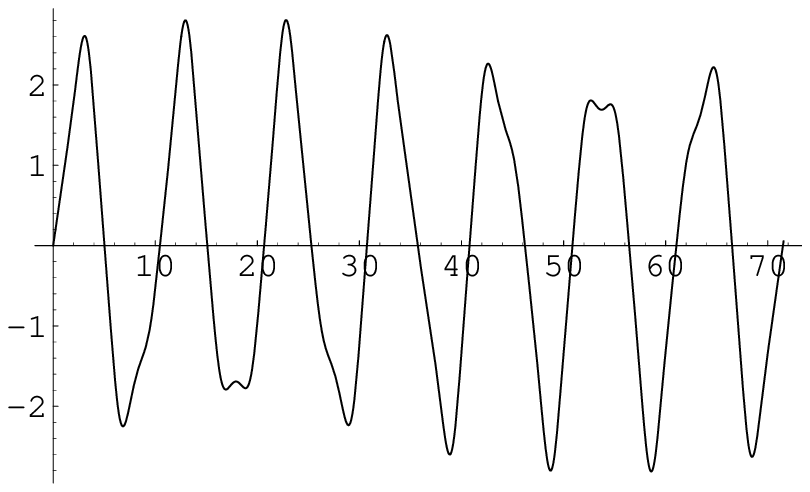}
  \includegraphics[scale=0.35]{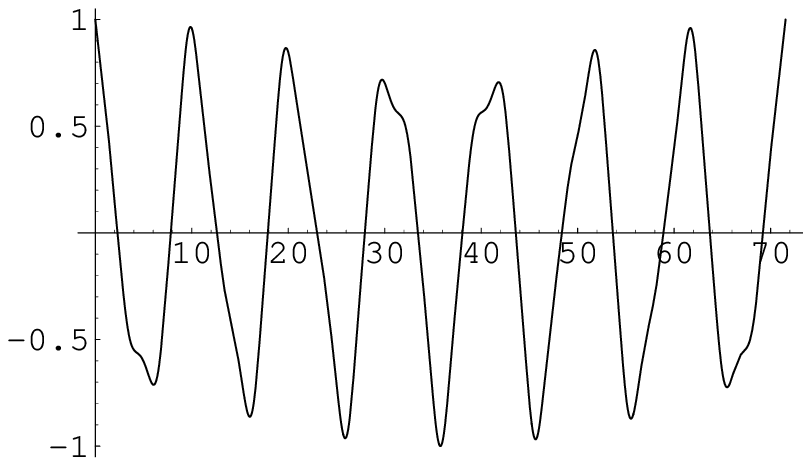}
  \includegraphics[scale=0.35]{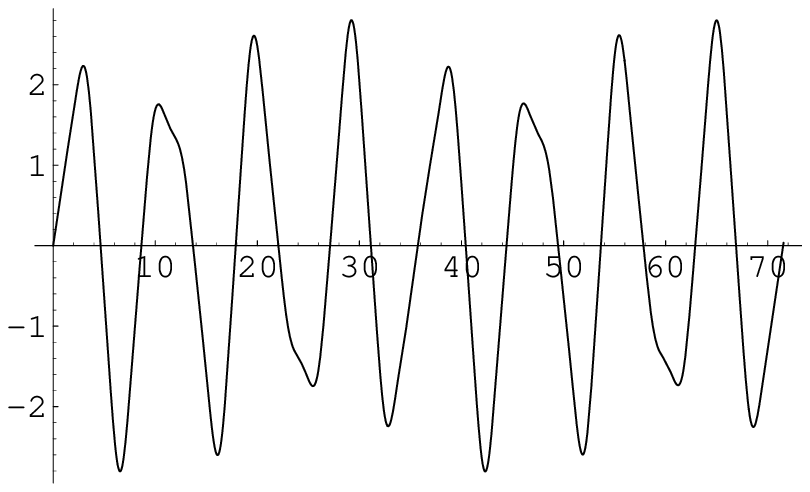}
    \includegraphics[scale=0.35]{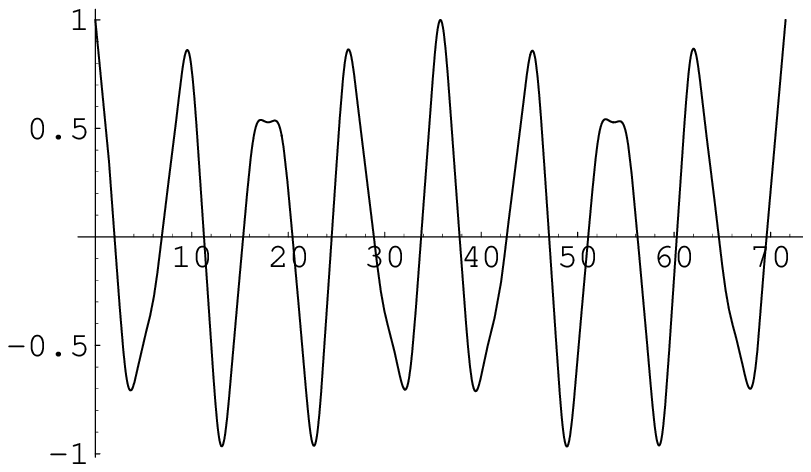}
  \includegraphics[scale=0.35]{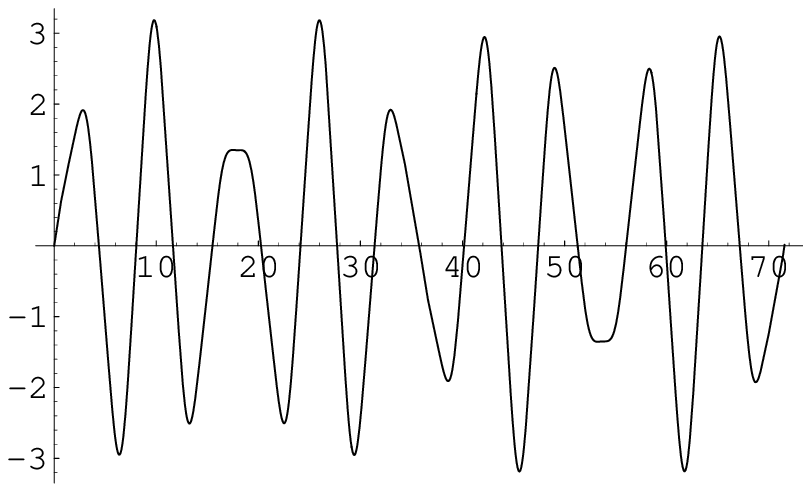}
  \includegraphics[scale=0.35]{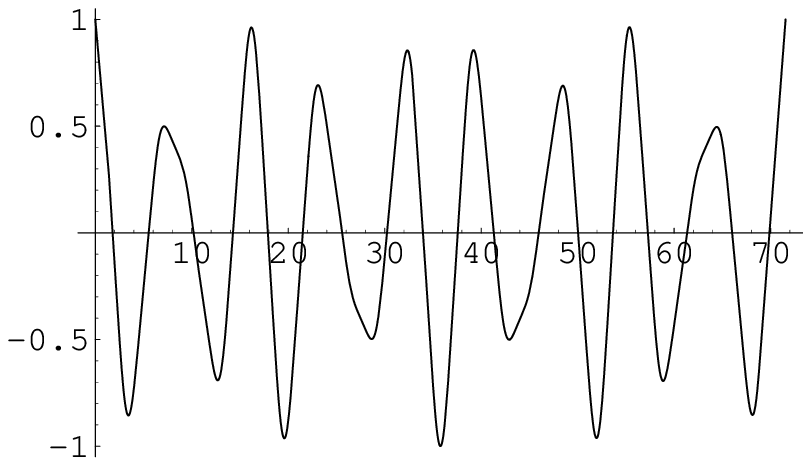}
     \includegraphics[scale=0.35]{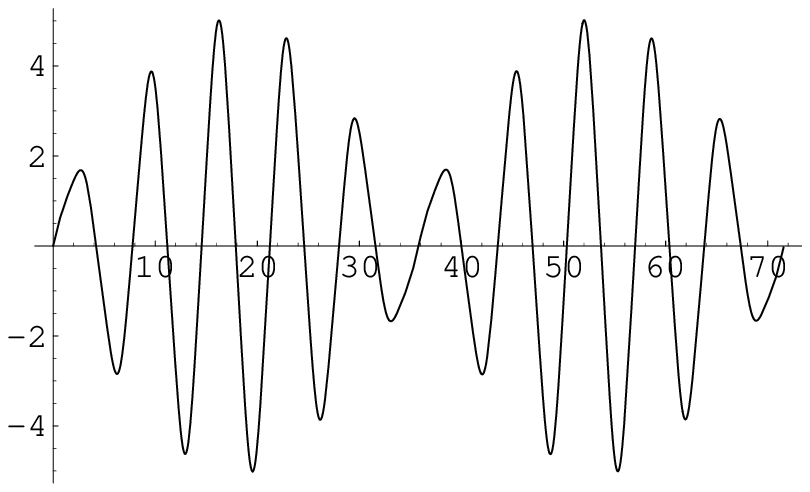}
   \includegraphics[scale=0.35]{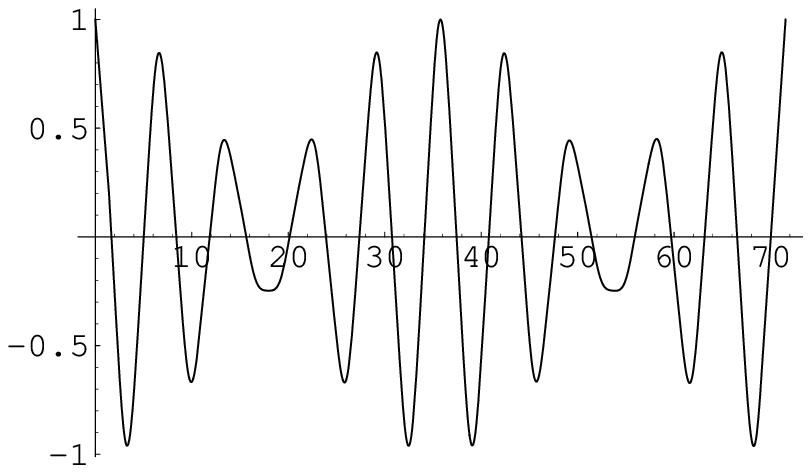}
    \includegraphics[scale=0.35]{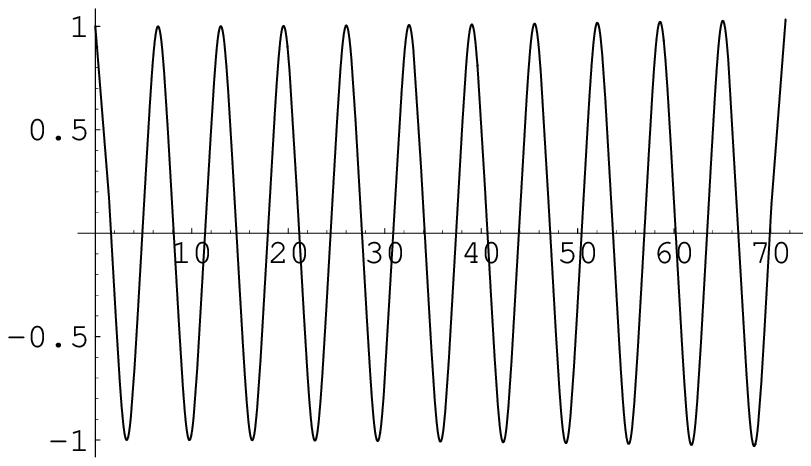}
        \includegraphics[scale=0.35]{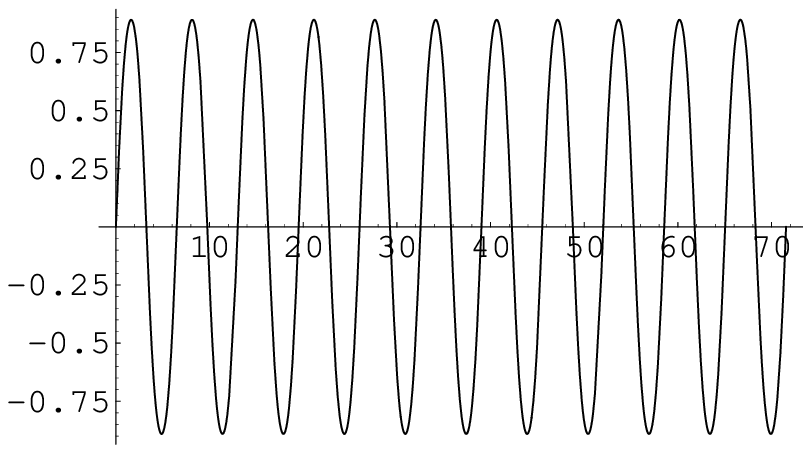}
 \caption{Eigenfunctions associated to the eigenvalues $\lambda_{1,0}$, ..., 
           $\lambda_{22,0}$, $\lambda_{23,0}=0$ of the surface $U_{17}$.}
  \label{R}
\end{figure}

\begin{figure}[phbt]
  \centering
  \includegraphics[scale=0.3]{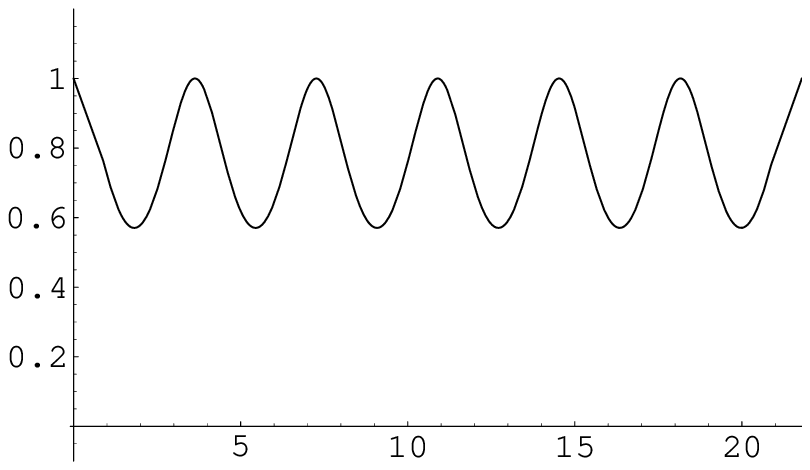}
  \includegraphics[scale=0.3]{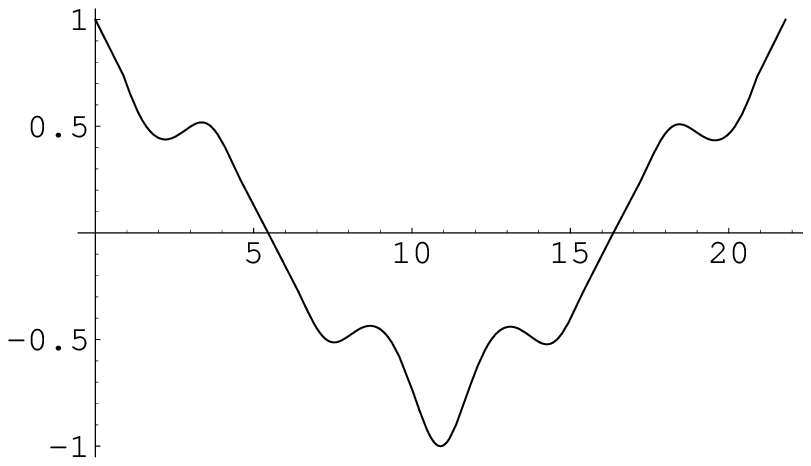}
  \includegraphics[scale=0.3]{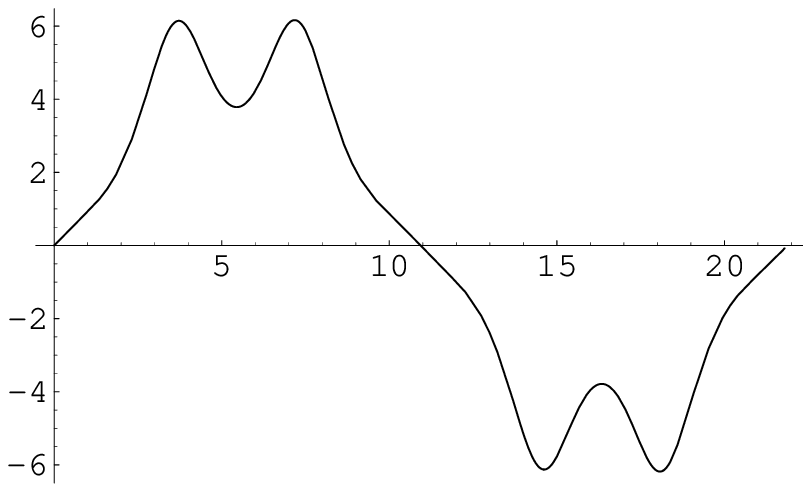}
   \includegraphics[scale=0.3]{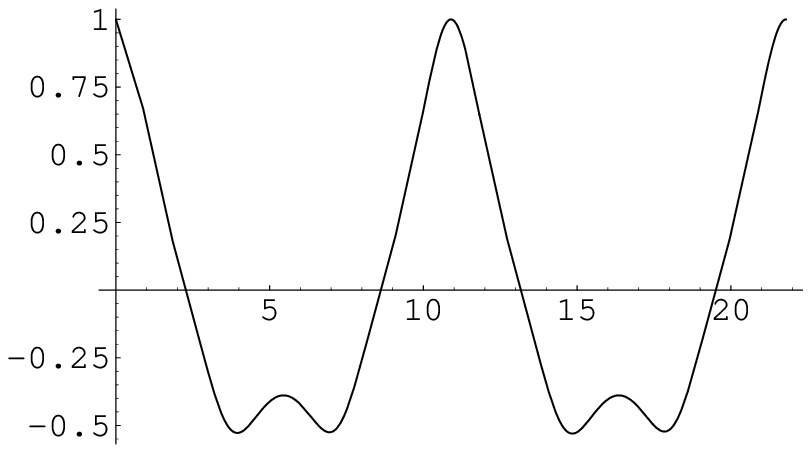}
   \includegraphics[scale=0.3]{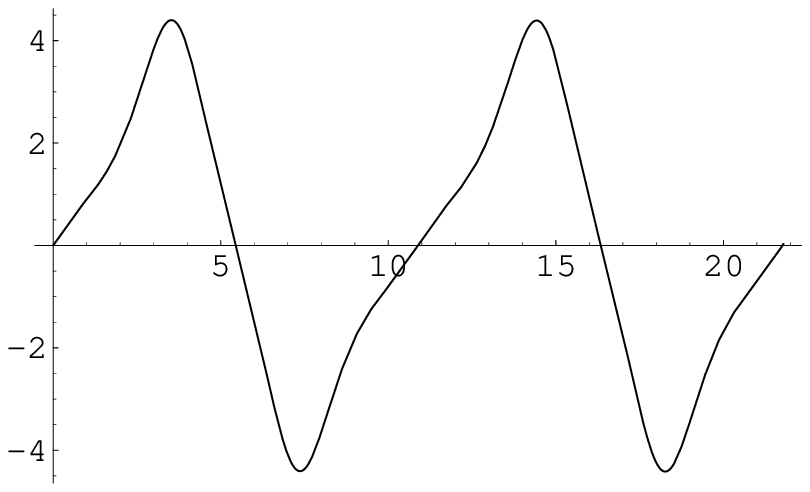}
  \includegraphics[scale=0.3]{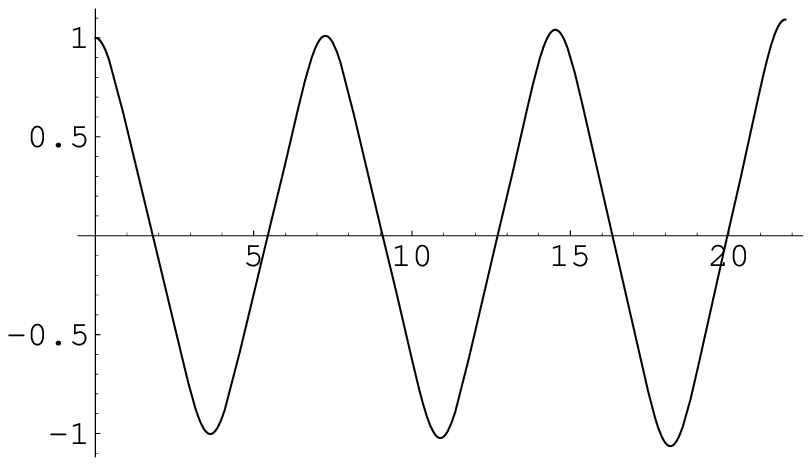}
  \includegraphics[scale=0.3]{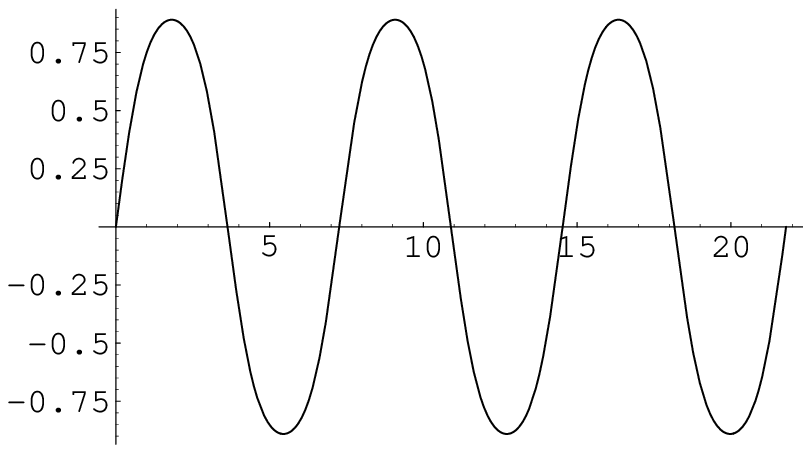}
 \caption{Eigenfunctions associated to the eigenvalues $\lambda_{1,0}$, ..., 
           $\lambda_{6,0}$, $\lambda_{7,0}=0$ of the surface $N_{1}$.}
  \label{S}
\end{figure}

\begin{figure}[phbt]
  \centering
 \includegraphics[scale=0.656]{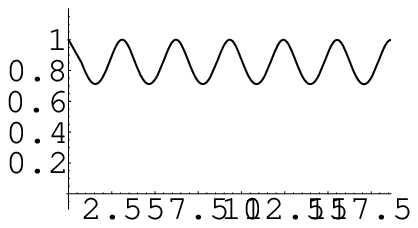}
  \includegraphics[scale=0.656]{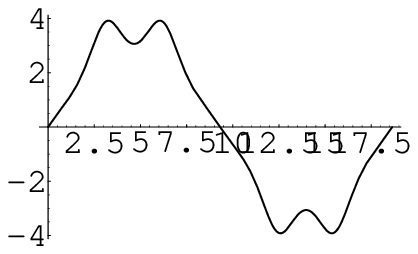}
  \includegraphics[scale=0.656]{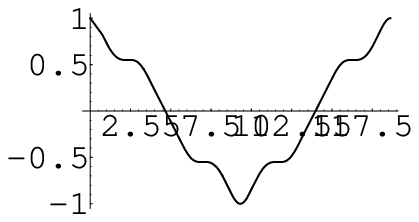}
  \includegraphics[scale=0.656]{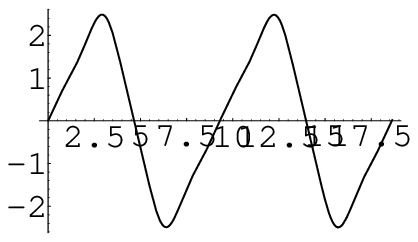}
      \includegraphics[scale=0.656]{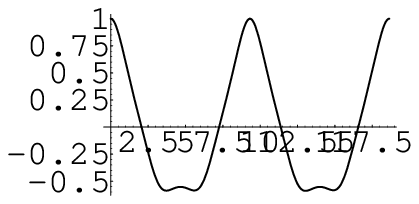}
  \includegraphics[scale=0.656]{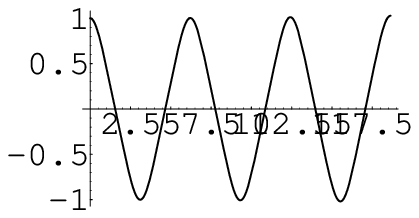}
  \includegraphics[scale=0.656]{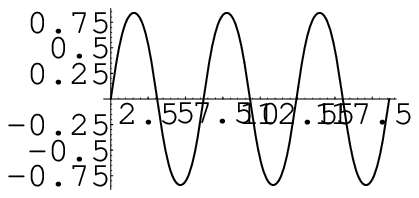}
 \caption{Eigenfunctions associated to the eigenvalues $\lambda_{1,0}$, ..., 
           $\lambda_{6,0}$, $\lambda_{7,0}=0$ of the surface $N_{2}$.}
  \label{T}
\end{figure}

\begin{figure}[phbt]
  \centering
  \includegraphics[scale=0.656]{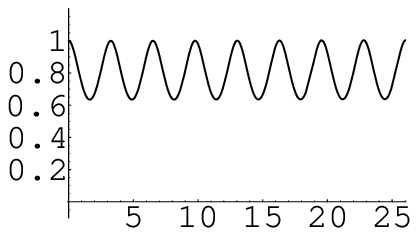}
  \includegraphics[scale=0.656]{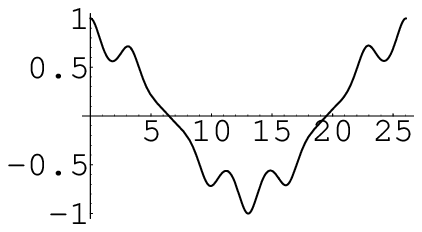}
  \includegraphics[scale=0.656]{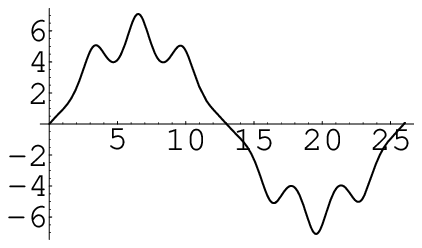}
  \includegraphics[scale=0.656]{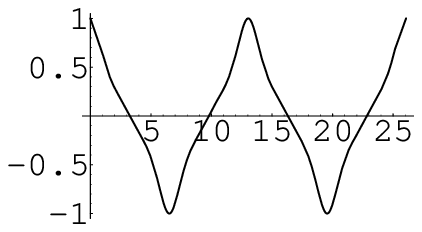}
      \includegraphics[scale=0.656]{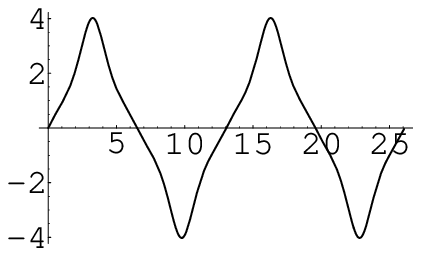}
   \includegraphics[scale=0.656]{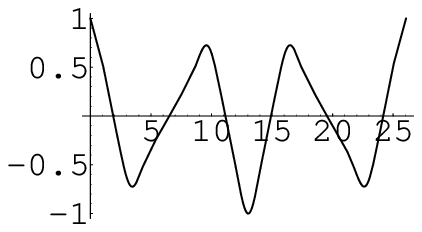}
  \includegraphics[scale=0.656]{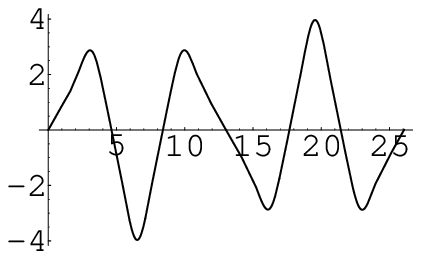}
  \includegraphics[scale=0.656]{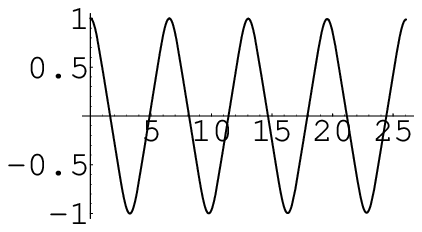}
  \includegraphics[scale=0.656]{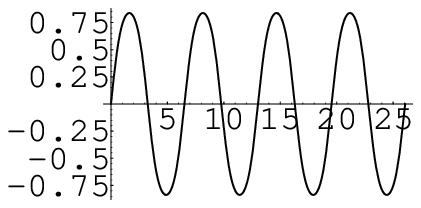}
 \caption{Eigenfunctions associated to the eigenvalues $\lambda_{1,0}$, ..., 
           $\lambda_{8,0}$, $\lambda_{9,0}=0$ of the surface $N_{3}$.}
  \label{U}
\end{figure}

\begin{figure}[phbt]
  \centering
  \includegraphics[scale=0.656]{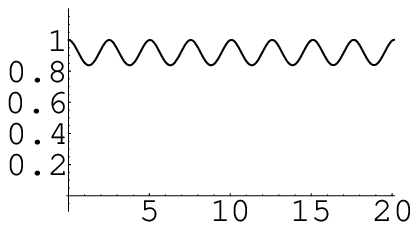}
  \includegraphics[scale=0.656]{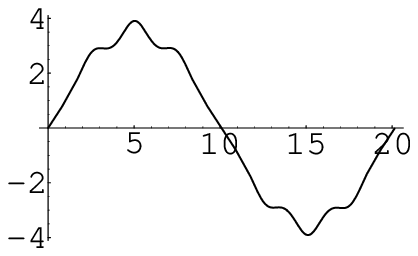}
  \includegraphics[scale=0.656]{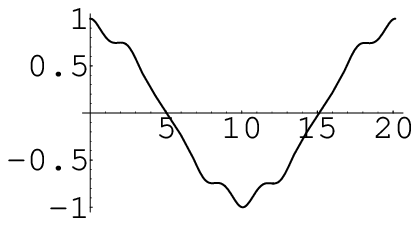}
  \includegraphics[scale=0.656]{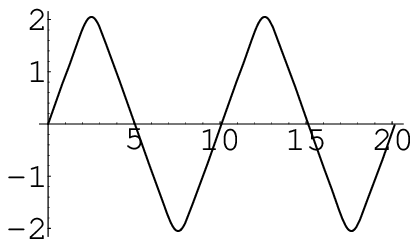}
  \includegraphics[scale=0.656]{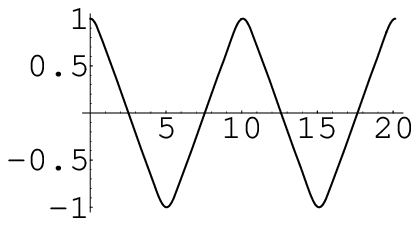}
  \includegraphics[scale=0.656]{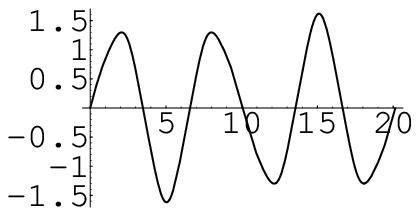}
  \includegraphics[scale=0.656]{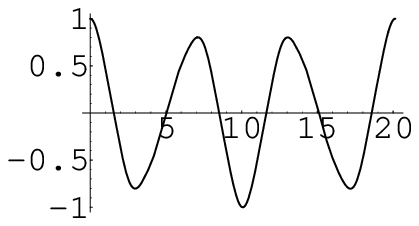}
  \includegraphics[scale=0.656]{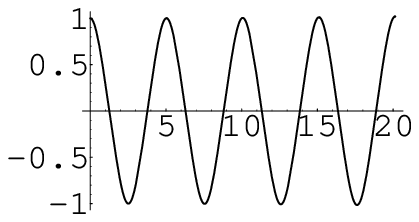}
  \includegraphics[scale=0.656]{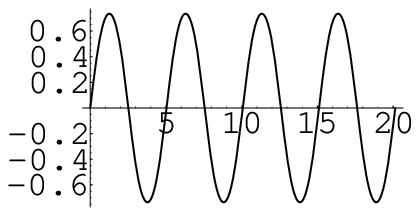}
 \caption{Eigenfunctions associated to the eigenvalues $\lambda_{1,0}$, ..., 
           $\lambda_{8,0}$, $\lambda_{9,0}=0$ of the surface $N_{4}$.}
  \label{V}
\end{figure}

\begin{figure}[phbt]
  \centering
  \includegraphics[scale=0.6]{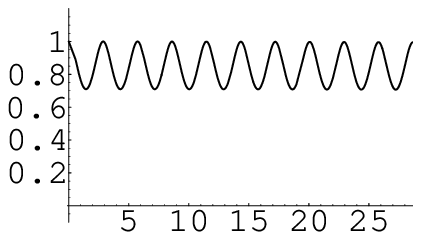}
  \includegraphics[scale=0.6]{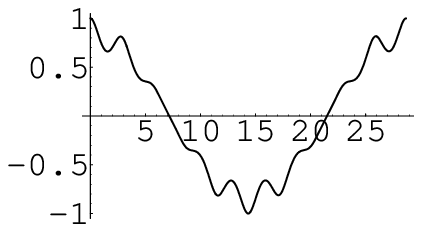}
  \includegraphics[scale=0.6]{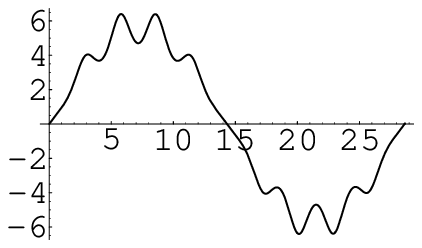}
  \includegraphics[scale=0.6]{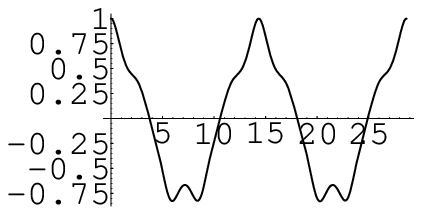}
  \includegraphics[scale=0.6]{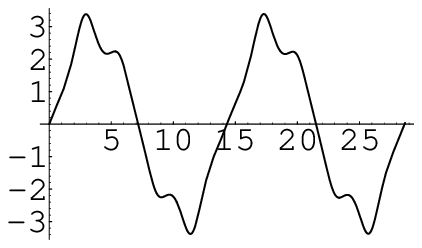}
  \includegraphics[scale=0.6]{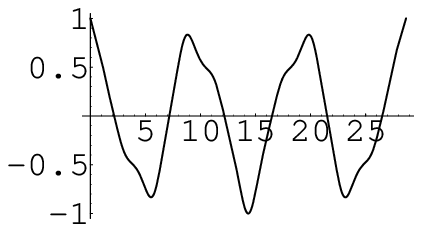}
  \includegraphics[scale=0.6]{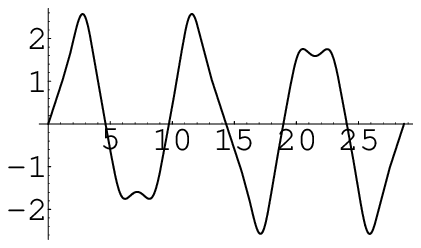}
  \includegraphics[scale=0.6]{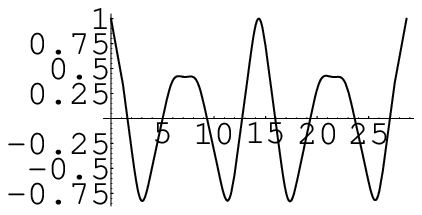}
    \includegraphics[scale=0.6]{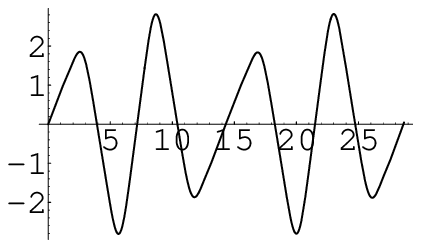}
  \includegraphics[scale=0.6]{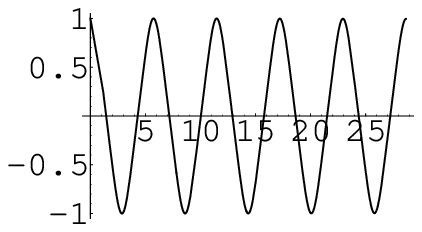}
  \includegraphics[scale=0.6]{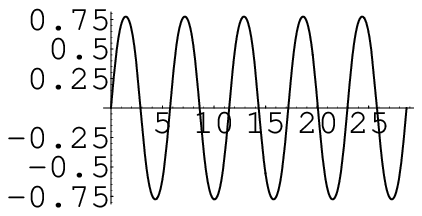}
 \caption{Eigenfunctions associated to the eigenvalues $\lambda_{1,0}$, ..., 
           $\lambda_{10,0}$, $\lambda_{11,0}=0$ of the surface $N_{5}$.}
  \label{W}
\end{figure}

\begin{figure}[phbt]
  \centering
  \includegraphics[scale=0.735]{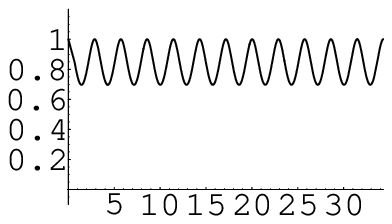}
  \includegraphics[scale=0.735]{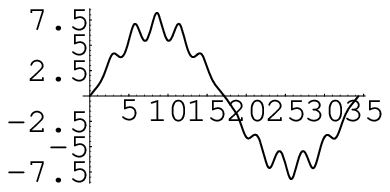}
  \includegraphics[scale=0.735]{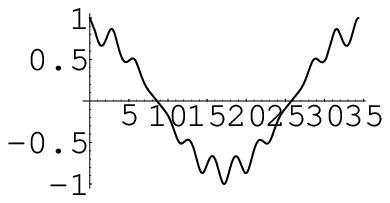}
  \includegraphics[scale=0.735]{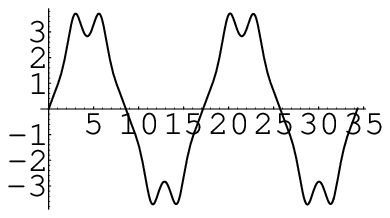}
   \includegraphics[scale=0.735]{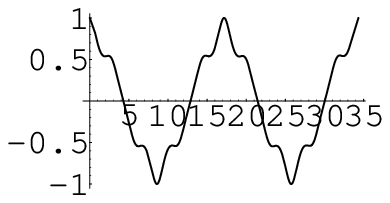}
  \includegraphics[scale=0.735]{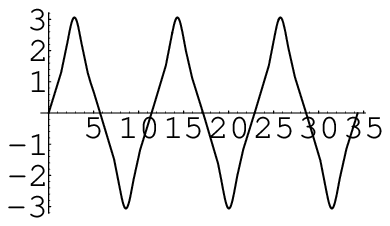}
  \includegraphics[scale=0.735]{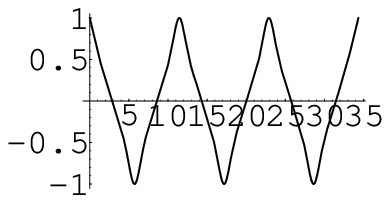}
  \includegraphics[scale=0.735]{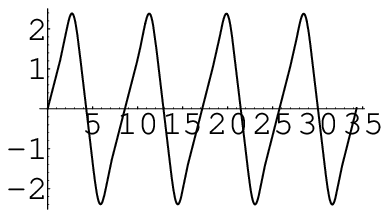}
  \includegraphics[scale=0.735]{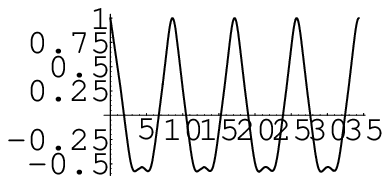}
  \includegraphics[scale=0.735]{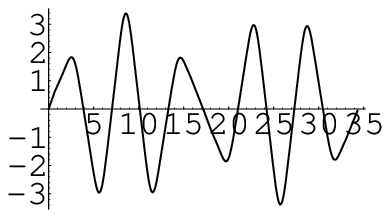}
  \includegraphics[scale=0.735]{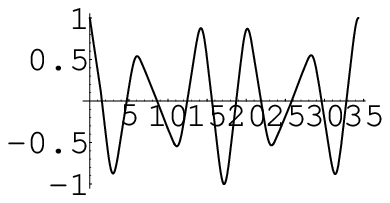}
  \includegraphics[scale=0.735]{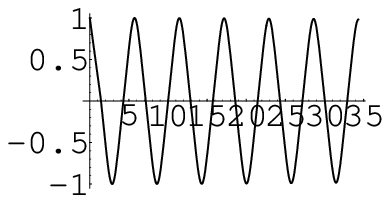}
   \includegraphics[scale=0.735]{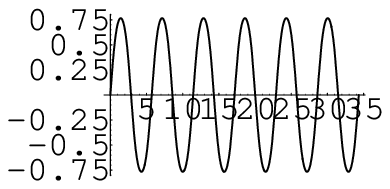}
 \caption{Eigenfunctions associated to the eigenvalues $\lambda_{1,0}$, ..., 
           $\lambda_{12,0}$, $\lambda_{13,0}=0$ of the surface $N_{6}$.}
  \label{X}
\end{figure}

\begin{figure}[phbt]
  \centering
  \includegraphics[scale=0.35]{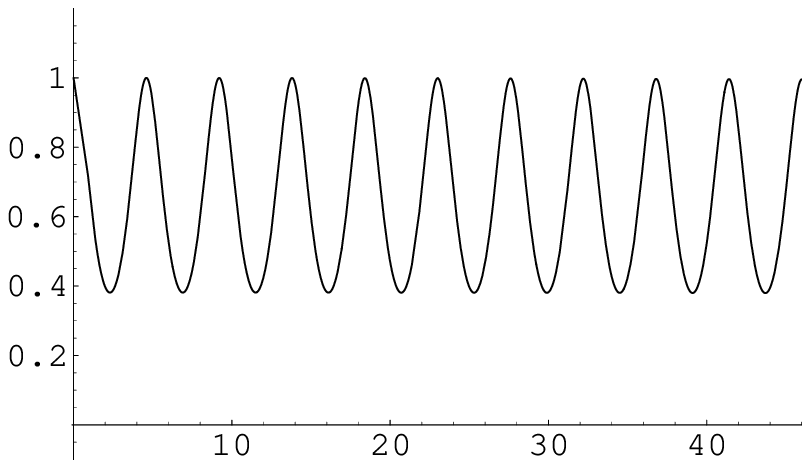}
  \includegraphics[scale=0.35]{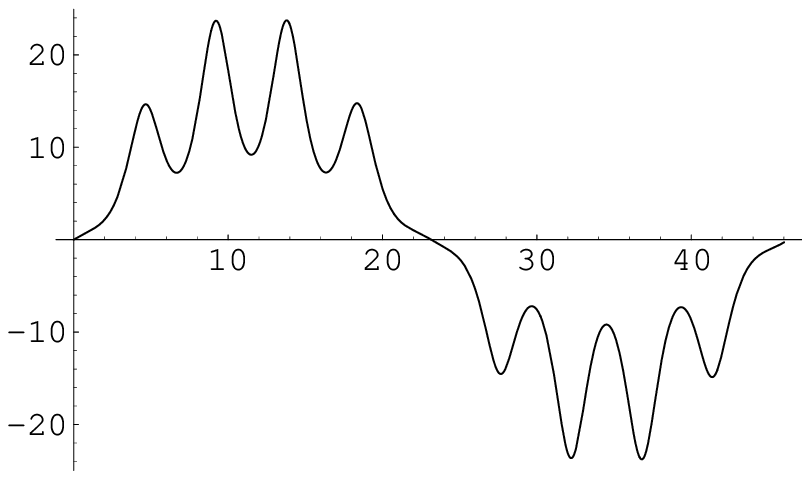}
  \includegraphics[scale=0.35]{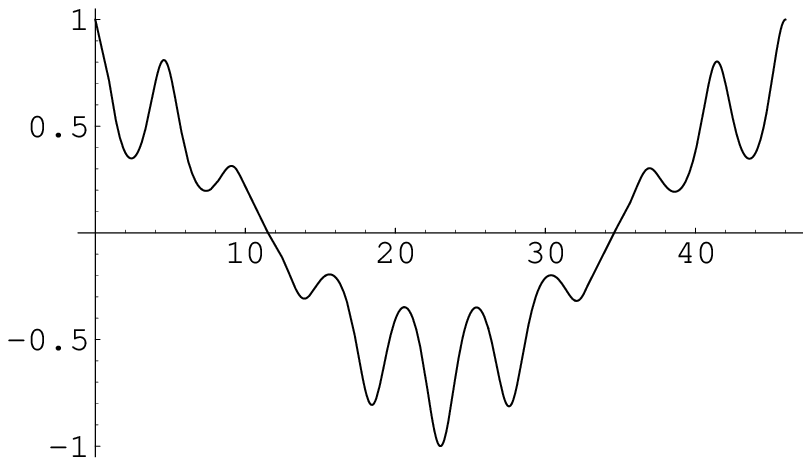}
  \includegraphics[scale=0.35]{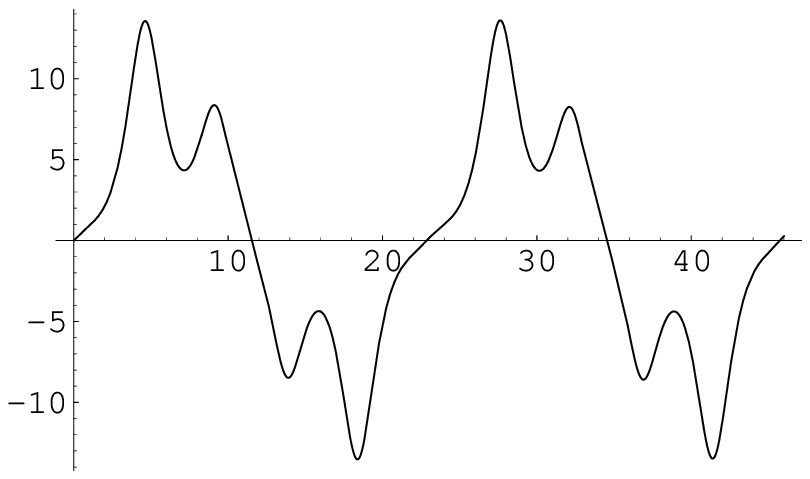}
   \includegraphics[scale=0.35]{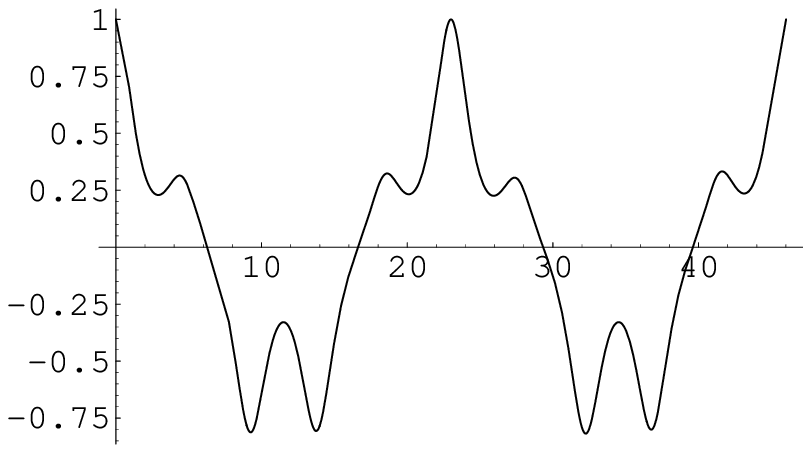}
  \includegraphics[scale=0.35]{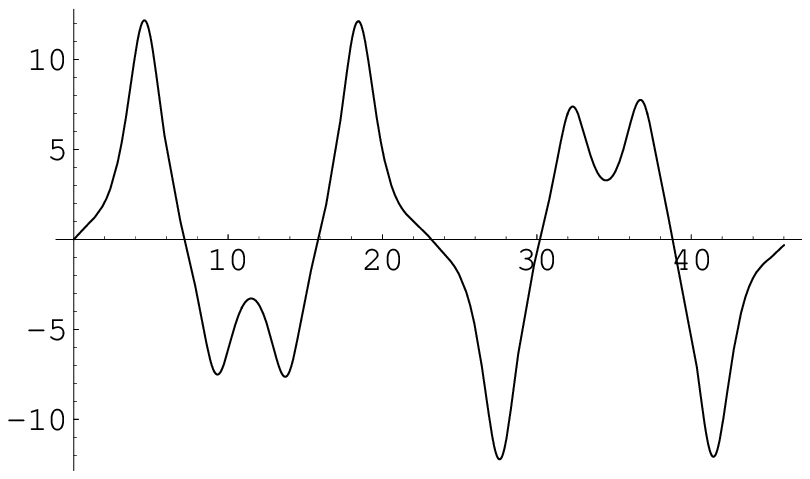}
  \includegraphics[scale=0.35]{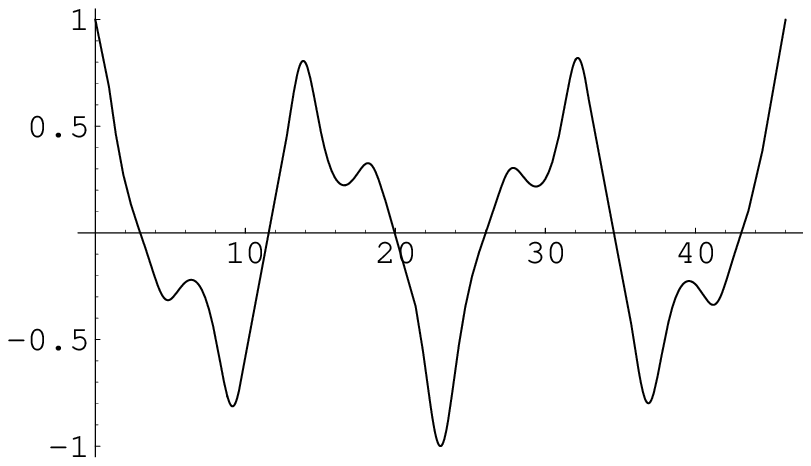}
  \includegraphics[scale=0.35]{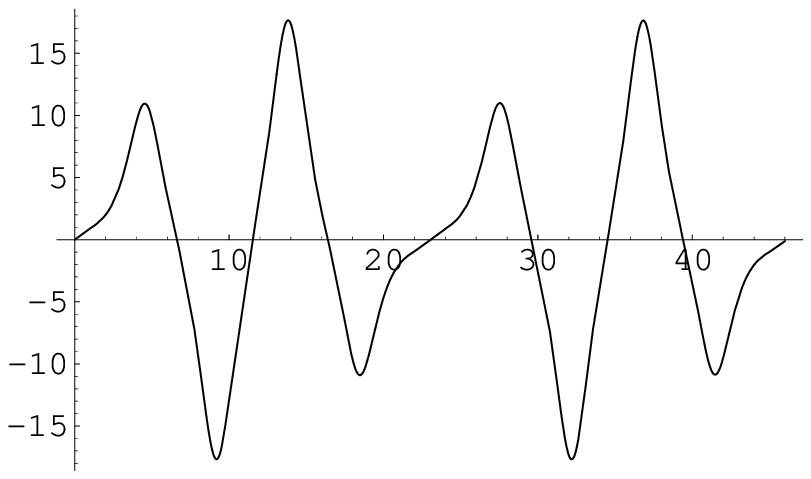}
   \includegraphics[scale=0.35]{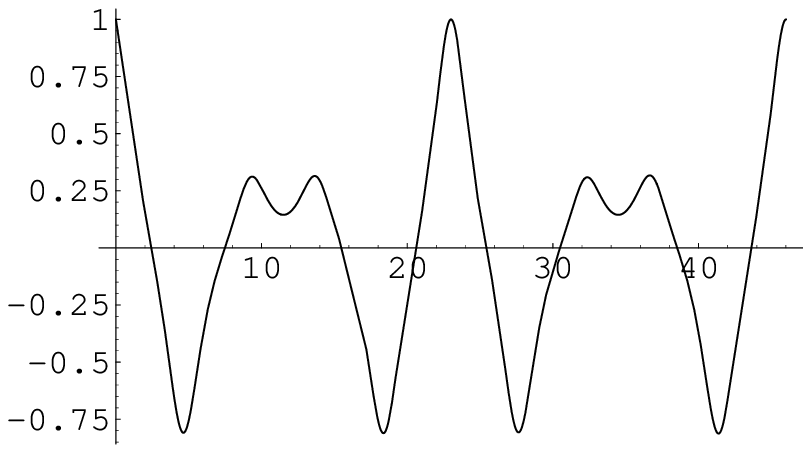}
  \includegraphics[scale=0.35]{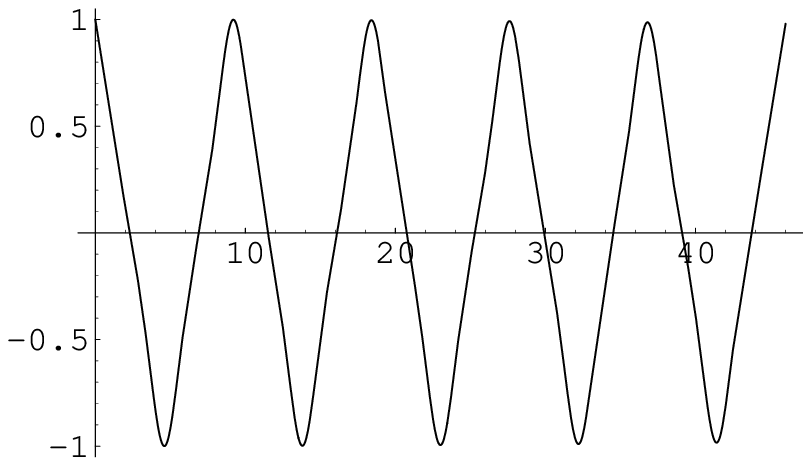}
  \includegraphics[scale=0.35]{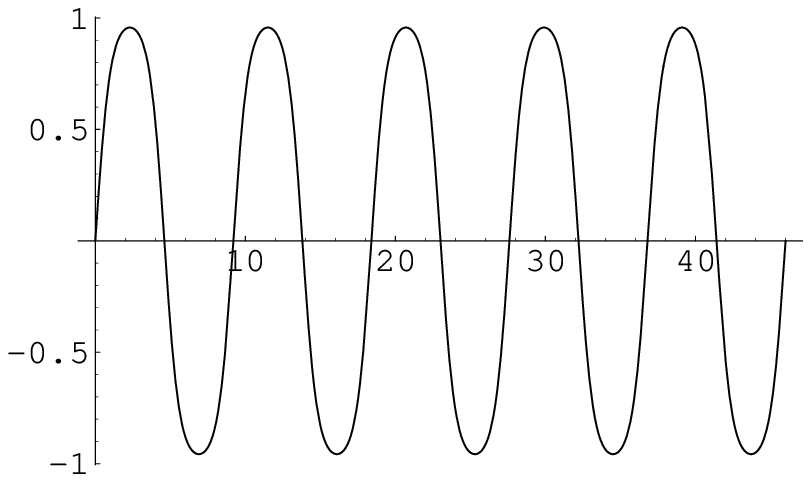}
  \caption{Eigenfunctions associated to the eigenvalues $\lambda_{1,0}$, ..., 
           $\lambda_{10,0}$, $\lambda_{11,0}=0$ of the surface $N_{7}$.}
  \label{Y}
  \end{figure}

\begin{figure}[phbt]
  \centering
  \includegraphics[scale=0.65]{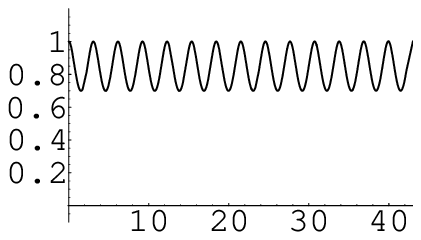}
  \includegraphics[scale=0.65]{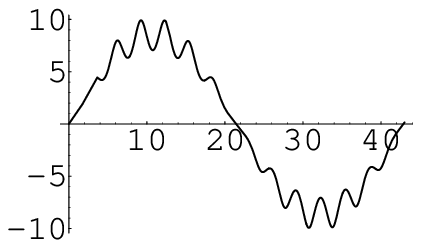}
  \includegraphics[scale=0.65]{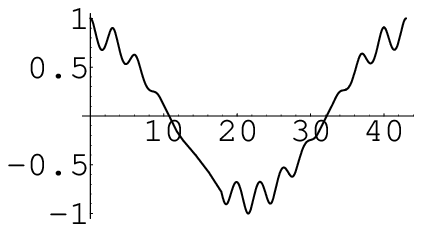}
  \includegraphics[scale=0.65]{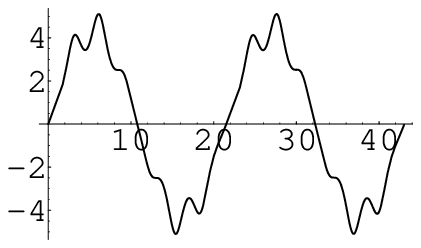}
    \includegraphics[scale=0.65]{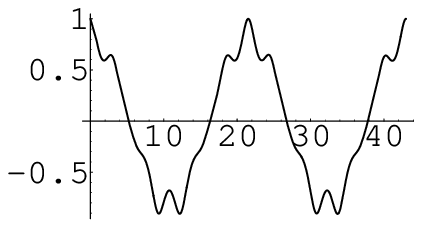}
  \includegraphics[scale=0.65]{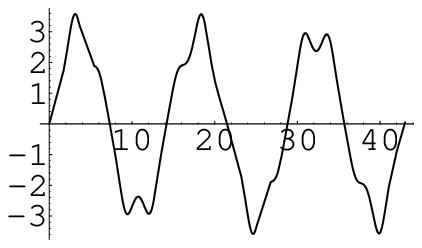}
  \includegraphics[scale=0.65]{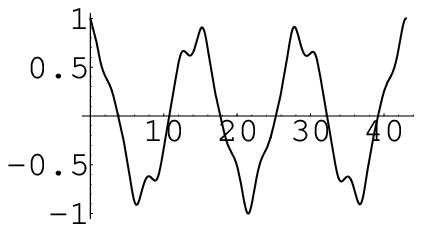}
  \includegraphics[scale=0.65]{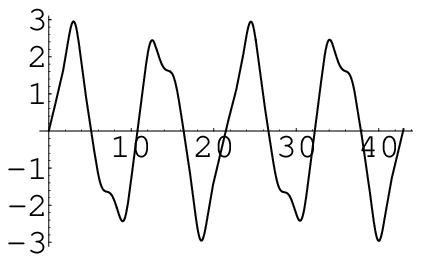}
 \includegraphics[scale=0.65]{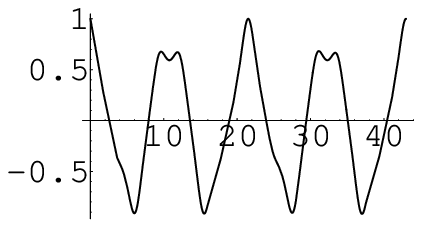}
  \includegraphics[scale=0.65]{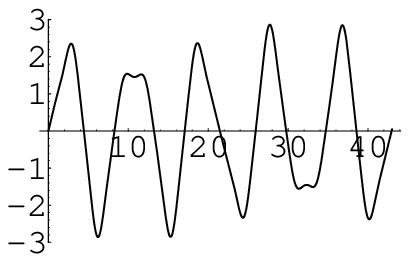}
  \includegraphics[scale=0.65]{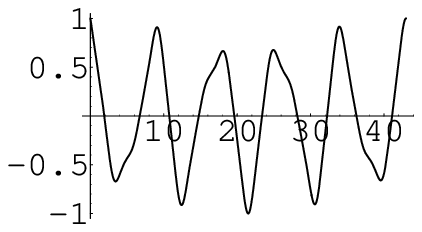}
  \includegraphics[scale=0.65]{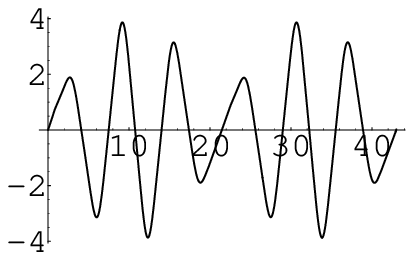}
   \includegraphics[scale=0.65]{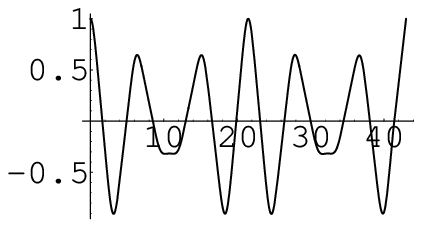}
  \includegraphics[scale=0.65]{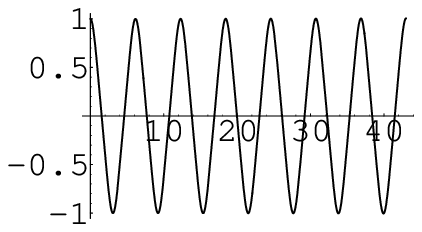}
  \includegraphics[scale=0.65]{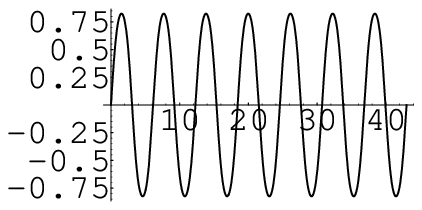}
 \caption{Eigenfunctions associated to the eigenvalues $\lambda_{1,0}$, ..., 
           $\lambda_{14,0}$, $\lambda_{15,0}=0$ of the surface $N_{8}$.}
  \label{Z}
\end{figure}

\begin{figure}[phbt]
  \centering
  \includegraphics[scale=0.35]{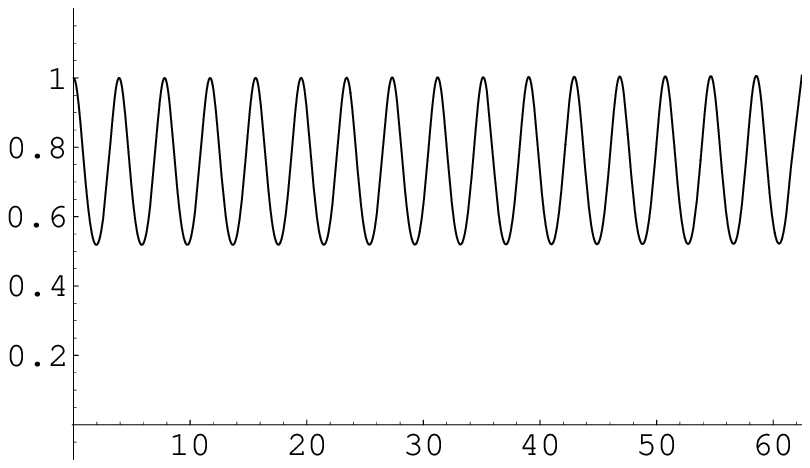}
  \includegraphics[scale=0.35]{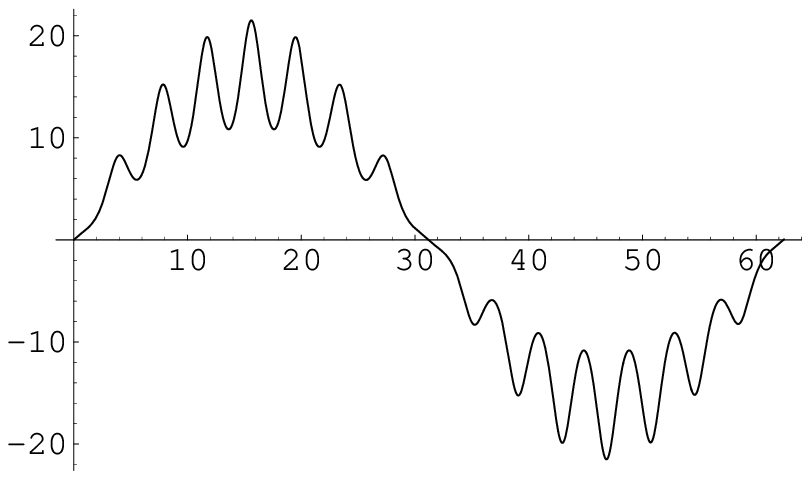}
  \includegraphics[scale=0.35]{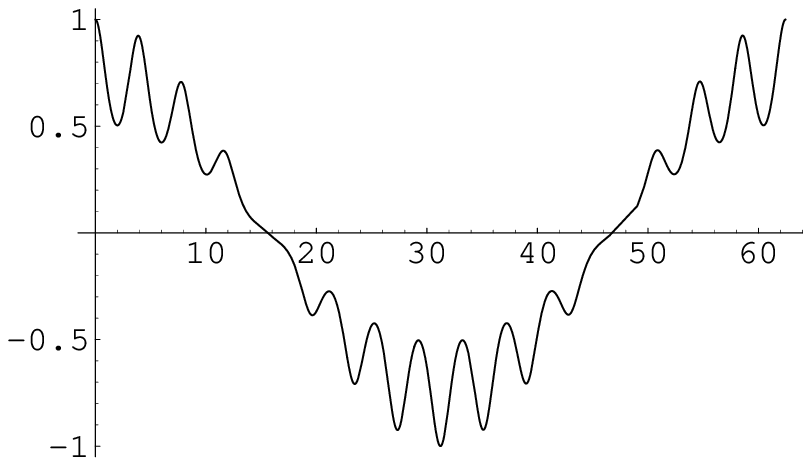}
  \includegraphics[scale=0.35]{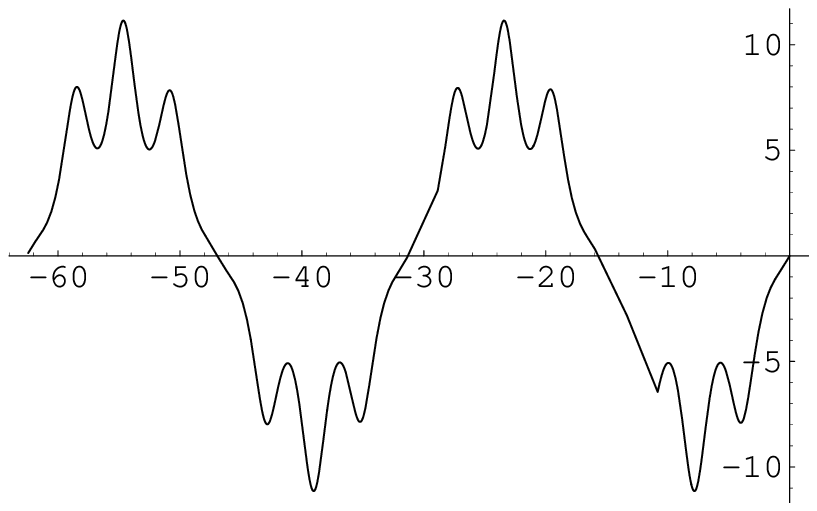}
   \includegraphics[scale=0.35]{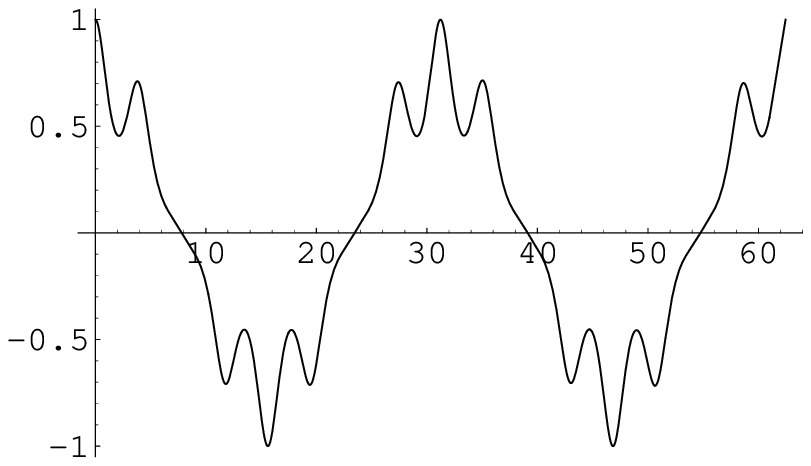}
  \includegraphics[scale=0.35]{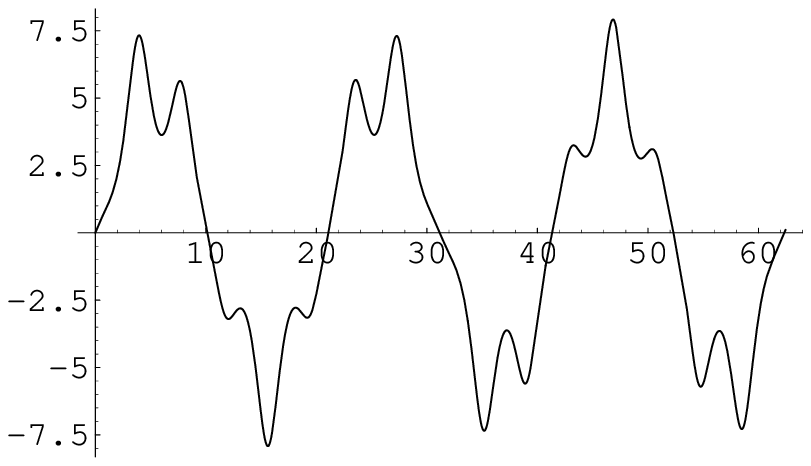}
  \includegraphics[scale=0.35]{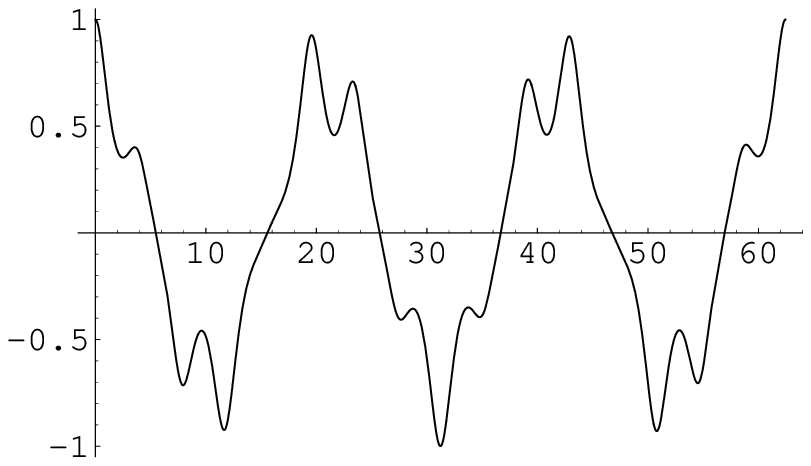}
  \includegraphics[scale=0.35]{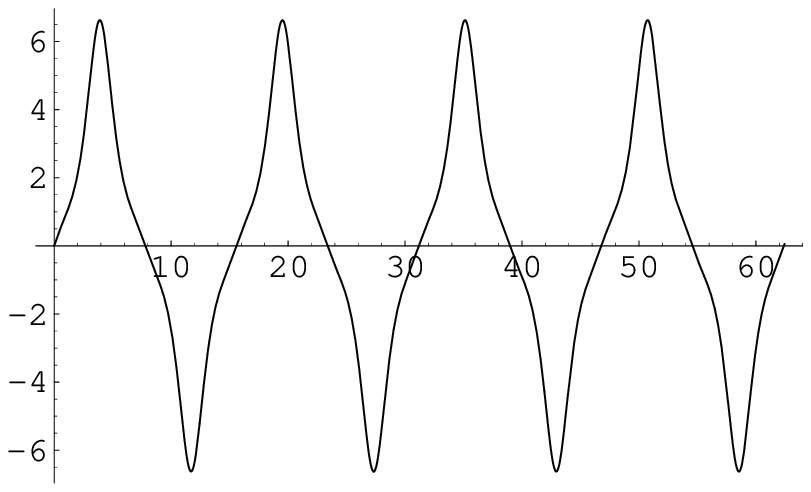}
   \includegraphics[scale=0.35]{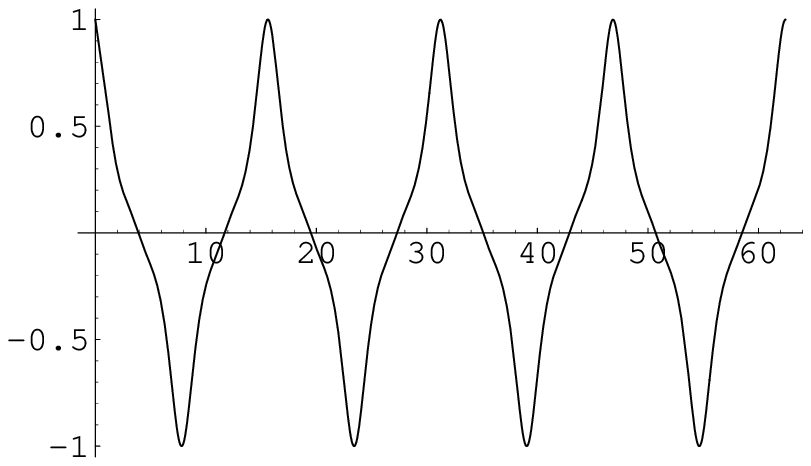}
  \includegraphics[scale=0.35]{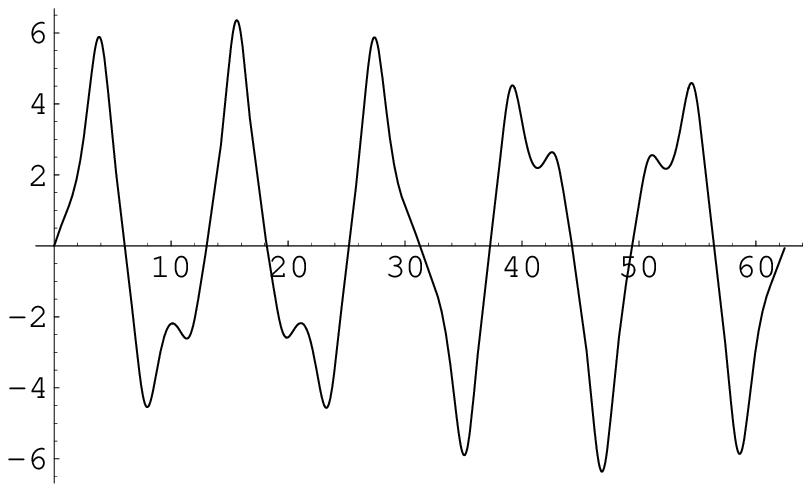}
  \includegraphics[scale=0.35]{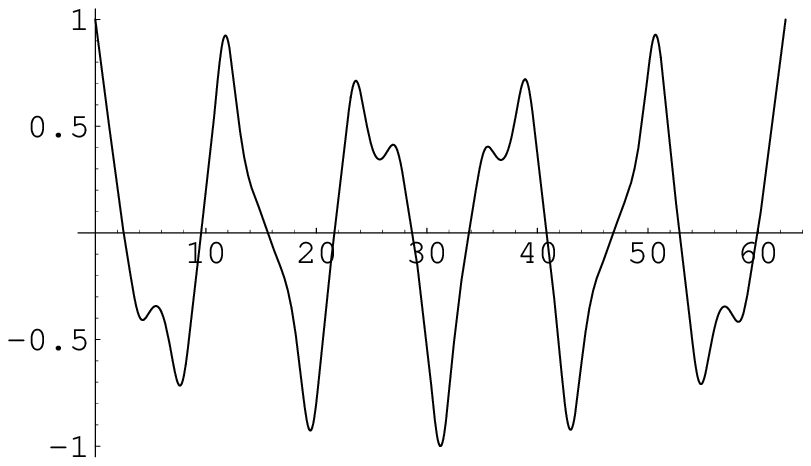}
  \includegraphics[scale=0.35]{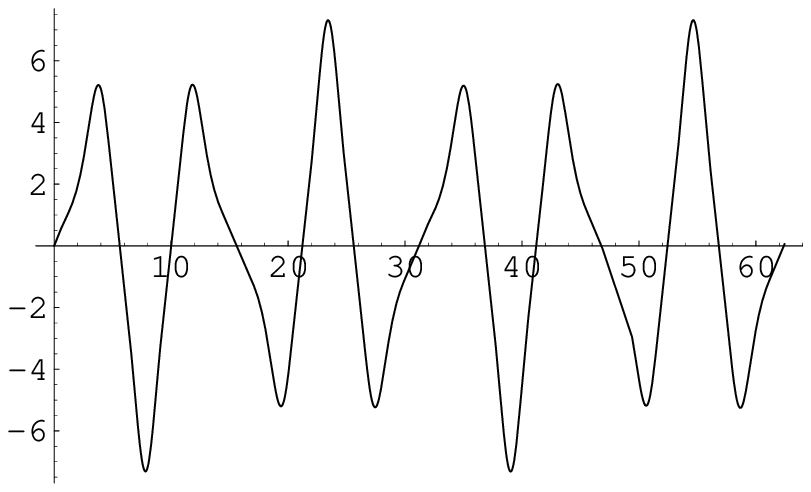}
   \includegraphics[scale=0.35]{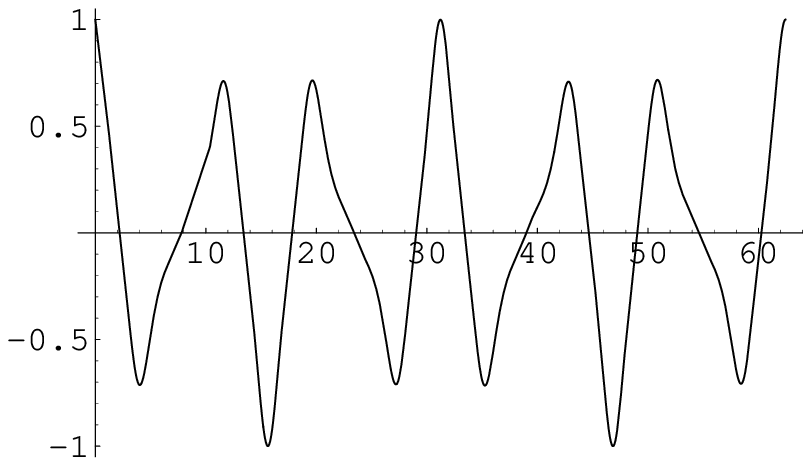}
   \includegraphics[scale=0.35]{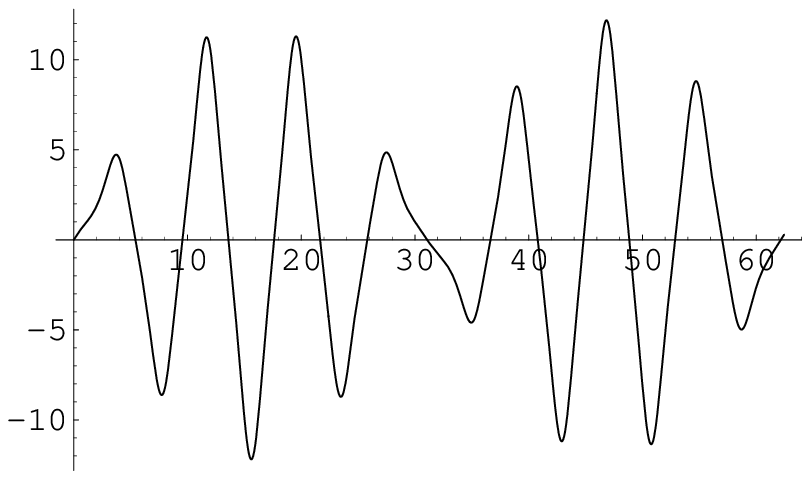}
   \includegraphics[scale=0.35]{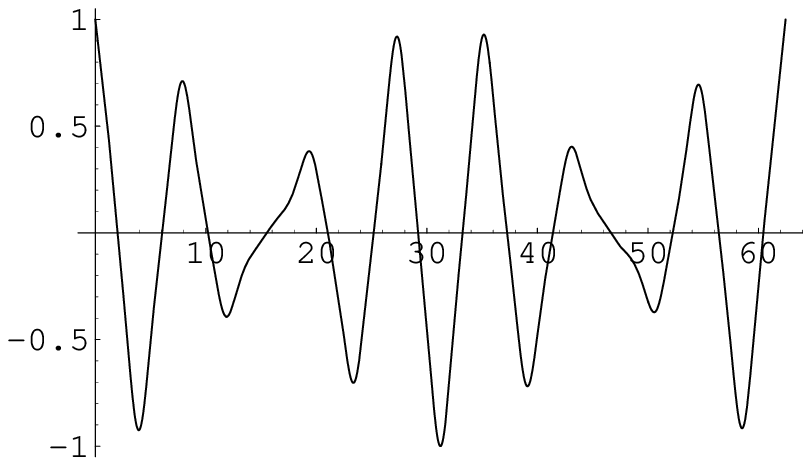}
   \includegraphics[scale=0.35]{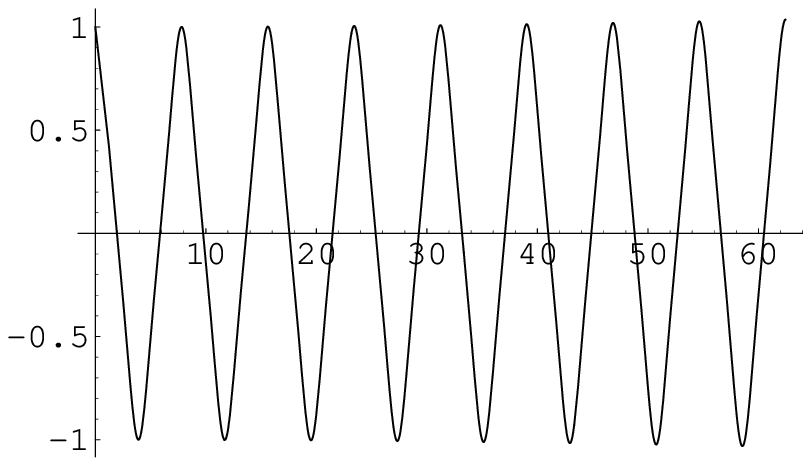}
   \includegraphics[scale=0.35]{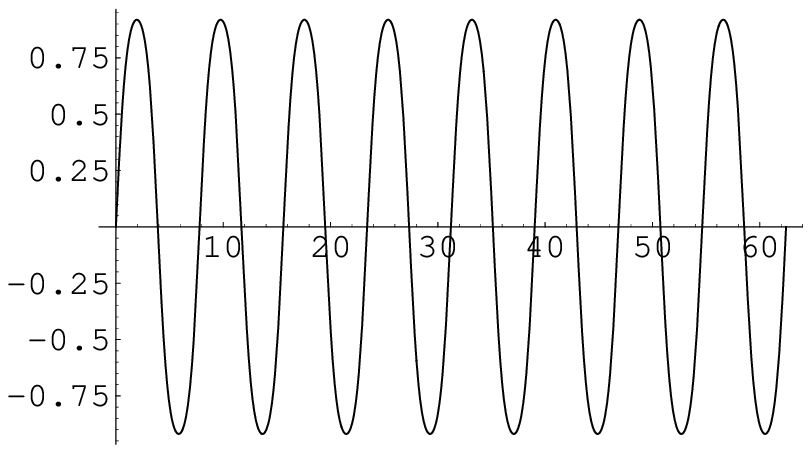}
  \caption{Eigenfunctions associated to the eigenvalues $\lambda_{1,0}$, ..., 
           $\lambda_{16,0}$, $\lambda_{17,0}=0$ of the surface $N_{9}$.}
  \label{Za}
  \end{figure}

\begin{figure}[phbt]
  \centering
  \includegraphics[scale=0.65]{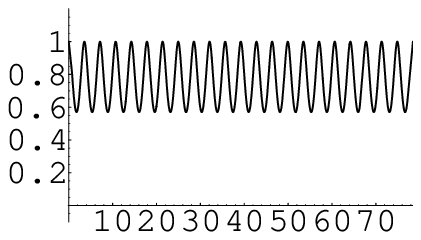}
  \includegraphics[scale=0.65]{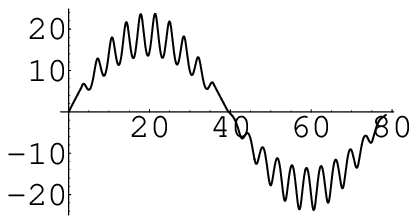}
  \includegraphics[scale=0.65]{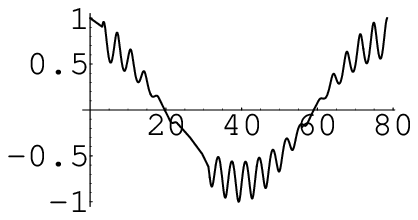}
  \includegraphics[scale=0.65]{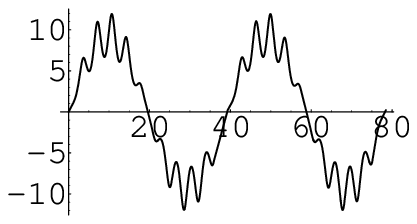}
   \includegraphics[scale=0.65]{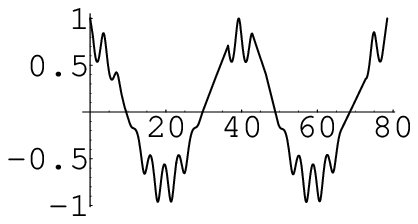}
  \includegraphics[scale=0.65]{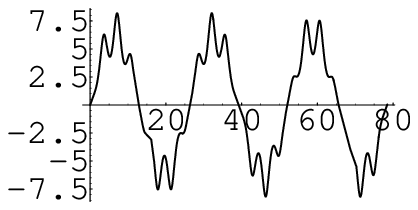}
  \includegraphics[scale=0.65]{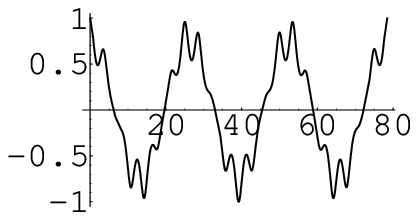}
  \includegraphics[scale=0.65]{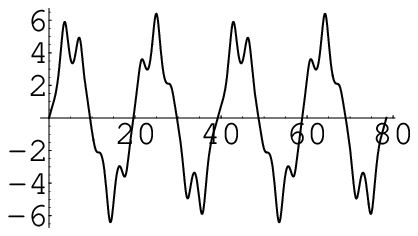}
  \includegraphics[scale=0.65]{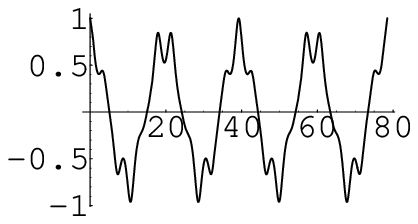}
  \includegraphics[scale=0.65]{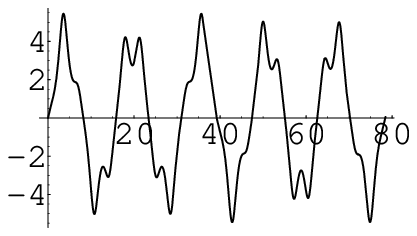}
  \includegraphics[scale=0.65]{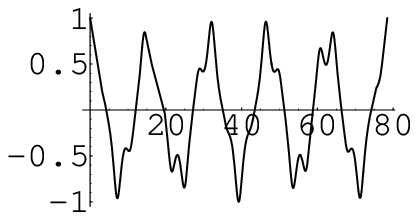}
  \includegraphics[scale=0.65]{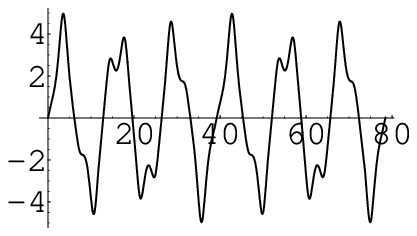}
   \includegraphics[scale=0.65]{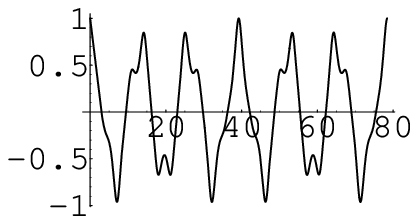}
  \includegraphics[scale=0.65]{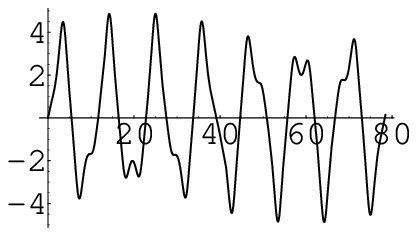}
  \includegraphics[scale=0.65]{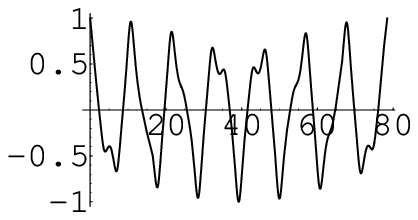}
  \includegraphics[scale=0.65]{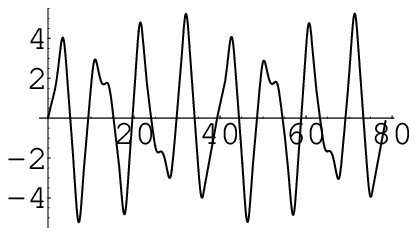}
  \includegraphics[scale=0.65]{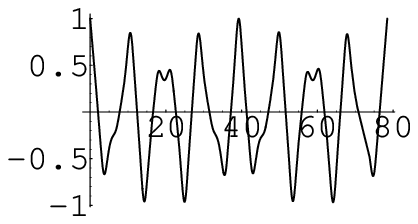}
  \includegraphics[scale=0.65]{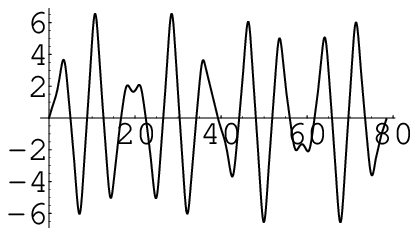}
  \includegraphics[scale=0.65]{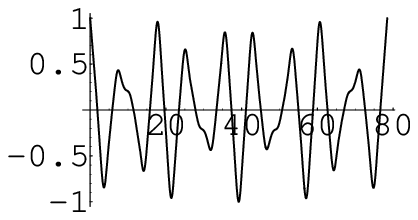}
  \includegraphics[scale=0.65]{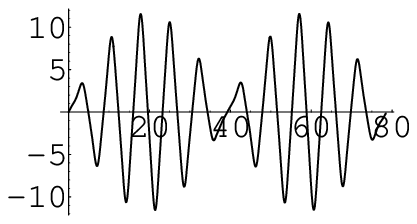}
  \includegraphics[scale=0.65]{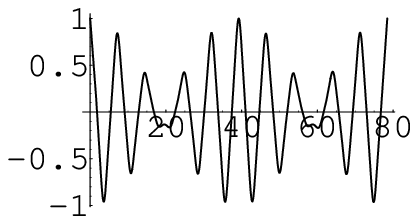}
  \includegraphics[scale=0.65]{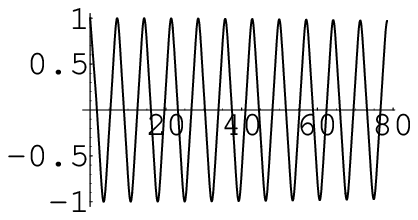}
  \includegraphics[scale=0.65]{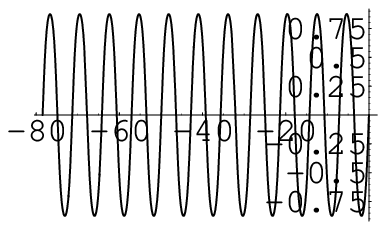}
 \caption{Eigenfunctions associated to the eigenvalues $\lambda_{1,0}$, ..., 
           $\lambda_{22,0}$, $\lambda_{23,0}=0$ of the surface $N_{10}$.}
  \label{Zb}
\end{figure}

\begin{figure}[phbt]
  \centering
  \includegraphics[scale=0.35]{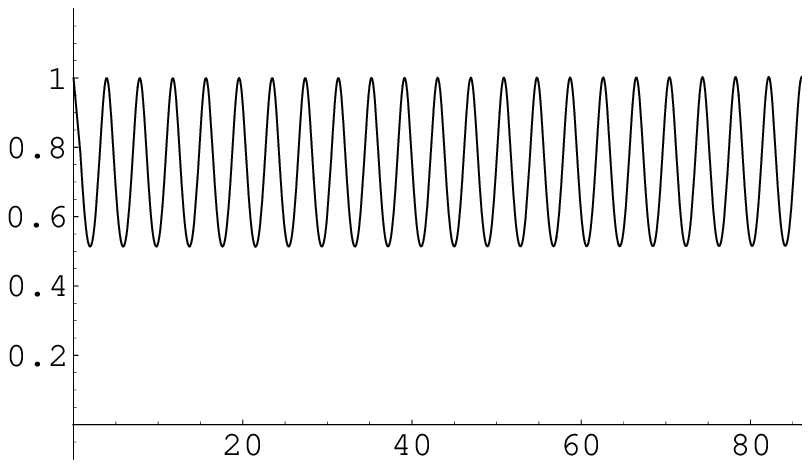}
  \includegraphics[scale=0.35]{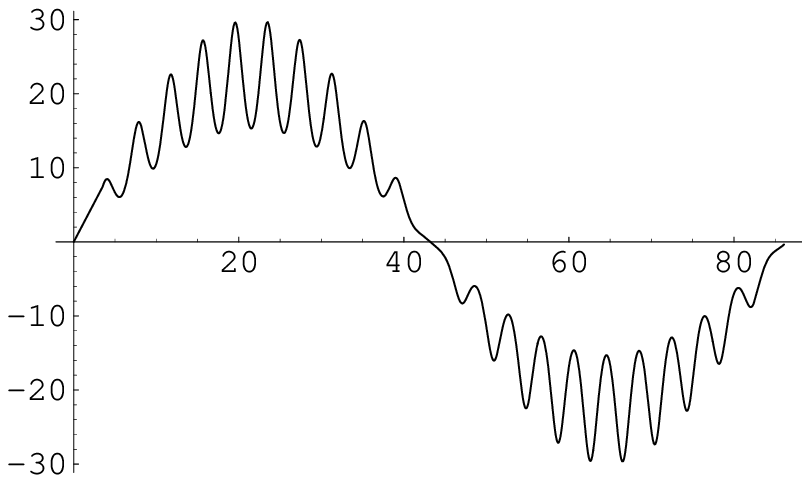}
  \includegraphics[scale=0.35]{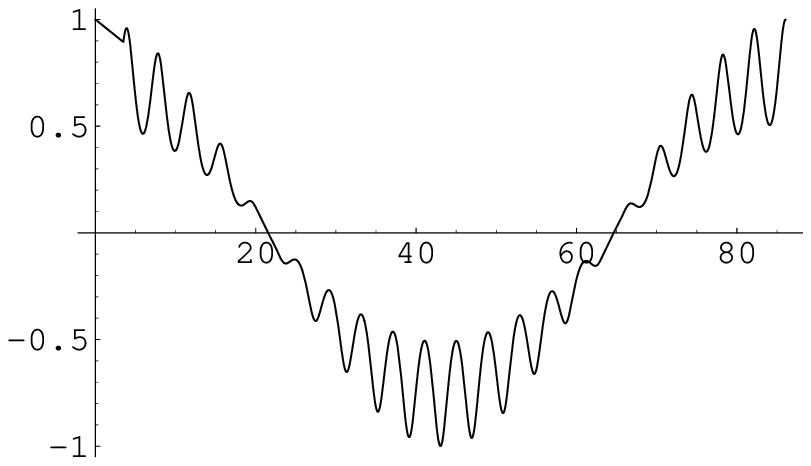}
  \includegraphics[scale=0.35]{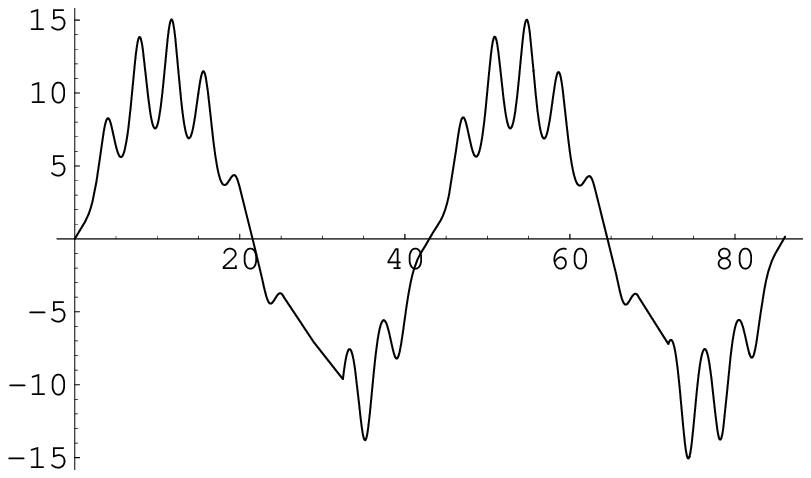}
   \includegraphics[scale=0.35]{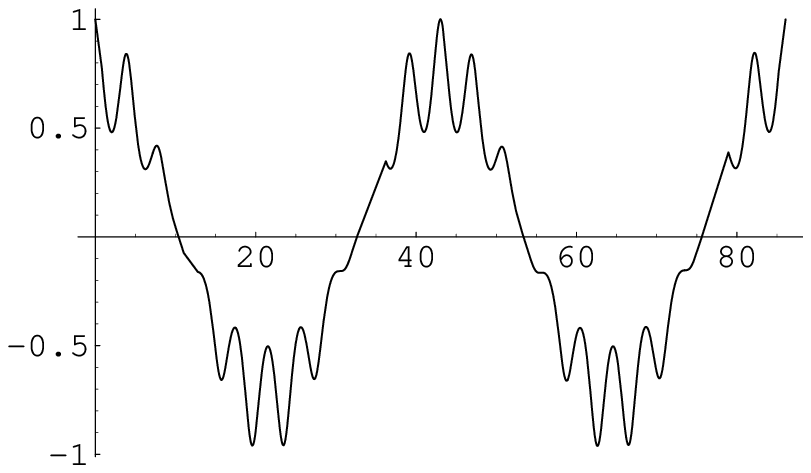}
  \includegraphics[scale=0.35]{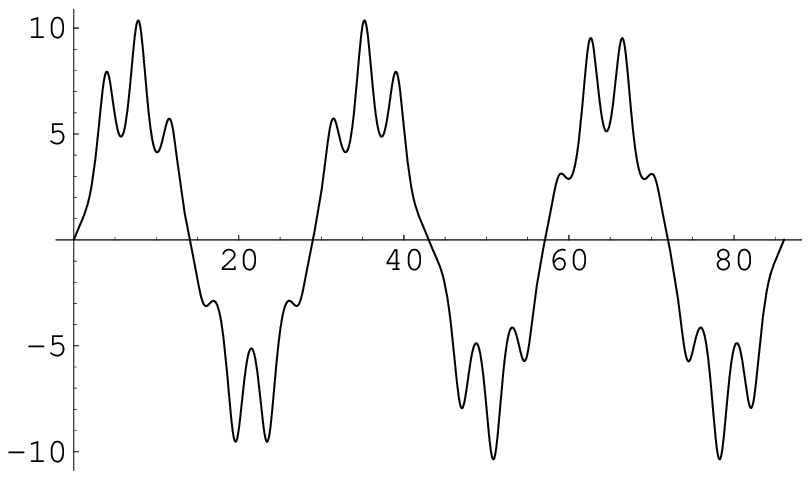}
  \includegraphics[scale=0.35]{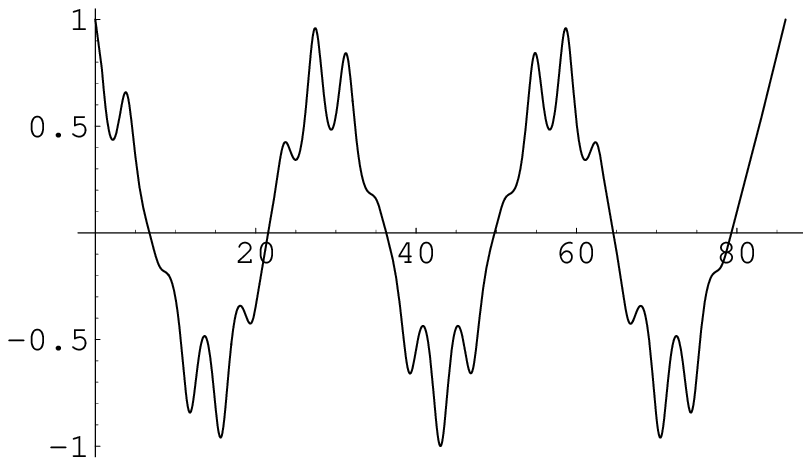}
  \includegraphics[scale=0.35]{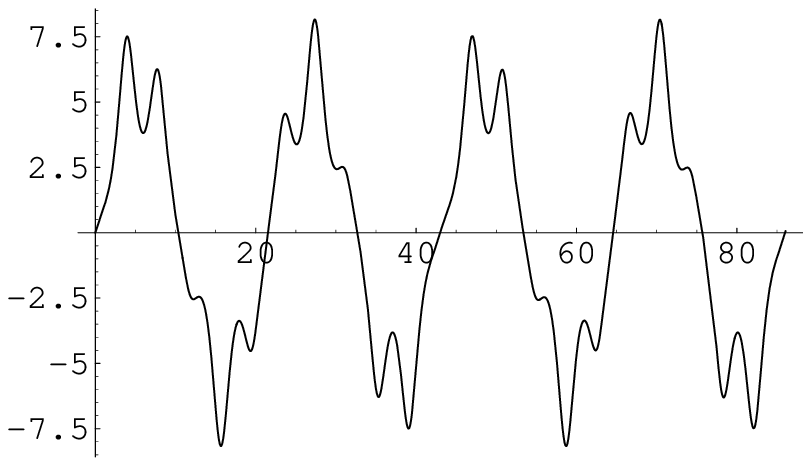}
   \includegraphics[scale=0.35]{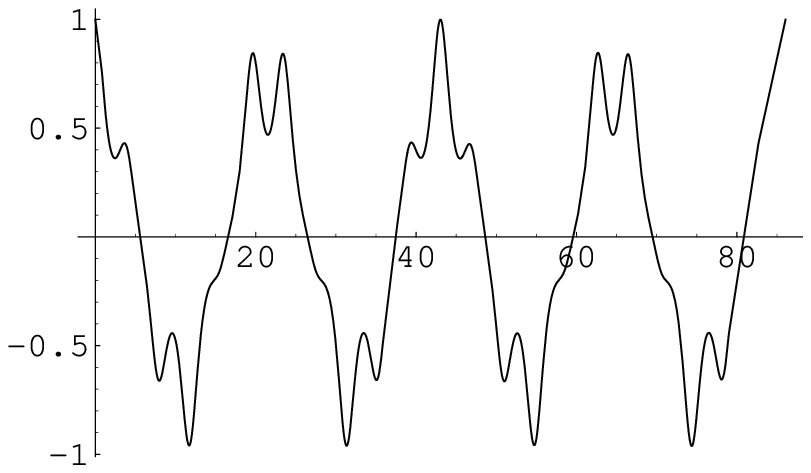}
  \includegraphics[scale=0.35]{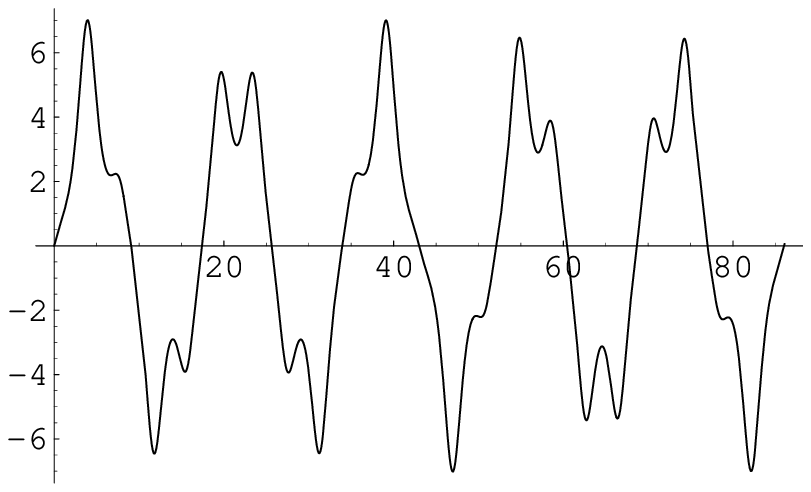}
  \includegraphics[scale=0.35]{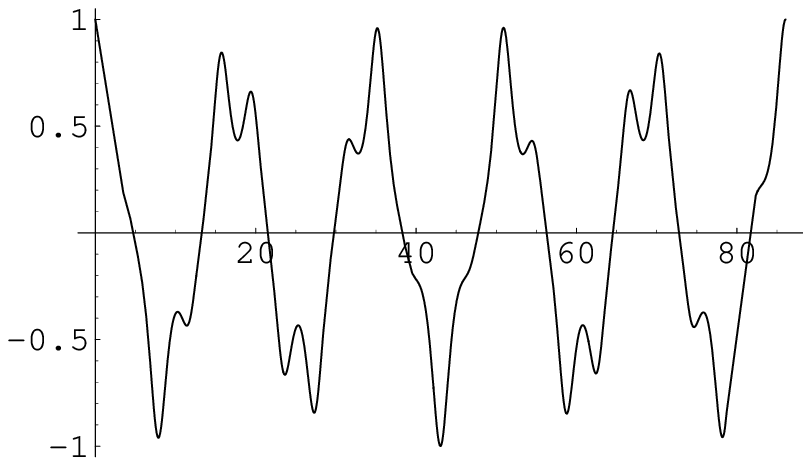}
  \includegraphics[scale=0.35]{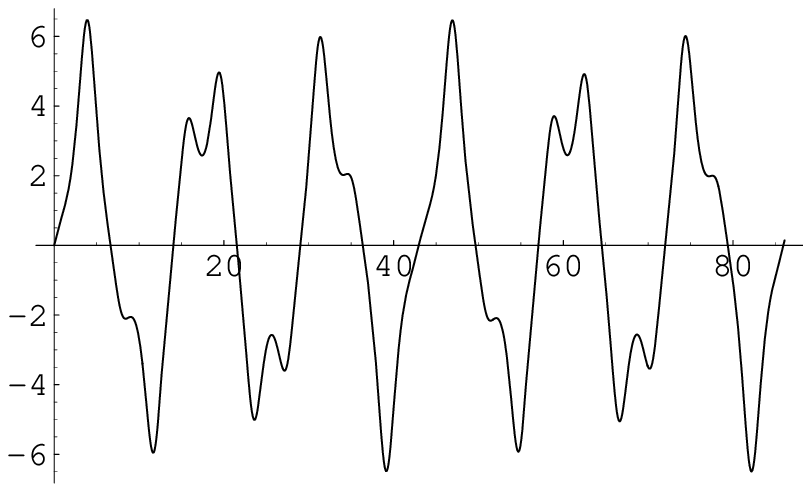}
   \includegraphics[scale=0.35]{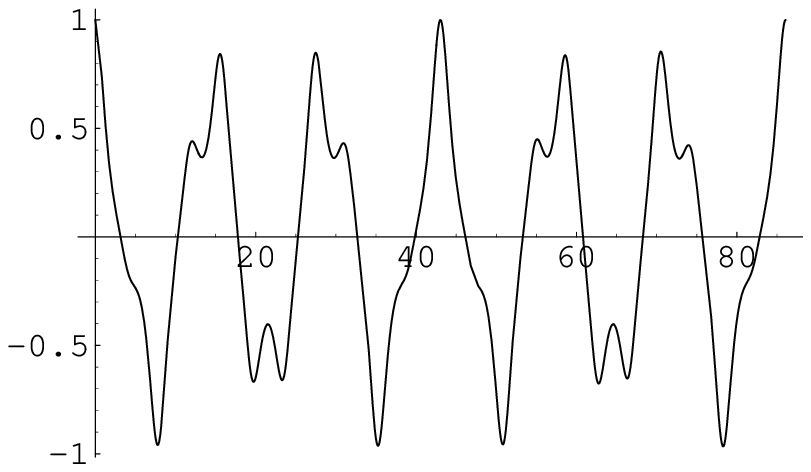}
   \includegraphics[scale=0.35]{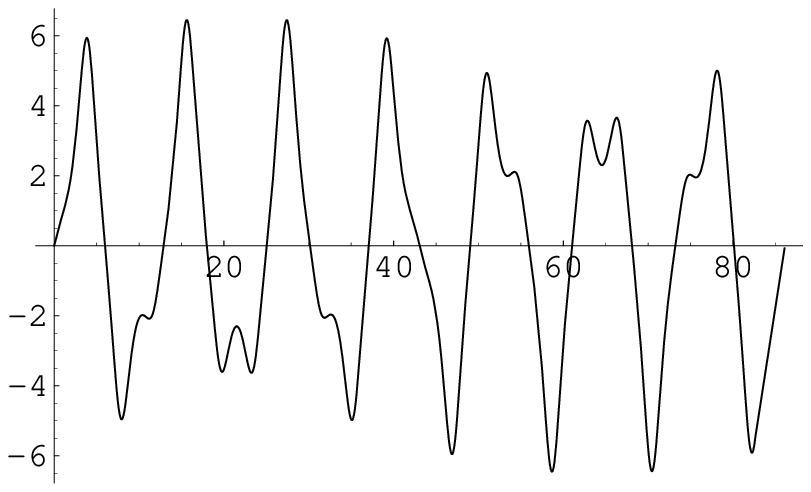}
   \includegraphics[scale=0.35]{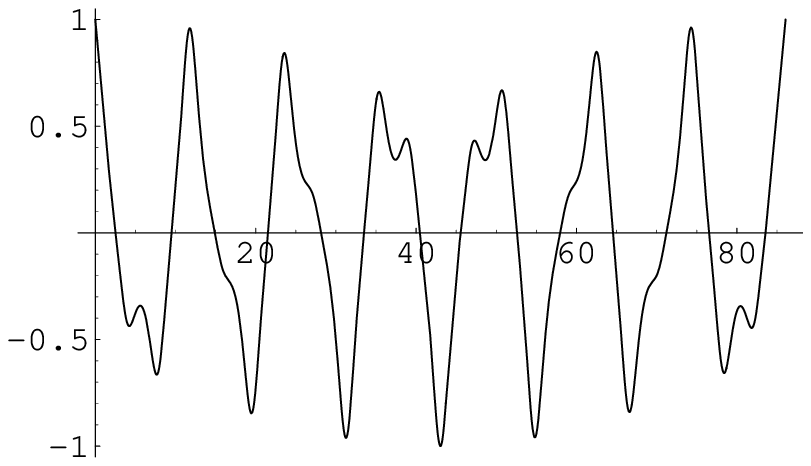}
   \includegraphics[scale=0.35]{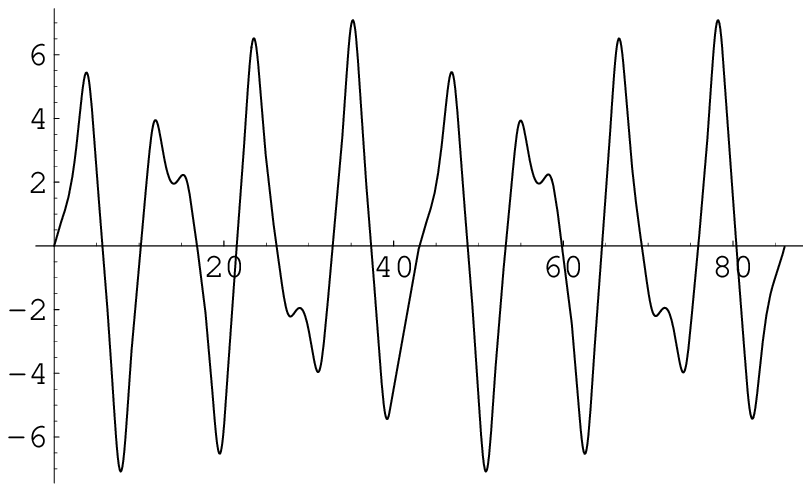}
   \includegraphics[scale=0.35]{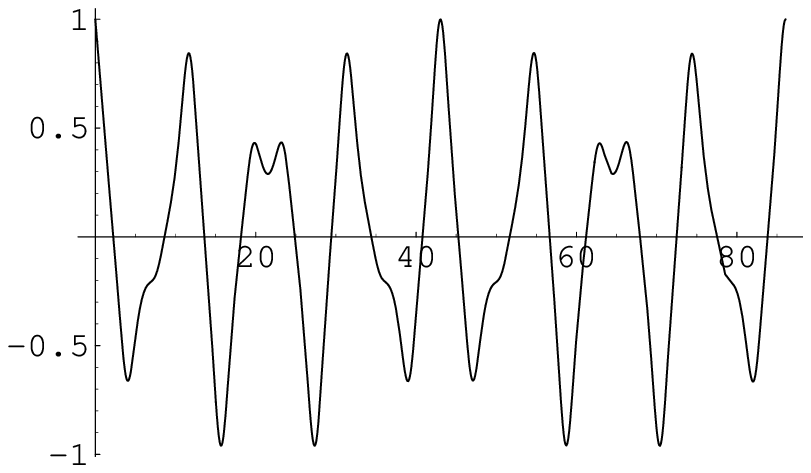}
   \includegraphics[scale=0.35]{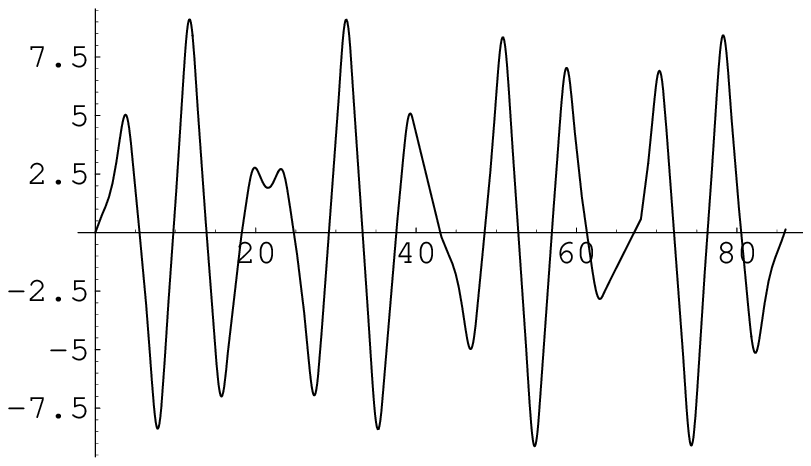}
   \includegraphics[scale=0.35]{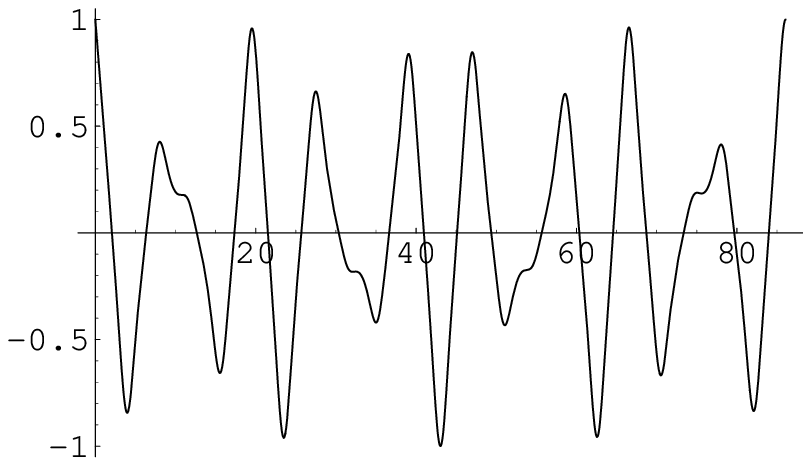}
   \includegraphics[scale=0.35]{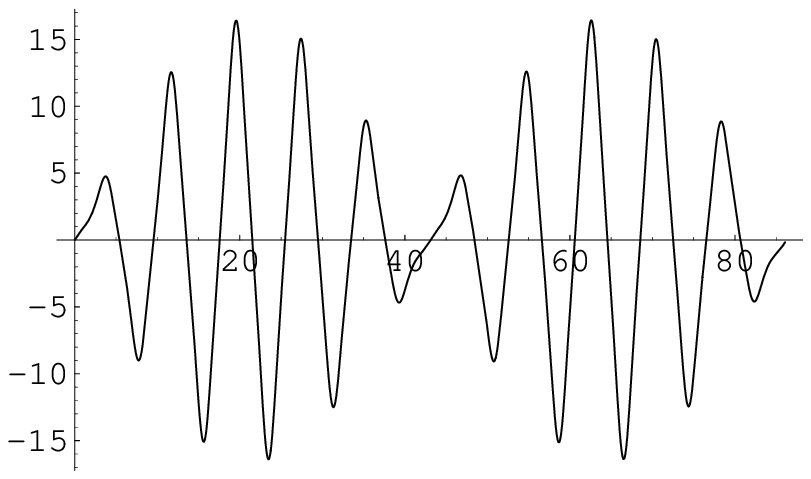}
   \includegraphics[scale=0.35]{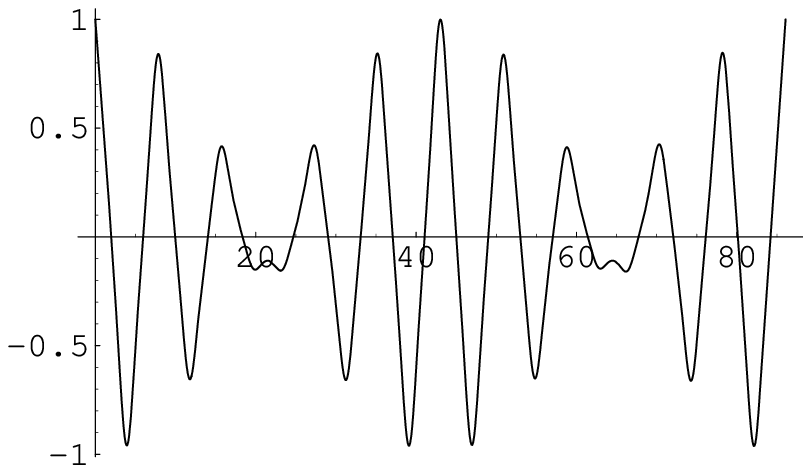}
   \includegraphics[scale=0.35]{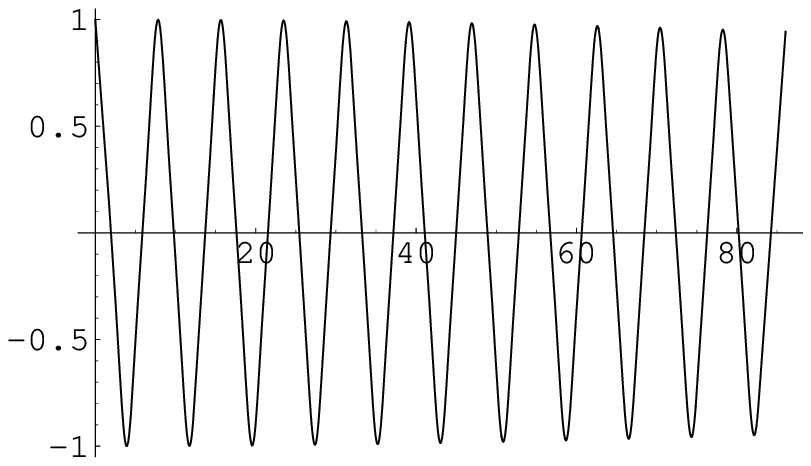}
   \includegraphics[scale=0.35]{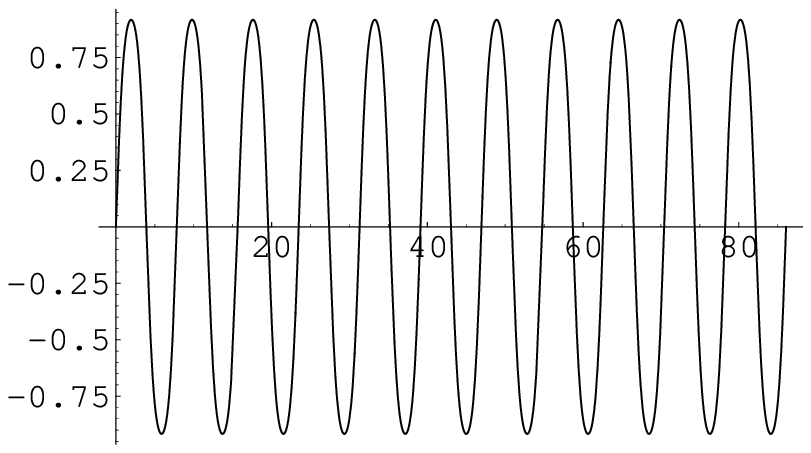}
  \caption{Eigenfunctions associated to the eigenvalues $\lambda_{1,0}$, ..., 
           $\lambda_{22,0}$, $\lambda_{23,0}=0$ of the surface $N_{11}$.}
  \label{Zc}
  \end{figure}

Tables \ref{table2}, \ref{table1} and Figures \ref{A}-\ref{Zc} show results 
of our numerical computations.  

\end{document}